\documentclass{amsart}
\usepackage{amssymb}
\usepackage[hyperref]{degt}

\def\2{\color{red}}
\def\itemrefform#1{$(#1)$}

\input NL.def

\def\sizemax{144}
\def\sizesub{132}

\def\minusform#1#2{\frac{\bar{#1}}{#2}}

\def\ltitle{\subsection{The lattice \lname{\thelattice}}}
\def\ltext{%
 There are \lcount{\thelattice} configurations to be considered,
 and the maximal na\"{\i}ve bound is $\bnd(\Orb)=\lmax{\thelattice}$%
 \iffullversion.\else\ (see \autoref{tab.\thelattice}).\fi
}
\def\latticetext{\iffullversion\ltitle\ltext\else
 \ifnum\lmax{\thelattice}>\mincount\ltitle\ltext\fi\fi
}
\def\atbeginorbits{\iffullversion\else\table[t]
 \caption{The lattice \lname{\thelattice}}\label{tab.\thelattice}%
 \hrule height0pt\relax\vskip-\abovedisplayskip\fi}
\def\atendorbits{\iffullversion\else\vskip-\belowdisplayskip\endtable\fi}

\def\specvec#1{\bar{#1}}

\let\rbar\rt

\def\bh{\specvec{h}}
\def\ortu{\specvec{u}}
\def\ortv{\specvec{v}}
\def\ortw{\specvec{w}}
\def\vv#1{\specvec1_{#1}}

\def\cw#1{[#1]}

\def\units{^\times}
\def\sdif{\mathbin\vartriangle}
\def\IS{\Cal I}
\def\CL{\Cal L}
\def\fB{\frak B}

\def\minitab#1{\vcenter\bgroup\rm
 \def\-##1{\setbox0\hbox{$00$}\hbox to\wd0{\hss$##1$\hss}}%
 \let\\\cr
 \halign\bgroup\strut\hss$##$&&#1\hss$##$\hss\cr}
\def\endminitab{\crcr\egroup\egroup}

\def\realcurve{\hyperref[item.real.curve]{^*}}

\title{Tritangents to smooth sextic curves}

\author{Alex Degtyarev}

\address{%
Department of Mathematics\\
Bilkent University\\
06800 Ankara, TURKEY}

\email{
degt@fen.bilkent.edu.tr}

\thanks{%
The author was partially supported by the T\"{U}B\DOTaccent{I}TAK grant 118F413}



\begin{document}

\begin{abstract}
We prove that a smooth plane sextic curve can have at most 72 tritangents,
whereas a smooth real sextic may have at most 66 real tritangents.
\end{abstract}

\maketitle

\section{Introduction}\label{S.intro}

All algebraic varieties considered in the paper are over~$\C$. On a few occasions,
we discuss subfields $\Bbbk\subset\C$ (most notably, $\Bbbk=\R$:
by definition,
a \emph{real} variety is a complex one equipped with an anti-holomorphic
involution), but we \emph{never} consider fields of positive characteristic.

\subsection{Principal results}\label{s.results}
This paper concludes the study of the maximal
number of straight lines in a
smooth polarized $K3$-surface. The most classical case, \viz. that of spatial
quartics $X\subset\Cp3$, goes as far back as to
A.~Clebsch~\cite{Clebsch} (the upper bound of $m(11m-24)$ lines in a smooth
degree~$m$ surface $X\subset\Cp3$, yielding at most $80$ lines in a quartic) and
F.~Schur~\cite{Schur:quartics}
(an example of a smooth quartic with $64$ lines). The sharp upper bound of
$64$ lines was established in B.~Segre~\cite{Segre}.
Two recent papers aroused new interest to this classical problem: first,
a minor gap in Segre's
proof was discovered and corrected by S.~Rams and
M.~Sch\"{u}tt~\cite{rams.schuett}, where the argument was also extended to
all characteristics $p\ge5$ of the ground field,
and second, an alternative, purely arithmetical (or rather lattice theoretic)
proof was given in~\cite{DIS}; this approach allowed us to obtain,
in addition to mere upper bounds, a complete classification of all large
(\ie, more than~$52$) configurations of lines, prove the uniqueness of
the line maximizing quartic, and establish the sharp upper bound of
$56$ \emph{real} lines in a \emph{real}
smooth quartic surface (also realized by a unique real quartic).
A great deal of other works have appeared almost immediately and, for the
moment,
the case of spatial quartic surfaces still
remains the best studied one: there are sharp upper bounds on the
number of lines over algebraically closed fields of positive characteristic
(see
\cite{degt:supersingular,rams.schuett:char3,rams.schuett,rams.schuett:char2})
and over~$\R$ (see~\cite{DIS}), partial bounds over~$\Q$
(see~\cite{degt:singular.K3}), upper bounds for singular quartics, both $K3$
(see \cite{Veniani:char2,Veniani}) and not (see~\cite{Rams.Gonzales}),
explicit equations of quartics with many lines (see
\cite{DIS,Shimada:X56,Veniani:equations}), \etc.

Lines in smooth polarized $K3$-surfaces $X\to\Cp{d+1}$
of all degrees $2d\ge4$, both birational and hyperelliptic
(\cf.~\cite{Saint-Donat}),
were studied
in~\cite{degt:lines},
using an arithmetical reduction similar to~\cite{DIS}
and an appropriate taxonomy of prospective Fano graphs.
(It appears that, so far,
the more conventional geometric arguments have failed to
produce even reasonable bounds on the number of lines.)
Among other results, found in~\cite{degt:lines} are sharp
upper bounds on the number of lines, both over~$\C$ and~$\R$,
and a complete description
of all large configurations of lines, especially in the two most ``classical''
cases, \viz. sextics in~$\Cp4$ and octics (most notably, triquadrics)
in~$\Cp5$, which give rise to a number of interesting Fano graphs.
An unexpected discovery is the fact that the configurations of lines simplify
dramatically when the degree grows: asymptotically, for $2d\gg0$, all lines
are either linearly independent in $H_2(X)$ or, else, among the
fiber components of a certain fixed elliptic pencil;
in either case, their number does not
exceed~$24$.
(The true sharp bound oscillates between $21$, $22$, and~$24$, periodically
in the degree $2d\gg0$, whereas its real counterpart oscillates between $19$, $20$,
and~$21$, with a larger period.)
The other side of the coin is the fact that in the remaining ``classical''
case, the smallest degree $2d=2$ (double planes),
the dual adjacency graph of lines may be too large: the star of a single
vertex is more complicated than the whole Fano graph of a quartic.
For this reason, the case $2d=2$ was left out
as not feasible in~\cite{degt:lines}; it is
treated in the present paper by means of
considerably different arithmetical techniques
(see \autoref{s.contents}).

In~\cite{degt:singular.K3} it was conjectured that the maximal number of
lines in a smooth $2$-polarized $K3$-surface is $144$, with the maximum
realized by the double plane $X\to\Cp2$ ramified over the sextic curve
\[
z_0^6+z_1^6+z_2^6=10(z_0^3z_1^3+z_1^3z_2^3+z_2^3z_0^3).
\label{eq.max.sextic}
\]
(This equation is borrowed from Sh.~Mukai~\cite{Mukai}, as the surface in
question admits a faithful action of the Mukai group~$M_9$; explicit
equations of the predicted $144$ lines were found independently by D.~Festi and
Y.~Zaytman, private communication.)
The conjecture is motivated by the fact that, like Schur's
quartic~\cite{Schur:quartics} and some line maximizing sextics in~$\Cp4$ and
octics in~$\Cp5$ (see \cite{degt:lines,degt:singular.K3}), this surface
minimizes the discriminant of a singular $K3$-surface admitting a smooth
model of a given degree.
In the present paper we settle (in the affirmative) and extend the
conjecture, see \autoref{th.main}, \autoref{ad.main}, and \autoref{th.real}.
We state our principal results in terms of tritangents to the ramification
locus $C\subset\Cp2$ (a smooth sextic curve) rather than lines in the surface
$X\to\Cp2$, dividing the numbers by~$2$
(see \autoref{s.K3} below for further details).
Certainly, when speaking about tritangents, we allow the collision of some of
the tangency points; in other words, a \emph{tritangent} to a smooth sextic
$C\subset\Cp2$ is merely a line $L\subset\Cp2$ such that the local
intersection index $(L\circ C)_P$ at each
intersection point $P\in L\cap C$ is even.

\theorem[see \autoref{proof.main} and \autoref{proof.main.2}]\label{th.main}
Let $t(C)$ denote the number of tritangents to a smooth sextic
$C\subset\Cp2$. Then either
\roster*
\item
$t(C)=72$, and then $C$ is the sextic given by~\eqref{eq.max.sextic}, or
\item
$t(C)=66$, and then $C$ is one of the two sextics described in
\autoref{proof.main}\iref{132.1}, \iref{132.2},
\endroster
or $t(C)\le65$.
\endtheorem

Previously known bounds are $t(C)\le76$ in
N.~Elkies~\cite{Elkies}
(\cf. \autoref{cor.Elkies} below) and $t(C)\le108$ given
by Pl\"{u}cker's formulas. Note that, unlike the $28$ bitangents to any smooth
quartic curve (and like the case of other polarized $K3$-surfaces),
counting tritangents
to a sextic is not an enumerative problem: tritangents are not stable under
deformation and a typical sextic has no tritangents at all.

\addendum[see \autoref{proof.main}]\label{ad.main}
The number $t(C)$ as in \autoref{th.main} takes all values in the set
$\{0,1,\ldots,65,66,72\}$ except, possibly, $61$.
\endaddendum

Twelve sextics (six configurations of lines) with $62\le t(C)\le65$ are
described in \autoref{proof.main}\iref{130}--\iref{124},
but we do not assert the completeness of this list. In spite of extensive,
although not exhaustive, search, we could not find a sextic with $61$
tritangents.
There are reasons (\eg, \autoref{cor.Elkies} below or the large number of
sextics with $60$ tritangents) to believe that $61$ is a natural threshold in
the problem, but taking the classification down to $61$
tritangents would require too much computing power.

As a by-product of the partial classification given by \autoref{th.main}, we
obtain a sharp upper bound on the number of real tritangents to a real
sextic.

\theorem[see \autoref{proof.real}]\label{th.real}
The number of real tritangents to a real smooth
\rom(over~$\C$\rom)
sextic $C\subset\Cp2$ does
not exceed $66$. Up to real projective transformation, a smooth
real sextic with $66$ real tritangents is unique,
see \autoref{proof.main}\iref{132.1}.
\endtheorem

\remark
At present, I do not know what other values are taken by the number of real
tritangents to a real smooth sextic. In the range $61\le t_\R\le65$,
\emph{among the known examples}, there is but one other configuration, with
$t_\R=63$ tritangents, see \autoref{proof.main}\iref{126.1} (and
\autoref{proof.real} for the explanation and further remarks).
\endremark

\subsection{Contents of the paper}\label{s.contents}
As in~\cite{degt:lines,DIS}, the line counting problem has a simple
arithmetical reduction (see \autoref{th.Fano.graph}): one can effectively
decide whether a given graph~$\Gamma$ can serve as the Fano graph of a
polarized $K3$-surface.
The candidates~$\Gamma$ to be tried were constructed in \cite{degt:lines,DIS}
line by line, starting from a sufficiently large and sufficiently simple
graph. Unfortunately, this straightforward approach seems to diverge in the
case of degree~$2$, and we choose another one, \viz. we replant the prospective
N\'{e}ron--Severi lattice $\NS:=\Z\Gamma/\ker$ to an appropriate Niemeier
lattice.
(The idea of embedding $h^\perp\subset\NS$ to a Niemeier lattice
is not new, \cf.
Kond\=o~\cite{Kondo},
Nikulin~\cite{Nikulin:degenerations},
Nishiyama~\cite{Nishiyama},
\etc.
The novelty is the fact that, as we need to keep track of the
polarization~$h$,
we have to rebuild the \emph{hyperbolic} lattice $\NS$ to embed it
to a \emph{definite} Niemeier lattice~$N$.
As a result, instead of counting roots in $\NS$, we work with
square~$4$ vectors in a certain root-free sublattice $S\subset N$;
unlike~\cite{degt:singular.K3}, this lattice~$S$ or
even its genus is not assumed fixed.
This construction is explained in \autoref{s.S},
see \autoref{prop.Niemeier}.) Then, instead of dealing with abstract graphs
of \latin{a priori} unbounded complexity, we merely need to consider
subsets~$\fL$ of
several finite sets $\fF(\hbar)$ \emph{known in advance}.
The precise arithmetical conditions on the subsets~$\fL$ that may serve as
Fano graphs are stated in \autoref{s.admissible} and \autoref{s.geometric}.

This approach has a number of advantages. First, for most
$6$-polarized Niemeier lattices $\bN\ni\hbar$ we have an immediate bound
$\ls|\fL|\le130$ (often even $\ls|\fF(\hbar)|\le130$)
obtained as explained in \autoref{S.bounds}. Second, the sets $\fF(\hbar)$
have rich intrinsic structure, splitting into orbits and combinatorial orbits
(see \autoref{s.notation}), which can be used in the construction of large
geometric subsets: instead of building them line-by-line from
scratch, we
try to patch together precomputed close to maximal intersections with the
combinatorial orbits. These algorithms are described in \autoref{S.approach}.
Finally, since we are working with known sets, all symmetry groups can be
expressed in terms of permutations, which makes the computation in
\GAP~\cite{GAP4} extremely effective.

In \autoref{S.few}--\autoref{S.24A1} we treat, one by one, the $23$ Niemeier
lattices rationally generated by roots, outlining the details of the
computation in those few cases where the \latin{a priori} upper bound
$\ls|\fL|\le130$ fails.
In \autoref{S.proofs}, we draw a formal punch-line, collecting together our
findings for individual Niemeier lattices and completing the proofs of the
principal results of the paper.

\subsection{Acknowledgements}
I would like to express my gratitude to
Noam Elkies, Dino Festi, Dmitrii Pasechnik, Ichiro Shimada,
and Davide Veniani for a number of fruitful discussions
concerning the subject. This paper was
completed during my research stay at the \emph{Max--Planck--Institut f\"{u}r
Mathematik}, Bonn; I am grateful to this institution
for its hospitality and financial support.

\section{The reduction}\label{S.reduction}

The tritangent problem is reduced to an arithmetical question about
the N\'{e}ron--Severi lattice $\NS(X)$ of a smooth $2$-polarized
$K3$-surface~$X$, see
\autoref{th.Fano.graph}.
The construction of \autoref{s.S}, combined with
\autoref{prop.Niemeier}, replants $\NS(X)$ to a
Niemeier lattice. The invertibility of this construction is discussed in
\autoref{s.geometric}.

\subsection{Lattices\pdfstr{}{
 {\rm(see~\cite{Nikulin:forms})}}}\label{s.lattice}
The principal goal of this section is fixing the terminology and notation.
A \emph{lattice} is a free abelian group~$L$ of finite rank equipped with a
symmetric bilinear form $b\:L\otimes L\to\Z$.
Since $b$ is assumed fixed (and omitted from the notation), we abbreviate
$x\cdot y:=b(x,y)$ and $x^2:=b(x,x)$. A lattice~$L$ is \emph{even} if
$x^2=0\bmod2$ for all $x\in L$; otherwise, $L$ is \emph{odd}.
The
\emph{determinant} $\det L\in\Z$ is the determinant of the Gram matrix of~$b$
in any integral basis; $L$ is called \emph{nondegenerate} (\emph{unimodular})
if $\det L\ne0$ (respectively, $\det L=\pm1$). The \emph{inertia indices}
$\Gs_\pm L$ are those of $L\otimes\R$. A nondegenerate lattice~$L$ is called
\emph{hyperbolic} if $\Gs_+L=1$.

The \emph{hyperbolic plane} is the only unimodular even lattice of rank~$2$.
Explicitly,
$\bU=\Z a+\Z b$, where $a^2=b^2=0$ and $a\cdot b=1$. One has
$\Gs_+\bU=\Gs_-\bU=1$.

A nondegenerate lattice~$L$ admits a canonical inclusion
\[*
L\into L\dual:=\bigl\{x\in L\otimes\Q\bigm|
 \mbox{$x\cdot y\in\Z$ for all $y\in L$}\bigr\}
\]
to the dual group~$L\dual$.
The finite abelian group $\CL:=\discr L:=L\dual\!/L$ ($q_L$
in~\cite{Nikulin:forms}) is called the \emph{discriminant group} of~$L$.
Clearly, $\ls|\CL|=(-1)^{\Gs_-L}\det L$.
This group is equipped with the nondegenerate symmetric bilinear form
\[*
\CL\otimes\CL\to\Q/\Z,\quad
(x\bmod L)\otimes(y\bmod L)\mapsto(x\cdot y)\bmod\Z,
\]
and, if $L$ is even, its quadratic extension
\[*
\CL\to\Q/2\Z,\quad x\bmod L\mapsto x^2\bmod2\Z.
\]
We denote by $\CL_p:=\discr_pL:=\CL\otimes\Z_p$ the $p$-primary
components of $\discr L$. The $2$-primary component~$\CL_2$ is called \emph{even}
if $x^2\in\Z$ for all order~$2$ elements $x\in\CL_2$; otherwise, $\CL_2$ is
\emph{odd}. The \emph{determinant} $\det\CL_p$ is the determinant of the
``Gram matrix'' of the quadratic form in any minimal set of generators.
(This is equivalent to the alternative definition given
in~\cite{Nikulin:forms}.)
Unless $p=2$ and $\CL_2$ is odd (in which case
the determinant is not defined or used),
we have $\det\CL_p=u_p/\ls|\CL_p|$, where $u_p$ is a well-defined element of
$\Z_p\units/(\Z_p\units)^2$.

The \emph{length} $\ell(\Cal A)$ of a finite abelian group~$\Cal A$ is the
minimal number of generators of~$\Cal A$. We abbreviate $\ell_p(\Cal
A):=\ell(\Cal A\otimes\Z_p)$ for a prime~$p$.

Given a lattice~$L$ and $q\in\Q$, we use the notation $L(q)$ for the same
abelian group with the form $x\otimes y\mapsto q(x\cdot y)$, assuming that it
is still a lattice. We abbreviate $-L:=L(-1)$, and this notation applies to
discriminant forms as well. The notation $nL$, $n\in\Z_+$, is used for the
orthogonal direct
sum of $n$ copies of~$L$.

A \emph{root} in an even lattice~$L$ is a vector of square~$\pm2$.
A \emph{root system} is a positive definite lattice generated by roots. Any
root system has a unique decomposition into orthogonal direct sum of
irreducible components, which are of types $\bA_n$, $n\ge1$, $\bD_n$,
$n\ge4$, $\bE_6$, $\bE_7$, or~$\bE_8$ (see, \eg, \cite{Bourbaki:Lie}),
according to their \emph{Dynkin diagrams}.

A \emph{Niemeier lattice} is a positive definite unimodular even lattice of
rank~$24$. Up to isomorphism, there are $24$ Niemeier lattices
(see~\cite{Niemeier}): the \emph{Leech lattice}~$\Lambda$, which is root free,
and $23$ lattices \emph{rationally} generated by roots.
In the latter case, the isomorphism class of a lattice $N:=\N(\bR)$ is uniquely
determined by that of its
maximal root system~$\bR$. For more details, see~\cite{Conway.Sloane}.

\subsection{The covering $K3$-surface}\label{s.K3}
Given a smooth sextic curve $C\subset\Cp2$,
the double covering $\Gf\:X\to\Cp2$ ramified over~$C$ is a $K3$-surface.
The ``hyperplane section'' $\Gf^*\CO_{\Cp2}(1)$ is
a \emph{$2$-polarization} of~$X$, \ie, a complete fixed
point free degree~$2$ linear system; it is viewed as an element
\[*
h\in\operatorname{Pic}X=\NS(X)\subset H_2(X;\Z)\cong-2\bE_8\oplus3\bU.
\]
Here, the group $H_2(X;\Z)=H^2(X;\Z)$ is regarded as a lattice
\via\ the intersection form;
it can be characterized as the only unimodular even lattice of rank~$22$ and
signature $\Gs_+-\Gs_-=-16$.
The \emph{N\'{e}ron--Severi lattice} $\NS(X)=H^{1,1}(X)\cap H_2(X;\Z)$ is
a primitive hyperbolic sublattice; in
particular, $\Gr(X):=\rank\NS(X)\le20$.

Conversely, any $2$-polarization~$h$ of a $K3$-surface $X$ gives rise to a
degree~$2$
map $\Gf_h\:X\to\Cp2$ ramified over a sextic curve $C\subset\Cp2$ (see
\cite{Persson:sextics,Saint-Donat}). This curve is smooth if and only if no
$(-2)$-curve is contracted by~$\Gf_h$, or, equivalently, there is no class
$e\in\NS(X)$ such that $e^2=-2$ and $e\cdot h=0$.
With the ramification locus in mind, a $2$-polarized $K3$-surface
$(X,h)$ with this
extra property is called \emph{smooth}.

A \emph{line} in a $2$-polarized $K3$-surface $(X,h)$ is a smooth rational curve
$L\subset X$ such that $L\cdot h=1$. Any two distinct lines
$L_1,L_2\subset X$ either are disjoint, $L_1\cdot L_2=0$,
or intersect at a single point, $L_1\cdot L_2=1$,
or intersect at three points, $L_1\cdot L_2=3$, the latter being the case if
and only if $L_1$, $L_2$ are interchanged by the deck translation of the
covering $\Gf_h\:X\to\Cp2$. Since, on the other hand, $L^2=-2$, each line is
unique in its homology class $[L]\in\NS(X)$.
Each $2$-polarized $K3$-surface has finitely many lines (typically none).
The \emph{Fano graph} $\Fn(X,h)$ is the set of lines in~$X$ in which each pair of
lines $L_1,L_2$ (regarded as vertices of the graph)
is connected by an edge of multiplicity
$L_1\cdot L_2$ (\ie, no edge, simple edge, or triple edge).

Let $C\subset\Cp2$ be a smooth sextic and $\Gf\:X\to\Cp2$ the covering
$K3$-surface. If $L\subset\Cp2$ is a tritangent to~$C$, its pull-back
$\Gf\1(L)$ splits into two lines $L_1,L_2$; they intersect at the three
points of tangency of~$L$ and~$C$ (possibly, infinitely near) and are
interchanged by the deck translation~$\tau$ of~$\Gf$. Conversely, any line
in~$X$ projects to a tritangent to~$C$. Thus, the set of tritangents to~$C$
is identified with $\Fn X/*$, where the free involution
${*}\:\Fn X\to\Fn X$ induced by~$\tau$
is intrinsic to the graph: it sends a vertex~$L$ to the only vertex connected
to~$L$ by a triple edge.

\subsection{The arithmetic reduction of the tritangent problem}\label{s.NS}
Throughout this paper, by a \emph{$2$-polarized lattice} we mean a hyperbolic
even
lattice~$\NS$ equipped with a distinguished class $h\in\NS$
of square $h^2=2$.
The \emph{Fano graph} of a $2$-polarized lattice $\NS\ni h$ is the set
\[*
\Fn(\NS,h):=\bigl\{l\in\NS\bigm|\text{$l^2=-2$, $l\cdot h=1$}\bigr\}
\]
with two points (vertices) $l_1$, $l_2$ connected by an edge of multiplicity
$l_1\cdot l_2$. This graph is equipped with a natural involution
\[*
l\mapsto l^*:=h-l;
\]
the vertex $l^*$, called the \emph{dual} of~$l$, is connected to~$l$ by a
triple edge.

Usually, we assume, in addition, that the orthogonal complement
$h^\perp\subset\NS$ is root free. Under this additional assumption,
\[*
\text{for $l_1,l_2\in\Fn(\NS,h)$, one has $l_1\cdot l_2=3$ (iff $l_1=l_2^*$),
 $1$, $0$, or $-2$ (iff $l_1=l_2$)};
\]
hence, all edges of $\Fn(\NS,h)$ other than $(l,l^*)$ are simple.

The following statement is well known:
it follows from the global Torelli theorem for
$K3$-surfaces~\cite{Pjatecki-Shapiro.Shafarevich}, surjectivity of the period
map~\cite{Kulikov:periods}, and Saint-Donat's results on projective
$K3$-surfaces~\cite{Saint-Donat} (\cf. also
\cite[Theorem 3.11]{DIS} or \cite[Theorem 7.3]{degt:singular.K3}).

\theorem\label{th.Fano.graph}
A graph $\Gamma$ is the Fano graph of a smooth $2$-polarized
$K3$-surface if and only if $\Gamma\cong\Fn(\NS,h)$ for some $2$-polarized
lattice $\NS\ni h$ admitting a primitive embedding
$\NS\into-2\bE_8\oplus3\bU$ and
such that $h^\perp\subset\NS$ is root free.
\done
\endtheorem

\subsection{Embedding to a Niemeier lattice}\label{s.S}
Let $\NS\ni h$ be a $2$-polarized lattice.
Consider the orthogonal complement $h^\perp\subset\NS$.
Each vector $l\in\Fn(\NS,h)$ projects to
$l':=l-\frac12h\in(h^\perp)\dual$, and, assuming
$\Fn(\NS,h)\ne\varnothing$,
there is a unique index~$2$ extension
\[
{-S}\supset h^\perp\oplus\Z\hbar,\quad \hbar^2=-6,
\label{eq.S}
\]
containing all vectors $l'+\frac12\hbar$, $l\in\Fn(\NS,h)$.
The lattice $S:=S(\NS,h)$ obtained from $-S$ by reverting the sign of the
binary form is
positive definite, and there is an obvious canonical bijection between $\Fn(\NS,h)$
and the set
\[*
\fL=\fL(S,\hbar):=\bigl\{l\in S\bigm|
 \mbox{$l^2=4$ and $l\cdot\hbar=3$}\bigr\};
\]
the elements of~$\fL$ are called \emph{lines} in~$S$. Furthermore,
the sublattice
$h^\perp\subset\NS$ is root free if and only if so is $\hbar^\perp\subset S$;
in this case, we call $S\ni\hbar$ \emph{admissible}.

For the images $l_1,l_2\in\fL$ of two lines $l_1',l_2'\in\Fn(\NS,h)$ one has
$l_1\cdot l_2=2-l_1'\cdot l_2'$. Hence, if $S\ni\hbar$ is admissible,
then
\[
\text{for $l_1,l_2\in\fL$, one has $l_1\cdot l_2=-1$ (iff $l_1=l_2^*$),
 $1$, $2$, or $4$ (iff $l_1=l_2$)}.
\label{eq.intersection}
\]
We say that $l_1$, $l_2$ \emph{intersect} (are \emph{disjoint}) if
$l_1\cdot l_2=1$ (respectively, $l_1\cdot l_2=2$). Accordingly, we regard
$\fL$ as a graph, with two distinct vertices $l_1$, $l_2$ connected by a
simple (triple) edge whenever $l_1\cdot l_2=1$
(respectively, $l_1\cdot l_2=-1$.)

\proposition\label{prop.Niemeier}
Let $\NS\ni h$ be a primitive $2$-polarized sublattice of
$-2\bE_8\oplus3\bU$,
$\Fn(\NS,h)\ne\varnothing$,
and let $S:=S(\NS,h)$ be the lattice constructed as
in~\eqref{eq.S}. Then
\roster
\item\label{Niemeier.primitive}
$S$ admits a primitive embedding to a Niemeier lattice~$N$\rom;
\item\label{Niemeier.root}
$S$ admits an embedding $S\into N$ to a Niemeier lattice such that the
torsion of $N/S$ is a $3$-group and $S$ is orthogonal to a root $\rbar\in N$.
\endroster
\endproposition

\proof
Denote $\Gr:=\rank\NS$ and $\CN:=\discr\NS$,
so that we have $\ell(\CN)\le22-\Gr$ by Theorem 1.12.2
in~\cite{Nikulin:forms}.
Since $h\notin2\NS\dual$ (by the assumption that $\Fn(\NS,h)\ne\varnothing$),
we have
\[*
\discr h^\perp=\bigl\<\tfrac12h\bigr\>\oplus\CN,
\quad\bigl(\tfrac12h\bigr)^2=\tfrac32\bmod2\Z,
\]
and the construction changes this to
\[*
\discr S=\bigl\<\tfrac12\hbar\bigr\>\oplus(-\CN),
\quad\bigl(\tfrac12\hbar\bigr)^2=\tfrac23\bmod2\Z.
\]
In particular, $\ell(\discr S)\le\ell(\CN)+1<24-\Gr$, and
Theorem 1.12.2 in~\cite{Nikulin:forms} implies the existence of a primitive
embedding $S\into N$.
For the second statement, we compute
\[*
\CS:=\discr(S\oplus\Z\rbar)
 =\bigl\<\tfrac12\hbar\bigr\>\oplus\bigl\<\tfrac12\rbar\bigr\>\oplus(-\CN),
\quad\bigl(\tfrac12\rbar\bigr)^2=\tfrac12\bmod2\Z.
\]
This time we have $\ell(\CS_p)=\ell(\CN_p)<23-\Gr=24-\rank(S\oplus\Z\rbar)$
for each prime $p>3$,
whereas $\ell(\CS_p)=\ell(\CN_p)+1\le24-\rank(S\oplus\Z\rbar)$ for $p=2,3$.
Since $\CS_2$ is odd, the possible equality does not impose any extra
restriction at $p=2$. For $p=3$, in the case of equality, the
``wrong'' determinant $\det(-\CS_3)=-\ls|\CS|\bmod(\Z_3\units)^2$
does inhibit the existence of a primitive embedding. However, since
$\ell(\CS_3)\ge3$ in this case, we may pass to an iterated index~$3$
extension and reduce the length.
\endproof

\subsection{Admissible sets}\label{s.admissible}
In the rest of the paper, we mainly use
statement~\iref{Niemeier.root} of
\autoref{prop.Niemeier}:
it lets us avoid the Leech lattice, although at the expense of the possible
imprimitivity (which makes some statements somewhat weaker and
more complicated, see, \eg, \autoref{prop.geometric} below).
The idea is to construct a lattice~$S$ (or, rather, its set of
lines) directly inside a Niemeier lattice. Thus, we fix a
Niemeier lattice $\bN$, a square~$6$ vector
$\hbar\in\bN$, and, optionally, a root $\rbar\in\hbar^\perp$
(which is typically omitted from the notation).
Consider the set
\[*
\fF:=\fF(\hbar):=\bigl\{l\in\bN\bigm|
 \text{$l^2=4$, $l\cdot\hbar=3$ (and $l\cdot\rbar=0$)}\bigr\}.
\]
It is equipped with the involution
\[*
*\:l\mapsto l^*:=\hbar-l.
\]
The elements of~$\fF(\hbar)$ are called \emph{lines}. The \emph{span} of
a subset $\fL\subset\fF(\hbar)$ is the lattice
\[*
\spn\fL:=(\Z_3\fL+\Z_3\hbar)\cap\bN\subset\bN.
\]
If $\fL$ is symmetric, $\fL^*=\fL$, the summation with $\Z_3\hbar$
is redundant as $\hbar\in\Z\fL$.
On a few occasions, we also consider the \emph{integral} and
\emph{rational span}
\[*
\spn_\Z\fL:=(\Z\fL+\Z\hbar)\cap\bN\subset\spn\fL\subset
 \spn_\Q\fL:=(\Q\fL+\Q\hbar)\cap\bN.
\]
(The latter is primitive in~$\bN$.)
\latin{Via} $\spn$, we extend to subsets $\fL\subset\fF(\hbar)$
much of the terminology applied to lattices. Thus,
the \emph{rank} of~$\fL$ is
$\rank\fL:=\rank\spn\fL$, and
we say that $\fL$ is \emph{generated} by a subset $\fL'\subset\fL$ if
$\fL=\fF(\hbar)\cap\spn\fL'$.

By definition, the torsion of $N/\spn\fL$ is a $3$-group and
$\hbar\in3(\spn\fL)\dual$.
A finite index extension $S\supset\spn\fL$
is called \emph{mild}
if $S\subset\{v\in\bN\,|\,v\cdot\hbar=0\bmod3\}$
(\ie, $S\subset\bN$ and still $\hbar\in3S\dual$) and $S$
contains no roots $r\in\hbar^\perp\subset\bN$.
Note that the latter condition is equivalent to the requirement that $S$ itself should
be root free. Indeed, since $S$ is positive definite and $\hbar\in3S\dual$,
we have $r\cdot\hbar=0$ or~$\pm3$ for any root $r\in S$, and in the latter
case $\hbar\mp r$ is a root in~$\hbar^\perp$.

\definition\label{def.admissible}
A subset $\fL\subset\fF(\hbar)$ is called \emph{admissible} if
\roster
\item\label{admissible.1}
$\fL$ is \emph{symmetric} (or $*$-invariant), \ie, $\fL^*=\fL$, and
\item\label{admissible.2}
the sublattice $\hbar^\perp\cap\spn\fL$ contains no roots.
\endroster
A subset $\fL\subset\fF(\hbar)$ is \emph{complete}
if $\fL=\fF(\hbar)\cap\spn\fL$.
A subset $\fL$ is \emph{saturated} if
the identity $\fL=\fF(\hbar)\cap S$ holds for
any mild extension $S\supset\spn\fL$.
Finally, we say that $\fL$ is \emph{$\Q$-complete} if
$\fL=\fF(\hbar)\cap\spn_\Q\fL$.
\enddefinition

Often, it is easier to check~\eqref{eq.intersection}, which follows
from~\iref{admissible.1}, \iref{admissible.2} above.
Indeed, since $S$ is definite, we have
$-1\le l_1\cdot l_2\le4$. Thus, forbidden are
$l_1\cdot l_2=3$ or~$0$, as then $l_1-l_2$
or $l_1-l_2^*=l_1+l_2-\hbar$, respectively, would be a root in $\hbar^\perp$.

The following bound is due to N.~Elkies.

\theorem[N.~Elkies, \cite{Elkies}]\label{th.Elkies}
Let $V$ be a Euclidean vector space, $\dim V=n$, and let
$v_1,\ldots,v_N\in V$ be a collection of unit vectors such that the
products $v_i\cdot v_j$, $i\ne j$, take but two values $\tau_1$, $\tau_2$.
Assume that $\tau_1+\tau_2\le0$ and $1+\tau_1\tau_2n>0$. Then
\[*
N\le\frac{(1-\tau_1)(1-\tau_2)n}{1+\tau_1\tau_2n}.
\]
\endtheorem

Selecting a single vector from each pair $l,l^*\in\fL$ and applying
\autoref{th.Elkies} to the normalized projections to $\hbar^\perp\subset\spn\fL$,
we arrive at the following corollary.

\corollary[N.~Elkies \cite{Elkies}]\label{cor.Elkies}
The size of an admissible set~$\fL$ is bounded \via
\[*
\ls|\fL|\le\frac{48(\rank\fL-1)}{26-\rank\fL}.
\]
\endcorollary

Since $\ls|\fL|$ is even,
this gives us $\ls|\fL|\le152$ or $122$ for $\rank\fL=20$ or~$19$,
respectively.

\subsection{Geometric sets}\label{s.geometric}
According to \autoref{th.Fano.graph} and \autoref{prop.Niemeier},
the Fano graph of any smooth $2$-polarized $K3$-surface~$X$
can be represented as
a complete admissible subset $\fL\subset\fF(\hbar)$ for an appropriate pair
$\hbar,\rbar\in\bN$ as in \autoref{s.admissible}.

For some lattices (those with few roots), the admissibility condition is not
enough to eliminate large sets of lines, and we need to use the full range of
restrictions.

Recall that we start with the N\'{e}ron--Severi lattice $\NS(X)\ni h$, and as
long as the configurations of lines are concerned, we can assume this lattice
rationally generated by lines. Indeed, let $N:=\Q\Fn(X)\cap\NS(X)$;
clearly, $h\in N$.
We can pick
a vector $\Go\in(N^\perp)\otimes\C$, $\Go^2=0$, in the same component of the
positive cone as the period
(class of a holomorphic $2$-form) $\Go_X$ of~$X$,
and such that $\Go^\perp\cap H_2(X)=N$.
(This condition merely means that $\Go$ is generic; it can be chosen
arbitrary close to~$\Go_X$.)
By
the surjectivity of the period map~\cite{Kulikov:periods}, there is a
$K3$-surface $X'\to\Cp2$ such that $\Go_{X'}=\Go$, so that $\NS(X')=N$. (The
fact that $h\in N$ defines a map $X'\to\Cp2$ with a smooth ramification locus
follows from that for~$X$, \cf.~\cite{Saint-Donat} or \autoref{s.K3}.)
Clearly, we have $\Fn(X')=\Fn(X)$, \cf. \autoref{s.K3}, and the lines in~$X'$
generate (over~$\Q$) its N\'{e}ron--Severi lattice~$N$ by the construction.

Thus, assuming that $\NS(X)\ni h$ is rationally generated by lines,
we can pass to the positive definite lattice $S\ni\hbar$ as
in~\eqref{eq.S} and embed the latter to a Niemeier lattice~$\bN$, mapping
$\Fn(X)$ bijectively onto the admissible set $\fL=\fF(\hbar)\cap S$.
Most steps of this construction are invertible. However,
starting from an admissible set $\fL\subset\fF(\hbar)$, we may have to take
for~$S$ a mild extension of $\spn\fL$ rather than $\spn\fL$ itself and,
still, we cannot guarantee that the lattice $\NS$ obtained from~$S$ by the
backward construction admits a primitive embedding to
$H_2(X)\cong-2\bE_8\oplus3\bU$.
This discussion motivates the following definition.

\definition\label{def.geometric}
An admissible set $\fL\subset\fF(\hbar)$ is called \emph{geometric} if $\fL$
is complete in some mild extension $S\supset\spn\fL$ such that the lattice
$\NS$ obtained from $S\ni\hbar$ by the inverse
of construction~\eqref{eq.S} admits a
primitive embedding to $-2\bE_8\oplus3\bU$.
\enddefinition

Using Theorem 1.12.2
in~\cite{Nikulin:forms}, one can recast this property as follows.
(For a \emph{mild} extension $S\supset\spn\fL$ there is a splitting
$\discr S=\bigl\<\frac12\hbar\bigr\>\oplus\Cal T$, and we merely restate the
restrictions on~$\Cal T\cong-\discr\NS$ in terms of $\discr S$.)

\proposition\label{prop.geometric}
Let $N\ni\hbar$ be as above. An admissible set $\fL\subset\fF(\hbar)$ is
geometric if and only if\rom:
\roster
\item\label{geometric.rank}
$\rank\fL\le20$\rom; we denote $\Gd:=22-\rank\fL\ge2$, and
\endroster
there is a mild extension $S\supset\spn\fL$ in which $\fL$ is a complete
subset and such that the discriminant $\CS:=\discr S$ has the following
properties at each prime~$p$\rom:
\roster[\lastitem]
\item\label{geometric.p}
if $p>3$, then either $\ell(\CS_p)<\Gd$ or $\ell(\CS_p)=\Gd$ and
$\det\CS_p=3\ls|\CS|\bmod(\Q_p\units)^2$\rom;
\item\label{geometric.2}
either $\ell(\CS_2)<\Gd$ or $\ell(\CS_2)=\Gd$ and $\CS_2$ is odd or
$\det\CS_2=\pm3\ls|\CS|\bmod(\Q_2\units)^2$\rom;
\item\label{geometric.3}
either $\ell(\CS_3)\le\Gd$ or $\ell(\CS_3)=\Gd+1$ and
$\det\CS_3=\ls|\CS|\bmod(\Q_3\units)^2$.
\endroster
\endproposition

\remark\label{rem.geometric}
In practice, when eliminating large admissible sets, we use just a few
simple consequences of \autoref{prop.geometric}. The main r\^{o}le is played by
condition~\iref{geometric.rank}, see \autoref{ss.maximal} below.
Then, conditions~\iref{geometric.p} and~\iref{geometric.2} are used, as they
apply directly to the original discriminant $\discr_p(\spn\fL)=\CS_p$,
$p\ne3$.
Condition~\iref{geometric.3} is typically used when there is an
obvious maximal mild extension, and we never insist that $\fL$ should be
complete in~$S$, thus eliminating both~$\fL$ itself and all its oversets.
\endremark

\section{The approach}\label{S.approach}

Throughout this section, we consider
a Niemeier lattice $\bN:=\N(\bR)$ generated over~$\Q$ by a
fixed root
system $\bR=\bigoplus_k\bR_k$, $k\in\Omega$,
where $\bR_k$ are the irreducible components
(\emph{aka} Dynkin diagrams) and $\Omega$ is the index set.
We construct~$\bN$ as a subgroup of $\bigoplus_i\bR_i\dual$; the vectors
in
\[*\textstyle
\discr\bR:=\bR\dual\!/\bR=\bigoplus_k\discr\bR_k
\]
that are declared ``integral'' are as described
in~\cite[Table 16.1]{Conway.Sloane}.
(We also use the convention of~\cite{Conway.Sloane} for the numbering of the
discriminant classes of irreducible root systems.)
We denote by $\OG:=\OG(\bN)$ the full orthogonal
group of~$\bN$, and by $\RG:=\RG(\bN)\subset\OG(\bN)$ its subgroup generated by
reflections. Both groups preserve~$\bR$; the reflection group
$\RG(\bN)$ preserves each~$\bR_k$ and acts identically on $\discr\bR$.

\subsection{Notation}\label{s.notation}
We fix a square~$6$ vector $\hbar\in\bN$ and, sometimes, a root $\rbar\in\bR$
orthogonal to~$\hbar$. (This root is usually omitted from the notation.)
We denote by $\OG_\hbar(\bN)\subset\OG(\bN)$ and
$\RG_\hbar(\bN)\subset\RG(\bN)$ the subgroups stabilizing~$\hbar$
(and~$\rbar$).
Let
\[*
\fF=\fF(\hbar)=\fF(\hbar,\rbar):=\bigl\{l\in\bN\bigm|
 \mbox{$l^2=4$, $l\cdot\hbar=3$ (and $l\cdot\rbar=0$)}\bigr\}
\]
be the set of lines. This set splits into a number of $\OG_\hbar(\bN)$-orbits
$\borb_n$, which split further into $\RG_\hbar(\bN)$-orbits
$\orb\subset\borb_n$; the latter are called \emph{combinatorial orbits}.
It is immediate that the duality $l\mapsto l^*$ preserves orbits and
combinatorial orbits; hence, we can speak about the dual orbits $\borb_n^*$
and $\orb^*$. The number of combinatorial orbits in an orbit $\borb_n$ is
denoted by~$m(\borb_n)$. The set of all combinatorial orbits is denoted by
$\Orb:=\Orb(\hbar)$. This set inherits a natural action of the group
\[*
\stab\hbar:=\OG_\hbar(\bN)/\!\RG_\hbar(\bN),
\]
which preserves each orbit~$\borb_n$.
(By an obvious abuse of notation, occasionally we treat $\borb_n$ as
a subset of $\Orb$; likewise, subsets of $\Orb$ are sometimes treated
as sets of lines.)
We denote by $\Orbit_m(\borb_n,k)$ the length~$m$ orbit of the action of
$\stab\hbar$ on the set of unordered
$*$-invariant (if so is~$\borb_n$) $k$-tuples
of
combinatorial orbits $\orb\subset\borb_n$.
The usage of this
notation implies implicitly that such an orbit is unique.

The \emph{support} of a vector $v\in\bN=\N(\bR)$ is the subset
\[*
\supp v:=\bigl\{k\in\Omega\bigm|v_k\ne0\in\bR_k\dual\}\subset\Omega.
\]
The support is invariant under reflections; hence, we can speak about the
\emph{support}
$\supp\orb$ of a combinatorial orbit~$\orb$.

The \emph{count} and \emph{bound} of a combinatorial orbit~$\orb$ are defined
\via
\[
\cnt(\orb):=\ls|\orb|,\qquad
\bnd(\orb):=\max\bigl\{\ls|\fL\cap\orb|\bigm|
 \mbox{$\fL\subset\fF$ is geometric}\bigr\}.
\label{eq.count.bound}
\]
Clearly, $\cnt$ and $\bnd$ are constant within each orbit~$\borb_n$ and
invariant under duality. In some cases, we replace $\bnd(\orb)$ by rough
bounds, see \autoref{s.computing} below for details.
We extend these notions to subsets $\Cluster\subset\Orb$
\emph{by additivity}:
\[*
\cnt(\Cluster):=\sum_{\orb\in\Cluster}\cnt(\orb),
\qquad
\bnd(\Cluster):=\sum_{\orb\in\Cluster}\bnd(\orb).
\]
Thus, we have a na\"{\i}ve \latin{a priori} bound
\[
\ls|\fL|\le\bnd(\Orb)=\sum m(\borb_n)\bnd(\orb),\quad\orb\subset\borb_n.
\label{eq.naive}
\]
Clearly, the true count $\ls|\fL\cap\Cluster|$ is genuinely additive, whereas
the sharp bound on $\ls|\fL\cap\Cluster|$ is only subadditive; thus, our
proof of \autoref{th.main} will essentially consist in
reducing~\eqref{eq.naive} down to a preset goal. To this end, we will
consider the set
\[*
\Bnd=\Bnd(\fF):=
 \bigl\{\fL\subset\fF\bigm|\mbox{$\fL$ is geometric}\}/\!\OG_\hbar(\bN)
\]
and,
for a collection of orbits $\Cluster=\borb_1\cup\ldots$ and
integer $d\in\NN$, let
\[*
\Bnd_d(\Cluster):=\bigl\{[\fL]\in\Bnd\bigm|
 \mbox{$\fL$ is generated by $\fL\cap\Cluster$ and
 $\ls|\fL\cap\Cluster|\ge\bnd(\Cluster)-d$}\bigr\}.
\]
Unless specified otherwise, the sets $\Bnd_d(\Cluster)$
(for reasonably small values of~$d$) are computed by brute force, using
patterns (see \autoref{s.patterns} below).

\subsection{Idea of the proof}\label{s.idea}
To prove \autoref{th.main}, we consider, one by one, all $23$ Niemeier
lattices generated by roots. For each lattice~$N$, we set a goal
\[
\ls|\fL|\ge M:=122\ \mbox{or}\ 132
\label{eq.goal}
\]
and try to find all geometric subsets $\fL\subset N$ satisfying this
inequality. First, we list all $\OG(N)$-orbits of square~$6$ vectors
$\hbar\in N$, compute the na\"{\i}ve bounds $\bnd(\Orb)$ given
by~\eqref{eq.naive}, and disregard those vectors for which $\bnd(\Orb)<M$.
In the remaining cases, we list all $\OG_\hbar(N)$-orbits of
roots~$\rbar$ orthogonal to~$\hbar$ and repeat the procedure.
This leaves us with relatively few triples $\hbar,\rbar\in N$, which are
treated on a case-by-case basis in \autoref{S.few}--\autoref{S.24A1} below.

A typical argument runs as follows. We choose a self-dual union~$\Cluster$ of
orbits~$\borb_n$ and use \emph{patterns} (see \autoref{s.patterns} below) to
compute the set $\Bnd_{\bnd(\Orb)-M}(\Cluster)$.
(As a modification, we take $\Cluster$ \emph{disjoint}
from its dual $\Cluster^*$ and use the obvious relation
$\Bnd_d(\Cluster)=\Bnd_{2d}(\Cluster\cup\Cluster^*)$.)
More generally, we can
consider several
pairwise disjoint self-dual unions of orbits $\Cluster_1,\ldots,\Cluster_m$
and compute the sets
$\Bnd_{d_i}(\Cluster_i)$ for appropriately chosen integers $d_i\ge0$
such that
\[*
d_1+\ldots+d_m+2(m-1)\ge\bnd(\Orb)-M.
\]
As a result of this procedure,
we can assert that, apart from a few explicitly listed exceptions
$\fL_1,\ldots,\fL_s$ contained in the above sets $\Bnd_{d_i}(\Cluster_i)$,
we have $\ls|\fL|<M$ for any geometric set $\fL\subset\fF$.
In each case, we manage to choose the unions $\Cluster_i$ and
goals~$d_i$ so that the
exceptional sets $\fL_k$ are sufficiently large, so that they can be
analysed further as explained below.

\subsubsection{Maximal sets}\label{ss.maximal}
The best case scenario is that of a \emph{maximal} (with respect to inclusion,
in the class of geometric sets) geometric set~$\fL$.
Such a set admits no geometric extensions; hence, it can be
either discarded, if $\ls|\fL_k|<M$, or listed as an exception in the
respective
statement. Besides, maximal sets can be discarded at
early stages of the computation, without completing the whole pattern; however,
we only use this approach in \autoref{s.24A1-3}, where intermediate lists
grow too large.

An obvious sufficient condition of maximality is given by
\autoref{prop.geometric}.

\lemma\label{lem.maximal}
Any maximal geometric set is saturated. Conversely, any saturated
geometric set~$\fL$ of the
maximal rank $\rank\fL=20$ is maximal.
\done
\endlemma

\subsubsection{Extension by a maximal orbit}\label{ss.ext.max}
If a set $\fL\in\Bnd_d(\Cluster)$ is not maximal, we try to list its
geometric extensions $\fL'\supset\fL$ satisfying~\eqref{eq.goal}. Clearly, it
suffices to consider \emph{$\Cluster$-proper} extensions only, \ie, those
with the property that $\fL'\cap\Cluster=\fL\cap\Cluster$. Thus, we merely
extend the partial pattern
\[*
\pat\:\Cluster\to\N,\quad \orb\mapsto\ls|\fL\cap\orb|,
\]
(see \autoref{s.patterns} below) used to construct~$\fL$ by a few
(usually one or at most two)
extra values $\pat(\orb)$, $\orb\in\Orb\sminus\Cluster$.

In many cases, a geometric set $\fL\in\Bnd_d(\Cluster)$
has the property that
\[*
\sum\bigl(\bnd(\orb)-\bnd'(\orb)\bigr)
 \ge\bnd(\Orb)-M,\qquad
 \orb\in\Orb_\Gd:=\bigl\{\orb\in\Orb\bigm|\ls|\fL\cap\orb|<\bnd(\orb)\bigr\},
\]
where $\bnd'(\orb)$ is the second largest value taken by
$\ls|\fL'\cap\orb|$ with $\fL'$ admissible.
This implies that any admissible extension $\fL'\supset\fL$
satisfying~\eqref{eq.goal} must have maximal intersection,
$\ls|\fL'\cap\orb|=\bnd(\orb)$, for at least one orbit $\orb\in\Orb_\Gd$.
Trying these orbits one by one
(\ie, extending the pattern \via\ $\pat(\orb)=\bnd(\orb)$),
we obtain larger sets, which are usually
maximal.
(This computation uses patterns, see \autoref{s.patterns} below, and
takes into account the symmetry of~$\fL$.)


\subsubsection{Other extensions}\label{ss.ext.other}
In the few remaining cases, we either analyze the lines contained in
$\spn_\Q\fL$
(if $\rank\fL=20$) or obtain maximal $\Cluster$-proper extensions
$\fL'\supset\fL$ by adding one or, rarely, two extra lines.

\subsection{Patterns}\label{s.patterns}
Since we are interested in large geometric sets, we construct them
orbit-by-orbit, by piling together maximal or close to maximal intersections
$\fL\cap\orb$.
This process is guided by \emph{patterns}, \ie, $*$-invariant functions
\[*
\pat\:\Orb\to\NN,\quad\orb\mapsto\ls|\fL\cap\orb|.
\]
Having $\hbar,\rbar\in\bN$ fixed, we start with precomputing all geometric
sets $\fL\subset\orb$ in each combinatorial orbit~$\orb$. (Certainly, it
suffices to consider one representative in each orbit $\borb_n$; the rest is
obtained by translations.) Then, in order to compute one of the sets
$\Bnd_d(\Cluster)$ in \autoref{s.notation}, we list all $(\stab\hbar)$-orbits
of restricted patterns $\pat\:\Cluster\to\NN$ satisfying the inequality
$\sum\pat(\orb)\ge\bnd(\Cluster)-d$, $\orb\in\Cluster$, order the orbits
appropriately (typically, by the decreasing of $\pat(\orb)$),
and construct a
geometric set~$\fL$ by adding one orbit at a time, as a sequence
$\varnothing=\fL_0\subset\fL_1\subset\fL_2\subset\ldots$. At each step~$k$
and for each set~$\fL_{k-1}$ constructed at the previous step,
we proceed as follows:
\roster
\item
compute the stabilizer~$G$ of~$\fL_{k-1}$ under the action
of~$\RG_\hbar(\bN)$;
\item
compute the $G$-orbits of the geometric sets $\fL'\subset\orb_k$
of size $\ls|\fL'|=\pat(\orb_k)$;
\item
for a representative~$\fL'$ of each $G$-orbit, consider the set $\fL_k$ generated
by the union $\fL_{k-1}\cup\fL'$; then, select those sets $\fL_k$ that are
geometric;
\item
to reduce the overcounting,
select, for the next step, those sets~$\fL_k$ that satisfy the equality
$\ls|\fL_k\cap\orb_i|=\pat(\orb_i)$ for each $i\le k$.
\endroster
If the defect~$d$ is not too large, this procedure works reasonably fast and
results in a reasonably small collection of sets that are to be analyzed
further.

\remark\label{rem.pattern.unique}
Although it is not obvious \latin{a priori}, it turns out that large
geometric sets are often determined by their patterns uniquely up
to~$\RG_\hbar(\bN)$. Furthermore, a
large set is easily reconstructed from its pattern, as the algorithm above
converges very fast. For this reason, we often
describe large geometric sets, especially those that are \emph{not}
$\Q$-complete (see \autoref{def.admissible}),
by their patterns.

A pattern~$\pat$ taking a constant value~$v_n$ on each orbit~$\borb_n$
is described \via
\[*
\pat=\patternform{v_1,v_2,\ldots}.
\]
Sometimes, we use a ``double value'' $v_n=a|b$;
this means that a cluster $\cluster_n\subset\borb_n$ is fixed (and
described elsewhere) so that the restriction of~$\pat$ to $\borb_n$
takes two values:
$\pat(\orb)=a$ for $\orb\subset\cluster_n$ and $\pi(\orb)=b$ for
$\orb\subset\borb_n\sminus\cluster_n$.
\endremark

\remark\label{s.check}
In some cases, where $\bnd(\Orb)$ exceeds the goal by just a few units, we
use patterns to show directly that $\Bnd_{\bnd(\Orb)-M}(\Orb)=\varnothing$.
These cases are
marked with a \checkmark\ in the tables, and
any further explanation is typically omitted%
\iffullversion\else\
(and so usually is the list of orbits)\fi.
\endremark

\subsection{Clusters}\label{ss.clusters}
Sometimes, the number of combinatorial orbits in an orbit~$\borb$ is too
large, making it difficult to compute all patterns. In these cases, we
subdivide $\borb$ into a number of \emph{clusters}
$\cluster_k\subset\borb$, not
necessarily disjoint, and compute patterns and, then, geometric sets cluster
by cluster. The subdivision is chosen so that $\stab\hbar$ acts transitively
on the set of clusters.
To reduce the overcounting, we assume that the clusters are ordered
lexicographically, by the decreasing of the sequence
\[*
\bigl(\ls|\fL\cap\cluster_k|,\dif_0(\cluster_k),\dif_1(\cluster_k),\ldots\bigr),
\qquad
\dif_i(\cluster_k):=
 \#\bigl\{\orb\subset\cluster_k\bigm|\ls|\fL\cap\orb|=\bnd(\orb)-i\bigr\}.
\]
In particular, this convention implies that, when computing the set
$\Bnd_d(\borb)$, for the first cluster $\cluster_1$ one must have
$\ls|\fL\cap\cluster_1|\ge\bnd(\cluster_1)-md/n$, where $n$ is the
total number of
clusters and the \emph{multiplicity} $m$ is the number of clusters containing
any fixed orbit $\orb\subset\borb$. More generally, extending a geometric
set~$\fL$ from $\cluster_1,\ldots,\cluster_{k}$ to the next cluster
$\cluster_{k+1}$, one must have
$\ls|\fL\cap\cluster_{k+1}|\le\ls|\fL\cap\cluster_{k}|$ and
\[*
\ls|\fL\cap\cluster_{k+1}|\ge\bnd(\cluster_{k+1})-
 \frac1{n-k}\left(md-\sum_{i=1}^k\bigl(\bnd(\cluster_i)-\ls|\fL\cap\cluster_i|\bigr)\right).
\]
Certainly, if the clusters are not disjoint, we also take into account the
intersections $\cluster_{k+1}\cap\cluster_i$, $i=1,\ldots,k$, when computing
the restricted patterns $\pat\:\cluster_{k+1}\to\NN$.

\section{Counts and bounds}\label{S.bounds}

In this section, we explain the computation of the bounds $\bnd(\orb)$ on the
number of lines in an admissible set
within a combinatorial orbit~$\orb$, see \eqref{eq.count.bound}.

\subsection{Blocks}\label{s.orbits}
Consider a combinatorial
orbit~$\orb$. In order to estimate the count $\cnt(\orb)$ and
bound $\bnd(\orb)$, we break the root system~$\bR$ into \emph{blocks},
$\bR=\bB_1\oplus\bB_2\oplus\ldots$, each block~$\bB_k$ consisting of whole
components~$\bR_i$. Then, $\hbar$ and $l\in\fF(\hbar)\cap\orb$ split into
$\bigoplus_k\hbar_k$ and $\bigoplus_kl_k$, respectively, with
$\hbar_k,l_k\in\bB_k\dual$. We denote by
$\orb|_k:=\orb|_{\bB_k}\subset\bB_k\dual$ the
\emph{restriction} of~$\orb$ to~$\bB_k$
(which, in fact, is nothing but the
orthogonal projection of $\orb$ to $\bB_k\dual$). This restriction
consists of a whole $\RG_{\hbar_k\!}(\bB_k)$-orbit of vectors; in particular,
we have
a well defined square $l_k^2\in\Q$, product $l_k\cdot\hbar_k\in\Q$,
and discriminant
class $l_k\bmod\bB_k\in\discr\bB_k$. Usually, these data determine an
irreducible
block up to isomorphism, the reason being the following simple observation
(which follows from the fact that \emph{all} roots in~$\bN$ are assumed to
lie in~$\bR$):
\roster*
\item
each vector $l_k\in\orb|_k$ is either integral, $l_k\in\bB_k$
(and then $l_k^2=0$, $2$,
or~$4$) or \emph{shortest vector} in its discriminant class;
\item
each vector $\hbar_k$ is either integral, $\hbar_k\in\bB_k$ (and then
$\hbar_k^2=0$, $2$, $4$, or~$6$) or shortest or \emph{second shortest}
vector in its discriminant class.
\endroster
Here, \emph{shortest} are the vectors minimizing the square within a given
discriminant class, whereas \emph{second shortest} are those of square
$(\mbox{minimum}+2)$.
In fact, $\hbar_k$ can be a second shortest vector in at most one
block~$\bB_k$.

The \emph{count} of a block~$\bB$ is defined in the obvious way:
$\cnt(\bB)=\bigl|{\orb|_\bB}\bigr|$. The \emph{bound} is defined \via\
$\bnd(\bB)=\max\ls|\fB|$, where $\fB\subset\orb|_\bB$ is a
$*$-invariant (if $\orb^*=\orb$) subset satisfying
the following condition:
\[
\text{for $l',l''\in\fB$, one has $l^{\prime2}-l'\cdot l''=0$ (iff $l'=l''$),
 $2$, $3$, or $5$ (iff $l'=(l'')^*$)}.
\label{eq.intersection.block}
\]
In other words, we bound the cardinality of subsets $\fL\subset\orb$
satisfying~\eqref{eq.intersection} and such that all lines $l\in\fL$ have the
same fixed restriction to all other blocks $\bB'\ne\bB$.

If $\bR$ is broken into two blocks, $\bB_1\oplus\bB_2$, we obviously have
\[
\cnt(\orb)=\cnt(\bB_1)\cnt(\bB_2),\qquad
\bnd(\orb)\le\min\bigl\{\cnt(\bB_1)\bnd(\bB_2),\bnd(\bB_1)\cnt(\bB_2)\bigr\}.
\label{eq.bound.2}
\]
By induction, for any number of blocks~$\bB_k$, this implies
\[
\cnt(\orb)=\prod_k\cnt(\bB_k),\qquad
\bnd(\orb)\le\cnt(\orb)\min_k\frac{\bnd(\bB_k)}{\cnt(\bB_k)}.
\label{eq.bound}
\]

This bound (with $\bB_k=\bR_k$ the irreducible components of~$\bR$)
and corresponding bound on $\bnd(\Orb)$ given
by~\eqref{eq.naive} are always listed first in the tables below.
If $\bnd(\Orb)\ge M$, we try to improve the bounds
$\bnd(\orb)$ using one of the following
arguments:
\roster
\item\label{bnd.self-dual}
\autoref{lem.3} below applied to an appropriate splitting into two blocks;
\item\label{bnd.pattern.like}
a computation using larger blocks, see \autoref{ss.brute.force.blocks} below;
\item\label{bnd.brute.force}
a brute force enumeration of admissible subsets $\fL\subset\orb$;
the bounds whose sharpness is confirmed by this
computation are underlined.
\endroster
In the tables, we refer to this list for the reasons for the improved
bounds.

\newcs{reason.2}{\noexpand\iref{bnd.self-dual}}
\newcs{reason.3}{\noexpand\iref{bnd.pattern.like}}
\newcs{reason.5}{\noexpand\iref{bnd.brute.force}}

\subsection{Brute force \via\ blocks}\label{ss.brute.force.blocks}
For some large combinatorial orbits~$\orb$, the exact computation of
$\bnd(\orb)$ by brute force is not feasible, and we improve the original
bound given by~\eqref{eq.bound} by using larger blocks.
Typically, we consider two blocks $\bB_1$ (one of the irreducible components
of~$\bR$) and $\bB_2$ (the sum of all other components).
Then, we compute all \emph{admissible} (rather than just
satisfying~\eqref{eq.intersection.block}) sets
$\fL(l_1)\subset\orb$ with a fixed
restriction $l_1\in\bB_1\dual$, replacing~\eqref{eq.bound.2} with
\[*
\bnd(\orb)\le\cnt(\bB_1)\max\ls|\fL(l_1)|.
\]
If this bound is still not good enough, we vary $l_1\in\bB_1\dual$ and try to
construct an admissible set $\fL\subset\orb$ by packing together
precomputed \emph{large}
(usually maximal or submaximal) sets $\fL(l_1')$, $\fL(l_1'')$, \etc.,
obtaining a better bound and, if necessary, a complete list of large
admissible sets in~$\orb$.

\subsection{Self-dual combinatorial orbits}\label{s.self-dual}
Let $\orb$ be a self-dual combinatorial orbit, $\orb^*=\orb$,
and break~$\bR$ into
blocks~$\bB_k$. Each block is also self-dual: $\bar l_k:=\hbar_k-l_k\in\orb|_k$
whenever $l_k\in\orb|_k$. In particular, $\bar l_k=l_k\bmod\bB_k$.
Hence, we have
\begin{align*}
  & 2l_k\cdot\hbar_k=\hbar_k^2\quad(\text{since $(\bar l_k^2=l_k^2$}), \\
  & l_k\cdot\bar l_k=l_k^2-\Gd_k\quad\text{for some $\Gd_k\in\Z$}.
\end{align*}
The integer $\Gd(\bB_k):=\Gd_k=2l_k^2-l_k\cdot\hbar_k$,
constant throughout the block, is called the \emph{defect} of
the block~$\bB_k$; it
takes values in the range $0\le\Gd_k\le5$, and the defects of all blocks sum
up to $5=l^2-l\cdot l^*$. Furthermore, for any pair of vectors
$l',l''\in\orb|_k$, the difference $l_k^2-l'\cdot l''$ is an integer
taking values in
\[
0\le l_k^2-l'\cdot l''\le\Gd_k,
\label{eq.self-dual}
\]
the two extreme values corresponding to $l''=l'$ and $l''=\bar l'$,
respectively. As a consequence, we have $\bnd(\bB_k)\le1$ if $\Gd(\bB_k)=1$ and
$\bnd(\bB_k)\le2$ if $\Gd(\bB_k)=2$; in the latter case, all maximal
admissible subsets
are of the form $\{l_k,\bar l_k\}$.

\lemma\label{lem.3}
Assume that a self-dual orbit~$\orb$ is broken into two blocks,
$B_2$ and~$B_3$, of defects~$2$ and~$3$, respectively. Then
\[*
\bnd(\orb)\le\max\bigl\{4u + \min\{\cnt_3 - 2u, (\cnt_2 - 2u)\bnd_3\}\bigm|
 u=0,\ldots,\tfrac12\min\{\cnt_2, \cnt_3\}\bigr\},
\]
where we abbreviate $\cnt_\Gd:=\cnt(\bB_\Gd)$ and $\bnd_\Gd:=\bnd(\bB_\Gd)$,
$\Gd=2,3$.
\endlemma

\proof
Let $\fL\subset\orb$ be an admissible set, and let $l_2\oplus l_3\in\fL$.
There is a dichotomy: either $\bar l_2\oplus l_3$ is in~$\fL$ or it is not.
In the former case, we have
\[*
\{l_2\oplus l_3,\bar l_2\oplus l_3,l_2\oplus\bar l_3,\bar l_2\oplus\bar l_3\}
 \subset\fL
\]
and, by~\eqref{eq.intersection} and~\eqref{eq.self-dual}, no other vector
$l_2\oplus l_3'$ or $\bar l_2\oplus l_3'$ with $l_3'\ne l_3,\bar l_3$ is
in~$\fL$. Each $4$-element subset of this form consumes two vectors
from~$\orb|_{3}$, and all these vectors are pairwise distinct.
Let $U\subset\orb|_2$ be the set of vectors~$l_2$ as above, and
denote $u:=\ls|U|$; clearly,
$0\le2u\le\min\{\cnt_2,\cnt_3\}$.

Otherwise, in the obvious notation, we have
\[*
l_2\oplus S(l_2)\subset\fL,\qquad \bar l_2\oplus S(\bar l_2)\subset\fL,
\]
where $S(l_2)\subset\orb|_3$ is a certain subset and
$S(\bar l_2)=\overline{S(l_2)}$.
Since $S(l_2)\cap S(\bar l_2)=\varnothing$ by the assumption, all subsets
$S(l_2)$, $l_2\in\orb|_2\sminus U$, are pairwise disjoint and do not
contain any of the $2u$ vectors~$l_3$ coupled with $l_2\in U$; hence, their total
cardinality does not exceed $\cnt_3-2u$. On the other hand, since
$\ls|S(l_2)|\le\bnd_3$ for each $l_2\in\orb|_2$, this cardinality does
not exceed $(\cnt_2-2u)\bnd_2$. Taking the minimum and maximizing over all
values of~$u$, we arrive at the bound in the statement.
\endproof

\subsection{Computing counts and bounds}\label{s.computing}
For ``small'' blocks $\bB_k\cong\bA_{\le7}$, $\bD_{\le7}$, $\bE_6$, $\bE_7$,
$\bE_8$, the counts $\cnt(\bB_k)$ and bounds $\bnd(\bB_k)$
used in~\eqref{eq.bound} are obtained by a
direct computation. For larger blocks, we use the standard combinatorial
description of the $\bA$- and $\bD$-type root systems as sublattices of the
odd unimodular lattice
\[*
\bH_{n}:=\bigoplus\Z\e_i,
\quad e^2_i=1,
\quad i\in\IS:=\{1,\ldots,n\}.
\]
(When working with this lattice, we let $\vv{o}:=\sum_{i\in o}\e_i$
for a subset $o\subset\IS$.)
Then, given a vector $\hbar_k=\sum_i\Ga_ie_i\in\bH_n\otimes\Q$,
we subdivide the block $\bB_k\dual\subset\bH_n\otimes\Q$ into ``subblocks''
\[*
\textstyle
\bB_k(\Ga):=\bigl\{\sum_i\Gb_ie_i\bigm|i\in\supp(\Ga)\bigr\},
\quad\supp(\Ga):=\bigl\{i\in\IS\bigm|\Ga_i=\Ga\bigr\},
\]
on which $\hbar_k$ is constant.
We obtain counts and bounds, in the sense of~\eqref{eq.intersection.block},
for each subblock and use an obvious analogue of~\eqref{eq.bound} to
estimate $\bnd(\bB_k)$. The technical details are outlined in the
next two sections.

\subsection{Root systems $\bA\sb{n}$}\label{s.An}
A block~$\bB_k$ of type $\bA_n$ is $\vv{\IS}^\perp\subset\bH_{n+1}$:
\[*
\textstyle
\bA_n=\bigl\{\sum_i\Ga_i\e_i\in\bH_{n+1}\bigm|\sum_i\Ga_i=0\bigr\}.
\]
One has $\discr\bA_n=\Z/(n+1)$, with a generator of square $n/(n+1)\bmod2\Z$,
and the shortest representatives of the
discriminant classes are vectors of the form
\[*
\be_{o}:=\dfrac1{n+1}
 \bigl(\ls|\bar o|\vv{o}-\ls|o|\vv{\bar o}\bigr),\qquad
 \be_{o}^2=\dfrac{\ls|o|\ls|\bar o|}{n+1},
\]
where $o\subset\IS$ and $\bar o$ is the complement. We have
$\be_{\bar o}=-\be_o$ and
\[*
\be_{r}\cdot\be_{s}=\ls|r\cap s|-\frac{\ls|r|\ls|s|}{n+1}.
\]
If $\ls|r|=\ls|s|$, or, equivalently,
$\e_r$ and $\e_s$ are in the same
discriminant class, then
\[
\be_{r}^2-\be_{r}\cdot\be_{s}=\frac12\ls|r\sdif s|,
\label{eq.A.intr}
\]
where $\sdif$ is the symmetric difference. Hence, in the case where $l_k$ is a
shortest vector in its (nonzero) discriminant class, the bound
$\bnd(\bB_k(\Ga))$
can be estimated by the following lemma, applied to $S=\supp(\Ga)$.

\lemma\label{lem.max.set.A}
Consider a finite set~$S$, $\ls|S|=n$, and let $\fS$ be a collection of
subsets $s\subset S$ with the following properties\rom:
\roster
\item\label{sets.m}
all subsets $s\in\fS$ have the same fixed cardinality~$m$\rom;
\item\label{sets.dif}
if $r,s\in\fS$, then $\ls|r\sdif s|\in\{0,4,6,10\}$\rom;
\item\label{sets.compl}
in the case $(n,m)=(10,5)$, if $s\in\fS$, then also $\bar s\in\fS$.
\endroster
Then, for small $(n,m)$, the maximal cardinality $\ls|\fS|$
is as follows\rom:
\[*
\minitab\
(n,m):&    (n,1)&(n,2)             &(6,3)&(7,3)&(8,3)&(9,3)&(10,3)&(11,3)&(8,4)&(9,4)&(10,5)\cr
\max\ls|\fS|:& 1&\lfloor n/2\rfloor&    4&    7&    8&   12&    13&    17&    9&   12&    24\cr
\endminitab
\]
More generally, for $m=3$ one has
$\ls|\fS|\le\bigl\lfloor n\lfloor(n-1)/2\rfloor/3\bigr\rfloor$.
\endlemma

Note that, if a collection~$\fS$ is as in the lemma, then so
is the collection $\{\bar s\,|\,s\in\fS\}$.
Hence, we can always assume that $2m\le n$.

\proof[Proof of \autoref{lem.max.set.A}]
The first two values are obvious;
the others are obtained by listing all admissible
collections.
The general estimate for $m=3$ follows from the observation that any two
subsets in~$\fS$ have at most one common point and, hence, each point of~$S$ is
contained in at most $\lfloor(n-1)/2\rfloor$ subsets.
\endproof

There remains to consider a subblock $\bB_k(\Ga)$ of a block~$\bB_k$
containing vectors of the form $l_k=\vv{r}-\vv{s}$, where
$r,s\subset\IS$, $r\cap s=\varnothing$, and $\ls|r|=\ls|s|=1$ or~$2$.
In the latter case, one must have $l_k\cdot\hbar_k=3$, and it follows that
$\ls|(r\cup s)\cap\supp(\Ga)|\le2$ for each $\Ga\in\Q$.
The bounds are as follows:
\roster
\item\label{bound.An.1}
if $\ls|(r\cup s)\cap\supp(\Ga)|=1$, then, obviously, $\bnd(\bB_k(\Ga))=1$;
\item\label{bound.An.2}
if $\ls|r\cap\supp(\Ga)|=2$ (or $\ls|s\cap\supp(\Ga)|=2$),
the distinct sets $r\cap\supp(\Ga)$
must be
pairwise disjoint and, hence,
$\bnd(\bB_k(\Ga))=\bigl\lfloor\frac12\ls|\supp(\Ga)|\bigr\rfloor$;
\item\label{bound.An.roots}
if $\ls|r\cap\supp(\Ga)|=\ls|s\cap\supp(\Ga)|=1$,
then the distinct sets $r\cap\supp(\Ga)$ must also be pairwise disjoint and, hence,
$\bnd(\bB_k(\Ga))=\ls|\supp(\Ga)|$.
\endroster

\subsection{Root systems $\bD\sb{n}$}\label{s.Dn}
A block~$\bB_k$ of type~$\bD_n$
can be defined as the maximal even sublattice in~$\bH_n$:
\[*
\textstyle
\bD_n=\bigl\{\sum_i\Ga_i\e_i\in\bH_{n}\bigm|\sum_i\Ga_i=0\bmod2\bigr\}.
\label{eq.Dn}
\]
One has $\discr\bD_n=\Z/2\oplus\Z/2$ (if $n$ is even) or $\Z/4$ (if $n$ is
odd); the shortest vectors are
\[*
\e_i,\ i\in\IS,\quad\text{and}\quad
 \be_o:=\frac12(\vv{o}-\vv{\bar o}),\
 o\subset\IS,\quad
 \be_o^2=\frac{n}4
\]
(the class $\be_o\bmod\bD_n$ depends on the parity of~$\ls|o|$)
and we have a literal analogue of~\eqref{eq.A.intr} for any pair
$r,s\subset\IS$.
Thus, if $\bB_k\ni\be_o$,
the bounds $\bnd(\bB_k(\Ga))$ are estimated by
\autoref{lem.max.set.A} (if $\Ga\ne0$) or \autoref{lem.max.set.D} below (if
$\Ga=0$)
applied to $S=\supp(\Ga)$.

\lemma\label{lem.max.set.D}
For $n\le10$, the maximal cardinality
of a collection~$\fS$ satisfying conditions~\iref{sets.dif}
and~\iref{sets.compl} \rom(if $n=10$\rom) of \autoref{lem.max.set.A}
is bounded as follows\rom:
\[*
\minitab\quad
n=&         1&2&3&4&5&6&7& 8& 9&10\cr
\ls|\fS|\le&1&1&1&2&2&4&8&10&16&32\cr
\endminitab
\]
These bounds are sharp for $n\le8$.
\endlemma

\proof
If $n\le6$, the statement is easily proved by inspection, using
\autoref{lem.max.set.A}.

Let $n=8$. Represent a subset $s\in\fS$ as the root
$\be_s\in\bD_8\dual$. Then,
condition~\iref{sets.dif} implies that all subsets $s\in\fS$ have
cardinality of the same parity and, hence, all roots are in the same
discriminant class; thus, they
lie in an extension $\bE_8\supset\bD_8$.
By~\iref{sets.dif}, the roots~$\be_s$ constitute a union~$\Gamma$ of
(affine) Dynkin diagrams other than \smash{$\tA_1$}
admitting an isometry to~$\bE_8$,
which gives us a
bound $\ls|\fS|\le12$. Furthermore, the
roots~$\be_s$ are distinguished by the property $\be_s\cdot2\e_1=1\bmod2$.
Thus,
each affine component of~$\Gamma$ must have even degree.
The maximal graph with these properties is $2\tD_4$, resulting in the bound
$\ls|\fS|\le10$.

If $n=7$, we extend the ambient set~$S$ and each subset by an extra point and
argue as above, obtaining roots $\be_s\in\bE_8$ with the
property $\be_s\cdot2\e_1=1$. This time, the roots are linearly independent and
$\ls|\fS|\le\rank\bE_8=8$. This is realized by $2\bD_4$.

In general, represent $s\in\fS$ by the vector $\vv{s}\in\bH_n$. (If
$n=10$, select one subset~$s$ from each pair $s$, $\bar s$.)
Then
$\vv{s}^2=n$ and the products $\vv{r}\cdot\vv{s}=n-2\ls|r\sdif s|$, $r\ne s$,
take but two values $n-8$ or $n-12$.
Since $n\le10$, \autoref{th.Elkies} applies and bounds the number of vectors
by $16$. If $n=10$, this bound is to be doubled.
\endproof

The few remaining cases are listed below.
\roster
\item\label{bound.Dn.2e}
If $\bB_k(\Ga)\ni\pm2\e_i$, $i\in\supp(\Ga)$, then $\bnd(\bB_k(0))=1$.
\endroster
Assume that
$l_k=\sum(\pm\e_i)$, $i\in o\subset\CS$, $\ls|o|\le4$.
If $\Ga=0$, then
\roster[\lastitem]
\item\label{bound.Dn.zero}
$\ls|o\cap\supp(\Ga)|=0$, $1$, or $2$ and
$\bnd(\bB_k(\Ga))\le1$, $2$, or $\frac43\ls|\supp(\Ga)|$, respectively,
\endroster
similar to \autoref{s.An}.
(Here, the last number is a bound on the size of
a union of (affine) Dynkin diagrams other than \smash{$\tA_1$} admitting an
isometry to $\bD_{\ls|\supp(\Ga)|}$.)
If $\Ga\ne0$, the numbers of
signs $\pm$ within $\supp(\Ga)$ are also fixed, and the options are as
follows:
\roster[\lastitem]
\item\label{bound.Dn.same}
$m:=\ls|o\cap\supp(\Ga)|\le3$ and all signs are the same:
by an analogue of~\eqref{eq.A.intr},
a bound
on $\bnd(\bB_k(\Ga))$ is given by \autoref{lem.max.set.A} applied to
$S=\supp(\Ga)$;
\item\label{bound.Dn.roots}
$\ls|o\cap\supp(\Ga)|=2$ and the signs differ:
$\bnd(\bB_k(\Ga))=\ls|\supp(\Ga)|$ as in \autoref{s.An}\iref{bound.An.roots}.
\endroster

\remark\label{rem.Dn}
If $n\ge5$, the group $\OG(\bD_n)$ is an index~$2$ extension of~$\RG(\bD_n)$:
it is
generated by the reflection against the hyperplane orthogonal to any
of~$\e_i$. Hence, up to $\OG(\bD_n)$, we can assume that, in the expression
$\hbar_k=\sum_i\Ga_i\e_i$, all coefficients $\Ga_i\ge0$. We always make this
assumption (and adjust the results afterwards) when describing the orbits and
computing counts and bounds.
\endremark

\begingroup
\iffullversion\raggedbottom\fi

\section{Root systems with few components}\label{S.few}

In this section, we consider the $20$ Niemeier lattices generated
over~$\Q$ by
root systems with few (up to six) irreducible components.
We set the goal
\[*
\ls|\fL|\ge M:=122
\]
and prove the following theorem.

\theorem\label{th.<=6}
Fix a root system~$\bR$ with at most six irreducible components and
a configuration $(\hbar,\rbar)$ in the Niemeier lattice $\N(\bR)$.
Then, with
the exception of
\roster*
\item
$\ls|\fL|=144$ and $\fL$ is conjugate to $\Lmax1\subset\N(4\bA_5\oplus\bD_4)$,
see \eqref{eq.Lmax.1}, or
\item
$\ls|\fL|=130$ and $\fL$ is conjugate to $\Lmisc1\subset\N(6\bA_4)$,
see \eqref{eq.Lmisc.1},
\endroster
one has $\ls|\fL|\le120$ for each geometric set~$\fL$.
\endtheorem


\proof
For each configuration $(\hbar,\rbar)$ (or just vector $\hbar$), we list all
$\OG_\hbar(\bN)$-orbits $\borb_n$ and indicate the number
$m(\borb_n)$ of combinatorial orbits $\orb\subset\borb_n$, the count
$\cnt(\orb)$, and the na\"{\i}ve bound on $\ls|\fL\cap\orb|$ given
by~\eqref{eq.bound}. Sometimes, this bound is improved by one of the
arguments \iref{bnd.self-dual}--\iref{bnd.brute.force} in \autoref{s.orbits};
the best bound obtained is denoted by $\bnd(\orb)$.
The results are listed in several tables below.

\convention\label{conv.tables}
In the tables, the number
$m(\borb_n)$ is marked with a $^*$ if $\borb_n$ is
self-dual; it is marked with $^{**}$ if also each combinatorial orbit
$\orb\subset\borb$ is self-dual. If $\borb_n$ is not self-dual, then its dual
$\borb_n^*=\borb_{n+1}$ is omitted.


For the components $\hbar_k$ of~$\hbar$ we use the notation
$\hvec{\hbar_k^2}_d$, where $d$ is either the discriminant class
of~$\hbar_k$ or, if $\hbar_k\in\bR_k$, the
symbol $0$ (if $\hbar_k=0$), $\circ$ (if $\hbar_k^2=2$),
$\bullet$ (if $\hbar_k^2=4$), or $*$ (if $\hbar_k^2=6$).
If these data do not determine~$\hbar_k$, we use a superscript:
\roster*
\item
$+$ or $-$ to select a second shortest vector
(one of the form $\be_o+r$, where $r=\e_i-\e_j\in\bR_k$ is a root and
$i,j\in o$ or $i,j\in\bar o$,
respectively, see \autoref{s.An} and \autoref{s.Dn})
in a discriminant class
$d\ne0$;
\item
$+$, if $\hbar=\hbar_k\in\bD_n\subset\bH_n$ or $\bA_{n-1}\subset\bH_n$ is
of the form $2\e_1-\e_2-\e_3$ rather than $\e_1+\e_2+\e_3-\e_4-\e_5-\e_6$,
see \autoref{s.An} and \autoref{s.Dn};
\item
the discriminant class of $\frac12\hbar_k$, if
$\hbar_2\in\bD_n\cap2\bD_n\dual$.
\endroster
If $\bR_k$ contains the root~$\rbar$, this notation is changed to
$\makeatletter\r@@ttrue\hvec{\hbar_k^2}_d$.

For the components $l_k$ of a line, we use the notation
$\lvec{l_k\cdot\hbar_k}_d$, where $d$ and an occasional superscript have
the same meaning as for~$\hbar$.
\endconvention

Also shown in the tables is the na\"{\i}ve \latin{a priori} estimate $\bnd(\Orb)$
given by~\eqref{eq.naive}.
For the vast majority of configurations we
have
$\bnd(\Orb)\le M$\iffullversion. \else, and these configurations are omitted.
(The complete set of tables is available
in~\cite{degt:sextics.full}.) \fi
The few cases where $\bnd(\Orb)\ge M$
are shown in bold, and we treat them separately below, except those marked with a
\checkmark\ (see \autoref{s.check}).
\iffullversion\else
In these ``trivial'' cases marked with
a \checkmark\ we usually also omit the list of orbits.
\fi


\newcs{max-D24}{56}%
\newcs{count-D24}{3}%
\newcs{counts-D24}{\DO{1}{56}\DO{2}{40}\DO{3}{0}}%
\orbits
\orbit{1}{1}{56}
\h{1}&\int(6){\6}{\spec{1}}&&\cnt{1680}&\bounds{56}{56}{1}{false}&&\cr
\l{1}&\int(3){\4}{}&\ln{1\sdd}&\cnt{1680}&\bounds{56}{56}{1}{false}&&\cr
\endorbit
\orbit{1}{2}{40}
\h{2}&\int(6){\6}{}&&\cnt{720}&\bounds{40}{40}{1}{false}&&\cr
\l{1}&\int(3){\4}{}&\ln{1\sdd}&\cnt{720}&\bounds{40}{40}{1}{false}&&\cr
\endorbit
\orbit{1}{3}{0}
\h{3}&\int(6){1}{}&&\cnt{0}&\bounds{0}{0}{1}{false}&&\cr
\endorbit
\endorbits

\newcs{max-D16_E8}{76}%
\newcs{count-D16_E8}{8}%
\newcs{counts-D16_E8}{\DO{1}{76}\DO{2}{72}\DO{3}{72}\DO{4}{72}\DO{5}{60}\DO{6}{58}\DO{7}{40}\DO{8}{16}}%
\orbits
\orbit{2}{1}{76}
\h{1}&\int(4){\4}{2}&\int(2){\2}{}&&\cnt{1680}&\bounds{76}{112}{fail}{false}&&\cr
\l{1}&\int(2){\2}{}&\int(1){\2}{}&\ln{1\sdd}&\cnt{1680}&\bounds{76}{112}{2}{false}&&\cr
\endorbit
\orbit{2}{2}{72}
\h{2}&\int(2){\2}{}&\int(4){\4}{}&&\cnt{912}&\bounds{72}{72}{1}{false}&&\cr
\l{1}&\int(1){\2}{}&\int(2){\2}{}&\ln{1\sdd}&\cnt{784}&\bounds{56}{56}{1}{false}&&\cr
\l{2}&\int(2){\2}{}&\int(1){\2}{}&\ln{1}&\cnt{64}&\bounds{8}{8}{1}{false}&&\cr
\endorbit
\orbit{2}{3}{72}
\h{3}&\int(6){\6}{}&\int(0){\0}{}&&\cnt{912}&\bounds{72}{72}{1}{false}&&\cr
\l{1}&\int(3){1}{}&\int(0){\0}{}&\ln{1\sdd}&\cnt{512}&\bounds{32}{32}{1}{false}&&\cr
\l{2}&\int(3){\4}{}&\int(0){\0}{}&\ln{1\sdd}&\cnt{400}&\bounds{40}{40}{1}{false}&&\cr
\endorbit
\orbit{2}{4}{72}
\h{4}&\int(6){1}{}&\int(0){\0}{}&&\cnt{912}&\bounds{72}{72}{1}{false}&&\cr
\l{1}&\int(3){1}{\spec{1}}&\int(0){\0}{}&\ln{1}&\cnt{455}&\bounds{35}{35}{1}{false}&&\cr
\l{3}&\int(3){1}{\spec{2}}&\int(0){\0}{}&\ln{1}&\cnt{1}&\bounds{1}{1}{1}{false}&&\cr
\endorbit
\orbit{2}{5}{60}
\h{5}&\int(0){\0}{}&\int(6){\6}{}&&\cnt{1104}&\bounds{60}{60}{1}{false}&&\cr
\l{1}&\int(0){\2}{}&\int(3){\2}{}&\ln{1\sdd}&\cnt{960}&\bounds{42}{42}{1}{false}&&\cr
\l{2}&\int(0){\0}{}&\int(3){\4}{}&\ln{1\sdd}&\cnt{144}&\bounds{18}{18}{1}{false}&&\cr
\endorbit
\orbit{2}{6}{58}
\h{6}&\int(6){\6}{\spec{1}}&\int(0){\0}{}&&\cnt{1104}&\bounds{58}{58}{1}{false}&&\cr
\l{1}&\int(3){\4}{}&\int(0){\0}{}&\ln{1\sdd}&\cnt{624}&\bounds{34}{34}{1}{false}&&\cr
\l{2}&\int(3){\2}{}&\int(0){\2}{}&\ln{1\sdd}&\cnt{480}&\bounds{24}{24}{1}{false}&&\cr
\endorbit
\orbit{2}{7}{40}
\h{7}&\int(4){\4}{}&\int(2){\2}{}&&\cnt{528}&\bounds{40}{40}{1}{false}&&\cr
\l{1}&\int(2){\2}{}&\int(1){\2}{}&\ln{1\sdd}&\cnt{336}&\bounds{24}{24}{1}{false}&&\cr
\l{2}&\int(3){\4}{}&\int(0){\0}{}&\ln{1}&\cnt{96}&\bounds{8}{8}{1}{false}&&\cr
\endorbit
\orbit{2}{8}{16}
\h{8}&\int(4){1}{}&\int(2){\2}{}&&\cnt{240}&\bounds{16}{16}{1}{false}&&\cr
\l{1}&\int(3){1}{}&\int(0){\0}{}&\ln{1}&\cnt{120}&\bounds{8}{8}{1}{false}&&\cr
\endorbit
\endorbits

\newcs{max-3E8}{72}%
\newcs{count-3E8}{3}%
\newcs{counts-3E8}{\DO{1}{72}\DO{2}{66}\DO{3}{24}}%
\orbits
\orbit{3}{1}{72}
\h{1}&\int(4){\4}{}&\int(2){\2}{}&\int(0){\0}{}&&\cnt{912}&\bounds{72}{72}{1}{false}&&\cr
\l{1}&\int(2){\2}{}&\int(1){\2}{}&\int(0){\0}{}&\ln{1\sdd}&\cnt{784}&\bounds{56}{56}{1}{false}&&\cr
\l{2}&\int(3){\4}{}&\int(0){\0}{}&\int(0){\0}{}&\ln{1}&\cnt{64}&\bounds{8}{8}{1}{false}&&\cr
\endorbit
\orbit{3}{2}{66}
\h{2}&\int(6){\6}{}&\int(0){\0}{}&\int(0){\0}{}&&\cnt{1104}&\bounds{66}{66}{1}{false}&&\cr
\l{1}&\int(3){\4}{}&\int(0){\0}{}&\int(0){\0}{}&\ln{1\sdd}&\cnt{144}&\bounds{18}{18}{1}{false}&&\cr
\l{2}&\int(3){\2}{}&\int(0){\2}{}&\int(0){\0}{}&\ln{2\sdd}&\cnt{480}&\bounds{24}{24}{1}{false}&&\cr
\endorbit
\orbit{3}{3}{24}
\h{3}&\int(2){\2}{}&\int(2){\2}{}&\int(2){\2}{}&&\cnt{336}&\bounds{24}{24}{1}{false}&&\cr
\l{1}&\int(2){\2}{}&\int(1){\2}{}&\int(0){\0}{}&\ln{6\sd}&\cnt{56}&\bounds{4}{4}{1}{false}&&\cr
\endorbit
\endorbits

\newcs{max-A24}{60}%
\newcs{count-A24}{5}%
\newcs{counts-A24}{\DO{1}{60}\DO{2}{44}\DO{3}{36}\DO{4}{30}\DO{5}{24}}%
\orbits
\orbit{1}{1}{60}
\h{1}&\int(6){20}{\spec{1}}&&\cnt{1260}&\bounds{60}{60}{1}{false}&&\cr
\l{1}&\int(3){20}{}&\ln{1}&\cnt{630}&\bounds{30}{30}{1}{false}&&\cr
\endorbit
\orbit{1}{2}{44}
\h{2}&\int(6){\6}{\spec{1}}&&\cnt{924}&\bounds{44}{44}{1}{false}&&\cr
\l{1}&\int(3){\4}{}&\ln{1\sdd}&\cnt{924}&\bounds{44}{44}{1}{false}&&\cr
\endorbit
\orbit{1}{3}{36}
\h{3}&\int(6){\6}{}&&\cnt{684}&\bounds{36}{36}{1}{false}&&\cr
\l{1}&\int(3){20}{}&\ln{2\sd}&\cnt{171}&\bounds{9}{9}{1}{false}&&\cr
\l{2}&\int(3){\4}{\spec{1}}&\ln{2\sd}&\cnt{171}&\bounds{9}{9}{1}{false}&&\cr
\endorbit
\orbit{1}{4}{30}
\h{4}&\int(6){20}{\spec{2}}&&\cnt{540}&\bounds{30}{30}{1}{false}&&\cr
\l{1}&\int(3){20}{}&\ln{1}&\cnt{270}&\bounds{15}{15}{1}{false}&&\cr
\endorbit
\orbit{1}{5}{24}
\h{5}&\int(6){15}{}&&\cnt{252}&\bounds{24}{24}{1}{false}&&\cr
\l{1}&\int(3){20}{}&\ln{1\sdd}&\cnt{252}&\bounds{24}{24}{1}{false}&&\cr
\endorbit
\endorbits

\newcs{max-2D12}{92}%
\newcs{count-2D12}{7}%
\newcs{counts-2D12}{\DO{1}{92}\DO{2}{88}\DO{3}{80}\DO{4}{60}\DO{5}{56}\DO{6}{40}\DO{7}{24}}%
\orbits
\orbit{2}{1}{92}
\h{1}&\int(6){\6}{}&\int(0){\0}{}&&\cnt{1008}&\bounds{92}{104}{fail}{false}&&\cr
\l{1}&\int(3){\4}{}&\int(0){\0}{}&\ln{1\sdd}&\cnt{240}&\bounds{40}{40}{1}{false}&&\cr
\l{2}&\int(3){1}{}&\int(0){2}{}&\ln{1\sdd}&\cnt{768}&\bounds{52}{64}{2}{false}&&\cr
\endorbit
\orbit{2}{2}{88}
\h{2}&\int(4){\4}{2}&\int(2){\2}{}&&\cnt{1392}&\bounds{88}{112}{fail}{false}&&\cr
\l{1}&\int(2){\2}{}&\int(1){\2}{}&\ln{1\sdd}&\cnt{880}&\bounds{56}{80}{2}{false}&&\cr
\l{2}&\int(2){2}{}&\int(1){1}{}&\ln{1\sdd}&\cnt{512}&\bounds{32}{32}{1}{false}&&\cr
\endorbit
\orbit{2}{3}{80}
\h{3}&\int(1){2}{}&\int(5){1}{}&&\cnt{816}&\bounds{80}{80}{1}{false}&&\cr
\l{1}&\int(1){\2}{}&\int(2){\2}{}&\ln{1}&\cnt{242}&\bounds{22}{22}{1}{false}&&\cr
\l{3}&\int(1){2}{}&\int(2){1}{\spec{1}}&\ln{1}&\cnt{165}&\bounds{17}{17}{1}{false}&&\cr
\l{5}&\int(1){2}{}&\int(2){1}{\spec{2}}&\ln{1}&\cnt{1}&\bounds{1}{1}{1}{false}&&\cr
\endorbit
\orbit{2}{4}{60}
\h{4}&\int(4){\4}{}&\int(2){\2}{}&&\cnt{624}&\bounds{60}{60}{1}{false}&&\cr
\l{1}&\int(2){1}{}&\int(1){2}{}&\ln{1\sdd}&\cnt{256}&\bounds{20}{20}{1}{false}&&\cr
\l{2}&\int(2){\2}{}&\int(1){\2}{}&\ln{1\sdd}&\cnt{240}&\bounds{24}{24}{1}{false}&&\cr
\l{3}&\int(3){\4}{}&\int(0){\0}{}&\ln{1}&\cnt{64}&\bounds{8}{8}{1}{false}&&\cr
\endorbit
\orbit{2}{5}{56}
\h{5}&\int(6){\6}{\spec{1}}&\int(0){\0}{}&&\cnt{816}&\bounds{56}{56}{1}{false}&&\cr
\l{1}&\int(3){\4}{}&\int(0){\0}{}&\ln{1\sdd}&\cnt{288}&\bounds{24}{24}{1}{false}&&\cr
\l{2}&\int(3){\2}{}&\int(0){\2}{}&\ln{1\sdd}&\cnt{528}&\bounds{32}{32}{1}{false}&&\cr
\endorbit
\orbit{2}{6}{40}
\h{6}&\int(3){2}{}&\int(3){1}{}&&\cnt{432}&\bounds{40}{40}{1}{false}&&\cr
\l{1}&\int(3){\4}{}&\int(0){\0}{}&\ln{1}&\cnt{18}&\bounds{2}{2}{1}{false}&&\cr
\l{3}&\int(2){\2}{}&\int(1){\2}{}&\ln{1}&\cnt{198}&\bounds{18}{18}{1}{false}&&\cr
\endorbit
\orbit{2}{7}{24}
\h{7}&\int(3){3}{}&\int(3){3}{}&&\cnt{288}&\bounds{24}{24}{1}{false}&&\cr
\l{1}&\rat(5/2){1}{}&\rat(1/2){2}{}&\ln{2\sd}&\cnt{144}&\bounds{12}{12}{1}{false}&&\cr
\endorbit
\endorbits

\newcs{max-A17_E7}{96}%
\newcs{count-A17_E7}{13}%
\newcs{counts-A17_E7}{\DO{1}{96}\DO{2}{72}\DO{3}{70}\DO{4}{66}\DO{5}{66}\DO{6}{60}\DO{7}{48}\DO{8}{48}\DO{9}{48}\DO{10}{48}\DO{11}{48}\DO{12}{48}\DO{13}{40}}%
\orbits
\orbit{2}{1}{96}
\h{1}&\int(0){\0}{}&\int(6){\6}{\spec{2}}&&\cnt{1632}&\bounds{96}{96}{1}{false}&&\cr
\l{1}&\int(0){15}{}&\int(3){1}{}&\ln{2\sd}&\cnt{816}&\bounds{48}{48}{1}{false}&&\cr
\endorbit
\orbit{2}{2}{72}
\h{2}&\int(6){12}{\spec{1}}&\int(0){\0}{}&&\cnt{960}&\bounds{72}{72}{1}{false}&&\cr
\l{1}&\int(3){15}{}&\int(0){1}{}&\ln{1\sdd}&\cnt{336}&\bounds{24}{24}{1}{false}&&\cr
\l{2}&\int(3){12}{}&\int(0){\0}{}&\ln{1}&\cnt{312}&\bounds{24}{24}{1}{false}&&\cr
\endorbit
\orbit{2}{3}{70}
\h{3}&\rat(9/2){15}{\spec{1}}&\rat(3/2){1}{}&&\cnt{1104}&\bounds{70}{70}{1}{false}&&\cr
\l{1}&\int(3){\4}{}&\int(0){\0}{}&\ln{1}&\cnt{120}&\bounds{8}{8}{1}{false}&&\cr
\l{3}&\rat(5/2){15}{}&\rat(1/2){1}{}&\ln{1}&\cnt{432}&\bounds{27}{27}{1}{false}&&\cr
\endorbit
\orbit{2}{4}{66}
\h{4}&\int(6){\6}{}&\int(0){\0}{}&&\cnt{768}&\bounds{66}{66}{1}{false}&&\cr
\l{1}&\int(3){12}{}&\int(0){\0}{}&\ln{2\sd}&\cnt{220}&\bounds{20}{20}{1}{false}&&\cr
\l{2}&\int(3){\4}{\spec{1}}&\int(0){\0}{}&\ln{2\sd}&\cnt{108}&\bounds{9}{9}{1}{false}&&\cr
\l{3}&\int(3){15}{}&\int(0){1}{}&\ln{2\sd}&\cnt{56}&\bounds{4}{4}{1}{false}&&\cr
\endorbit
\orbit{2}{5}{66}
\h{5}&\int(6){12}{\spec{2}}&\int(0){\0}{}&&\cnt{672}&\bounds{66}{66}{1}{false}&&\cr
\l{1}&\int(3){6}{}&\int(0){\0}{}&\ln{1\sdd}&\cnt{252}&\bounds{24}{24}{1}{false}&&\cr
\l{2}&\int(3){12}{}&\int(0){\0}{}&\ln{1}&\cnt{210}&\bounds{21}{21}{1}{false}&&\cr
\endorbit
\orbit{2}{6}{60}
\h{6}&\int(2){\2}{}&\int(4){\4}{}&&\cnt{864}&\bounds{60}{60}{1}{false}&&\cr
\l{1}&\int(1){15}{}&\int(2){1}{}&\ln{2\sd}&\cnt{240}&\bounds{16}{16}{1}{false}&&\cr
\l{2}&\int(1){\2}{\spec{1}}&\int(2){\2}{}&\ln{2\sd}&\cnt{160}&\bounds{10}{10}{1}{false}&&\cr
\l{3}&\int(2){\2}{}&\int(1){\2}{}&\ln{1}&\cnt{32}&\bounds{4}{4}{1}{false}&&\cr
\endorbit
\orbit{2}{7}{48}
\h{7}&\int(0){\0}{}&\int(6){\6}{\spec{1}}&&\cnt{672}&\bounds{48}{48}{1}{false}&&\cr
\l{1}&\int(0){\2}{}&\int(3){\2}{}&\ln{1\sdd}&\cnt{612}&\bounds{36}{36}{1}{false}&&\cr
\l{2}&\int(0){\0}{}&\int(3){\4}{}&\ln{1\sdd}&\cnt{60}&\bounds{12}{12}{1}{false}&&\cr
\endorbit
\orbit{2}{8}{48}
\h{8}&\int(4){\4}{}&\int(2){\2}{}&&\cnt{576}&\bounds{48}{48}{1}{false}&&\cr
\l{1}&\int(2){15}{}&\int(1){1}{}&\ln{2\sd}&\cnt{168}&\bounds{12}{12}{1}{false}&&\cr
\l{2}&\int(2){\2}{}&\int(1){\2}{}&\ln{1\sdd}&\cnt{128}&\bounds{16}{16}{1}{false}&&\cr
\l{3}&\int(3){\4}{\spec{1}}&\int(0){\0}{}&\ln{2}&\cnt{28}&\bounds{2}{2}{1}{false}&&\cr
\endorbit
\orbit{2}{9}{48}
\h{9}&\int(6){\6}{\spec{1}}&\int(0){\0}{}&&\cnt{672}&\bounds{48}{48}{1}{false}&&\cr
\l{1}&\int(3){\4}{}&\int(0){\0}{}&\ln{1\sdd}&\cnt{420}&\bounds{30}{30}{1}{false}&&\cr
\l{2}&\int(3){\2}{}&\int(0){\2}{}&\ln{1\sdd}&\cnt{252}&\bounds{18}{18}{1}{false}&&\cr
\endorbit
\orbit{2}{10}{48}
\h{10}&\rat(5/2){15}{}&\rat(7/2){1}{}&&\cnt{672}&\bounds{48}{48}{1}{false}&&\cr
\l{1}&\rat(5/2){15}{}&\rat(1/2){1}{}&\ln{1}&\cnt{21}&\bounds{3}{3}{1}{false}&&\cr
\l{3}&\rat(3/2){15}{}&\rat(3/2){1}{}&\ln{1}&\cnt{315}&\bounds{21}{21}{1}{false}&&\cr
\endorbit
\orbit{2}{11}{48}
\h{11}&\rat(9/2){15}{\spec{2}}&\rat(3/2){1}{}&&\cnt{528}&\bounds{48}{48}{1}{false}&&\cr
\l{1}&\int(3){12}{}&\int(0){\0}{}&\ln{1}&\cnt{78}&\bounds{6}{6}{1}{false}&&\cr
\l{3}&\int(3){\4}{}&\int(0){\0}{}&\ln{1}&\cnt{78}&\bounds{6}{6}{1}{false}&&\cr
\l{5}&\rat(5/2){15}{}&\rat(1/2){1}{}&\ln{1}&\cnt{108}&\bounds{12}{12}{1}{false}&&\cr
\endorbit
\orbit{2}{12}{48}
\h{12}&\rat(9/2){9}{}&\rat(3/2){1}{}&&\cnt{336}&\bounds{48}{48}{1}{false}&&\cr
\l{1}&\int(3){12}{}&\int(0){\0}{}&\ln{2}&\cnt{84}&\bounds{12}{12}{1}{false}&&\cr
\endorbit
\orbit{2}{13}{40}
\h{13}&\int(4){12}{}&\int(2){\2}{}&&\cnt{384}&\bounds{40}{52}{fail}{false}&&\cr
\l{1}&\int(2){15}{}&\int(1){1}{}&\ln{1\sdd}&\cnt{240}&\bounds{28}{40}{2}{false}&&\cr
\l{2}&\int(3){12}{}&\int(0){\0}{}&\ln{1}&\cnt{72}&\bounds{6}{6}{1}{false}&&\cr
\endorbit
\endorbits

\newcs{max-D10_2E7}{96}%
\newcs{count-D10_2E7}{14}%
\newcs{counts-D10_2E7}{\DO{1}{96}\DO{2}{96}\DO{3}{92}\DO{4}{88}\DO{5}{84}\DO{6}{78}\DO{7}{76}\DO{8}{76}\DO{9}{76}\DO{10}{56}\DO{11}{56}\DO{12}{56}\DO{13}{54}\DO{14}{40}}%
\orbits
\orbit{3}{1}{96}
\h{1}&\int(0){\0}{}&\int(6){\6}{\spec{2}}&\int(0){\0}{}&&\cnt{1632}&\bounds{96}{112}{fail}{false}&&\cr
\l{1}&\int(0){1}{}&\int(3){1}{}&\int(0){\0}{}&\ln{1\sdd}&\cnt{512}&\bounds{32}{32}{1}{false}&&\cr
\l{2}&\int(0){2}{}&\int(3){1}{}&\int(0){1}{}&\ln{1\sdd}&\cnt{1120}&\bounds{64}{80}{2}{false}&&\cr
\endorbit
\orbit{3}{2}{96}
\h{2}&\int(6){\6}{}&\int(0){\0}{}&\int(0){\0}{}&&\cnt{1056}&\bounds{96}{96}{1}{false}&&\cr
\l{1}&\int(3){\4}{}&\int(0){\0}{}&\int(0){\0}{}&\ln{1\sdd}&\cnt{160}&\bounds{32}{32}{1}{false}&&\cr
\l{2}&\int(3){3}{}&\int(0){\0}{}&\int(0){1}{}&\ln{2\sdd}&\cnt{448}&\bounds{32}{32}{1}{false}&&\cr
\endorbit
\orbit{3}{3}{92}
\h{3}&\int(4){\4}{2}&\int(2){\2}{}&\int(0){\0}{}&&\cnt{1248}&\bounds{92}{112}{fail}{false}&&\cr
\l{1}&\int(2){\2}{}&\int(1){\2}{}&\int(0){\0}{}&\ln{1\sdd}&\cnt{576}&\bounds{44}{64}{2}{false}&&\cr
\l{2}&\int(2){2}{}&\int(1){1}{}&\int(0){1}{}&\ln{1\sdd}&\cnt{672}&\bounds{48}{48}{1}{false}&&\cr
\endorbit
\orbit{3}{4}{88}
\h{4}&\rat(9/2){3}{}&\int(0){\0}{}&\rat(3/2){1}{}&&\cnt{768}&\bounds{88}{88}{1}{false}&&\cr
\l{1}&\int(3){\4}{}&\int(0){\0}{}&\int(0){\0}{}&\ln{1}&\cnt{84}&\bounds{12}{12}{1}{false}&&\cr
\l{3}&\int(3){\4}{2}&\int(0){\0}{}&\int(0){\0}{}&\ln{1}&\cnt{1}&\bounds{1}{1}{1}{false}&&\cr
\l{5}&\int(3){1}{}&\int(0){1}{}&\int(0){\0}{}&\ln{1}&\cnt{56}&\bounds{4}{4}{1}{false}&&\cr
\l{7}&\rat(5/2){3}{}&\int(0){\0}{}&\rat(1/2){1}{}&\ln{1}&\cnt{243}&\bounds{27}{27}{1}{false}&&\cr
\endorbit
\orbit{3}{5}{84}
\h{5}&\int(0){\0}{}&\int(4){\4}{}&\int(2){\2}{}&&\cnt{864}&\bounds{84}{96}{fail}{false}&&\cr
\l{1}&\int(0){2}{}&\int(2){1}{}&\int(1){1}{}&\ln{1\sdd}&\cnt{480}&\bounds{42}{48}{2}{false}&&\cr
\l{2}&\int(0){\0}{}&\int(2){\2}{}&\int(1){\2}{}&\ln{1\sdd}&\cnt{320}&\bounds{34}{40}{2}{false}&&\cr
\l{3}&\int(0){\0}{}&\int(3){\4}{}&\int(0){\0}{}&\ln{1}&\cnt{32}&\bounds{4}{4}{1}{false}&&\cr
\endorbit
\orbit{3}{6}{78}
\h{6}&\int(2){\2}{}&\int(4){\4}{}&\int(0){\0}{}&&\cnt{864}&\bounds{78}{84}{fail}{false}&&\cr
\l{1}&\int(1){1}{}&\int(2){1}{}&\int(0){\0}{}&\ln{1\sdd}&\cnt{256}&\bounds{20}{20}{1}{false}&&\cr
\l{2}&\int(1){\2}{}&\int(2){\2}{}&\int(0){\0}{}&\ln{1\sdd}&\cnt{320}&\bounds{34}{40}{2}{false}&&\cr
\l{3}&\int(1){2}{}&\int(2){1}{}&\int(0){1}{}&\ln{1\sdd}&\cnt{224}&\bounds{16}{16}{1}{false}&&\cr
\l{4}&\int(2){\2}{}&\int(1){\2}{}&\int(0){\0}{}&\ln{1}&\cnt{32}&\bounds{4}{4}{1}{false}&&\cr
\endorbit
\orbit{3}{7}{76}
\h{7}&\int(4){\4}{}&\int(2){\2}{}&\int(0){\0}{}&&\cnt{672}&\bounds{76}{88}{fail}{false}&&\cr
\l{1}&\int(2){1}{}&\int(1){1}{}&\int(0){\0}{}&\ln{1\sdd}&\cnt{384}&\bounds{36}{48}{2}{false}&&\cr
\l{2}&\int(2){\2}{}&\int(1){\2}{}&\int(0){\0}{}&\ln{1\sdd}&\cnt{192}&\bounds{24}{24}{1}{false}&&\cr
\l{3}&\int(3){\4}{}&\int(0){\0}{}&\int(0){\0}{}&\ln{1}&\cnt{48}&\bounds{8}{8}{1}{false}&&\cr
\endorbit
\orbit{3}{8}{76}
\h{8}&\rat(5/2){3}{}&\int(0){\0}{}&\rat(7/2){1}{}&&\cnt{672}&\bounds{76}{76}{1}{false}&&\cr
\l{1}&\rat(5/2){3}{}&\int(0){\0}{}&\rat(1/2){1}{}&\ln{1}&\cnt{21}&\bounds{3}{3}{1}{false}&&\cr
\l{3}&\rat(3/2){3}{}&\int(0){\0}{}&\rat(3/2){1}{}&\ln{1}&\cnt{315}&\bounds{35}{35}{1}{false}&&\cr
\endorbit
\orbit{3}{9}{76}
\h{9}&\int(1){2}{}&\rat(7/2){1}{}&\rat(3/2){1}{}&&\cnt{672}&\bounds{76}{76}{1}{false}&&\cr
\l{1}&\int(1){\2}{}&\int(2){\2}{}&\int(0){\0}{}&\ln{1}&\cnt{126}&\bounds{14}{14}{1}{false}&&\cr
\l{3}&\int(1){2}{}&\rat(3/2){1}{}&\rat(1/2){1}{}&\ln{1}&\cnt{189}&\bounds{21}{21}{1}{false}&&\cr
\l{5}&\int(1){2}{}&\rat(1/2){1}{}&\rat(3/2){1}{}&\ln{1}&\cnt{21}&\bounds{3}{3}{1}{false}&&\cr
\endorbit
\orbit{3}{10}{56}
\h{10}&\int(0){\0}{}&\int(6){\6}{\spec{1}}&\int(0){\0}{}&&\cnt{672}&\bounds{56}{56}{1}{false}&&\cr
\l{1}&\int(0){\2}{}&\int(3){\2}{}&\int(0){\0}{}&\ln{1\sdd}&\cnt{360}&\bounds{26}{26}{1}{false}&&\cr
\l{2}&\int(0){\0}{}&\int(3){\4}{}&\int(0){\0}{}&\ln{1\sdd}&\cnt{60}&\bounds{12}{12}{1}{false}&&\cr
\l{3}&\int(0){\0}{}&\int(3){\2}{}&\int(0){\2}{}&\ln{1\sdd}&\cnt{252}&\bounds{18}{18}{1}{false}&&\cr
\endorbit
\orbit{3}{11}{56}
\h{11}&\int(2){\2}{}&\int(2){\2}{}&\int(2){\2}{}&&\cnt{480}&\bounds{56}{72}{fail}{false}&&\cr
\l{1}&\int(1){\2}{}&\int(2){\2}{}&\int(0){\0}{}&\ln{2\sd}&\cnt{32}&\bounds{4}{4}{1}{false}&&\cr
\l{2}&\int(1){2}{}&\int(1){1}{}&\int(1){1}{}&\ln{1\sdd}&\cnt{288}&\bounds{32}{48}{2}{false}&&\cr
\l{3}&\int(2){\2}{}&\int(1){\2}{}&\int(0){\0}{}&\ln{2}&\cnt{32}&\bounds{4}{4}{1}{false}&&\cr
\endorbit
\orbit{3}{12}{56}
\h{12}&\int(3){2}{}&\rat(3/2){1}{}&\rat(3/2){1}{}&&\cnt{480}&\bounds{56}{56}{1}{false}&&\cr
\l{1}&\rat(3/2){3}{}&\int(0){\0}{}&\rat(3/2){1}{}&\ln{2\sd}&\cnt{64}&\bounds{8}{8}{1}{false}&&\cr
\l{2}&\int(3){\4}{}&\int(0){\0}{}&\int(0){\0}{}&\ln{1}&\cnt{14}&\bounds{2}{2}{1}{false}&&\cr
\l{4}&\int(2){\2}{}&\int(1){\2}{}&\int(0){\0}{}&\ln{2}&\cnt{81}&\bounds{9}{9}{1}{false}&&\cr
\endorbit
\orbit{3}{13}{54}
\h{13}&\int(6){\6}{\spec{1}}&\int(0){\0}{}&\int(0){\0}{}&&\cnt{672}&\bounds{54}{54}{1}{false}&&\cr
\l{1}&\int(3){\4}{}&\int(0){\0}{}&\int(0){\0}{}&\ln{1\sdd}&\cnt{168}&\bounds{18}{18}{1}{false}&&\cr
\l{2}&\int(3){\2}{}&\int(0){\2}{}&\int(0){\0}{}&\ln{2\sdd}&\cnt{252}&\bounds{18}{18}{1}{false}&&\cr
\endorbit
\orbit{3}{14}{40}
\h{14}&\rat(5/2){3}{}&\int(2){\2}{}&\rat(3/2){1}{}&&\cnt{384}&\bounds{40}{40}{1}{false}&&\cr
\l{1}&\rat(5/2){3}{}&\int(0){\0}{}&\rat(1/2){1}{}&\ln{1}&\cnt{27}&\bounds{3}{3}{1}{false}&&\cr
\l{3}&\int(2){1}{}&\int(1){1}{}&\int(0){\0}{}&\ln{1}&\cnt{120}&\bounds{12}{12}{1}{false}&&\cr
\l{5}&\rat(3/2){3}{}&\int(0){\0}{}&\rat(3/2){1}{}&\ln{1}&\cnt{45}&\bounds{5}{5}{1}{false}&&\cr
\endorbit
\endorbits

\newcs{max-A15_D9}{88}%
\newcs{count-A15_D9}{16}%
\newcs{counts-A15_D9}{\DO{1}{88}\DO{2}{88}\DO{3}{84}\DO{4}{78}\DO{5}{76}\DO{6}{72}\DO{7}{66}\DO{8}{64}\DO{9}{64}\DO{10}{60}\DO{11}{52}\DO{12}{52}\DO{13}{50}\DO{14}{48}\DO{15}{48}\DO{16}{40}}%
\orbits
\orbit{2}{1}{88}
\h{1}&\int(0){\0}{}&\int(6){\6}{}&&\cnt{1080}&\bounds{88}{88}{1}{false}&&\cr
\l{1}&\int(0){14}{}&\int(3){3}{}&\ln{2\sd}&\cnt{480}&\bounds{32}{32}{1}{false}&&\cr
\l{2}&\int(0){\0}{}&\int(3){\4}{}&\ln{1\sdd}&\cnt{120}&\bounds{24}{24}{1}{false}&&\cr
\endorbit
\orbit{2}{2}{88}
\h{2}&\int(2){\2}{}&\int(4){\4}{2}&&\cnt{1176}&\bounds{88}{88}{1}{false}&&\cr
\l{1}&\int(1){12}{}&\int(2){2}{}&\ln{2\sd}&\cnt{364}&\bounds{28}{28}{1}{false}&&\cr
\l{2}&\int(1){\2}{\spec{1}}&\int(2){\2}{}&\ln{2\sd}&\cnt{224}&\bounds{16}{16}{1}{false}&&\cr
\endorbit
\orbit{2}{3}{84}
\h{3}&\int(6){8}{\spec{1}}&\int(0){\0}{}&&\cnt{792}&\bounds{84}{88}{fail}{false}&&\cr
\l{1}&\int(3){12}{}&\int(0){2}{}&\ln{1\sdd}&\cnt{360}&\bounds{36}{40}{2}{false}&&\cr
\l{2}&\int(3){8}{}&\int(0){\0}{}&\ln{1}&\cnt{216}&\bounds{24}{24}{1}{false}&&\cr
\endorbit
\orbit{2}{4}{78}
\h{4}&\int(6){\6}{}&\int(0){\0}{}&&\cnt{792}&\bounds{78}{78}{1}{false}&&\cr
\l{1}&\int(3){8}{}&\int(0){\0}{}&\ln{1\sdd}&\cnt{252}&\bounds{24}{24}{1}{false}&&\cr
\l{2}&\int(3){\4}{\spec{1}}&\int(0){\0}{}&\ln{2\sd}&\cnt{90}&\bounds{9}{9}{1}{false}&&\cr
\l{3}&\int(3){12}{}&\int(0){2}{}&\ln{2\sd}&\cnt{180}&\bounds{18}{18}{1}{false}&&\cr
\endorbit
\orbit{2}{5}{76}
\h{5}&\int(5){12}{\spec{1}}&\int(1){2}{}&&\cnt{984}&\bounds{76}{76}{1}{false}&&\cr
\l{1}&\rat(5/2){14}{}&\rat(1/2){3}{}&\ln{1\sdd}&\cnt{256}&\bounds{20}{20}{1}{false}&&\cr
\l{2}&\int(3){\4}{}&\int(0){\0}{}&\ln{1}&\cnt{156}&\bounds{12}{12}{1}{false}&&\cr
\l{4}&\int(3){12}{}&\int(0){2}{}&\ln{1}&\cnt{208}&\bounds{16}{16}{1}{false}&&\cr
\endorbit
\orbit{2}{6}{72}
\h{6}&\rat(15/4){14}{\spec{1}}&\rat(9/4){3}{}&&\cnt{1080}&\bounds{72}{72}{1}{false}&&\cr
\l{1}&\int(2){\2}{}&\int(1){\2}{}&\ln{1}&\cnt{540}&\bounds{36}{36}{1}{false}&&\cr
\endorbit
\orbit{2}{7}{66}
\h{7}&\int(5){12}{\spec{2}}&\int(1){2}{}&&\cnt{600}&\bounds{66}{66}{1}{false}&&\cr
\l{1}&\int(3){8}{}&\int(0){\0}{}&\ln{1}&\cnt{120}&\bounds{13}{13}{1}{false}&&\cr
\l{3}&\int(3){\4}{}&\int(0){\0}{}&\ln{1}&\cnt{100}&\bounds{10}{10}{1}{false}&&\cr
\l{5}&\int(3){12}{}&\int(0){2}{}&\ln{1}&\cnt{80}&\bounds{10}{10}{1}{false}&&\cr
\endorbit
\orbit{2}{8}{64}
\h{8}&\int(4){\4}{}&\int(2){\2}{}&&\cnt{600}&\bounds{64}{64}{1}{false}&&\cr
\l{1}&\int(2){12}{}&\int(1){2}{}&\ln{2\sd}&\cnt{132}&\bounds{12}{12}{1}{false}&&\cr
\l{2}&\int(2){\2}{}&\int(1){\2}{}&\ln{1\sdd}&\cnt{112}&\bounds{16}{16}{1}{false}&&\cr
\l{3}&\int(2){14}{}&\int(1){3}{}&\ln{2\sd}&\cnt{64}&\bounds{8}{8}{1}{false}&&\cr
\l{4}&\int(3){\4}{\spec{1}}&\int(0){\0}{}&\ln{2}&\cnt{24}&\bounds{2}{2}{1}{false}&&\cr
\endorbit
\orbit{2}{9}{64}
\h{9}&\rat(7/4){14}{}&\rat(17/4){3}{}&&\cnt{744}&\bounds{64}{64}{1}{false}&&\cr
\l{1}&\rat(3/2){12}{}&\rat(3/2){2}{}&\ln{1}&\cnt{91}&\bounds{7}{7}{1}{false}&&\cr
\l{3}&\int(1){\2}{}&\int(2){\2}{}&\ln{1}&\cnt{224}&\bounds{16}{16}{1}{false}&&\cr
\l{5}&\rat(7/4){14}{}&\rat(5/4){3}{\spec{1}}&\ln{1}&\cnt{56}&\bounds{8}{8}{1}{false}&&\cr
\l{7}&\rat(7/4){14}{}&\rat(5/4){3}{\spec{2}}&\ln{1}&\cnt{1}&\bounds{1}{1}{1}{false}&&\cr
\endorbit
\orbit{2}{10}{60}
\h{10}&\int(2){\2}{}&\int(4){\4}{}&&\cnt{696}&\bounds{60}{60}{1}{false}&&\cr
\l{1}&\int(1){\2}{\spec{1}}&\int(2){\2}{}&\ln{2\sd}&\cnt{84}&\bounds{6}{6}{1}{false}&&\cr
\l{2}&\int(1){14}{}&\int(2){3}{}&\ln{2\sd}&\cnt{224}&\bounds{16}{16}{1}{false}&&\cr
\l{3}&\int(2){\2}{}&\int(1){\2}{}&\ln{1}&\cnt{40}&\bounds{8}{8}{1}{false}&&\cr
\endorbit
\orbit{2}{11}{52}
\h{11}&\int(3){12}{}&\int(3){2}{}&&\cnt{504}&\bounds{52}{52}{1}{false}&&\cr
\l{1}&\rat(3/2){14}{}&\rat(3/2){3}{}&\ln{1\sdd}&\cnt{192}&\bounds{24}{24}{1}{false}&&\cr
\l{2}&\int(3){12}{}&\int(0){2}{}&\ln{1}&\cnt{12}&\bounds{2}{2}{1}{false}&&\cr
\l{4}&\int(2){12}{}&\int(1){2}{}&\ln{1}&\cnt{144}&\bounds{12}{12}{1}{false}&&\cr
\endorbit
\orbit{2}{12}{52}
\h{12}&\int(4){8}{}&\int(2){\2}{}&&\cnt{408}&\bounds{52}{52}{1}{false}&&\cr
\l{1}&\int(2){12}{}&\int(1){2}{}&\ln{2\sdd}&\cnt{140}&\bounds{18}{18}{1}{false}&&\cr
\l{2}&\int(3){8}{}&\int(0){\0}{}&\ln{1}&\cnt{64}&\bounds{8}{8}{1}{false}&&\cr
\endorbit
\orbit{2}{13}{50}
\h{13}&\int(6){\6}{\spec{1}}&\int(0){\0}{}&&\cnt{600}&\bounds{50}{50}{1}{false}&&\cr
\l{1}&\int(3){\4}{}&\int(0){\0}{}&\ln{1\sdd}&\cnt{312}&\bounds{26}{26}{1}{false}&&\cr
\l{2}&\int(3){\2}{}&\int(0){\2}{}&\ln{1\sdd}&\cnt{288}&\bounds{24}{24}{1}{false}&&\cr
\endorbit
\orbit{2}{14}{48}
\h{14}&\int(0){\0}{}&\int(6){\6}{\spec{1}}&&\cnt{600}&\bounds{48}{48}{1}{false}&&\cr
\l{1}&\int(0){\2}{}&\int(3){\2}{}&\ln{1\sdd}&\cnt{480}&\bounds{32}{32}{1}{false}&&\cr
\l{2}&\int(0){\0}{}&\int(3){\4}{}&\ln{1\sdd}&\cnt{120}&\bounds{16}{16}{1}{false}&&\cr
\endorbit
\orbit{2}{15}{48}
\h{15}&\rat(15/4){14}{\spec{2}}&\rat(9/4){3}{}&&\cnt{504}&\bounds{48}{48}{1}{false}&&\cr
\l{1}&\int(3){\4}{}&\int(0){\0}{}&\ln{1}&\cnt{36}&\bounds{3}{3}{1}{false}&&\cr
\l{3}&\rat(5/2){12}{}&\rat(1/2){2}{}&\ln{1}&\cnt{108}&\bounds{9}{9}{1}{false}&&\cr
\l{5}&\int(2){\2}{}&\int(1){\2}{}&\ln{1}&\cnt{108}&\bounds{12}{12}{1}{false}&&\cr
\endorbit
\orbit{2}{16}{40}
\h{16}&\rat(15/4){10}{}&\rat(9/4){1}{}&&\cnt{360}&\bounds{40}{40}{1}{false}&&\cr
\l{1}&\int(3){8}{}&\int(0){\0}{}&\ln{1}&\cnt{45}&\bounds{5}{5}{1}{false}&&\cr
\l{3}&\rat(5/2){12}{}&\rat(1/2){2}{}&\ln{1}&\cnt{135}&\bounds{15}{15}{1}{false}&&\cr
\endorbit
\endorbits

\newcs{max-3D8}{104}%
\newcs{count-3D8}{10}%
\newcs{counts-3D8}{\DO{1}{104}\DO{2}{104}\DO{3}{104}\DO{4}{98}\DO{5}{96}\DO{6}{72}\DO{7}{72}\DO{8}{56}\DO{9}{52}\DO{10}{48}}%
\orbits
\orbit{3}{1}{104}
\h{1}&\int(4){\4}{}&\int(2){\2}{}&\int(0){\0}{}&&\cnt{720}&\bounds{104}{104}{1}{false}&&\cr
\l{1}&\int(2){3}{}&\int(1){3}{}&\int(0){\0}{}&\ln{1\sdd}&\cnt{256}&\bounds{32}{32}{1}{false}&&\cr
\l{2}&\int(2){1}{}&\int(1){2}{}&\int(0){2}{}&\ln{1\sdd}&\cnt{256}&\bounds{32}{32}{1}{false}&&\cr
\l{3}&\int(2){\2}{}&\int(1){\2}{}&\int(0){\0}{}&\ln{1\sdd}&\cnt{144}&\bounds{24}{24}{1}{false}&&\cr
\l{4}&\int(3){\4}{}&\int(0){\0}{}&\int(0){\0}{}&\ln{1}&\cnt{32}&\bounds{8}{8}{1}{false}&&\cr
\endorbit
\orbit{3}{2}{104}
\h{2}&\int(4){3}{}&\int(2){3}{}&\int(0){\0}{}&&\cnt{720}&\bounds{104}{104}{1}{false}&&\cr
\l{1}&\int(3){\4}{}&\int(0){\0}{}&\int(0){\0}{}&\ln{1}&\cnt{35}&\bounds{7}{7}{1}{false}&&\cr
\l{3}&\int(3){\4}{2}&\int(0){\0}{}&\int(0){\0}{}&\ln{1}&\cnt{1}&\bounds{1}{1}{1}{false}&&\cr
\l{5}&\rat(5/2){1}{}&\rat(1/2){2}{}&\int(0){2}{}&\ln{1}&\cnt{128}&\bounds{16}{16}{1}{false}&&\cr
\l{7}&\int(2){3}{}&\int(1){3}{}&\int(0){\0}{}&\ln{1}&\cnt{196}&\bounds{28}{28}{1}{false}&&\cr
\endorbit
\orbit{3}{3}{104}
\h{3}&\int(1){2}{}&\int(1){2}{}&\int(4){1}{}&&\cnt{720}&\bounds{104}{104}{1}{false}&&\cr
\l{1}&\int(1){\2}{}&\int(0){\0}{}&\int(2){\2}{}&\ln{2}&\cnt{98}&\bounds{14}{14}{1}{false}&&\cr
\l{3}&\rat(1/2){3}{}&\int(0){\0}{}&\rat(5/2){3}{}&\ln{2}&\cnt{64}&\bounds{8}{8}{1}{false}&&\cr
\l{5}&\int(1){2}{}&\int(1){2}{}&\int(1){1}{\spec{1}}&\ln{1}&\cnt{35}&\bounds{7}{7}{1}{false}&&\cr
\l{7}&\int(1){2}{}&\int(1){2}{}&\int(1){1}{\spec{2}}&\ln{1}&\cnt{1}&\bounds{1}{1}{1}{false}&&\cr
\endorbit
\orbit{3}{4}{98}
\h{4}&\int(6){\6}{}&\int(0){\0}{}&\int(0){\0}{}&&\cnt{1104}&\bounds{98}{120}{fail}{false}&&\cr
\l{1}&\int(3){\4}{}&\int(0){\0}{}&\int(0){\0}{}&\ln{1\sdd}&\cnt{80}&\bounds{16}{16}{1}{false}&&\cr
\l{2}&\int(3){3}{}&\int(0){3}{}&\int(0){\0}{}&\ln{2\sdd}&\cnt{256}&\bounds{20}{20}{1}{false}&&\cr
\l{3}&\int(3){1}{}&\int(0){2}{}&\int(0){2}{}&\ln{1\sdd}&\cnt{512}&\bounds{42}{64}{2}{false}&&\cr
\endorbit
\orbit{3}{5}{96}
\h{5}&\int(4){\4}{2}&\int(2){\2}{}&\int(0){\0}{}&&\cnt{1104}&\bounds{96}{132}{fail}{false}&&\cr
\l{1}&\int(2){\2}{}&\int(1){\2}{}&\int(0){\0}{}&\ln{1\sdd}&\cnt{336}&\bounds{34}{48}{2}{false}&&\cr
\l{2}&\int(2){2}{}&\int(1){1}{}&\int(0){2}{}&\ln{1\sdd}&\cnt{512}&\bounds{42}{64}{2}{false}&&\cr
\l{3}&\int(2){2}{}&\int(1){2}{}&\int(0){1}{}&\ln{1\sdd}&\cnt{256}&\bounds{20}{20}{1}{false}&&\cr
\endorbit
\orbit{3}{6}{72}
\h{6}&\int(2){\2}{}&\int(2){\2}{}&\int(2){\2}{}&&\cnt{528}&\bounds{72}{72}{1}{false}&&\cr
\l{1}&\int(2){\2}{}&\int(1){\2}{}&\int(0){\0}{}&\ln{6\sd}&\cnt{24}&\bounds{4}{4}{1}{false}&&\cr
\l{2}&\int(1){1}{}&\int(1){2}{}&\int(1){2}{}&\ln{3\sdd}&\cnt{128}&\bounds{16}{16}{1}{false}&&\cr
\endorbit
\orbit{3}{7}{72}
\h{7}&\int(3){2}{}&\int(1){2}{}&\int(2){1}{}&&\cnt{528}&\bounds{72}{72}{1}{false}&&\cr
\l{1}&\int(3){\4}{}&\int(0){\0}{}&\int(0){\0}{}&\ln{1}&\cnt{10}&\bounds{2}{2}{1}{false}&&\cr
\l{3}&\int(2){\2}{}&\int(1){\2}{}&\int(0){\0}{}&\ln{1}&\cnt{42}&\bounds{6}{6}{1}{false}&&\cr
\l{5}&\int(2){\2}{}&\int(0){\0}{}&\int(1){\2}{}&\ln{1}&\cnt{84}&\bounds{12}{12}{1}{false}&&\cr
\l{7}&\rat(3/2){3}{}&\int(0){\0}{}&\rat(3/2){3}{}&\ln{1}&\cnt{128}&\bounds{16}{16}{1}{false}&&\cr
\endorbit
\orbit{3}{8}{56}
\h{8}&\int(2){3}{}&\int(2){3}{}&\int(2){\2}{}&&\cnt{432}&\bounds{56}{56}{1}{false}&&\cr
\l{1}&\int(2){3}{}&\int(0){\0}{}&\int(1){3}{}&\ln{2\sd}&\cnt{32}&\bounds{4}{4}{1}{false}&&\cr
\l{2}&\rat(3/2){1}{}&\rat(1/2){2}{}&\int(1){2}{}&\ln{2\sd}&\cnt{128}&\bounds{16}{16}{1}{false}&&\cr
\l{3}&\int(2){3}{}&\int(1){3}{}&\int(0){\0}{}&\ln{2}&\cnt{28}&\bounds{4}{4}{1}{false}&&\cr
\endorbit
\orbit{3}{9}{52}
\h{9}&\int(6){\6}{\spec{1}}&\int(0){\0}{}&\int(0){\0}{}&&\cnt{528}&\bounds{52}{52}{1}{false}&&\cr
\l{1}&\int(3){\4}{}&\int(0){\0}{}&\int(0){\0}{}&\ln{1\sdd}&\cnt{80}&\bounds{12}{12}{1}{false}&&\cr
\l{2}&\int(3){\2}{}&\int(0){\2}{}&\int(0){\0}{}&\ln{2\sdd}&\cnt{224}&\bounds{20}{20}{1}{false}&&\cr
\endorbit
\orbit{3}{10}{48}
\h{10}&\int(2){1}{}&\int(2){1}{}&\int(2){1}{}&&\cnt{384}&\bounds{48}{48}{1}{false}&&\cr
\l{1}&\int(2){1}{}&\rat(1/2){2}{}&\rat(1/2){2}{}&\ln{3}&\cnt{64}&\bounds{8}{8}{1}{false}&&\cr
\endorbit
\endorbits

\newcs{max-2A12}{88}%
\newcs{count-2A12}{11}%
\newcs{counts-2A12}{\DO{1}{88}\DO{2}{84}\DO{3}{78}\DO{4}{72}\DO{5}{72}\DO{6}{68}\DO{7}{62}\DO{8}{58}\DO{9}{58}\DO{10}{56}\DO{11}{46}}%
\orbits
\orbit{2}{1}{88}
\h{1}&\rat(22/13){11}{}&\rat(56/13){3}{\spec{2}}&&\cnt{924}&\bounds{88}{88}{1}{false}&&\cr
\l{1}&\rat(16/13){8}{}&\rat(23/13){1}{}&\ln{1}&\cnt{165}&\bounds{17}{17}{1}{false}&&\cr
\l{3}&\int(1){\2}{}&\int(2){\2}{}&\ln{1}&\cnt{242}&\bounds{22}{22}{1}{false}&&\cr
\l{5}&\rat(22/13){11}{}&\rat(17/13){3}{}&\ln{1}&\cnt{55}&\bounds{5}{5}{1}{false}&&\cr
\endorbit
\orbit{2}{2}{84}
\h{2}&\rat(48/13){11}{\spec{1}}&\rat(30/13){3}{}&&\cnt{972}&\bounds{84}{84}{1}{false}&&\cr
\l{1}&\rat(24/13){12}{}&\rat(15/13){8}{}&\ln{1\sdd}&\cnt{252}&\bounds{24}{24}{1}{false}&&\cr
\l{2}&\int(2){\2}{}&\int(1){\2}{}&\ln{1}&\cnt{360}&\bounds{30}{30}{1}{false}&&\cr
\endorbit
\orbit{2}{3}{78}
\h{3}&\rat(12/13){12}{}&\rat(66/13){8}{\spec{1}}&&\cnt{828}&\bounds{78}{78}{1}{false}&&\cr
\l{1}&\rat(10/13){10}{}&\rat(29/13){11}{}&\ln{1}&\cnt{198}&\bounds{18}{18}{1}{false}&&\cr
\l{3}&\int(1){\2}{}&\int(2){\2}{}&\ln{1}&\cnt{108}&\bounds{9}{9}{1}{false}&&\cr
\l{5}&\rat(12/13){12}{}&\rat(27/13){8}{}&\ln{1}&\cnt{108}&\bounds{12}{12}{1}{false}&&\cr
\endorbit
\orbit{2}{4}{72}
\h{4}&\int(6){\6}{}&\int(0){\0}{}&&\cnt{828}&\bounds{72}{72}{1}{false}&&\cr
\l{1}&\int(3){\4}{\spec{1}}&\int(0){\0}{}&\ln{2\sd}&\cnt{63}&\bounds{9}{9}{1}{false}&&\cr
\l{2}&\int(3){8}{}&\int(0){1}{}&\ln{2\sd}&\cnt{273}&\bounds{21}{21}{1}{false}&&\cr
\l{3}&\int(3){10}{}&\int(0){11}{}&\ln{2\sd}&\cnt{78}&\bounds{6}{6}{1}{false}&&\cr
\endorbit
\orbit{2}{5}{72}
\h{5}&\rat(12/13){12}{}&\rat(66/13){8}{\spec{2}}&&\cnt{684}&\bounds{72}{72}{1}{false}&&\cr
\l{1}&\int(1){\2}{}&\int(2){\2}{}&\ln{1}&\cnt{72}&\bounds{6}{6}{1}{false}&&\cr
\l{3}&\rat(11/13){11}{}&\rat(28/13){3}{}&\ln{1}&\cnt{180}&\bounds{15}{15}{1}{false}&&\cr
\l{5}&\rat(12/13){12}{}&\rat(27/13){8}{}&\ln{1}&\cnt{90}&\bounds{15}{15}{1}{false}&&\cr
\endorbit
\orbit{2}{6}{68}
\h{6}&\int(4){\4}{}&\int(2){\2}{}&&\cnt{636}&\bounds{68}{68}{1}{false}&&\cr
\l{1}&\int(2){8}{}&\int(1){1}{}&\ln{2\sd}&\cnt{84}&\bounds{12}{12}{1}{false}&&\cr
\l{2}&\int(2){10}{}&\int(1){11}{}&\ln{2\sd}&\cnt{99}&\bounds{9}{9}{1}{false}&&\cr
\l{3}&\int(2){\2}{}&\int(1){\2}{\spec{1}}&\ln{2\sd}&\cnt{44}&\bounds{4}{4}{1}{false}&&\cr
\l{4}&\int(2){11}{}&\int(1){3}{}&\ln{2\sd}&\cnt{55}&\bounds{5}{5}{1}{false}&&\cr
\l{5}&\int(3){\4}{\spec{1}}&\int(0){\0}{}&\ln{2}&\cnt{18}&\bounds{2}{2}{1}{false}&&\cr
\endorbit
\orbit{2}{7}{62}
\h{7}&\rat(22/13){11}{}&\rat(56/13){3}{\spec{1}}&&\cnt{588}&\bounds{62}{62}{1}{false}&&\cr
\l{1}&\rat(11/13){12}{}&\rat(28/13){8}{}&\ln{1\sdd}&\cnt{140}&\bounds{18}{18}{1}{false}&&\cr
\l{2}&\rat(20/13){10}{}&\rat(19/13){11}{}&\ln{1}&\cnt{88}&\bounds{8}{8}{1}{false}&&\cr
\l{4}&\int(1){\2}{}&\int(2){\2}{}&\ln{1}&\cnt{88}&\bounds{8}{8}{1}{false}&&\cr
\l{6}&\rat(22/13){11}{}&\rat(17/13){3}{}&\ln{1}&\cnt{48}&\bounds{6}{6}{1}{false}&&\cr
\endorbit
\orbit{2}{8}{58}
\h{8}&\rat(38/13){12}{\spec{2}}&\rat(40/13){8}{}&&\cnt{492}&\bounds{58}{58}{1}{false}&&\cr
\l{1}&\int(3){\4}{}&\int(0){\0}{}&\ln{1}&\cnt{10}&\bounds{1}{1}{1}{false}&&\cr
\l{3}&\rat(23/13){10}{}&\rat(16/13){11}{}&\ln{1}&\cnt{100}&\bounds{10}{10}{1}{false}&&\cr
\l{5}&\int(2){\2}{}&\int(1){\2}{}&\ln{1}&\cnt{80}&\bounds{10}{10}{1}{false}&&\cr
\l{7}&\rat(24/13){11}{}&\rat(15/13){3}{}&\ln{1}&\cnt{56}&\bounds{8}{8}{1}{false}&&\cr
\endorbit
\orbit{2}{9}{58}
\h{9}&\rat(48/13){11}{\spec{2}}&\rat(30/13){3}{}&&\cnt{540}&\bounds{58}{58}{1}{false}&&\cr
\l{1}&\int(3){\4}{}&\int(0){\0}{}&\ln{1}&\cnt{27}&\bounds{3}{3}{1}{false}&&\cr
\l{3}&\rat(29/13){8}{}&\rat(10/13){1}{}&\ln{1}&\cnt{108}&\bounds{12}{12}{1}{false}&&\cr
\l{5}&\rat(33/13){10}{}&\rat(6/13){11}{}&\ln{1}&\cnt{45}&\bounds{5}{5}{1}{false}&&\cr
\l{7}&\int(2){\2}{}&\int(1){\2}{}&\ln{1}&\cnt{90}&\bounds{9}{9}{1}{false}&&\cr
\endorbit
\orbit{2}{10}{56}
\h{10}&\rat(36/13){9}{}&\rat(42/13){6}{}&&\cnt{396}&\bounds{56}{60}{fail}{false}&&\cr
\l{1}&\rat(18/13){11}{}&\rat(21/13){3}{}&\ln{1\sdd}&\cnt{120}&\bounds{20}{24}{2}{false}&&\cr
\l{2}&\rat(32/13){8}{}&\rat(7/13){1}{}&\ln{1}&\cnt{54}&\bounds{6}{6}{1}{false}&&\cr
\l{4}&\rat(27/13){10}{}&\rat(12/13){11}{}&\ln{1}&\cnt{84}&\bounds{12}{12}{1}{false}&&\cr
\endorbit
\orbit{2}{11}{46}
\h{11}&\int(6){\6}{\spec{1}}&\int(0){\0}{}&&\cnt{492}&\bounds{46}{46}{1}{false}&&\cr
\l{1}&\int(3){\4}{}&\int(0){\0}{}&\ln{1\sdd}&\cnt{180}&\bounds{20}{20}{1}{false}&&\cr
\l{2}&\int(3){\2}{}&\int(0){\2}{}&\ln{1\sdd}&\cnt{312}&\bounds{26}{26}{1}{false}&&\cr
\endorbit
\endorbits

\newcs{max-A11_D7_E6}{110}%
\newcs{count-A11_D7_E6}{33}%
\newcs{counts-A11_D7_E6}{\DO{1}{110}\DO{2}{102}\DO{3}{102}\DO{4}{100}\DO{5}{96}\DO{6}{96}\DO{7}{96}\DO{8}{94}\DO{9}{94}\DO{10}{90}\DO{11}{90}\DO{12}{88}\DO{13}{88}\DO{14}{88}\DO{15}{88}\DO{16}{86}\DO{17}{86}\DO{18}{84}\DO{19}{80}\DO{20}{80}\DO{21}{76}\DO{22}{76}\DO{23}{74}\DO{24}{72}\DO{25}{68}\DO{26}{68}\DO{27}{68}\DO{28}{64}\DO{29}{64}\DO{30}{60}\DO{31}{54}\DO{32}{54}\DO{33}{52}}%
\orbits
\orbit{3}{1}{110}
\h{1}&\int(0){\0}{}&\int(6){\6}{}&\int(0){\0}{}&&\cnt{1128}&\bounds{110}{110}{1}{false}&&\cr
\l{1}&\int(0){9}{}&\int(3){1}{}&\int(0){\0}{}&\ln{2\sd}&\cnt{220}&\bounds{20}{20}{1}{false}&&\cr
\l{2}&\int(0){11}{}&\int(3){3}{}&\int(0){2}{}&\ln{2\sd}&\cnt{324}&\bounds{27}{27}{1}{false}&&\cr
\l{3}&\int(0){\0}{}&\int(3){\4}{\spec{1}}&\int(0){\0}{}&\ln{2\sd}&\cnt{20}&\bounds{8}{8}{1}{false}&&\cr
\endorbit
\orbit{3}{2}{102}
\h{2}&\int(2){\2}{}&\int(4){\4}{2}&\int(0){\0}{}&&\cnt{1032}&\bounds{102}{102}{1}{false}&&\cr
\l{1}&\int(1){6}{}&\int(2){2}{}&\int(0){\0}{}&\ln{1\sdd}&\cnt{252}&\bounds{24}{24}{1}{false}&&\cr
\l{2}&\int(1){\2}{\spec{1}}&\int(2){\2}{}&\int(0){\0}{}&\ln{2\sd}&\cnt{120}&\bounds{12}{12}{1}{false}&&\cr
\l{3}&\int(1){10}{}&\int(2){2}{}&\int(0){1}{}&\ln{2\sd}&\cnt{270}&\bounds{27}{27}{1}{false}&&\cr
\endorbit
\orbit{3}{3}{102}
\h{3}&\int(6){\6}{}&\int(0){\0}{}&\int(0){\0}{}&&\cnt{840}&\bounds{102}{110}{fail}{false}&&\cr
\l{1}&\int(3){\4}{\spec{1}}&\int(0){\0}{}&\int(0){\0}{}&\ln{2\sd}&\cnt{54}&\bounds{9}{9}{1}{false}&&\cr
\l{2}&\int(3){6}{}&\int(0){2}{}&\int(0){\0}{}&\ln{1\sdd}&\cnt{280}&\bounds{32}{40}{2}{false}&&\cr
\l{3}&\int(3){8}{}&\int(0){\0}{}&\int(0){2}{}&\ln{2\sd}&\cnt{162}&\bounds{18}{18}{1}{false}&&\cr
\l{4}&\int(3){9}{}&\int(0){1}{}&\int(0){\0}{}&\ln{2\sd}&\cnt{64}&\bounds{8}{8}{1}{false}&&\cr
\endorbit
\orbit{3}{4}{100}
\h{4}&\int(0){\0}{}&\int(4){\4}{2}&\int(2){\2}{}&&\cnt{1032}&\bounds{100}{112}{fail}{false}&&\cr
\l{1}&\int(0){10}{}&\int(2){2}{}&\int(1){1}{}&\ln{2\sd}&\cnt{396}&\bounds{36}{36}{1}{false}&&\cr
\l{2}&\int(0){\0}{}&\int(2){\2}{}&\int(1){\2}{}&\ln{1\sdd}&\cnt{240}&\bounds{28}{40}{2}{false}&&\cr
\endorbit
\orbit{3}{5}{96}
\h{5}&\int(5){6}{\spec{1}}&\int(1){2}{}&\int(0){\0}{}&&\cnt{744}&\bounds{96}{96}{1}{false}&&\cr
\l{1}&\rat(5/2){9}{}&\rat(1/2){1}{}&\int(0){\0}{}&\ln{1\sdd}&\cnt{192}&\bounds{24}{24}{1}{false}&&\cr
\l{2}&\int(3){\4}{}&\int(0){\0}{}&\int(0){\0}{}&\ln{1}&\cnt{84}&\bounds{12}{12}{1}{false}&&\cr
\l{4}&\int(3){6}{}&\int(0){2}{}&\int(0){\0}{}&\ln{1}&\cnt{84}&\bounds{12}{12}{1}{false}&&\cr
\l{6}&\int(3){8}{}&\int(0){\0}{}&\int(0){2}{}&\ln{1}&\cnt{108}&\bounds{12}{12}{1}{false}&&\cr
\endorbit
\orbit{3}{6}{96}
\h{6}&\int(4){\4}{}&\int(0){\0}{}&\int(2){\2}{}&&\cnt{648}&\bounds{96}{96}{1}{false}&&\cr
\l{1}&\int(2){8}{}&\int(0){\0}{}&\int(1){2}{}&\ln{2\sd}&\cnt{168}&\bounds{24}{24}{1}{false}&&\cr
\l{2}&\int(2){\2}{}&\int(0){\0}{}&\int(1){\2}{}&\ln{1\sdd}&\cnt{80}&\bounds{16}{16}{1}{false}&&\cr
\l{3}&\int(2){10}{}&\int(0){2}{}&\int(1){1}{}&\ln{2\sd}&\cnt{84}&\bounds{12}{12}{1}{false}&&\cr
\l{4}&\int(3){\4}{\spec{1}}&\int(0){\0}{}&\int(0){\0}{}&\ln{2}&\cnt{16}&\bounds{2}{2}{1}{false}&&\cr
\endorbit
\orbit{3}{7}{96}
\h{7}&\rat(14/3){8}{\spec{1}}&\int(0){\0}{}&\rat(4/3){2}{}&&\cnt{840}&\bounds{96}{104}{fail}{false}&&\cr
\l{1}&\rat(7/3){10}{}&\int(0){2}{}&\rat(2/3){1}{}&\ln{1\sdd}&\cnt{280}&\bounds{32}{40}{2}{false}&&\cr
\l{2}&\int(3){\4}{}&\int(0){\0}{}&\int(0){\0}{}&\ln{1}&\cnt{72}&\bounds{8}{8}{1}{false}&&\cr
\l{4}&\rat(8/3){8}{}&\int(0){\0}{}&\rat(1/3){2}{}&\ln{1}&\cnt{144}&\bounds{16}{16}{1}{false}&&\cr
\l{6}&\int(3){9}{}&\int(0){1}{}&\int(0){\0}{}&\ln{1}&\cnt{64}&\bounds{8}{8}{1}{false}&&\cr
\endorbit
\orbit{3}{8}{94}
\h{8}&\int(2){\2}{}&\int(0){\0}{}&\int(4){\4}{}&&\cnt{840}&\bounds{94}{94}{1}{false}&&\cr
\l{1}&\int(1){8}{}&\int(0){\0}{}&\int(2){2}{}&\ln{2\sd}&\cnt{120}&\bounds{13}{13}{1}{false}&&\cr
\l{2}&\int(1){\2}{\spec{1}}&\int(0){\0}{}&\int(2){\2}{}&\ln{2\sd}&\cnt{80}&\bounds{8}{8}{1}{false}&&\cr
\l{3}&\int(1){10}{}&\int(0){2}{}&\int(2){1}{}&\ln{2\sd}&\cnt{140}&\bounds{14}{14}{1}{false}&&\cr
\l{4}&\int(1){11}{}&\int(0){3}{}&\int(2){2}{}&\ln{2\sd}&\cnt{64}&\bounds{8}{8}{1}{false}&&\cr
\l{5}&\int(2){\2}{}&\int(0){\0}{}&\int(1){\2}{\spec{1}}&\ln{2}&\cnt{8}&\bounds{2}{2}{1}{false}&&\cr
\endorbit
\orbit{3}{9}{94}
\h{9}&\rat(17/4){9}{\spec{1}}&\rat(7/4){1}{}&\int(0){\0}{}&&\cnt{888}&\bounds{94}{94}{1}{false}&&\cr
\l{1}&\int(3){\4}{}&\int(0){\0}{}&\int(0){\0}{}&\ln{1}&\cnt{45}&\bounds{5}{5}{1}{false}&&\cr
\l{3}&\rat(9/4){9}{}&\rat(3/4){1}{}&\int(0){\0}{}&\ln{1}&\cnt{210}&\bounds{21}{21}{1}{false}&&\cr
\l{5}&\rat(5/2){10}{}&\rat(1/2){2}{}&\int(0){1}{}&\ln{1}&\cnt{189}&\bounds{21}{21}{1}{false}&&\cr
\endorbit
\orbit{3}{10}{90}
\h{10}&\rat(14/3){8}{\spec{2}}&\int(0){\0}{}&\rat(4/3){2}{}&&\cnt{648}&\bounds{90}{104}{fail}{false}&&\cr
\l{1}&\rat(7/3){4}{}&\int(0){\0}{}&\rat(2/3){1}{}&\ln{1\sdd}&\cnt{200}&\bounds{26}{40}{2}{false}&&\cr
\l{2}&\int(3){\4}{}&\int(0){\0}{}&\int(0){\0}{}&\ln{1}&\cnt{60}&\bounds{10}{10}{1}{false}&&\cr
\l{4}&\int(3){6}{}&\int(0){2}{}&\int(0){\0}{}&\ln{1}&\cnt{84}&\bounds{12}{12}{1}{false}&&\cr
\l{6}&\rat(8/3){8}{}&\int(0){\0}{}&\rat(1/3){2}{}&\ln{1}&\cnt{80}&\bounds{10}{10}{1}{false}&&\cr
\endorbit
\orbit{3}{11}{90}
\h{11}&\rat(11/3){10}{\spec{1}}&\int(1){2}{}&\rat(4/3){1}{}&&\cnt{936}&\bounds{90}{96}{fail}{false}&&\cr
\l{1}&\rat(11/6){11}{}&\rat(1/2){3}{}&\rat(2/3){2}{}&\ln{1\sdd}&\cnt{320}&\bounds{34}{40}{2}{false}&&\cr
\l{2}&\int(2){\2}{}&\int(1){\2}{}&\int(0){\0}{}&\ln{1}&\cnt{132}&\bounds{12}{12}{1}{false}&&\cr
\l{4}&\int(2){\2}{}&\int(0){\0}{}&\int(1){\2}{}&\ln{1}&\cnt{176}&\bounds{16}{16}{1}{false}&&\cr
\endorbit
\orbit{3}{12}{88}
\h{12}&\int(0){\0}{}&\int(2){\2}{}&\int(4){\4}{}&&\cnt{840}&\bounds{88}{96}{fail}{false}&&\cr
\l{1}&\int(0){10}{}&\int(1){2}{}&\int(2){1}{}&\ln{2\sd}&\cnt{132}&\bounds{12}{12}{1}{false}&&\cr
\l{2}&\int(0){11}{}&\int(1){3}{}&\int(2){2}{}&\ln{2\sd}&\cnt{192}&\bounds{16}{16}{1}{false}&&\cr
\l{3}&\int(0){\0}{}&\int(1){\2}{}&\int(2){\2}{}&\ln{1\sdd}&\cnt{160}&\bounds{24}{32}{2}{false}&&\cr
\l{4}&\int(0){\0}{}&\int(2){\2}{}&\int(1){\2}{\spec{1}}&\ln{2}&\cnt{8}&\bounds{2}{2}{1}{false}&&\cr
\endorbit
\orbit{3}{13}{88}
\h{13}&\int(3){6}{}&\int(3){2}{}&\int(0){\0}{}&&\cnt{552}&\bounds{88}{104}{fail}{false}&&\cr
\l{1}&\rat(3/2){9}{}&\rat(3/2){1}{}&\int(0){\0}{}&\ln{2\sdd}&\cnt{160}&\bounds{24}{32}{2}{false}&&\cr
\l{2}&\int(3){6}{}&\int(0){2}{}&\int(0){\0}{}&\ln{1}&\cnt{8}&\bounds{2}{2}{1}{false}&&\cr
\l{4}&\int(2){6}{}&\int(1){2}{}&\int(0){\0}{}&\ln{1}&\cnt{108}&\bounds{18}{18}{1}{false}&&\cr
\endorbit
\orbit{3}{14}{88}
\h{14}&\int(2){\2}{}&\int(4){\4}{}&\int(0){\0}{}&&\cnt{744}&\bounds{88}{88}{1}{false}&&\cr
\l{1}&\int(1){9}{}&\int(2){1}{}&\int(0){\0}{}&\ln{2\sd}&\cnt{180}&\bounds{20}{20}{1}{false}&&\cr
\l{2}&\int(1){\2}{\spec{1}}&\int(2){\2}{}&\int(0){\0}{}&\ln{2\sd}&\cnt{60}&\bounds{6}{6}{1}{false}&&\cr
\l{3}&\int(1){11}{}&\int(2){3}{}&\int(0){2}{}&\ln{2\sd}&\cnt{108}&\bounds{12}{12}{1}{false}&&\cr
\l{4}&\int(2){\2}{}&\int(1){\2}{}&\int(0){\0}{}&\ln{1}&\cnt{24}&\bounds{6}{6}{1}{false}&&\cr
\endorbit
\orbit{3}{15}{88}
\h{15}&\rat(9/4){9}{}&\rat(15/4){1}{}&\int(0){\0}{}&&\cnt{696}&\bounds{88}{88}{1}{false}&&\cr
\l{1}&\rat(3/2){6}{}&\rat(3/2){2}{}&\int(0){\0}{}&\ln{1}&\cnt{84}&\bounds{12}{12}{1}{false}&&\cr
\l{3}&\rat(9/4){9}{}&\rat(3/4){1}{\spec{1}}&\int(0){\0}{}&\ln{1}&\cnt{1}&\bounds{1}{1}{1}{false}&&\cr
\l{5}&\rat(9/4){9}{}&\rat(3/4){1}{\spec{2}}&\int(0){\0}{}&\ln{1}&\cnt{20}&\bounds{4}{4}{1}{false}&&\cr
\l{7}&\rat(5/4){9}{}&\rat(7/4){1}{}&\int(0){\0}{}&\ln{1}&\cnt{162}&\bounds{18}{18}{1}{false}&&\cr
\l{9}&\rat(3/2){10}{}&\rat(3/2){2}{}&\int(0){1}{}&\ln{1}&\cnt{81}&\bounds{9}{9}{1}{false}&&\cr
\endorbit
\orbit{3}{16}{86}
\h{16}&\int(4){\4}{}&\int(2){\2}{}&\int(0){\0}{}&&\cnt{648}&\bounds{86}{86}{1}{false}&&\cr
\l{1}&\int(2){6}{}&\int(1){2}{}&\int(0){\0}{}&\ln{1\sdd}&\cnt{140}&\bounds{18}{18}{1}{false}&&\cr
\l{2}&\int(2){9}{}&\int(1){1}{}&\int(0){\0}{}&\ln{2\sd}&\cnt{128}&\bounds{16}{16}{1}{false}&&\cr
\l{3}&\int(2){\2}{}&\int(1){\2}{}&\int(0){\0}{}&\ln{1\sdd}&\cnt{80}&\bounds{16}{16}{1}{false}&&\cr
\l{4}&\int(2){10}{}&\int(1){2}{}&\int(0){1}{}&\ln{2\sd}&\cnt{54}&\bounds{6}{6}{1}{false}&&\cr
\l{5}&\int(3){\4}{\spec{1}}&\int(0){\0}{}&\int(0){\0}{}&\ln{2}&\cnt{16}&\bounds{2}{2}{1}{false}&&\cr
\endorbit
\orbit{3}{17}{86}
\h{17}&\rat(8/3){8}{}&\int(0){\0}{}&\rat(10/3){2}{}&&\cnt{648}&\bounds{86}{86}{1}{false}&&\cr
\l{1}&\rat(4/3){4}{}&\int(0){\0}{}&\rat(5/3){1}{}&\ln{1\sdd}&\cnt{140}&\bounds{18}{18}{1}{false}&&\cr
\l{2}&\rat(4/3){10}{}&\int(0){2}{}&\rat(5/3){1}{}&\ln{1\sdd}&\cnt{168}&\bounds{24}{24}{1}{false}&&\cr
\l{3}&\rat(8/3){8}{}&\int(0){\0}{}&\rat(1/3){2}{}&\ln{1}&\cnt{10}&\bounds{2}{2}{1}{false}&&\cr
\l{5}&\rat(5/3){8}{}&\int(0){\0}{}&\rat(4/3){2}{}&\ln{1}&\cnt{160}&\bounds{20}{20}{1}{false}&&\cr
\endorbit
\orbit{3}{18}{84}
\h{18}&\rat(17/4){9}{\spec{2}}&\rat(7/4){1}{}&\int(0){\0}{}&&\cnt{600}&\bounds{84}{84}{1}{false}&&\cr
\l{1}&\int(3){\4}{}&\int(0){\0}{}&\int(0){\0}{}&\ln{1}&\cnt{42}&\bounds{6}{6}{1}{false}&&\cr
\l{3}&\rat(5/2){6}{}&\rat(1/2){2}{}&\int(0){\0}{}&\ln{1}&\cnt{147}&\bounds{21}{21}{1}{false}&&\cr
\l{5}&\int(3){8}{}&\int(0){\0}{}&\int(0){2}{}&\ln{1}&\cnt{27}&\bounds{3}{3}{1}{false}&&\cr
\l{7}&\rat(9/4){9}{}&\rat(3/4){1}{}&\int(0){\0}{}&\ln{1}&\cnt{84}&\bounds{12}{12}{1}{false}&&\cr
\endorbit
\orbit{3}{19}{80}
\h{19}&\int(0){\0}{}&\int(4){\4}{}&\int(2){\2}{}&&\cnt{744}&\bounds{80}{84}{fail}{false}&&\cr
\l{1}&\int(0){11}{}&\int(2){3}{}&\int(1){2}{}&\ln{2\sd}&\cnt{288}&\bounds{24}{24}{1}{false}&&\cr
\l{2}&\int(0){\0}{}&\int(2){\2}{}&\int(1){\2}{}&\ln{1\sdd}&\cnt{120}&\bounds{20}{24}{2}{false}&&\cr
\l{3}&\int(0){\0}{}&\int(3){\4}{}&\int(0){\0}{}&\ln{1}&\cnt{24}&\bounds{6}{6}{1}{false}&&\cr
\endorbit
\orbit{3}{20}{80}
\h{20}&\rat(5/3){10}{}&\int(1){2}{}&\rat(10/3){1}{}&&\cnt{648}&\bounds{80}{80}{1}{false}&&\cr
\l{1}&\rat(5/6){11}{}&\rat(1/2){3}{}&\rat(5/3){2}{}&\ln{1\sdd}&\cnt{128}&\bounds{16}{16}{1}{false}&&\cr
\l{2}&\rat(4/3){8}{}&\int(0){\0}{}&\rat(5/3){2}{}&\ln{1}&\cnt{90}&\bounds{10}{10}{1}{false}&&\cr
\l{4}&\int(1){\2}{}&\int(0){\0}{}&\int(2){\2}{}&\ln{1}&\cnt{100}&\bounds{10}{10}{1}{false}&&\cr
\l{6}&\rat(5/3){10}{}&\int(1){2}{}&\rat(1/3){1}{}&\ln{1}&\cnt{10}&\bounds{2}{2}{1}{false}&&\cr
\l{8}&\rat(5/3){10}{}&\int(0){2}{}&\rat(4/3){1}{}&\ln{1}&\cnt{60}&\bounds{10}{10}{1}{false}&&\cr
\endorbit
\orbit{3}{21}{76}
\h{21}&\int(3){6}{}&\int(1){2}{}&\int(2){\2}{}&&\cnt{456}&\bounds{76}{76}{1}{false}&&\cr
\l{1}&\int(3){6}{}&\int(0){2}{}&\int(0){\0}{}&\ln{1}&\cnt{12}&\bounds{2}{2}{1}{false}&&\cr
\l{3}&\int(2){6}{}&\int(1){2}{}&\int(0){\0}{}&\ln{1}&\cnt{36}&\bounds{6}{6}{1}{false}&&\cr
\l{5}&\int(2){8}{}&\int(0){\0}{}&\int(1){2}{}&\ln{2}&\cnt{90}&\bounds{15}{15}{1}{false}&&\cr
\endorbit
\orbit{3}{22}{76}
\h{22}&\rat(11/12){11}{}&\rat(15/4){3}{}&\rat(4/3){2}{}&&\cnt{696}&\bounds{76}{76}{1}{false}&&\cr
\l{1}&\rat(3/4){9}{}&\rat(9/4){1}{}&\int(0){\0}{}&\ln{1}&\cnt{55}&\bounds{5}{5}{1}{false}&&\cr
\l{3}&\int(1){\2}{}&\int(2){\2}{}&\int(0){\0}{}&\ln{1}&\cnt{66}&\bounds{6}{6}{1}{false}&&\cr
\l{5}&\rat(5/6){10}{}&\rat(3/2){2}{}&\rat(2/3){1}{}&\ln{1}&\cnt{110}&\bounds{10}{10}{1}{false}&&\cr
\l{7}&\rat(11/12){11}{}&\rat(7/4){3}{}&\rat(1/3){2}{}&\ln{1}&\cnt{96}&\bounds{12}{12}{1}{false}&&\cr
\l{9}&\rat(11/12){11}{}&\rat(3/4){3}{\spec{1}}&\rat(4/3){2}{}&\ln{1}&\cnt{1}&\bounds{1}{1}{1}{false}&&\cr
\l{11}&\rat(11/12){11}{}&\rat(3/4){3}{\spec{2}}&\rat(4/3){2}{}&\ln{1}&\cnt{20}&\bounds{4}{4}{1}{false}&&\cr
\endorbit
\orbit{3}{23}{74}
\h{23}&\rat(11/3){10}{\spec{2}}&\int(1){2}{}&\rat(4/3){1}{}&&\cnt{552}&\bounds{74}{74}{1}{false}&&\cr
\l{1}&\int(3){\4}{}&\int(0){\0}{}&\int(0){\0}{}&\ln{1}&\cnt{24}&\bounds{3}{3}{1}{false}&&\cr
\l{3}&\int(2){6}{}&\int(1){2}{}&\int(0){\0}{}&\ln{1}&\cnt{56}&\bounds{8}{8}{1}{false}&&\cr
\l{5}&\rat(7/3){8}{}&\int(0){\0}{}&\rat(2/3){2}{}&\ln{1}&\cnt{80}&\bounds{10}{10}{1}{false}&&\cr
\l{7}&\rat(5/2){9}{}&\rat(1/2){1}{}&\int(0){\0}{}&\ln{1}&\cnt{32}&\bounds{4}{4}{1}{false}&&\cr
\l{9}&\int(2){\2}{}&\int(1){\2}{}&\int(0){\0}{}&\ln{1}&\cnt{36}&\bounds{6}{6}{1}{false}&&\cr
\l{11}&\int(2){\2}{}&\int(0){\0}{}&\int(1){\2}{}&\ln{1}&\cnt{48}&\bounds{6}{6}{1}{false}&&\cr
\endorbit
\orbit{3}{24}{72}
\h{24}&\rat(11/12){11}{}&\rat(7/4){3}{}&\rat(10/3){2}{}&&\cnt{648}&\bounds{72}{72}{1}{false}&&\cr
\l{1}&\int(1){\2}{}&\int(0){\0}{}&\int(2){\2}{}&\ln{1}&\cnt{55}&\bounds{5}{5}{1}{false}&&\cr
\l{3}&\rat(5/6){10}{}&\rat(1/2){2}{}&\rat(5/3){1}{}&\ln{1}&\cnt{154}&\bounds{14}{14}{1}{false}&&\cr
\l{5}&\rat(11/12){11}{}&\rat(7/4){3}{}&\rat(1/3){2}{}&\ln{1}&\cnt{10}&\bounds{2}{2}{1}{false}&&\cr
\l{7}&\rat(11/12){11}{}&\rat(3/4){3}{}&\rat(4/3){2}{}&\ln{1}&\cnt{105}&\bounds{15}{15}{1}{false}&&\cr
\endorbit
\orbit{3}{25}{68}
\h{25}&\rat(5/3){10}{}&\int(3){2}{}&\rat(4/3){1}{}&&\cnt{552}&\bounds{68}{76}{fail}{false}&&\cr
\l{1}&\rat(5/6){11}{}&\rat(3/2){3}{}&\rat(2/3){2}{}&\ln{1\sdd}&\cnt{160}&\bounds{24}{32}{2}{false}&&\cr
\l{2}&\rat(3/2){9}{}&\rat(3/2){1}{}&\int(0){\0}{}&\ln{1}&\cnt{80}&\bounds{8}{8}{1}{false}&&\cr
\l{4}&\int(1){\2}{}&\int(2){\2}{}&\int(0){\0}{}&\ln{1}&\cnt{60}&\bounds{6}{6}{1}{false}&&\cr
\l{6}&\rat(5/3){10}{}&\int(1){2}{}&\rat(1/3){1}{}&\ln{1}&\cnt{48}&\bounds{6}{6}{1}{false}&&\cr
\l{8}&\rat(5/3){10}{}&\int(0){2}{}&\rat(4/3){1}{}&\ln{1}&\cnt{8}&\bounds{2}{2}{1}{false}&&\cr
\endorbit
\orbit{3}{26}{68}
\h{26}&\rat(35/12){7}{}&\rat(7/4){3}{}&\rat(4/3){1}{}&&\cnt{408}&\bounds{68}{68}{1}{false}&&\cr
\l{1}&\rat(5/2){6}{}&\rat(1/2){2}{}&\int(0){\0}{}&\ln{1}&\cnt{49}&\bounds{7}{7}{1}{false}&&\cr
\l{3}&\rat(7/3){8}{}&\int(0){\0}{}&\rat(2/3){2}{}&\ln{1}&\cnt{50}&\bounds{10}{10}{1}{false}&&\cr
\l{5}&\rat(5/3){4}{}&\int(0){\0}{}&\rat(4/3){1}{}&\ln{1}&\cnt{35}&\bounds{7}{7}{1}{false}&&\cr
\l{7}&\rat(7/4){9}{}&\rat(5/4){1}{}&\int(0){\0}{}&\ln{1}&\cnt{70}&\bounds{10}{10}{1}{false}&&\cr
\endorbit
\orbit{3}{27}{68}
\h{27}&\int(2){\2}{}&\int(2){\2}{}&\int(2){\2}{}&&\cnt{552}&\bounds{68}{68}{1}{false}&&\cr
\l{1}&\int(1){10}{}&\int(1){2}{}&\int(1){1}{}&\ln{2\sd}&\cnt{120}&\bounds{12}{12}{1}{false}&&\cr
\l{2}&\int(1){11}{}&\int(1){3}{}&\int(1){2}{}&\ln{2\sd}&\cnt{96}&\bounds{12}{12}{1}{false}&&\cr
\l{3}&\int(2){\2}{}&\int(1){\2}{}&\int(0){\0}{}&\ln{1}&\cnt{20}&\bounds{4}{4}{1}{false}&&\cr
\l{5}&\int(2){\2}{}&\int(0){\0}{}&\int(1){\2}{}&\ln{1}&\cnt{20}&\bounds{4}{4}{1}{false}&&\cr
\l{7}&\int(1){\2}{\spec{1}}&\int(2){\2}{}&\int(0){\0}{}&\ln{2}&\cnt{10}&\bounds{1}{1}{1}{false}&&\cr
\endorbit
\orbit{3}{28}{64}
\h{28}&\rat(35/12){11}{\spec{2}}&\rat(7/4){3}{}&\rat(4/3){2}{}&&\cnt{504}&\bounds{64}{64}{1}{false}&&\cr
\l{1}&\int(3){\4}{}&\int(0){\0}{}&\int(0){\0}{}&\ln{1}&\cnt{9}&\bounds{1}{1}{1}{false}&&\cr
\l{3}&\rat(5/3){8}{}&\int(0){\0}{}&\rat(4/3){2}{}&\ln{1}&\cnt{36}&\bounds{4}{4}{1}{false}&&\cr
\l{5}&\rat(7/4){9}{}&\rat(5/4){1}{}&\int(0){\0}{}&\ln{1}&\cnt{63}&\bounds{7}{7}{1}{false}&&\cr
\l{7}&\int(2){\2}{}&\int(1){\2}{}&\int(0){\0}{}&\ln{1}&\cnt{42}&\bounds{6}{6}{1}{false}&&\cr
\l{9}&\int(2){\2}{}&\int(0){\0}{}&\int(1){\2}{}&\ln{1}&\cnt{32}&\bounds{4}{4}{1}{false}&&\cr
\l{11}&\rat(11/6){10}{}&\rat(1/2){2}{}&\rat(2/3){1}{}&\ln{1}&\cnt{70}&\bounds{10}{10}{1}{false}&&\cr
\endorbit
\orbit{3}{29}{64}
\h{29}&\rat(8/3){8}{}&\int(2){\2}{}&\rat(4/3){2}{}&&\cnt{456}&\bounds{64}{68}{fail}{false}&&\cr
\l{1}&\rat(4/3){10}{}&\int(1){2}{}&\rat(2/3){1}{}&\ln{1\sdd}&\cnt{120}&\bounds{20}{24}{2}{false}&&\cr
\l{2}&\int(2){6}{}&\int(1){2}{}&\int(0){\0}{}&\ln{1}&\cnt{56}&\bounds{8}{8}{1}{false}&&\cr
\l{4}&\rat(8/3){8}{}&\int(0){\0}{}&\rat(1/3){2}{}&\ln{1}&\cnt{16}&\bounds{2}{2}{1}{false}&&\cr
\l{6}&\rat(5/3){8}{}&\int(0){\0}{}&\rat(4/3){2}{}&\ln{1}&\cnt{32}&\bounds{4}{4}{1}{false}&&\cr
\l{8}&\int(2){9}{}&\int(1){1}{}&\int(0){\0}{}&\ln{1}&\cnt{64}&\bounds{8}{8}{1}{false}&&\cr
\endorbit
\orbit{3}{30}{60}
\h{30}&\rat(9/4){9}{}&\rat(7/4){1}{}&\int(2){\2}{}&&\cnt{456}&\bounds{60}{60}{1}{false}&&\cr
\l{1}&\int(2){8}{}&\int(0){\0}{}&\int(1){2}{}&\ln{1}&\cnt{54}&\bounds{6}{6}{1}{false}&&\cr
\l{3}&\rat(9/4){9}{}&\rat(3/4){1}{}&\int(0){\0}{}&\ln{1}&\cnt{21}&\bounds{3}{3}{1}{false}&&\cr
\l{5}&\rat(5/4){9}{}&\rat(7/4){1}{}&\int(0){\0}{}&\ln{1}&\cnt{27}&\bounds{3}{3}{1}{false}&&\cr
\l{7}&\rat(3/2){10}{}&\rat(1/2){2}{}&\int(1){1}{}&\ln{1}&\cnt{126}&\bounds{18}{18}{1}{false}&&\cr
\endorbit
\orbit{3}{31}{54}
\h{31}&\int(6){\6}{\spec{1}}&\int(0){\0}{}&\int(0){\0}{}&&\cnt{456}&\bounds{54}{54}{1}{false}&&\cr
\l{1}&\int(3){\4}{}&\int(0){\0}{}&\int(0){\0}{}&\ln{1\sdd}&\cnt{144}&\bounds{18}{18}{1}{false}&&\cr
\l{2}&\int(3){\2}{}&\int(0){\2}{}&\int(0){\0}{}&\ln{1\sdd}&\cnt{168}&\bounds{18}{18}{1}{false}&&\cr
\l{3}&\int(3){\2}{}&\int(0){\0}{}&\int(0){\2}{}&\ln{1\sdd}&\cnt{144}&\bounds{18}{18}{1}{false}&&\cr
\endorbit
\orbit{3}{32}{54}
\h{32}&\int(0){\0}{}&\int(0){\0}{}&\int(6){\6}{}&&\cnt{456}&\bounds{54}{54}{1}{false}&&\cr
\l{1}&\int(0){\2}{}&\int(0){\0}{}&\int(3){\2}{}&\ln{1\sdd}&\cnt{264}&\bounds{24}{24}{1}{false}&&\cr
\l{2}&\int(0){\0}{}&\int(0){\2}{}&\int(3){\2}{}&\ln{1\sdd}&\cnt{168}&\bounds{18}{18}{1}{false}&&\cr
\l{3}&\int(0){\0}{}&\int(0){\0}{}&\int(3){\4}{\spec{1}}&\ln{2\sdd}&\cnt{12}&\bounds{6}{6}{1}{false}&&\cr
\endorbit
\orbit{3}{33}{52}
\h{33}&\int(0){\0}{}&\int(6){\6}{\spec{1}}&\int(0){\0}{}&&\cnt{456}&\bounds{52}{52}{1}{false}&&\cr
\l{1}&\int(0){\2}{}&\int(3){\2}{}&\int(0){\0}{}&\ln{1\sdd}&\cnt{264}&\bounds{24}{24}{1}{false}&&\cr
\l{2}&\int(0){\0}{}&\int(3){\4}{}&\int(0){\0}{}&\ln{1\sdd}&\cnt{48}&\bounds{10}{10}{1}{false}&&\cr
\l{3}&\int(0){\0}{}&\int(3){\2}{}&\int(0){\2}{}&\ln{1\sdd}&\cnt{144}&\bounds{18}{18}{1}{false}&&\cr
\endorbit
\endorbits

\newcs{max-4E6}{104}%
\newcs{count-4E6}{5}%
\newcs{counts-4E6}{\DO{1}{104}\DO{2}{96}\DO{3}{82}\DO{4}{72}\DO{5}{66}}%
\orbits
\orbit{4}{1}{104}
\h{1}&\int(4){\4}{}&\int(2){\2}{}&\int(0){\0}{}&\int(0){\0}{}&&\cnt{840}&\bounds{104}{112}{fail}{false}&&\cr
\l{1}&\int(2){\2}{}&\int(1){\2}{}&\int(0){\0}{}&\int(0){\0}{}&\ln{1\sdd}&\cnt{160}&\bounds{24}{32}{2}{false}&&\cr
\l{2}&\int(2){2}{}&\int(1){2}{}&\int(0){1}{}&\int(0){\0}{}&\ln{4\sd}&\cnt{162}&\bounds{18}{18}{1}{false}&&\cr
\l{3}&\int(3){\4}{\spec{1}}&\int(0){\0}{}&\int(0){\0}{}&\int(0){\0}{}&\ln{2}&\cnt{8}&\bounds{2}{2}{1}{false}&&\cr
\endorbit
\orbit{4}{2}{96}
\h{2}&\int(2){\2}{}&\int(2){\2}{}&\int(2){\2}{}&\int(0){\0}{}&&\cnt{552}&\bounds{96}{96}{1}{false}&&\cr
\l{1}&\int(2){\2}{}&\int(1){\2}{}&\int(0){\0}{}&\int(0){\0}{}&\ln{6\sd}&\cnt{20}&\bounds{4}{4}{1}{false}&&\cr
\l{2}&\int(1){2}{}&\int(1){2}{}&\int(1){1}{}&\int(0){\0}{}&\ln{2\sd}&\cnt{216}&\bounds{36}{36}{1}{false}&&\cr
\endorbit
\orbit{4}{3}{82}
\h{3}&\rat(10/3){2}{}&\rat(4/3){2}{}&\rat(4/3){1}{}&\int(0){\0}{}&&\cnt{648}&\bounds{82}{96}{fail}{false}&&\cr
\l{1}&\rat(5/3){1}{}&\rat(4/3){2}{}&\int(0){\0}{}&\int(0){1}{}&\ln{2\sd}&\cnt{54}&\bounds{6}{6}{1}{false}&&\cr
\l{2}&\rat(5/3){1}{}&\rat(2/3){1}{}&\rat(2/3){2}{}&\int(0){\0}{}&\ln{1\sdd}&\cnt{200}&\bounds{26}{40}{2}{false}&&\cr
\l{3}&\int(3){\4}{}&\int(0){\0}{}&\int(0){\0}{}&\int(0){\0}{}&\ln{1}&\cnt{10}&\bounds{2}{2}{1}{false}&&\cr
\l{5}&\int(2){\2}{}&\int(1){\2}{}&\int(0){\0}{}&\int(0){\0}{}&\ln{2}&\cnt{80}&\bounds{10}{10}{1}{false}&&\cr
\endorbit
\orbit{4}{4}{72}
\h{4}&\rat(4/3){2}{}&\rat(4/3){2}{}&\rat(4/3){1}{}&\int(2){\2}{}&&\cnt{456}&\bounds{72}{72}{1}{false}&&\cr
\l{1}&\rat(4/3){2}{}&\rat(2/3){1}{}&\int(0){\0}{}&\int(1){2}{}&\ln{6\sd}&\cnt{60}&\bounds{10}{10}{1}{false}&&\cr
\l{2}&\int(1){\2}{}&\int(0){\0}{}&\int(0){\0}{}&\int(2){\2}{}&\ln{3}&\cnt{16}&\bounds{2}{2}{1}{false}&&\cr
\endorbit
\orbit{4}{5}{66}
\h{5}&\int(6){\6}{}&\int(0){\0}{}&\int(0){\0}{}&\int(0){\0}{}&&\cnt{456}&\bounds{66}{66}{1}{false}&&\cr
\l{1}&\int(3){\4}{\spec{1}}&\int(0){\0}{}&\int(0){\0}{}&\int(0){\0}{}&\ln{2\sdd}&\cnt{12}&\bounds{6}{6}{1}{false}&&\cr
\l{2}&\int(3){\2}{}&\int(0){\2}{}&\int(0){\0}{}&\int(0){\0}{}&\ln{3\sdd}&\cnt{144}&\bounds{18}{18}{1}{false}&&\cr
\endorbit
\endorbits

\newcs{max-2A9_D6}{114}%
\newcs{count-2A9_D6}{27}%
\newcs{counts-2A9_D6}{\DO{1}{114}\DO{2}{106}\DO{3}{102}\DO{4}{100}\DO{5}{100}\DO{6}{100}\DO{7}{98}\DO{8}{96}\DO{9}{96}\DO{10}{96}\DO{11}{94}\DO{12}{94}\DO{13}{90}\DO{14}{88}\DO{15}{84}\DO{16}{84}\DO{17}{84}\DO{18}{84}\DO{19}{80}\DO{20}{80}\DO{21}{80}\DO{22}{78}\DO{23}{76}\DO{24}{76}\DO{25}{72}\DO{26}{50}\DO{27}{48}}%
\orbits
\orbit{3}{1}{114}
\h{1}&\int(0){\0}{}&\int(0){\0}{}&\int(6){\6}{3}&&\cnt{1152}&\bounds{114}{114}{1}{false}&&\cr
\l{1}&\int(0){8}{}&\int(0){1}{}&\int(3){3}{}&\ln{2\sd}&\cnt{450}&\bounds{45}{45}{1}{false}&&\cr
\l{2}&\int(0){\0}{}&\int(0){5}{}&\int(3){3}{}&\ln{1\sdd}&\cnt{252}&\bounds{24}{24}{1}{false}&&\cr
\endorbit
\orbit{3}{2}{106}
\h{2}&\int(6){\6}{}&\int(0){\0}{}&\int(0){\0}{}&&\cnt{864}&\bounds{106}{106}{1}{false}&&\cr
\l{1}&\int(3){\4}{\spec{1}}&\int(0){\0}{}&\int(0){\0}{}&\ln{2\sd}&\cnt{36}&\bounds{9}{9}{1}{false}&&\cr
\l{2}&\int(3){5}{}&\int(0){\0}{}&\int(0){1}{}&\ln{1\sdd}&\cnt{192}&\bounds{24}{24}{1}{false}&&\cr
\l{3}&\int(3){6}{}&\int(0){2}{}&\int(0){\0}{}&\ln{2\sd}&\cnt{180}&\bounds{20}{20}{1}{false}&&\cr
\l{4}&\int(3){7}{}&\int(0){9}{}&\int(0){2}{}&\ln{2\sd}&\cnt{120}&\bounds{12}{12}{1}{false}&&\cr
\endorbit
\orbit{3}{3}{102}
\h{3}&\int(2){\2}{}&\int(0){\0}{}&\int(4){\4}{2}&&\cnt{960}&\bounds{102}{102}{1}{false}&&\cr
\l{1}&\int(1){7}{}&\int(0){9}{}&\int(2){2}{}&\ln{2\sd}&\cnt{280}&\bounds{28}{28}{1}{false}&&\cr
\l{2}&\int(1){\2}{\spec{1}}&\int(0){\0}{}&\int(2){\2}{}&\ln{2\sd}&\cnt{80}&\bounds{10}{10}{1}{false}&&\cr
\l{3}&\int(1){9}{}&\int(0){3}{}&\int(2){2}{}&\ln{2\sd}&\cnt{120}&\bounds{13}{13}{1}{false}&&\cr
\endorbit
\orbit{3}{4}{100}
\h{4}&\rat(18/5){8}{\spec{1}}&\rat(12/5){6}{}&\int(0){\0}{}&&\cnt{864}&\bounds{100}{112}{fail}{false}&&\cr
\l{1}&\rat(9/5){9}{}&\rat(6/5){3}{}&\int(0){2}{}&\ln{1\sdd}&\cnt{240}&\bounds{28}{40}{2}{false}&&\cr
\l{2}&\rat(9/5){9}{}&\rat(6/5){8}{}&\int(0){1}{}&\ln{1\sdd}&\cnt{192}&\bounds{24}{24}{1}{false}&&\cr
\l{3}&\int(2){\2}{}&\int(1){\2}{}&\int(0){\0}{}&\ln{1}&\cnt{216}&\bounds{24}{24}{1}{false}&&\cr
\endorbit
\orbit{3}{5}{100}
\h{5}&\rat(8/5){8}{}&\rat(22/5){6}{\spec{2}}&\int(0){\0}{}&&\cnt{672}&\bounds{100}{100}{1}{false}&&\cr
\l{1}&\rat(4/5){9}{}&\rat(11/5){3}{}&\int(0){2}{}&\ln{1\sdd}&\cnt{144}&\bounds{24}{24}{1}{false}&&\cr
\l{2}&\rat(6/5){6}{}&\rat(9/5){2}{}&\int(0){\0}{}&\ln{1}&\cnt{112}&\bounds{16}{16}{1}{false}&&\cr
\l{4}&\int(1){\2}{}&\int(2){\2}{}&\int(0){\0}{}&\ln{1}&\cnt{80}&\bounds{10}{10}{1}{false}&&\cr
\l{6}&\rat(8/5){8}{}&\rat(7/5){6}{}&\int(0){\0}{}&\ln{1}&\cnt{40}&\bounds{8}{8}{1}{false}&&\cr
\l{8}&\rat(8/5){8}{}&\rat(7/5){1}{}&\int(0){3}{}&\ln{1}&\cnt{32}&\bounds{4}{4}{1}{false}&&\cr
\endorbit
\orbit{3}{6}{100}
\h{6}&\rat(9/10){9}{}&\rat(41/10){3}{\spec{2}}&\int(1){2}{}&&\cnt{816}&\bounds{100}{100}{1}{false}&&\cr
\l{1}&\rat(3/5){6}{}&\rat(12/5){2}{}&\int(0){\0}{}&\ln{1}&\cnt{84}&\bounds{12}{12}{1}{false}&&\cr
\l{3}&\int(1){\2}{}&\int(2){\2}{}&\int(0){\0}{}&\ln{1}&\cnt{72}&\bounds{8}{8}{1}{false}&&\cr
\l{5}&\rat(4/5){8}{}&\rat(17/10){1}{}&\rat(1/2){3}{}&\ln{1}&\cnt{144}&\bounds{16}{16}{1}{false}&&\cr
\l{7}&\rat(9/10){9}{}&\rat(21/10){3}{}&\int(0){2}{}&\ln{1}&\cnt{80}&\bounds{10}{10}{1}{false}&&\cr
\l{9}&\rat(9/10){9}{}&\rat(11/10){3}{}&\int(1){2}{}&\ln{1}&\cnt{28}&\bounds{4}{4}{1}{false}&&\cr
\endorbit
\orbit{3}{7}{98}
\h{7}&\rat(8/5){8}{}&\rat(22/5){6}{\spec{1}}&\int(0){\0}{}&&\cnt{768}&\bounds{98}{98}{1}{false}&&\cr
\l{1}&\rat(4/5){4}{}&\rat(11/5){8}{}&\int(0){\0}{}&\ln{1\sdd}&\cnt{140}&\bounds{18}{18}{1}{false}&&\cr
\l{2}&\rat(4/5){9}{}&\rat(11/5){8}{}&\int(0){1}{}&\ln{1\sdd}&\cnt{128}&\bounds{16}{16}{1}{false}&&\cr
\l{3}&\rat(7/5){7}{}&\rat(8/5){9}{}&\int(0){2}{}&\ln{1}&\cnt{96}&\bounds{12}{12}{1}{false}&&\cr
\l{5}&\int(1){\2}{}&\int(2){\2}{}&\int(0){\0}{}&\ln{1}&\cnt{112}&\bounds{14}{14}{1}{false}&&\cr
\l{7}&\rat(8/5){8}{}&\rat(7/5){6}{}&\int(0){\0}{}&\ln{1}&\cnt{42}&\bounds{6}{6}{1}{false}&&\cr
\endorbit
\orbit{3}{8}{96}
\h{8}&\rat(18/5){8}{\spec{2}}&\rat(12/5){6}{}&\int(0){\0}{}&&\cnt{576}&\bounds{96}{100}{fail}{false}&&\cr
\l{1}&\rat(9/5){4}{}&\rat(6/5){8}{}&\int(0){\0}{}&\ln{1\sdd}&\cnt{120}&\bounds{20}{24}{2}{false}&&\cr
\l{2}&\int(3){\4}{}&\int(0){\0}{}&\int(0){\0}{}&\ln{1}&\cnt{18}&\bounds{3}{3}{1}{false}&&\cr
\l{4}&\rat(11/5){6}{}&\rat(4/5){2}{}&\int(0){\0}{}&\ln{1}&\cnt{90}&\bounds{15}{15}{1}{false}&&\cr
\l{6}&\rat(12/5){7}{}&\rat(3/5){9}{}&\int(0){2}{}&\ln{1}&\cnt{48}&\bounds{8}{8}{1}{false}&&\cr
\l{8}&\int(2){\2}{}&\int(1){\2}{}&\int(0){\0}{}&\ln{1}&\cnt{72}&\bounds{12}{12}{1}{false}&&\cr
\endorbit
\orbit{3}{9}{96}
\h{9}&\rat(5/2){5}{}&\int(0){\0}{}&\rat(7/2){1}{}&&\cnt{672}&\bounds{96}{96}{1}{false}&&\cr
\l{1}&\rat(5/2){5}{}&\int(0){\0}{}&\rat(1/2){1}{\spec{1}}&\ln{1}&\cnt{1}&\bounds{1}{1}{1}{false}&&\cr
\l{3}&\rat(5/2){5}{}&\int(0){\0}{}&\rat(1/2){1}{\spec{2}}&\ln{1}&\cnt{10}&\bounds{2}{2}{1}{false}&&\cr
\l{5}&\rat(3/2){5}{}&\int(0){\0}{}&\rat(3/2){1}{}&\ln{1}&\cnt{125}&\bounds{25}{25}{1}{false}&&\cr
\l{7}&\rat(3/2){7}{}&\int(0){9}{}&\rat(3/2){2}{}&\ln{2}&\cnt{100}&\bounds{10}{10}{1}{false}&&\cr
\endorbit
\orbit{3}{10}{96}
\h{10}&\rat(9/10){9}{}&\rat(18/5){8}{\spec{1}}&\rat(3/2){1}{}&&\cnt{864}&\bounds{96}{96}{1}{false}&&\cr
\l{1}&\rat(7/10){7}{}&\rat(9/5){9}{}&\rat(1/2){2}{}&\ln{1}&\cnt{216}&\bounds{24}{24}{1}{false}&&\cr
\l{3}&\int(1){\2}{}&\int(2){\2}{}&\int(0){\0}{}&\ln{1}&\cnt{81}&\bounds{9}{9}{1}{false}&&\cr
\l{5}&\rat(9/10){9}{}&\rat(8/5){8}{}&\rat(1/2){1}{}&\ln{1}&\cnt{135}&\bounds{15}{15}{1}{false}&&\cr
\endorbit
\orbit{3}{11}{94}
\h{11}&\int(4){\4}{}&\int(2){\2}{}&\int(0){\0}{}&&\cnt{672}&\bounds{94}{94}{1}{false}&&\cr
\l{1}&\int(2){6}{}&\int(1){2}{}&\int(0){\0}{}&\ln{2\sd}&\cnt{120}&\bounds{15}{15}{1}{false}&&\cr
\l{2}&\int(2){7}{}&\int(1){9}{}&\int(0){2}{}&\ln{2\sd}&\cnt{72}&\bounds{12}{12}{1}{false}&&\cr
\l{3}&\int(2){\2}{}&\int(1){\2}{\spec{1}}&\int(0){\0}{}&\ln{2\sd}&\cnt{32}&\bounds{4}{4}{1}{false}&&\cr
\l{4}&\int(2){8}{}&\int(1){6}{}&\int(0){\0}{}&\ln{2\sd}&\cnt{56}&\bounds{8}{8}{1}{false}&&\cr
\l{5}&\int(2){8}{}&\int(1){1}{}&\int(0){3}{}&\ln{2\sd}&\cnt{32}&\bounds{4}{4}{1}{false}&&\cr
\l{6}&\int(3){\4}{\spec{1}}&\int(0){\0}{}&\int(0){\0}{}&\ln{2}&\cnt{12}&\bounds{2}{2}{1}{false}&&\cr
\endorbit
\orbit{3}{12}{94}
\h{12}&\rat(9/2){5}{\spec{1}}&\int(0){\0}{}&\rat(3/2){1}{}&&\cnt{720}&\bounds{94}{94}{1}{false}&&\cr
\l{1}&\int(3){\4}{}&\int(0){\0}{}&\int(0){\0}{}&\ln{1}&\cnt{45}&\bounds{9}{9}{1}{false}&&\cr
\l{3}&\rat(5/2){5}{}&\int(0){\0}{}&\rat(1/2){1}{}&\ln{1}&\cnt{90}&\bounds{15}{15}{1}{false}&&\cr
\l{5}&\int(3){6}{}&\int(0){2}{}&\int(0){\0}{}&\ln{1}&\cnt{45}&\bounds{5}{5}{1}{false}&&\cr
\l{7}&\rat(5/2){7}{}&\int(0){9}{}&\rat(1/2){2}{}&\ln{1}&\cnt{180}&\bounds{18}{18}{1}{false}&&\cr
\endorbit
\orbit{3}{13}{90}
\h{13}&\int(2){\2}{}&\int(0){\0}{}&\int(4){\4}{}&&\cnt{768}&\bounds{90}{90}{1}{false}&&\cr
\l{1}&\int(1){5}{}&\int(0){\0}{}&\int(2){1}{}&\ln{1\sdd}&\cnt{140}&\bounds{18}{18}{1}{false}&&\cr
\l{2}&\int(1){\2}{\spec{1}}&\int(0){\0}{}&\int(2){\2}{}&\ln{2\sd}&\cnt{48}&\bounds{6}{6}{1}{false}&&\cr
\l{3}&\int(1){8}{}&\int(0){1}{}&\int(2){3}{}&\ln{2\sd}&\cnt{160}&\bounds{16}{16}{1}{false}&&\cr
\l{4}&\int(1){9}{}&\int(0){8}{}&\int(2){1}{}&\ln{2\sd}&\cnt{90}&\bounds{10}{10}{1}{false}&&\cr
\l{5}&\int(2){\2}{}&\int(0){\0}{}&\int(1){\2}{}&\ln{1}&\cnt{16}&\bounds{4}{4}{1}{false}&&\cr
\endorbit
\orbit{3}{14}{88}
\h{14}&\int(4){\4}{}&\int(0){\0}{}&\int(2){\2}{}&&\cnt{672}&\bounds{88}{96}{fail}{false}&&\cr
\l{1}&\int(2){5}{}&\int(0){\0}{}&\int(1){1}{}&\ln{1\sdd}&\cnt{160}&\bounds{24}{32}{2}{false}&&\cr
\l{2}&\int(2){7}{}&\int(0){9}{}&\int(1){2}{}&\ln{2\sd}&\cnt{120}&\bounds{12}{12}{1}{false}&&\cr
\l{3}&\int(2){\2}{}&\int(0){\0}{}&\int(1){\2}{}&\ln{1\sdd}&\cnt{64}&\bounds{16}{16}{1}{false}&&\cr
\l{4}&\int(2){8}{}&\int(0){1}{}&\int(1){3}{}&\ln{2\sd}&\cnt{80}&\bounds{8}{8}{1}{false}&&\cr
\l{5}&\int(3){\4}{\spec{1}}&\int(0){\0}{}&\int(0){\0}{}&\ln{2}&\cnt{12}&\bounds{2}{2}{1}{false}&&\cr
\endorbit
\orbit{3}{15}{84}
\h{15}&\rat(9/10){9}{}&\rat(8/5){8}{}&\rat(7/2){1}{}&&\cnt{672}&\bounds{84}{84}{1}{false}&&\cr
\l{1}&\rat(7/10){7}{}&\rat(4/5){9}{}&\rat(3/2){2}{}&\ln{1}&\cnt{72}&\bounds{8}{8}{1}{false}&&\cr
\l{3}&\int(1){\2}{}&\int(0){\0}{}&\int(2){\2}{}&\ln{1}&\cnt{45}&\bounds{5}{5}{1}{false}&&\cr
\l{5}&\rat(4/5){8}{}&\rat(1/5){1}{}&\int(2){3}{}&\ln{1}&\cnt{72}&\bounds{8}{8}{1}{false}&&\cr
\l{7}&\rat(9/10){9}{}&\rat(3/5){3}{}&\rat(3/2){2}{}&\ln{1}&\cnt{56}&\bounds{8}{8}{1}{false}&&\cr
\l{9}&\rat(9/10){9}{}&\rat(8/5){8}{}&\rat(1/2){1}{\spec{1}}&\ln{1}&\cnt{1}&\bounds{1}{1}{1}{false}&&\cr
\l{11}&\rat(9/10){9}{}&\rat(8/5){8}{}&\rat(1/2){1}{\spec{2}}&\ln{1}&\cnt{10}&\bounds{2}{2}{1}{false}&&\cr
\l{13}&\rat(9/10){9}{}&\rat(3/5){8}{}&\rat(3/2){1}{}&\ln{1}&\cnt{80}&\bounds{10}{10}{1}{false}&&\cr
\endorbit
\orbit{3}{16}{84}
\h{16}&\rat(5/2){5}{}&\rat(5/2){5}{}&\int(1){2}{}&&\cnt{432}&\bounds{84}{84}{1}{false}&&\cr
\l{1}&\rat(5/2){5}{}&\int(0){\0}{}&\rat(1/2){1}{}&\ln{2\sd}&\cnt{16}&\bounds{2}{2}{1}{false}&&\cr
\l{2}&\int(2){6}{}&\int(1){2}{}&\int(0){\0}{}&\ln{4}&\cnt{50}&\bounds{10}{10}{1}{false}&&\cr
\endorbit
\orbit{3}{17}{84}
\h{17}&\rat(9/10){9}{}&\rat(18/5){8}{\spec{2}}&\rat(3/2){1}{}&&\cnt{576}&\bounds{84}{84}{1}{false}&&\cr
\l{1}&\int(1){\2}{}&\int(2){\2}{}&\int(0){\0}{}&\ln{1}&\cnt{27}&\bounds{3}{3}{1}{false}&&\cr
\l{3}&\rat(4/5){8}{}&\rat(11/5){6}{}&\int(0){\0}{}&\ln{1}&\cnt{54}&\bounds{6}{6}{1}{false}&&\cr
\l{5}&\rat(4/5){8}{}&\rat(6/5){1}{}&\int(1){3}{}&\ln{1}&\cnt{54}&\bounds{6}{6}{1}{false}&&\cr
\l{7}&\rat(9/10){9}{}&\rat(8/5){3}{}&\rat(1/2){2}{}&\ln{1}&\cnt{90}&\bounds{15}{15}{1}{false}&&\cr
\l{9}&\rat(9/10){9}{}&\rat(8/5){8}{}&\rat(1/2){1}{}&\ln{1}&\cnt{45}&\bounds{9}{9}{1}{false}&&\cr
\l{11}&\rat(9/10){9}{}&\rat(3/5){8}{}&\rat(3/2){1}{}&\ln{1}&\cnt{18}&\bounds{3}{3}{1}{false}&&\cr
\endorbit
\orbit{3}{18}{84}
\h{18}&\rat(9/10){9}{}&\rat(41/10){3}{\spec{1}}&\int(1){2}{}&&\cnt{624}&\bounds{84}{84}{1}{false}&&\cr
\l{1}&\rat(7/10){7}{}&\rat(13/10){9}{}&\int(1){2}{}&\ln{1}&\cnt{36}&\bounds{4}{4}{1}{false}&&\cr
\l{3}&\int(1){\2}{}&\int(2){\2}{}&\int(0){\0}{}&\ln{1}&\cnt{36}&\bounds{4}{4}{1}{false}&&\cr
\l{5}&\rat(4/5){8}{}&\rat(11/5){6}{}&\int(0){\0}{}&\ln{1}&\cnt{90}&\bounds{10}{10}{1}{false}&&\cr
\l{7}&\rat(9/10){9}{}&\rat(21/10){3}{}&\int(0){2}{}&\ln{1}&\cnt{40}&\bounds{8}{8}{1}{false}&&\cr
\l{9}&\rat(9/10){9}{}&\rat(11/10){3}{}&\int(1){2}{}&\ln{1}&\cnt{30}&\bounds{6}{6}{1}{false}&&\cr
\l{11}&\rat(9/10){9}{}&\rat(8/5){8}{}&\rat(1/2){1}{}&\ln{1}&\cnt{80}&\bounds{10}{10}{1}{false}&&\cr
\endorbit
\orbit{3}{19}{80}
\h{19}&\rat(29/10){9}{\spec{2}}&\rat(21/10){3}{}&\int(1){2}{}&&\cnt{528}&\bounds{80}{80}{1}{false}&&\cr
\l{1}&\int(3){\4}{}&\int(0){\0}{}&\int(0){\0}{}&\ln{1}&\cnt{7}&\bounds{1}{1}{1}{false}&&\cr
\l{3}&\rat(8/5){6}{}&\rat(7/5){2}{}&\int(0){\0}{}&\ln{1}&\cnt{63}&\bounds{9}{9}{1}{false}&&\cr
\l{5}&\rat(17/10){7}{}&\rat(3/10){9}{}&\int(1){2}{}&\ln{1}&\cnt{49}&\bounds{7}{7}{1}{false}&&\cr
\l{7}&\int(2){\2}{}&\int(1){\2}{}&\int(0){\0}{}&\ln{1}&\cnt{42}&\bounds{6}{6}{1}{false}&&\cr
\l{9}&\int(2){\2}{}&\int(0){\0}{}&\int(1){\2}{}&\ln{1}&\cnt{20}&\bounds{4}{4}{1}{false}&&\cr
\l{11}&\rat(9/5){8}{}&\rat(6/5){6}{}&\int(0){\0}{}&\ln{1}&\cnt{35}&\bounds{7}{7}{1}{false}&&\cr
\l{13}&\rat(9/5){8}{}&\rat(7/10){1}{}&\rat(1/2){3}{}&\ln{1}&\cnt{48}&\bounds{6}{6}{1}{false}&&\cr
\endorbit
\orbit{3}{20}{80}
\h{20}&\rat(29/10){9}{\spec{2}}&\rat(8/5){8}{}&\rat(3/2){1}{}&&\cnt{528}&\bounds{80}{80}{1}{false}&&\cr
\l{1}&\int(3){\4}{}&\int(0){\0}{}&\int(0){\0}{}&\ln{1}&\cnt{7}&\bounds{1}{1}{1}{false}&&\cr
\l{3}&\rat(3/2){5}{}&\int(0){\0}{}&\rat(3/2){1}{}&\ln{1}&\cnt{35}&\bounds{7}{7}{1}{false}&&\cr
\l{5}&\rat(17/10){7}{}&\rat(4/5){9}{}&\rat(1/2){2}{}&\ln{1}&\cnt{84}&\bounds{12}{12}{1}{false}&&\cr
\l{7}&\int(2){\2}{}&\int(1){\2}{}&\int(0){\0}{}&\ln{1}&\cnt{32}&\bounds{4}{4}{1}{false}&&\cr
\l{9}&\int(2){\2}{}&\int(0){\0}{}&\int(1){\2}{}&\ln{1}&\cnt{30}&\bounds{6}{6}{1}{false}&&\cr
\l{11}&\rat(9/5){8}{}&\rat(6/5){6}{}&\int(0){\0}{}&\ln{1}&\cnt{28}&\bounds{4}{4}{1}{false}&&\cr
\l{13}&\rat(9/5){8}{}&\rat(1/5){1}{}&\int(1){3}{}&\ln{1}&\cnt{48}&\bounds{6}{6}{1}{false}&&\cr
\endorbit
\orbit{3}{21}{80}
\h{21}&\int(2){\2}{}&\int(2){\2}{}&\int(2){\2}{}&&\cnt{576}&\bounds{80}{80}{1}{false}&&\cr
\l{1}&\int(1){7}{}&\int(1){9}{}&\int(1){2}{}&\ln{4\sd}&\cnt{56}&\bounds{8}{8}{1}{false}&&\cr
\l{2}&\int(2){\2}{}&\int(0){\0}{}&\int(1){\2}{}&\ln{2\sd}&\cnt{16}&\bounds{4}{4}{1}{false}&&\cr
\l{3}&\int(1){8}{}&\int(1){1}{}&\int(1){3}{}&\ln{4\sd}&\cnt{64}&\bounds{8}{8}{1}{false}&&\cr
\l{4}&\int(2){\2}{}&\int(1){\2}{\spec{1}}&\int(0){\0}{}&\ln{4}&\cnt{8}&\bounds{1}{1}{1}{false}&&\cr
\endorbit
\orbit{3}{22}{78}
\h{22}&\rat(8/5){8}{}&\rat(12/5){6}{}&\int(2){\2}{}&&\cnt{480}&\bounds{78}{84}{fail}{false}&&\cr
\l{1}&\rat(4/5){9}{}&\rat(6/5){3}{}&\int(1){2}{}&\ln{1\sdd}&\cnt{80}&\bounds{16}{16}{1}{false}&&\cr
\l{2}&\rat(4/5){9}{}&\rat(6/5){8}{}&\int(1){1}{}&\ln{1\sdd}&\cnt{96}&\bounds{18}{24}{2}{false}&&\cr
\l{3}&\rat(7/5){7}{}&\rat(3/5){9}{}&\int(1){2}{}&\ln{1}&\cnt{64}&\bounds{8}{8}{1}{false}&&\cr
\l{5}&\int(1){\2}{}&\int(0){\0}{}&\int(2){\2}{}&\ln{1}&\cnt{16}&\bounds{2}{2}{1}{false}&&\cr
\l{7}&\rat(8/5){8}{}&\rat(7/5){6}{}&\int(0){\0}{}&\ln{1}&\cnt{24}&\bounds{4}{4}{1}{false}&&\cr
\l{9}&\rat(8/5){8}{}&\rat(2/5){1}{}&\int(1){3}{}&\ln{1}&\cnt{48}&\bounds{8}{8}{1}{false}&&\cr
\endorbit
\orbit{3}{23}{76}
\h{23}&\rat(9/10){9}{}&\rat(21/10){3}{}&\int(3){2}{}&&\cnt{576}&\bounds{76}{76}{1}{false}&&\cr
\l{1}&\int(1){\2}{}&\int(0){\0}{}&\int(2){\2}{}&\ln{1}&\cnt{27}&\bounds{3}{3}{1}{false}&&\cr
\l{3}&\rat(4/5){8}{}&\rat(7/10){1}{}&\rat(3/2){3}{}&\ln{1}&\cnt{108}&\bounds{12}{12}{1}{false}&&\cr
\l{5}&\rat(9/10){9}{}&\rat(21/10){3}{}&\int(0){2}{}&\ln{1}&\cnt{6}&\bounds{2}{2}{1}{false}&&\cr
\l{7}&\rat(9/10){9}{}&\rat(11/10){3}{}&\int(1){2}{}&\ln{1}&\cnt{63}&\bounds{9}{9}{1}{false}&&\cr
\l{9}&\rat(9/10){9}{}&\rat(3/5){8}{}&\rat(3/2){1}{}&\ln{1}&\cnt{84}&\bounds{12}{12}{1}{false}&&\cr
\endorbit
\orbit{3}{24}{76}
\h{24}&\rat(5/2){5}{}&\int(2){\2}{}&\rat(3/2){1}{}&&\cnt{480}&\bounds{76}{76}{1}{false}&&\cr
\l{1}&\rat(5/2){5}{}&\int(0){\0}{}&\rat(1/2){1}{}&\ln{1}&\cnt{15}&\bounds{3}{3}{1}{false}&&\cr
\l{3}&\rat(3/2){5}{}&\int(0){\0}{}&\rat(3/2){1}{}&\ln{1}&\cnt{25}&\bounds{5}{5}{1}{false}&&\cr
\l{5}&\int(2){6}{}&\int(1){2}{}&\int(0){\0}{}&\ln{2}&\cnt{40}&\bounds{5}{5}{1}{false}&&\cr
\l{7}&\rat(3/2){7}{}&\int(1){9}{}&\rat(1/2){2}{}&\ln{2}&\cnt{60}&\bounds{10}{10}{1}{false}&&\cr
\endorbit
\orbit{3}{25}{72}
\h{25}&\rat(21/10){7}{}&\rat(12/5){4}{}&\rat(3/2){1}{}&&\cnt{432}&\bounds{72}{72}{1}{false}&&\cr
\l{1}&\rat(3/2){5}{}&\int(0){\0}{}&\rat(3/2){1}{}&\ln{1}&\cnt{21}&\bounds{3}{3}{1}{false}&&\cr
\l{3}&\rat(9/5){6}{}&\rat(6/5){2}{}&\int(0){\0}{}&\ln{1}&\cnt{42}&\bounds{6}{6}{1}{false}&&\cr
\l{5}&\rat(21/10){7}{}&\rat(2/5){9}{}&\rat(1/2){2}{}&\ln{1}&\cnt{36}&\bounds{6}{6}{1}{false}&&\cr
\l{7}&\rat(7/5){8}{}&\rat(8/5){6}{}&\int(0){\0}{}&\ln{1}&\cnt{45}&\bounds{9}{9}{1}{false}&&\cr
\l{9}&\rat(7/5){8}{}&\rat(3/5){1}{}&\int(1){3}{}&\ln{1}&\cnt{72}&\bounds{12}{12}{1}{false}&&\cr
\endorbit
\orbit{3}{26}{50}
\h{26}&\int(6){\6}{\spec{1}}&\int(0){\0}{}&\int(0){\0}{}&&\cnt{384}&\bounds{50}{50}{1}{false}&&\cr
\l{1}&\int(3){\4}{}&\int(0){\0}{}&\int(0){\0}{}&\ln{1\sdd}&\cnt{84}&\bounds{14}{14}{1}{false}&&\cr
\l{2}&\int(3){\2}{}&\int(0){\2}{}&\int(0){\0}{}&\ln{1\sdd}&\cnt{180}&\bounds{20}{20}{1}{false}&&\cr
\l{3}&\int(3){\2}{}&\int(0){\0}{}&\int(0){\2}{}&\ln{1\sdd}&\cnt{120}&\bounds{16}{16}{1}{false}&&\cr
\endorbit
\orbit{3}{27}{48}
\h{27}&\int(0){\0}{}&\int(0){\0}{}&\int(6){\6}{\spec{1}}&&\cnt{384}&\bounds{48}{48}{1}{false}&&\cr
\l{1}&\int(0){\2}{}&\int(0){\0}{}&\int(3){\2}{}&\ln{2\sdd}&\cnt{180}&\bounds{20}{20}{1}{false}&&\cr
\l{2}&\int(0){\0}{}&\int(0){\0}{}&\int(3){\4}{}&\ln{1\sdd}&\cnt{24}&\bounds{8}{8}{1}{false}&&\cr
\endorbit
\endorbits

\newcs{max-4D6}{120}%
\newcs{count-4D6}{10}%
\newcs{counts-4D6}{\DO{1}{120}\DO{2}{112}\DO{3}{110}\DO{4}{108}\DO{5}{104}\DO{6}{100}\DO{7}{88}\DO{8}{80}\DO{9}{72}\DO{10}{56}}%
\orbits
\orbit{4}{1}{120}
\h{1}&\int(4){\4}{2}&\int(2){\2}{}&\int(0){\0}{}&\int(0){\0}{}&&\cnt{960}&\bounds{120}{144}{fail}{false}&&\cr
\l{1}&\int(2){\2}{}&\int(1){\2}{}&\int(0){\0}{}&\int(0){\0}{}&\ln{1\sdd}&\cnt{160}&\bounds{24}{32}{2}{false}&&\cr
\l{2}&\int(2){2}{}&\int(1){3}{}&\int(0){\0}{}&\int(0){1}{}&\ln{2\sdd}&\cnt{256}&\bounds{32}{32}{1}{false}&&\cr
\l{3}&\int(2){2}{}&\int(1){2}{}&\int(0){2}{}&\int(0){2}{}&\ln{1\sdd}&\cnt{288}&\bounds{32}{48}{2}{false}&&\cr
\endorbit
\orbit{4}{2}{112}
\h{2}&\rat(7/2){3}{}&\int(1){2}{}&\rat(3/2){1}{}&\int(0){\0}{}&&\cnt{672}&\bounds{112}{112}{1}{false}&&\cr
\l{1}&\int(3){\4}{}&\int(0){\0}{}&\int(0){\0}{}&\int(0){\0}{}&\ln{1}&\cnt{10}&\bounds{2}{2}{1}{false}&&\cr
\l{3}&\int(3){\4}{2}&\int(0){\0}{}&\int(0){\0}{}&\int(0){\0}{}&\ln{1}&\cnt{1}&\bounds{1}{1}{1}{false}&&\cr
\l{5}&\int(2){\2}{}&\int(1){\2}{}&\int(0){\0}{}&\int(0){\0}{}&\ln{1}&\cnt{50}&\bounds{10}{10}{1}{false}&&\cr
\l{7}&\int(2){\2}{}&\int(0){\0}{}&\int(1){\2}{}&\int(0){\0}{}&\ln{1}&\cnt{75}&\bounds{15}{15}{1}{false}&&\cr
\l{9}&\int(2){1}{}&\rat(1/2){3}{}&\rat(1/2){2}{}&\int(0){\0}{}&\ln{1}&\cnt{96}&\bounds{12}{12}{1}{false}&&\cr
\l{11}&\int(2){1}{}&\int(1){2}{}&\int(0){\0}{}&\int(0){3}{}&\ln{1}&\cnt{32}&\bounds{4}{4}{1}{false}&&\cr
\l{13}&\int(2){1}{}&\int(0){\0}{}&\int(1){3}{}&\int(0){2}{}&\ln{1}&\cnt{72}&\bounds{12}{12}{1}{false}&&\cr
\endorbit
\orbit{4}{3}{110}
\h{3}&\int(4){\4}{}&\int(2){\2}{}&\int(0){\0}{}&\int(0){\0}{}&&\cnt{768}&\bounds{110}{128}{fail}{false}&&\cr
\l{1}&\int(2){\2}{}&\int(1){\2}{}&\int(0){\0}{}&\int(0){\0}{}&\ln{1\sdd}&\cnt{96}&\bounds{18}{24}{2}{false}&&\cr
\l{2}&\int(2){3}{}&\int(1){1}{}&\int(0){\0}{}&\int(0){2}{}&\ln{2\sdd}&\cnt{192}&\bounds{26}{32}{2}{false}&&\cr
\l{3}&\int(2){3}{}&\int(1){2}{}&\int(0){1}{}&\int(0){\0}{}&\ln{2\sdd}&\cnt{128}&\bounds{16}{16}{1}{false}&&\cr
\l{4}&\int(3){\4}{}&\int(0){\0}{}&\int(0){\0}{}&\int(0){\0}{}&\ln{1}&\cnt{16}&\bounds{4}{4}{1}{false}&&\cr
\endorbit
\orbit{4}{4}{108}
\h{4}&\int(6){\6}{3}&\int(0){\0}{}&\int(0){\0}{}&\int(0){\0}{}&&\cnt{1152}&\bounds{108}{144}{fail}{false}&&\cr
\l{1}&\int(3){3}{}&\int(0){1}{}&\int(0){\0}{}&\int(0){2}{}&\ln{3\sdd}&\cnt{384}&\bounds{36}{48}{2}{false}&&\cr
\endorbit
\orbit{4}{5}{104}
\h{5}&\rat(3/2){3}{}&\int(3){2}{}&\rat(3/2){1}{}&\int(0){\0}{}&&\cnt{576}&\bounds{104}{104}{1}{false}&&\cr
\l{1}&\rat(3/2){3}{}&\rat(3/2){1}{}&\int(0){\0}{}&\int(0){2}{}&\ln{2\sd}&\cnt{48}&\bounds{8}{8}{1}{false}&&\cr
\l{2}&\int(1){1}{}&\rat(3/2){3}{}&\rat(1/2){2}{}&\int(0){\0}{}&\ln{2\sd}&\cnt{144}&\bounds{24}{24}{1}{false}&&\cr
\l{3}&\int(1){\2}{}&\int(2){\2}{}&\int(0){\0}{}&\int(0){\0}{}&\ln{2}&\cnt{45}&\bounds{9}{9}{1}{false}&&\cr
\l{5}&\rat(3/2){3}{}&\int(0){2}{}&\rat(3/2){1}{}&\int(0){\0}{}&\ln{1}&\cnt{6}&\bounds{2}{2}{1}{false}&&\cr
\endorbit
\orbit{4}{6}{100}
\h{6}&\int(2){\2}{}&\int(2){\2}{}&\int(2){\2}{}&\int(0){\0}{}&&\cnt{576}&\bounds{100}{136}{fail}{false}&&\cr
\l{1}&\int(2){\2}{}&\int(1){\2}{}&\int(0){\0}{}&\int(0){\0}{}&\ln{6\sd}&\cnt{16}&\bounds{4}{4}{1}{false}&&\cr
\l{2}&\int(1){3}{}&\int(1){2}{}&\int(1){1}{}&\int(0){\0}{}&\ln{3\sdd}&\cnt{128}&\bounds{20}{32}{2}{false}&&\cr
\l{3}&\int(1){2}{}&\int(1){2}{}&\int(1){2}{}&\int(0){2}{}&\ln{1\sdd}&\cnt{96}&\bounds{16}{16}{1}{false}&&\cr
\endorbit
\orbit{4}{7}{88}
\h{7}&\int(3){2}{}&\int(1){2}{}&\int(1){2}{}&\int(1){2}{}&&\cnt{576}&\bounds{88}{88}{1}{false}&&\cr
\l{1}&\rat(3/2){3}{}&\rat(1/2){1}{}&\int(0){\0}{}&\int(1){2}{}&\ln{6\sd}&\cnt{64}&\bounds{8}{8}{1}{false}&&\cr
\l{2}&\int(3){\4}{}&\int(0){\0}{}&\int(0){\0}{}&\int(0){\0}{}&\ln{1}&\cnt{6}&\bounds{2}{2}{1}{false}&&\cr
\l{4}&\int(2){\2}{}&\int(1){\2}{}&\int(0){\0}{}&\int(0){\0}{}&\ln{3}&\cnt{30}&\bounds{6}{6}{1}{false}&&\cr
\endorbit
\orbit{4}{8}{80}
\h{8}&\rat(3/2){3}{}&\int(1){2}{}&\rat(3/2){1}{}&\int(2){\2}{}&&\cnt{480}&\bounds{80}{80}{1}{false}&&\cr
\l{1}&\rat(3/2){3}{}&\rat(1/2){1}{}&\int(0){\0}{}&\int(1){2}{}&\ln{2\sd}&\cnt{32}&\bounds{4}{4}{1}{false}&&\cr
\l{2}&\int(1){\2}{}&\int(0){\0}{}&\int(0){\0}{}&\int(2){\2}{}&\ln{2}&\cnt{15}&\bounds{3}{3}{1}{false}&&\cr
\l{4}&\rat(3/2){3}{}&\int(0){2}{}&\rat(3/2){1}{}&\int(0){\0}{}&\ln{1}&\cnt{10}&\bounds{2}{2}{1}{false}&&\cr
\l{6}&\rat(3/2){3}{}&\int(0){\0}{}&\rat(1/2){2}{}&\int(1){1}{}&\ln{2}&\cnt{48}&\bounds{8}{8}{1}{false}&&\cr
\l{8}&\int(1){1}{}&\int(0){\0}{}&\int(1){3}{}&\int(1){2}{}&\ln{1}&\cnt{72}&\bounds{12}{12}{1}{false}&&\cr
\endorbit
\orbit{4}{9}{72}
\h{9}&\rat(3/2){3}{}&\rat(3/2){3}{}&\rat(3/2){3}{}&\rat(3/2){3}{}&&\cnt{432}&\bounds{72}{72}{1}{false}&&\cr
\l{1}&\rat(3/2){3}{}&\int(1){1}{}&\int(0){\0}{}&\rat(1/2){2}{}&\ln{12\sd}&\cnt{36}&\bounds{6}{6}{1}{false}&&\cr
\endorbit
\orbit{4}{10}{56}
\h{10}&\int(6){\6}{\spec{1}}&\int(0){\0}{}&\int(0){\0}{}&\int(0){\0}{}&&\cnt{384}&\bounds{56}{56}{1}{false}&&\cr
\l{1}&\int(3){\4}{}&\int(0){\0}{}&\int(0){\0}{}&\int(0){\0}{}&\ln{1\sdd}&\cnt{24}&\bounds{8}{8}{1}{false}&&\cr
\l{2}&\int(3){\2}{}&\int(0){\2}{}&\int(0){\0}{}&\int(0){\0}{}&\ln{3\sdd}&\cnt{120}&\bounds{16}{16}{1}{false}&&\cr
\endorbit
\endorbits

\newcs{max-3A8}{120}%
\newcs{count-3A8}{15}%
\newcs{counts-3A8}{\DO{1}{120}\DO{2}{104}\DO{3}{100}\DO{4}{96}\DO{5}{96}\DO{6}{96}\DO{7}{92}\DO{8}{88}\DO{9}{88}\DO{10}{86}\DO{11}{84}\DO{12}{84}\DO{13}{82}\DO{14}{78}\DO{15}{48}}%
\orbits
\orbit{3}{1}{120}
\h{1}&\int(6){\6}{}&\int(0){\0}{}&\int(0){\0}{}&&\cnt{876}&\bounds{120}{120}{1}{false}&&\cr
\l{1}&\int(3){\4}{\spec{1}}&\int(0){\0}{}&\int(0){\0}{}&\ln{2\sd}&\cnt{27}&\bounds{9}{9}{1}{false}&&\cr
\l{2}&\int(3){5}{}&\int(0){8}{}&\int(0){8}{}&\ln{2\sd}&\cnt{243}&\bounds{27}{27}{1}{false}&&\cr
\l{3}&\int(3){6}{}&\int(0){3}{}&\int(0){\0}{}&\ln{4\sd}&\cnt{84}&\bounds{12}{12}{1}{false}&&\cr
\endorbit
\orbit{3}{2}{104}
\h{2}&\rat(8/9){8}{}&\rat(32/9){2}{\spec{2}}&\rat(14/9){2}{}&&\cnt{828}&\bounds{104}{104}{1}{false}&&\cr
\l{1}&\rat(4/9){4}{}&\rat(16/9){1}{}&\rat(7/9){1}{}&\ln{1\sdd}&\cnt{140}&\bounds{18}{18}{1}{false}&&\cr
\l{2}&\int(1){\2}{}&\int(2){\2}{}&\int(0){\0}{}&\ln{1}&\cnt{64}&\bounds{8}{8}{1}{false}&&\cr
\l{4}&\rat(7/9){7}{}&\rat(16/9){1}{}&\rat(4/9){7}{}&\ln{1}&\cnt{168}&\bounds{21}{21}{1}{false}&&\cr
\l{6}&\rat(8/9){8}{}&\rat(14/9){2}{}&\rat(5/9){2}{}&\ln{1}&\cnt{112}&\bounds{14}{14}{1}{false}&&\cr
\endorbit
\orbit{3}{3}{100}
\h{3}&\rat(8/9){8}{}&\rat(8/9){8}{}&\rat(38/9){5}{\spec{1}}&&\cnt{732}&\bounds{100}{100}{1}{false}&&\cr
\l{1}&\rat(7/9){7}{}&\rat(1/9){1}{}&\rat(19/9){7}{}&\ln{2\sd}&\cnt{128}&\bounds{16}{16}{1}{false}&&\cr
\l{2}&\rat(5/9){5}{}&\rat(8/9){8}{}&\rat(14/9){8}{}&\ln{2}&\cnt{56}&\bounds{8}{8}{1}{false}&&\cr
\l{4}&\int(1){\2}{}&\int(0){\0}{}&\int(2){\2}{}&\ln{2}&\cnt{48}&\bounds{6}{6}{1}{false}&&\cr
\l{6}&\rat(8/9){8}{}&\rat(8/9){8}{}&\rat(11/9){5}{}&\ln{1}&\cnt{30}&\bounds{6}{6}{1}{false}&&\cr
\endorbit
\orbit{3}{4}{96}
\h{4}&\int(2){\2}{}&\int(2){\2}{}&\int(2){\2}{}&&\cnt{588}&\bounds{96}{96}{1}{false}&&\cr
\l{1}&\int(1){5}{}&\int(1){8}{}&\int(1){8}{}&\ln{6\sd}&\cnt{35}&\bounds{7}{7}{1}{false}&&\cr
\l{2}&\int(2){\2}{}&\int(1){\2}{\spec{1}}&\int(0){\0}{}&\ln{12\sd}&\cnt{7}&\bounds{1}{1}{1}{false}&&\cr
\l{3}&\int(1){7}{}&\int(1){7}{}&\int(1){1}{}&\ln{6\sd}&\cnt{49}&\bounds{7}{7}{1}{false}&&\cr
\endorbit
\orbit{3}{5}{96}
\h{5}&\int(4){6}{\spec{1}}&\int(2){3}{}&\int(0){\0}{}&&\cnt{780}&\bounds{96}{96}{1}{false}&&\cr
\l{1}&\int(3){\4}{}&\int(0){\0}{}&\int(0){\0}{}&\ln{1}&\cnt{21}&\bounds{3}{3}{1}{false}&&\cr
\l{3}&\int(2){6}{}&\int(1){3}{}&\int(0){\0}{}&\ln{1}&\cnt{126}&\bounds{18}{18}{1}{false}&&\cr
\l{5}&\rat(7/3){7}{}&\rat(2/3){7}{}&\int(0){1}{}&\ln{1}&\cnt{135}&\bounds{15}{15}{1}{false}&&\cr
\l{7}&\rat(7/3){7}{}&\rat(2/3){1}{}&\int(0){7}{}&\ln{1}&\cnt{108}&\bounds{12}{12}{1}{false}&&\cr
\endorbit
\orbit{3}{6}{96}
\h{6}&\rat(8/9){8}{}&\rat(8/9){8}{}&\rat(38/9){5}{\spec{2}}&&\cnt{684}&\bounds{96}{96}{1}{false}&&\cr
\l{1}&\int(1){\2}{}&\int(0){\0}{}&\int(2){\2}{}&\ln{2}&\cnt{40}&\bounds{5}{5}{1}{false}&&\cr
\l{3}&\rat(2/3){6}{}&\int(0){\0}{}&\rat(7/3){3}{}&\ln{2}&\cnt{84}&\bounds{12}{12}{1}{false}&&\cr
\l{5}&\rat(7/9){7}{}&\rat(7/9){7}{}&\rat(13/9){1}{}&\ln{1}&\cnt{64}&\bounds{8}{8}{1}{false}&&\cr
\l{7}&\rat(8/9){8}{}&\rat(8/9){8}{}&\rat(11/9){5}{}&\ln{1}&\cnt{30}&\bounds{6}{6}{1}{false}&&\cr
\endorbit
\orbit{3}{7}{92}
\h{7}&\int(4){6}{\spec{2}}&\int(2){3}{}&\int(0){\0}{}&&\cnt{636}&\bounds{92}{96}{fail}{false}&&\cr
\l{1}&\int(2){3}{}&\int(1){6}{}&\int(0){\0}{}&\ln{1\sdd}&\cnt{120}&\bounds{20}{24}{2}{false}&&\cr
\l{2}&\int(3){\4}{}&\int(0){\0}{}&\int(0){\0}{}&\ln{1}&\cnt{24}&\bounds{6}{6}{1}{false}&&\cr
\l{4}&\rat(8/3){5}{}&\rat(1/3){8}{}&\int(0){8}{}&\ln{1}&\cnt{54}&\bounds{6}{6}{1}{false}&&\cr
\l{6}&\rat(7/3){4}{}&\rat(2/3){1}{}&\int(0){1}{}&\ln{1}&\cnt{108}&\bounds{12}{12}{1}{false}&&\cr
\l{8}&\int(2){6}{}&\int(1){3}{}&\int(0){\0}{}&\ln{1}&\cnt{72}&\bounds{12}{12}{1}{false}&&\cr
\endorbit
\orbit{3}{8}{88}
\h{8}&\int(4){\4}{}&\int(2){\2}{}&\int(0){\0}{}&&\cnt{684}&\bounds{88}{88}{1}{false}&&\cr
\l{1}&\int(2){5}{}&\int(1){8}{}&\int(0){8}{}&\ln{2\sd}&\cnt{90}&\bounds{10}{10}{1}{false}&&\cr
\l{2}&\int(2){6}{}&\int(1){3}{}&\int(0){\0}{}&\ln{2\sd}&\cnt{105}&\bounds{15}{15}{1}{false}&&\cr
\l{3}&\int(2){\2}{}&\int(1){\2}{\spec{1}}&\int(0){\0}{}&\ln{2\sd}&\cnt{28}&\bounds{4}{4}{1}{false}&&\cr
\l{4}&\int(2){7}{}&\int(1){7}{}&\int(0){1}{}&\ln{2\sd}&\cnt{63}&\bounds{7}{7}{1}{false}&&\cr
\l{5}&\int(2){7}{}&\int(1){1}{}&\int(0){7}{}&\ln{2\sd}&\cnt{36}&\bounds{4}{4}{1}{false}&&\cr
\l{6}&\int(3){\4}{\spec{1}}&\int(0){\0}{}&\int(0){\0}{}&\ln{2}&\cnt{10}&\bounds{2}{2}{1}{false}&&\cr
\endorbit
\orbit{3}{9}{88}
\h{9}&\rat(26/9){8}{\spec{2}}&\rat(8/9){8}{}&\rat(20/9){5}{}&&\cnt{540}&\bounds{88}{88}{1}{false}&&\cr
\l{1}&\int(3){\4}{}&\int(0){\0}{}&\int(0){\0}{}&\ln{1}&\cnt{6}&\bounds{1}{1}{1}{false}&&\cr
\l{3}&\rat(14/9){5}{}&\rat(8/9){8}{}&\rat(5/9){8}{}&\ln{1}&\cnt{60}&\bounds{12}{12}{1}{false}&&\cr
\l{5}&\int(2){\2}{}&\int(1){\2}{}&\int(0){\0}{}&\ln{1}&\cnt{16}&\bounds{2}{2}{1}{false}&&\cr
\l{7}&\int(2){\2}{}&\int(0){\0}{}&\int(1){\2}{}&\ln{1}&\cnt{40}&\bounds{8}{8}{1}{false}&&\cr
\l{9}&\rat(5/3){6}{}&\int(0){\0}{}&\rat(4/3){3}{}&\ln{1}&\cnt{60}&\bounds{10}{10}{1}{false}&&\cr
\l{11}&\rat(16/9){7}{}&\rat(7/9){7}{}&\rat(4/9){1}{}&\ln{1}&\cnt{40}&\bounds{5}{5}{1}{false}&&\cr
\l{13}&\rat(16/9){7}{}&\rat(1/9){1}{}&\rat(10/9){7}{}&\ln{1}&\cnt{48}&\bounds{6}{6}{1}{false}&&\cr
\endorbit
\orbit{3}{10}{86}
\h{10}&\rat(8/9){8}{}&\rat(32/9){2}{\spec{1}}&\rat(14/9){2}{}&&\cnt{588}&\bounds{86}{86}{1}{false}&&\cr
\l{1}&\int(1){\2}{}&\int(2){\2}{}&\int(0){\0}{}&\ln{1}&\cnt{24}&\bounds{3}{3}{1}{false}&&\cr
\l{3}&\rat(2/3){6}{}&\rat(7/3){3}{}&\int(0){\0}{}&\ln{1}&\cnt{28}&\bounds{4}{4}{1}{false}&&\cr
\l{5}&\rat(7/9){7}{}&\rat(13/9){7}{}&\rat(7/9){1}{}&\ln{1}&\cnt{80}&\bounds{10}{10}{1}{false}&&\cr
\l{7}&\rat(8/9){8}{}&\rat(17/9){5}{}&\rat(2/9){8}{}&\ln{1}&\cnt{70}&\bounds{10}{10}{1}{false}&&\cr
\l{9}&\rat(8/9){8}{}&\rat(14/9){2}{}&\rat(5/9){2}{}&\ln{1}&\cnt{42}&\bounds{6}{6}{1}{false}&&\cr
\l{11}&\rat(8/9){8}{}&\rat(5/9){2}{}&\rat(14/9){2}{}&\ln{1}&\cnt{15}&\bounds{3}{3}{1}{false}&&\cr
\l{13}&\rat(8/9){8}{}&\rat(11/9){8}{}&\rat(8/9){5}{}&\ln{1}&\cnt{35}&\bounds{7}{7}{1}{false}&&\cr
\endorbit
\orbit{3}{11}{84}
\h{11}&\int(2){6}{}&\int(2){3}{}&\int(2){\2}{}&&\cnt{492}&\bounds{84}{84}{1}{false}&&\cr
\l{1}&\rat(5/3){5}{}&\rat(1/3){8}{}&\int(1){8}{}&\ln{2\sd}&\cnt{36}&\bounds{6}{6}{1}{false}&&\cr
\l{2}&\int(2){6}{}&\int(0){\0}{}&\int(1){3}{}&\ln{2\sd}&\cnt{21}&\bounds{3}{3}{1}{false}&&\cr
\l{3}&\rat(4/3){7}{}&\rat(2/3){1}{}&\int(1){7}{}&\ln{2\sd}&\cnt{63}&\bounds{9}{9}{1}{false}&&\cr
\l{4}&\rat(4/3){4}{}&\rat(2/3){1}{}&\int(1){1}{}&\ln{2}&\cnt{45}&\bounds{9}{9}{1}{false}&&\cr
\l{6}&\int(2){6}{}&\int(1){3}{}&\int(0){\0}{}&\ln{2}&\cnt{18}&\bounds{3}{3}{1}{false}&&\cr
\endorbit
\orbit{3}{12}{84}
\h{12}&\rat(14/9){7}{}&\rat(20/9){4}{}&\rat(20/9){4}{}&&\cnt{444}&\bounds{84}{92}{fail}{false}&&\cr
\l{1}&\rat(7/9){8}{}&\rat(16/9){5}{}&\rat(4/9){8}{}&\ln{2\sd}&\cnt{50}&\bounds{10}{10}{1}{false}&&\cr
\l{2}&\rat(7/9){8}{}&\rat(10/9){2}{}&\rat(10/9){2}{}&\ln{1\sdd}&\cnt{72}&\bounds{16}{24}{2}{false}&&\cr
\l{3}&\rat(4/3){6}{}&\rat(5/3){3}{}&\int(0){\0}{}&\ln{2}&\cnt{28}&\bounds{4}{4}{1}{false}&&\cr
\l{5}&\rat(14/9){7}{}&\rat(8/9){7}{}&\rat(5/9){1}{}&\ln{2}&\cnt{40}&\bounds{8}{8}{1}{false}&&\cr
\endorbit
\orbit{3}{13}{82}
\h{13}&\rat(26/9){8}{\spec{2}}&\rat(14/9){2}{}&\rat(14/9){2}{}&&\cnt{540}&\bounds{82}{82}{1}{false}&&\cr
\l{1}&\rat(13/9){4}{}&\rat(7/9){1}{}&\rat(7/9){1}{}&\ln{1\sdd}&\cnt{80}&\bounds{16}{16}{1}{false}&&\cr
\l{2}&\int(3){\4}{}&\int(0){\0}{}&\int(0){\0}{}&\ln{1}&\cnt{6}&\bounds{1}{1}{1}{false}&&\cr
\l{4}&\int(2){\2}{}&\int(1){\2}{}&\int(0){\0}{}&\ln{2}&\cnt{28}&\bounds{4}{4}{1}{false}&&\cr
\l{6}&\rat(5/3){6}{}&\rat(4/3){3}{}&\int(0){\0}{}&\ln{2}&\cnt{42}&\bounds{6}{6}{1}{false}&&\cr
\l{8}&\rat(16/9){7}{}&\rat(4/9){7}{}&\rat(7/9){1}{}&\ln{2}&\cnt{42}&\bounds{6}{6}{1}{false}&&\cr
\endorbit
\orbit{3}{14}{78}
\h{14}&\int(2){6}{}&\int(2){6}{}&\int(2){6}{}&&\cnt{444}&\bounds{78}{78}{1}{false}&&\cr
\l{1}&\int(2){6}{}&\int(1){3}{}&\int(0){\0}{}&\ln{6\sd}&\cnt{20}&\bounds{4}{4}{1}{false}&&\cr
\l{2}&\rat(5/3){5}{}&\rat(2/3){8}{}&\rat(2/3){8}{}&\ln{3}&\cnt{54}&\bounds{9}{9}{1}{false}&&\cr
\endorbit
\orbit{3}{15}{48}
\h{15}&\int(6){\6}{\spec{1}}&\int(0){\0}{}&\int(0){\0}{}&&\cnt{348}&\bounds{48}{48}{1}{false}&&\cr
\l{1}&\int(3){\4}{}&\int(0){\0}{}&\int(0){\0}{}&\ln{1\sdd}&\cnt{60}&\bounds{12}{12}{1}{false}&&\cr
\l{2}&\int(3){\2}{}&\int(0){\2}{}&\int(0){\0}{}&\ln{2\sdd}&\cnt{144}&\bounds{18}{18}{1}{false}&&\cr
\endorbit
\endorbits

\newcs{max-2A7_2D5}{120}%
\newcs{count-2A7_2D5}{34}%
\newcs{counts-2A7_2D5}{\DO{1}{118}\DO{2}{110}\DO{3}{120}\DO{4}{118}\DO{5}{116}\DO{6}{112}\DO{7}{112}\DO{8}{112}\DO{9}{112}\DO{10}{112}\DO{11}{112}\DO{12}{110}\DO{13}{108}\DO{14}{108}\DO{15}{108}\DO{16}{108}\DO{17}{108}\DO{18}{108}\DO{19}{106}\DO{20}{106}\DO{21}{106}\DO{22}{106}\DO{23}{104}\DO{24}{102}\DO{25}{100}\DO{26}{100}\DO{27}{100}\DO{28}{96}\DO{29}{92}\DO{30}{92}\DO{31}{92}\DO{32}{88}\DO{33}{50}\DO{34}{46}}%
\orbits
\orbit{4}{1}{118}
\h{1}&\int(6){\6}{}&\int(0){\0}{}&\int(0){\0}{}&\int(0){\0}{}&&\cnt{888}&\bounds{118}{142}{fail}{false}&&\cr
\l{1}&\int(3){\4}{\spec{1}}&\int(0){\0}{}&\int(0){\0}{}&\int(0){\0}{}&\ln{2\sd}&\cnt{18}&\bounds{5}{10}{5}{true}&&\cr
\l{2}&\int(3){4}{}&\int(0){4}{}&\int(0){\0}{}&\int(0){\0}{}&\ln{1\sdd}&\cnt{140}&\bounds{18}{18}{1}{false}&&\cr
\l{3}&\int(3){4}{}&\int(0){\0}{}&\int(0){2}{}&\int(0){2}{}&\ln{1\sdd}&\cnt{200}&\bounds{26}{40}{2}{false}&&\cr
\l{4}&\int(3){5}{}&\int(0){7}{}&\int(0){\0}{}&\int(0){3}{}&\ln{4\sd}&\cnt{128}&\bounds{16}{16}{1}{false}&&\cr
\endorbit
\orbit{4}{2}{110}
\h{2}&\int(2){\2}{}&\int(0){\0}{}&\int(4){\4}{2}&\int(0){\0}{}&&\cnt{888}&\bounds{110}{136}{fail}{false}&&\cr
\l{1}&\int(1){4}{}&\int(0){\0}{}&\int(2){2}{}&\int(0){2}{}&\ln{1\sdd}&\cnt{200}&\bounds{26}{40}{2}{false}&&\cr
\l{2}&\int(1){\2}{\spec{1}}&\int(0){\0}{}&\int(2){\2}{}&\int(0){\0}{}&\ln{2\sd}&\cnt{48}&\bounds{8}{8}{1}{false}&&\cr
\l{3}&\int(1){6}{}&\int(0){6}{}&\int(2){2}{}&\int(0){\0}{}&\ln{2\sd}&\cnt{168}&\bounds{18}{24}{3}{false}&&\cr
\l{4}&\int(1){7}{}&\int(0){1}{}&\int(2){2}{}&\int(0){3}{}&\ln{2\sd}&\cnt{128}&\bounds{16}{16}{1}{false}&&\cr
\endorbit
\orbit{4}{3}{120}
\h{3}&\int(2){4}{}&\int(2){4}{}&\int(2){\2}{}&\int(0){\0}{}&&\cnt{504}&\bounds{120}{136}{fail}{false}&&\cr
\l{1}&\int(2){4}{}&\int(0){\0}{}&\int(1){2}{}&\int(0){2}{}&\ln{2\sd}&\cnt{20}&\bounds{4}{4}{1}{false}&&\cr
\l{2}&\rat(3/2){5}{}&\rat(1/2){1}{}&\int(1){3}{}&\int(0){\0}{}&\ln{4\sd}&\cnt{64}&\bounds{16}{16}{1}{false}&&\cr
\l{3}&\int(1){6}{}&\int(1){6}{}&\int(1){2}{}&\int(0){\0}{}&\ln{2\sdd}&\cnt{72}&\bounds{16}{24}{2}{false}&&\cr
\l{4}&\int(2){4}{}&\int(1){4}{}&\int(0){\0}{}&\int(0){\0}{}&\ln{2}&\cnt{16}&\bounds{4}{4}{1}{false}&&\cr
\endorbit
\orbit{4}{4}{118}
\h{4}&\int(0){\0}{}&\int(0){\0}{}&\int(4){\4}{2}&\int(2){\2}{}&&\cnt{888}&\bounds{118}{124}{fail}{false}&&\cr
\l{1}&\int(0){4}{}&\int(0){\0}{}&\int(2){2}{}&\int(1){2}{}&\ln{2\sdd}&\cnt{140}&\bounds{18}{18}{1}{false}&&\cr
\l{2}&\int(0){7}{}&\int(0){1}{}&\int(2){2}{}&\int(1){3}{}&\ln{2\sd}&\cnt{256}&\bounds{32}{32}{1}{false}&&\cr
\l{3}&\int(0){\0}{}&\int(0){\0}{}&\int(2){\2}{}&\int(1){\2}{}&\ln{1\sdd}&\cnt{96}&\bounds{18}{24}{2}{false}&&\cr
\endorbit
\orbit{4}{5}{116}
\h{5}&\int(0){\0}{}&\int(0){\0}{}&\int(4){\4}{}&\int(2){\2}{}&&\cnt{792}&\bounds{116}{124}{fail}{false}&&\cr
\l{1}&\int(0){6}{}&\int(0){\0}{}&\int(2){1}{}&\int(1){1}{}&\ln{4\sd}&\cnt{112}&\bounds{16}{16}{1}{false}&&\cr
\l{2}&\int(0){7}{}&\int(0){7}{}&\int(2){3}{}&\int(1){2}{}&\ln{2\sd}&\cnt{128}&\bounds{16}{16}{1}{false}&&\cr
\l{3}&\int(0){\0}{}&\int(0){\0}{}&\int(2){\2}{}&\int(1){\2}{}&\ln{1\sdd}&\cnt{72}&\bounds{16}{24}{2}{false}&&\cr
\l{4}&\int(0){\0}{}&\int(0){\0}{}&\int(3){\4}{\spec{1}}&\int(0){\0}{}&\ln{2}&\cnt{4}&\bounds{1}{1}{1}{false}&&\cr
\endorbit
\orbit{4}{6}{112}
\h{6}&\int(2){4}{}&\int(2){\2}{}&\int(1){2}{}&\int(1){2}{}&&\cnt{504}&\bounds{112}{112}{1}{false}&&\cr
\l{1}&\int(1){6}{}&\int(1){6}{}&\int(1){2}{}&\int(0){\0}{}&\ln{4\sd}&\cnt{36}&\bounds{6}{6}{1}{false}&&\cr
\l{2}&\int(2){4}{}&\int(1){4}{}&\int(0){\0}{}&\int(0){\0}{}&\ln{1}&\cnt{20}&\bounds{4}{4}{1}{false}&&\cr
\l{4}&\int(2){4}{}&\int(0){\0}{}&\int(1){2}{}&\int(0){2}{}&\ln{2}&\cnt{8}&\bounds{2}{2}{1}{false}&&\cr
\l{6}&\int(1){4}{}&\int(0){\0}{}&\int(1){2}{}&\int(1){2}{}&\ln{1}&\cnt{16}&\bounds{4}{4}{1}{false}&&\cr
\l{8}&\rat(3/2){5}{}&\int(1){7}{}&\int(0){\0}{}&\rat(1/2){3}{}&\ln{4}&\cnt{32}&\bounds{8}{8}{1}{false}&&\cr
\endorbit
\orbit{4}{7}{112}
\h{7}&\rat(7/2){6}{\spec{1}}&\rat(3/2){6}{}&\int(1){2}{}&\int(0){\0}{}&&\cnt{792}&\bounds{112}{128}{fail}{false}&&\cr
\l{1}&\rat(7/4){7}{}&\rat(3/4){3}{}&\rat(1/2){1}{}&\int(0){\0}{}&\ln{1\sdd}&\cnt{160}&\bounds{24}{32}{2}{false}&&\cr
\l{2}&\rat(7/4){7}{}&\rat(3/4){7}{}&\rat(1/2){3}{}&\int(0){2}{}&\ln{1\sdd}&\cnt{160}&\bounds{24}{32}{2}{false}&&\cr
\l{3}&\int(2){\2}{}&\int(1){\2}{}&\int(0){\0}{}&\int(0){\0}{}&\ln{1}&\cnt{84}&\bounds{12}{12}{1}{false}&&\cr
\l{5}&\int(2){\2}{}&\int(0){\0}{}&\int(1){\2}{}&\int(0){\0}{}&\ln{1}&\cnt{56}&\bounds{8}{8}{1}{false}&&\cr
\l{7}&\rat(7/4){7}{}&\rat(5/4){5}{}&\int(0){\0}{}&\int(0){1}{}&\ln{1}&\cnt{96}&\bounds{12}{12}{1}{false}&&\cr
\endorbit
\orbit{4}{8}{112}
\h{8}&\rat(7/2){6}{\spec{2}}&\rat(3/2){6}{}&\int(1){2}{}&\int(0){\0}{}&&\cnt{600}&\bounds{112}{118}{fail}{false}&&\cr
\l{1}&\rat(7/4){3}{}&\rat(3/4){7}{}&\rat(1/2){1}{}&\int(0){\0}{}&\ln{1\sdd}&\cnt{96}&\bounds{18}{24}{2}{false}&&\cr
\l{2}&\int(3){\4}{}&\int(0){\0}{}&\int(0){\0}{}&\int(0){\0}{}&\ln{1}&\cnt{12}&\bounds{3}{3}{1}{false}&&\cr
\l{4}&\int(2){4}{}&\int(1){4}{}&\int(0){\0}{}&\int(0){\0}{}&\ln{1}&\cnt{60}&\bounds{12}{12}{1}{false}&&\cr
\l{6}&\int(2){4}{}&\int(0){\0}{}&\int(1){2}{}&\int(0){2}{}&\ln{1}&\cnt{40}&\bounds{8}{8}{1}{false}&&\cr
\l{8}&\int(2){\2}{}&\int(1){\2}{}&\int(0){\0}{}&\int(0){\0}{}&\ln{1}&\cnt{36}&\bounds{6}{6}{1}{false}&&\cr
\l{10}&\int(2){\2}{}&\int(0){\0}{}&\int(1){\2}{}&\int(0){\0}{}&\ln{1}&\cnt{24}&\bounds{6}{6}{1}{false}&&\cr
\l{12}&\rat(9/4){5}{}&\rat(3/4){7}{}&\int(0){\0}{}&\int(0){3}{}&\ln{1}&\cnt{32}&\bounds{4}{4}{1}{false}&&\cr
\l{14}&\rat(9/4){5}{}&\rat(1/4){1}{}&\rat(1/2){3}{}&\int(0){\0}{}&\ln{1}&\cnt{48}&\bounds{8}{8}{1}{false}&&\cr
\endorbit
\orbit{4}{9}{112}
\h{9}&\rat(7/8){7}{}&\rat(31/8){5}{\spec{2}}&\int(0){\0}{}&\rat(5/4){1}{}&&\cnt{648}&\bounds{112}{112}{1}{false}&&\cr
\l{1}&\int(1){\2}{}&\int(2){\2}{}&\int(0){\0}{}&\int(0){\0}{}&\ln{1}&\cnt{28}&\bounds{4}{4}{1}{false}&&\cr
\l{3}&\rat(1/2){4}{}&\rat(5/2){4}{}&\int(0){\0}{}&\int(0){\0}{}&\ln{1}&\cnt{35}&\bounds{7}{7}{1}{false}&&\cr
\l{5}&\rat(3/4){6}{}&\rat(7/4){2}{}&\int(0){\0}{}&\rat(1/2){2}{}&\ln{1}&\cnt{105}&\bounds{15}{15}{1}{false}&&\cr
\l{7}&\rat(7/8){7}{}&\rat(17/8){3}{}&\int(0){1}{}&\int(0){\0}{}&\ln{1}&\cnt{48}&\bounds{6}{6}{1}{false}&&\cr
\l{9}&\rat(7/8){7}{}&\rat(15/8){5}{}&\int(0){\0}{}&\rat(1/4){1}{}&\ln{1}&\cnt{40}&\bounds{8}{8}{1}{false}&&\cr
\l{11}&\rat(7/8){7}{}&\rat(7/8){5}{}&\int(0){\0}{}&\rat(5/4){1}{}&\ln{1}&\cnt{18}&\bounds{6}{6}{1}{false}&&\cr
\l{13}&\rat(7/8){7}{}&\rat(11/8){1}{}&\int(0){2}{}&\rat(3/4){3}{}&\ln{1}&\cnt{50}&\bounds{10}{10}{1}{false}&&\cr
\endorbit
\orbit{4}{10}{112}
\h{10}&\int(2){\2}{}&\int(2){\2}{}&\int(2){\2}{}&\int(0){\0}{}&&\cnt{600}&\bounds{112}{112}{1}{false}&&\cr
\l{1}&\int(2){\2}{}&\int(0){\0}{}&\int(1){\2}{}&\int(0){\0}{}&\ln{2\sd}&\cnt{12}&\bounds{4}{4}{1}{false}&&\cr
\l{2}&\int(1){5}{}&\int(1){1}{}&\int(1){3}{}&\int(0){\0}{}&\ln{4\sd}&\cnt{60}&\bounds{12}{12}{1}{false}&&\cr
\l{3}&\int(1){6}{}&\int(1){6}{}&\int(1){2}{}&\int(0){\0}{}&\ln{2\sd}&\cnt{72}&\bounds{12}{12}{1}{false}&&\cr
\l{4}&\int(1){7}{}&\int(1){7}{}&\int(1){3}{}&\int(0){2}{}&\ln{2\sd}&\cnt{40}&\bounds{8}{8}{1}{false}&&\cr
\l{5}&\int(1){7}{}&\int(1){1}{}&\int(1){2}{}&\int(0){3}{}&\ln{2\sd}&\cnt{32}&\bounds{4}{4}{1}{false}&&\cr
\l{6}&\int(2){\2}{}&\int(1){\2}{\spec{1}}&\int(0){\0}{}&\int(0){\0}{}&\ln{4}&\cnt{6}&\bounds{1}{1}{1}{false}&&\cr
\endorbit
\orbit{4}{11}{112}
\h{11}&\int(4){4}{\spec{1}}&\int(2){4}{}&\int(0){\0}{}&\int(0){\0}{}&&\cnt{696}&\bounds{112}{120}{fail}{false}&&\cr
\l{1}&\int(2){6}{}&\int(1){6}{}&\int(0){2}{}&\int(0){\0}{}&\ln{2\sdd}&\cnt{120}&\bounds{20}{24}{2}{false}&&\cr
\l{2}&\int(3){\4}{}&\int(0){\0}{}&\int(0){\0}{}&\int(0){\0}{}&\ln{1}&\cnt{20}&\bounds{4}{4}{1}{false}&&\cr
\l{4}&\int(2){4}{}&\int(1){4}{}&\int(0){\0}{}&\int(0){\0}{}&\ln{1}&\cnt{80}&\bounds{16}{16}{1}{false}&&\cr
\l{6}&\rat(5/2){5}{}&\rat(1/2){7}{}&\int(0){\0}{}&\int(0){3}{}&\ln{2}&\cnt{64}&\bounds{8}{8}{1}{false}&&\cr
\endorbit
\orbit{4}{12}{110}
\h{12}&\rat(7/8){7}{}&\rat(7/8){7}{}&\rat(13/4){3}{}&\int(1){2}{}&&\cnt{648}&\bounds{110}{110}{1}{false}&&\cr
\l{1}&\int(1){\2}{}&\int(0){\0}{}&\int(2){\2}{}&\int(0){\0}{}&\ln{2}&\cnt{28}&\bounds{4}{4}{1}{false}&&\cr
\l{3}&\rat(1/2){4}{}&\int(0){\0}{}&\rat(3/2){2}{}&\int(1){2}{}&\ln{2}&\cnt{35}&\bounds{7}{7}{1}{false}&&\cr
\l{5}&\rat(3/4){6}{}&\rat(3/4){6}{}&\rat(3/2){2}{}&\int(0){\0}{}&\ln{1}&\cnt{49}&\bounds{7}{7}{1}{false}&&\cr
\l{7}&\rat(3/4){6}{}&\int(0){\0}{}&\rat(7/4){1}{}&\rat(1/2){1}{}&\ln{2}&\cnt{56}&\bounds{8}{8}{1}{false}&&\cr
\l{9}&\rat(7/8){7}{}&\rat(7/8){7}{}&\rat(5/4){3}{}&\int(0){2}{}&\ln{1}&\cnt{32}&\bounds{8}{8}{1}{false}&&\cr
\l{11}&\rat(7/8){7}{}&\rat(7/8){7}{}&\rat(1/4){3}{\spec{1}}&\int(1){2}{}&\ln{1}&\cnt{1}&\bounds{1}{1}{1}{false}&&\cr
\l{13}&\rat(7/8){7}{}&\rat(7/8){7}{}&\rat(1/4){3}{\spec{2}}&\int(1){2}{}&\ln{1}&\cnt{4}&\bounds{1}{1}{1}{false}&&\cr
\endorbit
\orbit{4}{13}{108}
\h{13}&\int(2){4}{}&\int(0){\0}{}&\int(3){2}{}&\int(1){2}{}&&\cnt{600}&\bounds{108}{120}{fail}{false}&&\cr
\l{1}&\int(1){6}{}&\int(0){\0}{}&\rat(3/2){1}{}&\rat(1/2){1}{}&\ln{2\sdd}&\cnt{96}&\bounds{18}{24}{2}{false}&&\cr
\l{2}&\int(2){4}{}&\int(0){\0}{}&\int(1){2}{}&\int(0){2}{}&\ln{1}&\cnt{24}&\bounds{6}{6}{1}{false}&&\cr
\l{4}&\int(2){4}{}&\int(0){\0}{}&\int(0){2}{}&\int(1){2}{}&\ln{1}&\cnt{4}&\bounds{2}{2}{1}{false}&&\cr
\l{6}&\int(1){4}{}&\int(0){\0}{}&\int(1){2}{}&\int(1){2}{}&\ln{1}&\cnt{48}&\bounds{12}{12}{1}{false}&&\cr
\l{8}&\rat(3/2){5}{}&\int(0){1}{}&\rat(3/2){3}{}&\int(0){\0}{}&\ln{2}&\cnt{64}&\bounds{8}{8}{1}{false}&&\cr
\endorbit
\orbit{4}{14}{108}
\h{14}&\int(4){\4}{}&\int(2){\2}{}&\int(0){\0}{}&\int(0){\0}{}&&\cnt{696}&\bounds{108}{112}{fail}{false}&&\cr
\l{1}&\int(2){4}{}&\int(1){4}{}&\int(0){\0}{}&\int(0){\0}{}&\ln{1\sdd}&\cnt{120}&\bounds{20}{24}{2}{false}&&\cr
\l{2}&\int(2){\2}{}&\int(1){\2}{\spec{1}}&\int(0){\0}{}&\int(0){\0}{}&\ln{2\sd}&\cnt{24}&\bounds{4}{4}{1}{false}&&\cr
\l{3}&\int(2){5}{}&\int(1){7}{}&\int(0){\0}{}&\int(0){3}{}&\ln{4\sd}&\cnt{64}&\bounds{8}{8}{1}{false}&&\cr
\l{4}&\int(2){6}{}&\int(1){6}{}&\int(0){2}{}&\int(0){\0}{}&\ln{4\sd}&\cnt{60}&\bounds{10}{10}{1}{false}&&\cr
\l{5}&\int(3){\4}{\spec{1}}&\int(0){\0}{}&\int(0){\0}{}&\int(0){\0}{}&\ln{2}&\cnt{8}&\bounds{2}{2}{1}{false}&&\cr
\endorbit
\orbit{4}{15}{108}
\h{15}&\int(4){4}{\spec{1}}&\int(0){\0}{}&\int(1){2}{}&\int(1){2}{}&&\cnt{696}&\bounds{108}{120}{fail}{false}&&\cr
\l{1}&\int(2){6}{}&\int(0){6}{}&\int(1){2}{}&\int(0){\0}{}&\ln{2\sd}&\cnt{56}&\bounds{8}{8}{1}{false}&&\cr
\l{2}&\int(2){6}{}&\int(0){\0}{}&\rat(1/2){1}{}&\rat(1/2){1}{}&\ln{1\sdd}&\cnt{128}&\bounds{20}{32}{2}{false}&&\cr
\l{3}&\int(3){\4}{}&\int(0){\0}{}&\int(0){\0}{}&\int(0){\0}{}&\ln{1}&\cnt{20}&\bounds{4}{4}{1}{false}&&\cr
\l{5}&\int(2){4}{}&\int(0){\0}{}&\int(1){2}{}&\int(0){2}{}&\ln{2}&\cnt{40}&\bounds{8}{8}{1}{false}&&\cr
\l{7}&\rat(5/2){5}{}&\int(0){7}{}&\int(0){\0}{}&\rat(1/2){3}{}&\ln{2}&\cnt{64}&\bounds{8}{8}{1}{false}&&\cr
\endorbit
\orbit{4}{16}{108}
\h{16}&\rat(7/8){7}{}&\rat(31/8){5}{\spec{1}}&\int(0){\0}{}&\rat(5/4){1}{}&&\cnt{744}&\bounds{108}{108}{1}{false}&&\cr
\l{1}&\int(1){\2}{}&\int(2){\2}{}&\int(0){\0}{}&\int(0){\0}{}&\ln{1}&\cnt{42}&\bounds{6}{6}{1}{false}&&\cr
\l{3}&\rat(5/8){5}{}&\rat(13/8){7}{}&\int(0){\0}{}&\rat(3/4){3}{}&\ln{1}&\cnt{105}&\bounds{15}{15}{1}{false}&&\cr
\l{5}&\rat(3/4){6}{}&\rat(9/4){6}{}&\int(0){2}{}&\int(0){\0}{}&\ln{1}&\cnt{70}&\bounds{10}{10}{1}{false}&&\cr
\l{7}&\rat(7/8){7}{}&\rat(15/8){5}{}&\int(0){\0}{}&\rat(1/4){1}{}&\ln{1}&\cnt{60}&\bounds{10}{10}{1}{false}&&\cr
\l{9}&\rat(7/8){7}{}&\rat(7/8){5}{}&\int(0){\0}{}&\rat(5/4){1}{}&\ln{1}&\cnt{15}&\bounds{3}{3}{1}{false}&&\cr
\l{11}&\rat(7/8){7}{}&\rat(13/8){7}{}&\int(0){3}{}&\rat(1/2){2}{}&\ln{1}&\cnt{80}&\bounds{10}{10}{1}{false}&&\cr
\endorbit
\orbit{4}{17}{108}
\h{17}&\rat(23/8){7}{\spec{2}}&\rat(15/8){5}{}&\int(0){\0}{}&\rat(5/4){1}{}&&\cnt{552}&\bounds{108}{108}{1}{false}&&\cr
\l{1}&\int(3){\4}{}&\int(0){\0}{}&\int(0){\0}{}&\int(0){\0}{}&\ln{1}&\cnt{5}&\bounds{1}{1}{1}{false}&&\cr
\l{3}&\int(2){\2}{}&\int(1){\2}{}&\int(0){\0}{}&\int(0){\0}{}&\ln{1}&\cnt{30}&\bounds{6}{6}{1}{false}&&\cr
\l{5}&\int(2){\2}{}&\int(0){\0}{}&\int(0){\0}{}&\int(1){\2}{}&\ln{1}&\cnt{20}&\bounds{4}{4}{1}{false}&&\cr
\l{7}&\rat(3/2){4}{}&\rat(3/2){4}{}&\int(0){\0}{}&\int(0){\0}{}&\ln{1}&\cnt{50}&\bounds{10}{10}{1}{false}&&\cr
\l{9}&\rat(13/8){5}{}&\rat(5/8){7}{}&\int(0){\0}{}&\rat(3/4){3}{}&\ln{1}&\cnt{75}&\bounds{15}{15}{1}{false}&&\cr
\l{11}&\rat(7/4){6}{}&\rat(5/4){6}{}&\int(0){2}{}&\int(0){\0}{}&\ln{1}&\cnt{30}&\bounds{6}{6}{1}{false}&&\cr
\l{13}&\rat(7/4){6}{}&\rat(3/4){2}{}&\int(0){\0}{}&\rat(1/2){2}{}&\ln{1}&\cnt{50}&\bounds{10}{10}{1}{false}&&\cr
\l{15}&\rat(7/4){6}{}&\int(0){\0}{}&\int(0){1}{}&\rat(5/4){1}{}&\ln{1}&\cnt{16}&\bounds{2}{2}{1}{false}&&\cr
\endorbit
\orbit{4}{18}{108}
\h{18}&\rat(3/2){6}{}&\rat(3/2){6}{}&\int(3){2}{}&\int(0){\0}{}&&\cnt{600}&\bounds{108}{108}{1}{false}&&\cr
\l{1}&\rat(5/4){5}{}&\rat(1/4){1}{}&\rat(3/2){3}{}&\int(0){\0}{}&\ln{2\sd}&\cnt{72}&\bounds{12}{12}{1}{false}&&\cr
\l{2}&\rat(3/4){3}{}&\rat(3/4){7}{}&\rat(3/2){1}{}&\int(0){\0}{}&\ln{2\sdd}&\cnt{80}&\bounds{16}{16}{1}{false}&&\cr
\l{3}&\rat(3/2){6}{}&\int(0){\0}{}&\rat(3/2){1}{}&\int(0){1}{}&\ln{2\sd}&\cnt{32}&\bounds{4}{4}{1}{false}&&\cr
\l{4}&\rat(3/4){7}{}&\rat(3/4){7}{}&\rat(3/2){3}{}&\int(0){2}{}&\ln{1\sdd}&\cnt{80}&\bounds{16}{16}{1}{false}&&\cr
\l{5}&\int(1){\2}{}&\int(0){\0}{}&\int(2){\2}{}&\int(0){\0}{}&\ln{2}&\cnt{36}&\bounds{6}{6}{1}{false}&&\cr
\l{7}&\rat(3/2){6}{}&\rat(3/2){6}{}&\int(0){2}{}&\int(0){\0}{}&\ln{1}&\cnt{4}&\bounds{2}{2}{1}{false}&&\cr
\endorbit
\orbit{4}{19}{106}
\h{19}&\rat(3/2){6}{}&\int(0){\0}{}&\rat(13/4){1}{}&\rat(5/4){1}{}&&\cnt{648}&\bounds{106}{106}{1}{false}&&\cr
\l{1}&\int(1){4}{}&\int(0){\0}{}&\rat(3/2){2}{}&\rat(1/2){2}{}&\ln{1}&\cnt{75}&\bounds{15}{15}{1}{false}&&\cr
\l{3}&\int(1){\2}{}&\int(0){\0}{}&\int(2){\2}{}&\int(0){\0}{}&\ln{1}&\cnt{48}&\bounds{8}{8}{1}{false}&&\cr
\l{5}&\rat(5/4){5}{}&\int(0){1}{}&\rat(7/4){3}{}&\int(0){\0}{}&\ln{1}&\cnt{48}&\bounds{6}{6}{1}{false}&&\cr
\l{7}&\rat(3/2){6}{}&\int(0){6}{}&\rat(3/2){2}{}&\int(0){\0}{}&\ln{1}&\cnt{28}&\bounds{4}{4}{1}{false}&&\cr
\l{9}&\rat(3/2){6}{}&\int(0){\0}{}&\rat(5/4){1}{}&\rat(1/4){1}{}&\ln{1}&\cnt{40}&\bounds{8}{8}{1}{false}&&\cr
\l{11}&\rat(3/2){6}{}&\int(0){\0}{}&\rat(1/4){1}{\spec{1}}&\rat(5/4){1}{}&\ln{1}&\cnt{1}&\bounds{1}{1}{1}{false}&&\cr
\l{13}&\rat(3/2){6}{}&\int(0){\0}{}&\rat(1/4){1}{\spec{2}}&\rat(5/4){1}{}&\ln{1}&\cnt{4}&\bounds{1}{1}{1}{false}&&\cr
\l{15}&\rat(3/4){7}{}&\int(0){7}{}&\rat(7/4){3}{}&\rat(1/2){2}{}&\ln{1}&\cnt{80}&\bounds{10}{10}{1}{false}&&\cr
\endorbit
\orbit{4}{20}{106}
\h{20}&\rat(7/2){6}{\spec{1}}&\int(0){\0}{}&\rat(5/4){1}{}&\rat(5/4){1}{}&&\cnt{792}&\bounds{106}{106}{1}{false}&&\cr
\l{1}&\rat(7/4){7}{}&\int(0){5}{}&\int(0){\0}{}&\rat(5/4){1}{}&\ln{2\sd}&\cnt{56}&\bounds{8}{8}{1}{false}&&\cr
\l{2}&\rat(7/4){7}{}&\int(0){7}{}&\rat(3/4){3}{}&\rat(1/2){2}{}&\ln{2\sd}&\cnt{200}&\bounds{25}{25}{1}{false}&&\cr
\l{3}&\int(2){\2}{}&\int(0){\0}{}&\int(1){\2}{}&\int(0){\0}{}&\ln{2}&\cnt{70}&\bounds{10}{10}{1}{false}&&\cr
\endorbit
\orbit{4}{21}{106}
\h{21}&\rat(7/8){7}{}&\rat(15/8){5}{}&\int(0){\0}{}&\rat(13/4){1}{}&&\cnt{648}&\bounds{106}{106}{1}{false}&&\cr
\l{1}&\int(1){\2}{}&\int(0){\0}{}&\int(0){\0}{}&\int(2){\2}{}&\ln{1}&\cnt{28}&\bounds{4}{4}{1}{false}&&\cr
\l{3}&\rat(5/8){5}{}&\rat(5/8){7}{}&\int(0){\0}{}&\rat(7/4){3}{}&\ln{1}&\cnt{63}&\bounds{9}{9}{1}{false}&&\cr
\l{5}&\rat(3/4){6}{}&\rat(3/4){2}{}&\int(0){\0}{}&\rat(3/2){2}{}&\ln{1}&\cnt{70}&\bounds{10}{10}{1}{false}&&\cr
\l{7}&\rat(7/8){7}{}&\rat(15/8){5}{}&\int(0){\0}{}&\rat(1/4){1}{\spec{1}}&\ln{1}&\cnt{1}&\bounds{1}{1}{1}{false}&&\cr
\l{9}&\rat(7/8){7}{}&\rat(15/8){5}{}&\int(0){\0}{}&\rat(1/4){1}{\spec{2}}&\ln{1}&\cnt{4}&\bounds{1}{1}{1}{false}&&\cr
\l{11}&\rat(7/8){7}{}&\rat(7/8){5}{}&\int(0){\0}{}&\rat(5/4){1}{}&\ln{1}&\cnt{60}&\bounds{12}{12}{1}{false}&&\cr
\l{13}&\rat(7/8){7}{}&\rat(5/8){7}{}&\int(0){3}{}&\rat(3/2){2}{}&\ln{1}&\cnt{48}&\bounds{6}{6}{1}{false}&&\cr
\l{15}&\rat(7/8){7}{}&\rat(3/8){1}{}&\int(0){2}{}&\rat(7/4){3}{}&\ln{1}&\cnt{50}&\bounds{10}{10}{1}{false}&&\cr
\endorbit
\orbit{4}{22}{106}
\h{22}&\int(2){\2}{}&\int(0){\0}{}&\int(4){\4}{}&\int(0){\0}{}&&\cnt{792}&\bounds{106}{106}{1}{false}&&\cr
\l{1}&\int(1){\2}{\spec{1}}&\int(0){\0}{}&\int(2){\2}{}&\int(0){\0}{}&\ln{2\sd}&\cnt{36}&\bounds{6}{6}{1}{false}&&\cr
\l{2}&\int(1){5}{}&\int(0){1}{}&\int(2){3}{}&\int(0){\0}{}&\ln{2\sd}&\cnt{120}&\bounds{15}{15}{1}{false}&&\cr
\l{3}&\int(1){6}{}&\int(0){\0}{}&\int(2){1}{}&\int(0){1}{}&\ln{2\sd}&\cnt{96}&\bounds{12}{12}{1}{false}&&\cr
\l{4}&\int(1){7}{}&\int(0){3}{}&\int(2){1}{}&\int(0){\0}{}&\ln{2\sd}&\cnt{56}&\bounds{8}{8}{1}{false}&&\cr
\l{5}&\int(1){7}{}&\int(0){7}{}&\int(2){3}{}&\int(0){2}{}&\ln{2\sd}&\cnt{80}&\bounds{10}{10}{1}{false}&&\cr
\l{6}&\int(2){\2}{}&\int(0){\0}{}&\int(1){\2}{\spec{1}}&\int(0){\0}{}&\ln{2}&\cnt{4}&\bounds{1}{1}{1}{false}&&\cr
\endorbit
\orbit{4}{23}{104}
\h{23}&\int(4){\4}{}&\int(0){\0}{}&\int(2){\2}{}&\int(0){\0}{}&&\cnt{696}&\bounds{104}{112}{fail}{false}&&\cr
\l{1}&\int(2){4}{}&\int(0){\0}{}&\int(1){2}{}&\int(0){2}{}&\ln{1\sdd}&\cnt{120}&\bounds{20}{24}{2}{false}&&\cr
\l{2}&\int(2){\2}{}&\int(0){\0}{}&\int(1){\2}{}&\int(0){\0}{}&\ln{1\sdd}&\cnt{48}&\bounds{12}{16}{2}{false}&&\cr
\l{3}&\int(2){5}{}&\int(0){1}{}&\int(1){3}{}&\int(0){\0}{}&\ln{2\sd}&\cnt{128}&\bounds{16}{16}{1}{false}&&\cr
\l{4}&\int(2){6}{}&\int(0){6}{}&\int(1){2}{}&\int(0){\0}{}&\ln{2\sd}&\cnt{56}&\bounds{8}{8}{1}{false}&&\cr
\l{5}&\int(2){6}{}&\int(0){\0}{}&\int(1){1}{}&\int(0){1}{}&\ln{2\sd}&\cnt{64}&\bounds{8}{8}{1}{false}&&\cr
\l{6}&\int(3){\4}{\spec{1}}&\int(0){\0}{}&\int(0){\0}{}&\int(0){\0}{}&\ln{2}&\cnt{8}&\bounds{2}{2}{1}{false}&&\cr
\endorbit
\orbit{4}{24}{102}
\h{24}&\rat(7/2){6}{\spec{2}}&\int(0){\0}{}&\rat(5/4){1}{}&\rat(5/4){1}{}&&\cnt{600}&\bounds{102}{102}{1}{false}&&\cr
\l{1}&\rat(7/4){3}{}&\int(0){7}{}&\rat(5/4){1}{}&\int(0){\0}{}&\ln{2\sd}&\cnt{48}&\bounds{6}{6}{1}{false}&&\cr
\l{2}&\int(3){\4}{}&\int(0){\0}{}&\int(0){\0}{}&\int(0){\0}{}&\ln{1}&\cnt{12}&\bounds{3}{3}{1}{false}&&\cr
\l{4}&\int(2){4}{}&\int(0){\0}{}&\rat(1/2){2}{}&\rat(1/2){2}{}&\ln{1}&\cnt{100}&\bounds{20}{20}{1}{false}&&\cr
\l{6}&\int(2){\2}{}&\int(0){\0}{}&\int(1){\2}{}&\int(0){\0}{}&\ln{2}&\cnt{30}&\bounds{6}{6}{1}{false}&&\cr
\l{8}&\rat(9/4){5}{}&\int(0){7}{}&\int(0){\0}{}&\rat(3/4){3}{}&\ln{2}&\cnt{40}&\bounds{5}{5}{1}{false}&&\cr
\endorbit
\orbit{4}{25}{100}
\h{25}&\rat(23/8){7}{\spec{2}}&\rat(7/8){7}{}&\rat(5/4){3}{}&\int(1){2}{}&&\cnt{552}&\bounds{100}{100}{1}{false}&&\cr
\l{1}&\int(3){\4}{}&\int(0){\0}{}&\int(0){\0}{}&\int(0){\0}{}&\ln{1}&\cnt{5}&\bounds{1}{1}{1}{false}&&\cr
\l{3}&\int(2){\2}{}&\int(1){\2}{}&\int(0){\0}{}&\int(0){\0}{}&\ln{1}&\cnt{14}&\bounds{2}{2}{1}{false}&&\cr
\l{5}&\int(2){\2}{}&\int(0){\0}{}&\int(1){\2}{}&\int(0){\0}{}&\ln{1}&\cnt{20}&\bounds{4}{4}{1}{false}&&\cr
\l{7}&\int(2){\2}{}&\int(0){\0}{}&\int(0){\0}{}&\int(1){\2}{}&\ln{1}&\cnt{16}&\bounds{4}{4}{1}{false}&&\cr
\l{9}&\rat(3/2){4}{}&\int(0){\0}{}&\rat(1/2){2}{}&\int(1){2}{}&\ln{1}&\cnt{50}&\bounds{10}{10}{1}{false}&&\cr
\l{11}&\rat(13/8){5}{}&\rat(7/8){7}{}&\int(0){\0}{}&\rat(1/2){3}{}&\ln{1}&\cnt{40}&\bounds{8}{8}{1}{false}&&\cr
\l{13}&\rat(13/8){5}{}&\rat(1/8){1}{}&\rat(5/4){3}{}&\int(0){\0}{}&\ln{1}&\cnt{35}&\bounds{5}{5}{1}{false}&&\cr
\l{15}&\rat(7/4){6}{}&\rat(3/4){6}{}&\rat(1/2){2}{}&\int(0){\0}{}&\ln{1}&\cnt{35}&\bounds{5}{5}{1}{false}&&\cr
\l{17}&\rat(7/4){6}{}&\rat(1/4){2}{}&\int(0){\0}{}&\int(1){2}{}&\ln{1}&\cnt{21}&\bounds{3}{3}{1}{false}&&\cr
\l{19}&\rat(7/4){6}{}&\int(0){\0}{}&\rat(3/4){1}{}&\rat(1/2){1}{}&\ln{1}&\cnt{40}&\bounds{8}{8}{1}{false}&&\cr
\endorbit
\orbit{4}{26}{100}
\h{26}&\rat(3/2){6}{}&\rat(3/2){6}{}&\int(1){2}{}&\int(2){\2}{}&&\cnt{504}&\bounds{100}{100}{1}{false}&&\cr
\l{1}&\rat(3/2){6}{}&\int(0){\0}{}&\rat(1/2){1}{}&\int(1){1}{}&\ln{2\sd}&\cnt{32}&\bounds{8}{8}{1}{false}&&\cr
\l{2}&\rat(3/4){7}{}&\rat(3/4){7}{}&\rat(1/2){3}{}&\int(1){2}{}&\ln{1\sdd}&\cnt{64}&\bounds{16}{16}{1}{false}&&\cr
\l{3}&\int(1){4}{}&\int(0){\0}{}&\int(1){2}{}&\int(1){2}{}&\ln{2}&\cnt{30}&\bounds{6}{6}{1}{false}&&\cr
\l{5}&\int(1){\2}{}&\int(0){\0}{}&\int(0){\0}{}&\int(2){\2}{}&\ln{2}&\cnt{12}&\bounds{2}{2}{1}{false}&&\cr
\l{7}&\rat(5/4){5}{}&\rat(3/4){7}{}&\int(0){\0}{}&\int(1){3}{}&\ln{2}&\cnt{48}&\bounds{8}{8}{1}{false}&&\cr
\l{9}&\rat(3/2){6}{}&\rat(3/2){6}{}&\int(0){2}{}&\int(0){\0}{}&\ln{1}&\cnt{8}&\bounds{2}{2}{1}{false}&&\cr
\endorbit
\orbit{4}{27}{100}
\h{27}&\int(2){\2}{}&\int(0){\0}{}&\int(2){\2}{}&\int(2){\2}{}&&\cnt{600}&\bounds{100}{100}{1}{false}&&\cr
\l{1}&\int(1){4}{}&\int(0){\0}{}&\int(1){2}{}&\int(1){2}{}&\ln{1\sdd}&\cnt{80}&\bounds{16}{16}{1}{false}&&\cr
\l{2}&\int(1){\2}{\spec{1}}&\int(0){\0}{}&\int(2){\2}{}&\int(0){\0}{}&\ln{4\sd}&\cnt{6}&\bounds{1}{1}{1}{false}&&\cr
\l{3}&\int(1){6}{}&\int(0){\0}{}&\int(1){1}{}&\int(1){1}{}&\ln{2\sd}&\cnt{96}&\bounds{16}{16}{1}{false}&&\cr
\l{4}&\int(1){7}{}&\int(0){7}{}&\int(1){3}{}&\int(1){2}{}&\ln{4\sd}&\cnt{64}&\bounds{8}{8}{1}{false}&&\cr
\l{5}&\int(2){\2}{}&\int(0){\0}{}&\int(1){\2}{}&\int(0){\0}{}&\ln{2}&\cnt{12}&\bounds{4}{4}{1}{false}&&\cr
\endorbit
\orbit{4}{28}{96}
\h{28}&\rat(15/8){5}{}&\rat(15/8){5}{}&\rat(5/4){1}{}&\int(1){2}{}&&\cnt{456}&\bounds{96}{96}{1}{false}&&\cr
\l{1}&\rat(3/2){4}{}&\rat(3/2){4}{}&\int(0){\0}{}&\int(0){\0}{}&\ln{1}&\cnt{25}&\bounds{5}{5}{1}{false}&&\cr
\l{3}&\rat(3/2){4}{}&\int(0){\0}{}&\rat(1/2){2}{}&\int(1){2}{}&\ln{2}&\cnt{25}&\bounds{5}{5}{1}{false}&&\cr
\l{5}&\rat(15/8){5}{}&\rat(5/8){7}{}&\int(0){\0}{}&\rat(1/2){3}{}&\ln{2}&\cnt{24}&\bounds{6}{6}{1}{false}&&\cr
\l{7}&\rat(9/8){3}{}&\rat(5/8){7}{}&\rat(5/4){1}{}&\int(0){\0}{}&\ln{2}&\cnt{30}&\bounds{6}{6}{1}{false}&&\cr
\l{9}&\rat(5/4){6}{}&\rat(5/4){6}{}&\rat(1/2){2}{}&\int(0){\0}{}&\ln{1}&\cnt{45}&\bounds{9}{9}{1}{false}&&\cr
\endorbit
\orbit{4}{29}{92}
\h{29}&\rat(3/2){6}{}&\int(2){\2}{}&\rat(5/4){1}{}&\rat(5/4){1}{}&&\cnt{504}&\bounds{92}{92}{1}{false}&&\cr
\l{1}&\rat(3/4){3}{}&\int(1){7}{}&\rat(5/4){1}{}&\int(0){\0}{}&\ln{2\sd}&\cnt{20}&\bounds{4}{4}{1}{false}&&\cr
\l{2}&\rat(3/4){7}{}&\int(1){5}{}&\int(0){\0}{}&\rat(5/4){1}{}&\ln{2\sd}&\cnt{30}&\bounds{6}{6}{1}{false}&&\cr
\l{3}&\rat(3/4){7}{}&\int(1){7}{}&\rat(3/4){3}{}&\rat(1/2){2}{}&\ln{2\sd}&\cnt{50}&\bounds{10}{10}{1}{false}&&\cr
\l{4}&\int(1){\2}{}&\int(2){\2}{}&\int(0){\0}{}&\int(0){\0}{}&\ln{1}&\cnt{12}&\bounds{2}{2}{1}{false}&&\cr
\l{6}&\rat(5/4){5}{}&\int(1){7}{}&\int(0){\0}{}&\rat(3/4){3}{}&\ln{2}&\cnt{30}&\bounds{5}{5}{1}{false}&&\cr
\l{8}&\rat(3/2){6}{}&\int(1){6}{}&\rat(1/2){2}{}&\int(0){\0}{}&\ln{2}&\cnt{30}&\bounds{5}{5}{1}{false}&&\cr
\l{10}&\rat(3/2){6}{}&\int(0){\0}{}&\rat(5/4){1}{}&\rat(1/4){1}{}&\ln{2}&\cnt{10}&\bounds{2}{2}{1}{false}&&\cr
\endorbit
\orbit{4}{30}{92}
\h{30}&\rat(7/8){7}{}&\rat(7/8){7}{}&\rat(5/4){3}{}&\int(3){2}{}&&\cnt{600}&\bounds{92}{92}{1}{false}&&\cr
\l{1}&\int(1){\2}{}&\int(0){\0}{}&\int(0){\0}{}&\int(2){\2}{}&\ln{2}&\cnt{21}&\bounds{3}{3}{1}{false}&&\cr
\l{3}&\rat(5/8){5}{}&\rat(7/8){7}{}&\int(0){\0}{}&\rat(3/2){3}{}&\ln{2}&\cnt{42}&\bounds{6}{6}{1}{false}&&\cr
\l{5}&\rat(3/4){6}{}&\int(0){\0}{}&\rat(3/4){1}{}&\rat(3/2){1}{}&\ln{2}&\cnt{70}&\bounds{10}{10}{1}{false}&&\cr
\l{7}&\rat(7/8){7}{}&\rat(7/8){7}{}&\rat(5/4){3}{}&\int(0){2}{}&\ln{1}&\cnt{4}&\bounds{2}{2}{1}{false}&&\cr
\l{9}&\rat(7/8){7}{}&\rat(7/8){7}{}&\rat(1/4){3}{}&\int(1){2}{}&\ln{1}&\cnt{30}&\bounds{6}{6}{1}{false}&&\cr
\endorbit
\orbit{4}{31}{92}
\h{31}&\rat(7/8){7}{}&\rat(15/8){5}{}&\int(2){\2}{}&\rat(5/4){1}{}&&\cnt{504}&\bounds{92}{92}{1}{false}&&\cr
\l{1}&\int(1){\2}{}&\int(0){\0}{}&\int(2){\2}{}&\int(0){\0}{}&\ln{1}&\cnt{7}&\bounds{1}{1}{1}{false}&&\cr
\l{3}&\rat(3/4){6}{}&\rat(5/4){6}{}&\int(1){2}{}&\int(0){\0}{}&\ln{1}&\cnt{42}&\bounds{6}{6}{1}{false}&&\cr
\l{5}&\rat(3/4){6}{}&\int(0){\0}{}&\int(1){1}{}&\rat(5/4){1}{}&\ln{1}&\cnt{28}&\bounds{4}{4}{1}{false}&&\cr
\l{7}&\rat(7/8){7}{}&\rat(15/8){5}{}&\int(0){\0}{}&\rat(1/4){1}{}&\ln{1}&\cnt{10}&\bounds{2}{2}{1}{false}&&\cr
\l{9}&\rat(7/8){7}{}&\rat(9/8){3}{}&\int(1){1}{}&\int(0){\0}{}&\ln{1}&\cnt{40}&\bounds{8}{8}{1}{false}&&\cr
\l{11}&\rat(7/8){7}{}&\rat(7/8){5}{}&\int(0){\0}{}&\rat(5/4){1}{}&\ln{1}&\cnt{15}&\bounds{3}{3}{1}{false}&&\cr
\l{13}&\rat(7/8){7}{}&\rat(5/8){7}{}&\int(1){3}{}&\rat(1/2){2}{}&\ln{1}&\cnt{60}&\bounds{12}{12}{1}{false}&&\cr
\l{15}&\rat(7/8){7}{}&\rat(3/8){1}{}&\int(1){2}{}&\rat(3/4){3}{}&\ln{1}&\cnt{50}&\bounds{10}{10}{1}{false}&&\cr
\endorbit
\orbit{4}{32}{88}
\h{32}&\rat(3/2){6}{}&\int(2){4}{}&\rat(5/4){3}{}&\rat(5/4){3}{}&&\cnt{456}&\bounds{88}{88}{1}{false}&&\cr
\l{1}&\int(1){4}{}&\int(2){4}{}&\int(0){\0}{}&\int(0){\0}{}&\ln{1}&\cnt{15}&\bounds{3}{3}{1}{false}&&\cr
\l{3}&\rat(5/4){5}{}&\rat(1/2){7}{}&\int(0){\0}{}&\rat(5/4){3}{}&\ln{2}&\cnt{24}&\bounds{4}{4}{1}{false}&&\cr
\l{5}&\rat(3/2){6}{}&\int(1){6}{}&\rat(1/2){2}{}&\int(0){\0}{}&\ln{2}&\cnt{30}&\bounds{6}{6}{1}{false}&&\cr
\l{7}&\rat(3/2){6}{}&\int(0){\0}{}&\rat(3/4){1}{}&\rat(3/4){1}{}&\ln{1}&\cnt{25}&\bounds{5}{5}{1}{false}&&\cr
\l{9}&\rat(3/4){7}{}&\rat(3/2){5}{}&\int(0){\0}{}&\rat(3/4){1}{}&\ln{2}&\cnt{40}&\bounds{8}{8}{1}{false}&&\cr
\endorbit
\orbit{4}{33}{50}
\h{33}&\int(6){\6}{\spec{1}}&\int(0){\0}{}&\int(0){\0}{}&\int(0){\0}{}&&\cnt{312}&\bounds{50}{50}{1}{false}&&\cr
\l{1}&\int(3){\4}{}&\int(0){\0}{}&\int(0){\0}{}&\int(0){\0}{}&\ln{1\sdd}&\cnt{40}&\bounds{10}{10}{1}{false}&&\cr
\l{2}&\int(3){\2}{}&\int(0){\2}{}&\int(0){\0}{}&\int(0){\0}{}&\ln{1\sdd}&\cnt{112}&\bounds{16}{16}{1}{false}&&\cr
\l{3}&\int(3){\2}{}&\int(0){\0}{}&\int(0){\2}{}&\int(0){\0}{}&\ln{2\sdd}&\cnt{80}&\bounds{12}{12}{1}{false}&&\cr
\endorbit
\orbit{4}{34}{46}
\h{34}&\int(0){\0}{}&\int(0){\0}{}&\int(6){\6}{\spec{1}}&\int(0){\0}{}&&\cnt{312}&\bounds{46}{46}{1}{false}&&\cr
\l{1}&\int(0){\2}{}&\int(0){\0}{}&\int(3){\2}{}&\int(0){\0}{}&\ln{2\sdd}&\cnt{112}&\bounds{16}{16}{1}{false}&&\cr
\l{2}&\int(0){\0}{}&\int(0){\0}{}&\int(3){\4}{\spec{1}}&\int(0){\0}{}&\ln{2\sdd}&\cnt{4}&\bounds{1}{1}{1}{false}&&\cr
\l{3}&\int(0){\0}{}&\int(0){\0}{}&\int(3){\2}{}&\int(0){\2}{}&\ln{1\sdd}&\cnt{80}&\bounds{12}{12}{1}{false}&&\cr
\endorbit
\endorbits

\newcs{max-4A6}{116}%
\newcs{count-4A6}{16}%
\newcs{counts-4A6}{\DO{1}{104}\DO{2}{106}\DO{3}{108}\DO{4}{108}\DO{5}{116}\DO{6}{114}\DO{7}{112}\DO{8}{112}\DO{9}{108}\DO{10}{108}\DO{11}{108}\DO{12}{106}\DO{13}{102}\DO{14}{100}\DO{15}{96}\DO{16}{50}}%
\orbits
\orbit{4}{1}{104}
\h{1}&\root\int(6){\6}{\spec{1}}&\int(0){\0}{}&\int(0){\0}{}&\int(0){\0}{}&&\cnt{888}&\bounds{104}{128}{fail}{false}&&\cr
\l{1}&\int(3){\4}{\spec{1}}&\int(0){\0}{}&\int(0){\0}{}&\int(0){\0}{}&\ln{1}&\cnt{3}&\bounds{1}{1}{1}{false}&&\cr
\l{3}&\int(3){4}{}&\int(0){2}{}&\int(0){6}{}&\int(0){\0}{}&\ln{3}&\cnt{147}&\bounds{17}{21}{3}{false}&&\cr
\endorbit
\orbit{4}{2}{106}
\h{2}&\int(6){\6}{}&\root\int(0){\0}{}&\int(0){\0}{}&\int(0){\0}{}&&\cnt{676}&\bounds{106}{106}{1}{false}&&\cr
\l{1}&\int(3){\4}{\spec{1}}&\int(0){\0}{}&\int(0){\0}{}&\int(0){\0}{}&\ln{2\sd}&\cnt{9}&\bounds{6}{6}{1}{false}&&\cr
\l{2}&\int(3){4}{}&\int(0){2}{\spec{1}}&\int(0){6}{}&\int(0){\0}{}&\ln{2\sd}&\cnt{7}&\bounds{1}{1}{1}{false}&&\cr
\l{3}&\int(3){4}{}&\int(0){2}{\spec{2}}&\int(0){6}{}&\int(0){\0}{}&\ln{2\sd}&\cnt{70}&\bounds{10}{10}{1}{false}&&\cr
\l{4}&\int(3){4}{}&\int(0){6}{}&\int(0){\0}{}&\int(0){2}{}&\ln{2\sd}&\cnt{105}&\bounds{15}{15}{1}{false}&&\cr
\l{5}&\int(3){4}{}&\int(0){\0}{}&\int(0){2}{}&\int(0){6}{}&\ln{2\sd}&\cnt{147}&\bounds{21}{21}{1}{false}&&\cr
\endorbit
\orbit{4}{3}{108}
\h{3}&\rat(6/7){6}{}&\rat(6/7){6}{}&\rat(6/7){1}{}&\rat(24/7){5}{\spec{1}}&&\cnt{756}&\bounds{108}{126}{fail}{false}&&\cr
\l{1}&\int(1){\2}{}&\int(0){\0}{}&\int(0){\0}{}&\int(2){\2}{}&\ln{3}&\cnt{36}&\bounds{6}{6}{1}{false}&&\cr
\l{3}&\rat(4/7){4}{}&\int(0){\0}{}&\rat(5/7){2}{}&\rat(12/7){6}{}&\ln{3}&\cnt{90}&\bounds{12}{15}{3}{false}&&\cr
\endorbit
\orbit{4}{4}{108}
\h{4}&\rat(6/7){6}{}&\rat(26/7){3}{\spec{2}}&\rat(10/7){2}{}&\int(0){\0}{}&&\cnt{708}&\bounds{108}{122}{fail}{false}&&\cr
\l{1}&\int(1){\2}{}&\int(2){\2}{}&\int(0){\0}{}&\int(0){\0}{}&\ln{1}&\cnt{30}&\bounds{5}{5}{1}{true}&&\cr
\l{3}&\rat(4/7){4}{}&\rat(15/7){2}{}&\rat(2/7){6}{}&\int(0){\0}{}&\ln{1}&\cnt{75}&\bounds{11}{15}{5}{true}&&\cr
\l{5}&\rat(5/7){5}{}&\rat(11/7){1}{}&\rat(5/7){1}{}&\int(0){1}{}&\ln{1}&\cnt{84}&\bounds{12}{12}{1}{false}&&\cr
\l{7}&\rat(6/7){6}{}&\rat(12/7){3}{}&\rat(3/7){2}{}&\int(0){\0}{}&\ln{1}&\cnt{50}&\bounds{9}{10}{5}{true}&&\cr
\l{9}&\rat(6/7){6}{}&\rat(5/7){3}{}&\rat(10/7){2}{}&\int(0){\0}{}&\ln{1}&\cnt{10}&\bounds{2}{2}{1}{true}&&\cr
\l{11}&\rat(6/7){6}{}&\rat(15/7){2}{}&\int(0){\0}{}&\int(0){3}{}&\ln{1}&\cnt{35}&\bounds{5}{7}{5}{true}&&\cr
\l{13}&\rat(6/7){6}{}&\rat(11/7){1}{}&\rat(4/7){5}{}&\int(0){6}{}&\ln{1}&\cnt{70}&\bounds{10}{10}{1}{true}&&\cr
\endorbit
\orbit{4}{5}{116}
\h{5}&\rat(6/7){6}{}&\rat(12/7){3}{}&\rat(24/7){2}{\spec{2}}&\int(0){\0}{}&&\cnt{756}&\bounds{116}{120}{fail}{false}&&\cr
\l{1}&\rat(3/7){3}{}&\rat(6/7){5}{}&\rat(12/7){1}{}&\int(0){\0}{}&\ln{1\sdd}&\cnt{120}&\bounds{20}{24}{2}{false}&&\cr
\l{2}&\int(1){\2}{}&\int(0){\0}{}&\int(2){\2}{}&\int(0){\0}{}&\ln{1}&\cnt{36}&\bounds{6}{6}{1}{false}&&\cr
\l{4}&\rat(5/7){5}{}&\rat(4/7){1}{}&\rat(12/7){1}{}&\int(0){1}{}&\ln{1}&\cnt{126}&\bounds{18}{18}{1}{false}&&\cr
\l{6}&\rat(6/7){6}{}&\rat(5/7){3}{}&\rat(10/7){2}{}&\int(0){\0}{}&\ln{1}&\cnt{72}&\bounds{12}{12}{1}{false}&&\cr
\l{8}&\rat(6/7){6}{}&\rat(3/7){6}{}&\rat(12/7){1}{}&\int(0){5}{}&\ln{1}&\cnt{84}&\bounds{12}{12}{1}{false}&&\cr
\endorbit
\orbit{4}{6}{114}
\h{6}&\rat(6/7){6}{}&\rat(6/7){6}{}&\rat(6/7){1}{}&\rat(24/7){5}{\spec{2}}&&\cnt{612}&\bounds{114}{114}{1}{false}&&\cr
\l{1}&\int(1){\2}{}&\int(0){\0}{}&\int(0){\0}{}&\int(2){\2}{}&\ln{3}&\cnt{18}&\bounds{3}{3}{1}{false}&&\cr
\l{3}&\rat(4/7){4}{}&\rat(6/7){6}{}&\int(0){\0}{}&\rat(11/7){2}{}&\ln{3}&\cnt{45}&\bounds{9}{9}{1}{false}&&\cr
\l{5}&\rat(5/7){5}{}&\rat(1/7){1}{}&\rat(6/7){1}{}&\rat(9/7){1}{}&\ln{3}&\cnt{36}&\bounds{6}{6}{1}{false}&&\cr
\l{7}&\rat(6/7){6}{}&\rat(6/7){6}{}&\rat(6/7){1}{}&\rat(3/7){5}{}&\ln{1}&\cnt{9}&\bounds{3}{3}{1}{false}&&\cr
\endorbit
\orbit{4}{7}{112}
\h{7}&\int(4){\4}{}&\int(2){\2}{}&\int(0){\0}{}&\int(0){\0}{}&&\cnt{708}&\bounds{112}{112}{1}{false}&&\cr
\l{1}&\int(2){\2}{}&\int(1){\2}{\spec{1}}&\int(0){\0}{}&\int(0){\0}{}&\ln{2\sd}&\cnt{20}&\bounds{4}{4}{1}{false}&&\cr
\l{2}&\int(2){4}{}&\int(1){2}{}&\int(0){6}{}&\int(0){\0}{}&\ln{2\sd}&\cnt{105}&\bounds{15}{15}{1}{false}&&\cr
\l{3}&\int(2){4}{}&\int(1){6}{}&\int(0){\0}{}&\int(0){2}{}&\ln{2\sd}&\cnt{63}&\bounds{9}{9}{1}{false}&&\cr
\l{4}&\int(2){5}{}&\int(1){4}{}&\int(0){\0}{}&\int(0){6}{}&\ln{2\sd}&\cnt{70}&\bounds{10}{10}{1}{false}&&\cr
\l{5}&\int(2){5}{}&\int(1){6}{}&\int(0){4}{}&\int(0){\0}{}&\ln{2\sd}&\cnt{35}&\bounds{7}{7}{1}{false}&&\cr
\l{6}&\int(2){5}{}&\int(1){1}{}&\int(0){1}{}&\int(0){1}{}&\ln{2\sd}&\cnt{49}&\bounds{7}{7}{1}{false}&&\cr
\l{7}&\int(3){\4}{\spec{1}}&\int(0){\0}{}&\int(0){\0}{}&\int(0){\0}{}&\ln{2}&\cnt{6}&\bounds{2}{2}{1}{false}&&\cr
\endorbit
\orbit{4}{8}{112}
\h{8}&\rat(6/7){6}{}&\rat(26/7){3}{\spec{1}}&\rat(10/7){2}{}&\int(0){\0}{}&&\cnt{660}&\bounds{112}{112}{1}{false}&&\cr
\l{1}&\rat(3/7){3}{}&\rat(13/7){5}{}&\rat(5/7){1}{}&\int(0){\0}{}&\ln{1\sdd}&\cnt{80}&\bounds{16}{16}{1}{false}&&\cr
\l{2}&\int(1){\2}{}&\int(2){\2}{}&\int(0){\0}{}&\int(0){\0}{}&\ln{1}&\cnt{24}&\bounds{4}{4}{1}{false}&&\cr
\l{4}&\rat(5/7){5}{}&\rat(16/7){4}{}&\int(0){\0}{}&\int(0){6}{}&\ln{1}&\cnt{42}&\bounds{6}{6}{1}{false}&&\cr
\l{6}&\rat(5/7){5}{}&\rat(10/7){6}{}&\rat(6/7){4}{}&\int(0){\0}{}&\ln{1}&\cnt{60}&\bounds{10}{10}{1}{false}&&\cr
\l{8}&\rat(6/7){6}{}&\rat(12/7){3}{}&\rat(3/7){2}{}&\int(0){\0}{}&\ln{1}&\cnt{40}&\bounds{8}{8}{1}{false}&&\cr
\l{10}&\rat(6/7){6}{}&\rat(5/7){3}{}&\rat(10/7){2}{}&\int(0){\0}{}&\ln{1}&\cnt{12}&\bounds{4}{4}{1}{false}&&\cr
\l{12}&\rat(6/7){6}{}&\rat(13/7){5}{}&\rat(2/7){6}{}&\int(0){1}{}&\ln{1}&\cnt{70}&\bounds{10}{10}{1}{false}&&\cr
\l{14}&\rat(6/7){6}{}&\rat(10/7){6}{}&\rat(5/7){1}{}&\int(0){5}{}&\ln{1}&\cnt{42}&\bounds{6}{6}{1}{false}&&\cr
\endorbit
\orbit{4}{9}{108}
\h{9}&\int(2){\2}{}&\int(2){\2}{}&\int(2){\2}{}&\int(0){\0}{}&&\cnt{612}&\bounds{108}{108}{1}{false}&&\cr
\l{1}&\int(1){4}{}&\int(1){2}{}&\int(1){6}{}&\int(0){\0}{}&\ln{6\sd}&\cnt{50}&\bounds{10}{10}{1}{false}&&\cr
\l{2}&\int(1){5}{}&\int(1){1}{}&\int(1){1}{}&\int(0){1}{}&\ln{6\sd}&\cnt{35}&\bounds{5}{5}{1}{false}&&\cr
\l{3}&\int(1){6}{}&\int(1){6}{}&\int(1){1}{}&\int(0){5}{}&\ln{2\sd}&\cnt{21}&\bounds{3}{3}{1}{false}&&\cr
\l{4}&\int(2){\2}{}&\int(1){\2}{\spec{1}}&\int(0){\0}{}&\int(0){\0}{}&\ln{6}&\cnt{5}&\bounds{1}{1}{1}{false}&&\cr
\endorbit
\orbit{4}{10}{108}
\h{10}&\rat(6/7){6}{}&\rat(12/7){3}{}&\rat(24/7){2}{\spec{1}}&\int(0){\0}{}&&\cnt{612}&\bounds{108}{108}{1}{false}&&\cr
\l{1}&\int(1){\2}{}&\int(0){\0}{}&\int(2){\2}{}&\int(0){\0}{}&\ln{1}&\cnt{18}&\bounds{3}{3}{1}{false}&&\cr
\l{3}&\rat(4/7){4}{}&\rat(8/7){2}{}&\rat(9/7){6}{}&\int(0){\0}{}&\ln{1}&\cnt{45}&\bounds{9}{9}{1}{false}&&\cr
\l{5}&\rat(5/7){5}{}&\rat(3/7){6}{}&\rat(13/7){4}{}&\int(0){\0}{}&\ln{1}&\cnt{72}&\bounds{12}{12}{1}{false}&&\cr
\l{7}&\rat(6/7){6}{}&\rat(12/7){3}{}&\rat(3/7){2}{}&\int(0){\0}{}&\ln{1}&\cnt{9}&\bounds{3}{3}{1}{false}&&\cr
\l{9}&\rat(6/7){6}{}&\rat(5/7){3}{}&\rat(10/7){2}{}&\int(0){\0}{}&\ln{1}&\cnt{36}&\bounds{9}{9}{1}{false}&&\cr
\l{11}&\rat(6/7){6}{}&\rat(6/7){5}{}&\rat(9/7){6}{}&\int(0){1}{}&\ln{1}&\cnt{42}&\bounds{6}{6}{1}{false}&&\cr
\l{13}&\rat(6/7){6}{}&\rat(4/7){1}{}&\rat(11/7){5}{}&\int(0){6}{}&\ln{1}&\cnt{63}&\bounds{9}{9}{1}{false}&&\cr
\l{15}&\rat(6/7){6}{}&\int(0){\0}{}&\rat(15/7){3}{}&\int(0){2}{}&\ln{1}&\cnt{21}&\bounds{3}{3}{1}{false}&&\cr
\endorbit
\orbit{4}{11}{108}
\h{11}&\rat(6/7){6}{}&\rat(12/7){4}{}&\rat(12/7){4}{}&\rat(12/7){4}{}&&\cnt{468}&\bounds{108}{108}{1}{false}&&\cr
\l{1}&\rat(5/7){5}{}&\rat(12/7){4}{}&\int(0){\0}{}&\rat(4/7){6}{}&\ln{3}&\cnt{18}&\bounds{3}{3}{1}{false}&&\cr
\l{3}&\rat(6/7){6}{}&\rat(9/7){3}{}&\rat(6/7){2}{}&\int(0){\0}{}&\ln{3}&\cnt{24}&\bounds{6}{6}{1}{false}&&\cr
\l{5}&\rat(6/7){6}{}&\rat(8/7){5}{}&\rat(4/7){6}{}&\rat(3/7){1}{}&\ln{3}&\cnt{36}&\bounds{9}{9}{1}{false}&&\cr
\endorbit
\orbit{4}{12}{106}
\h{12}&\rat(20/7){6}{\spec{2}}&\rat(12/7){3}{}&\rat(10/7){2}{}&\int(0){\0}{}&&\cnt{564}&\bounds{106}{114}{fail}{false}&&\cr
\l{1}&\rat(10/7){3}{}&\rat(6/7){5}{}&\rat(5/7){1}{}&\int(0){\0}{}&\ln{1\sdd}&\cnt{72}&\bounds{16}{24}{2}{false}&&\cr
\l{2}&\int(3){\4}{}&\int(0){\0}{}&\int(0){\0}{}&\int(0){\0}{}&\ln{1}&\cnt{4}&\bounds{1}{1}{1}{false}&&\cr
\l{4}&\int(2){\2}{}&\int(1){\2}{}&\int(0){\0}{}&\int(0){\0}{}&\ln{1}&\cnt{24}&\bounds{6}{6}{1}{false}&&\cr
\l{6}&\int(2){\2}{}&\int(0){\0}{}&\int(1){\2}{}&\int(0){\0}{}&\ln{1}&\cnt{20}&\bounds{4}{4}{1}{false}&&\cr
\l{8}&\rat(11/7){4}{}&\rat(8/7){2}{}&\rat(2/7){6}{}&\int(0){\0}{}&\ln{1}&\cnt{60}&\bounds{12}{12}{1}{false}&&\cr
\l{10}&\rat(11/7){4}{}&\int(0){\0}{}&\rat(10/7){2}{}&\int(0){6}{}&\ln{1}&\cnt{28}&\bounds{4}{4}{1}{false}&&\cr
\l{12}&\rat(12/7){5}{}&\rat(9/7){4}{}&\int(0){\0}{}&\int(0){6}{}&\ln{1}&\cnt{28}&\bounds{4}{4}{1}{false}&&\cr
\l{14}&\rat(12/7){5}{}&\rat(4/7){1}{}&\rat(5/7){1}{}&\int(0){1}{}&\ln{1}&\cnt{42}&\bounds{6}{6}{1}{false}&&\cr
\l{16}&\rat(12/7){5}{}&\rat(3/7){6}{}&\rat(6/7){4}{}&\int(0){\0}{}&\ln{1}&\cnt{40}&\bounds{8}{8}{1}{false}&&\cr
\endorbit
\orbit{4}{13}{102}
\h{13}&\rat(20/7){6}{\spec{2}}&\rat(6/7){6}{}&\rat(6/7){1}{}&\rat(10/7){5}{}&&\cnt{564}&\bounds{102}{102}{1}{false}&&\cr
\l{1}&\int(3){\4}{}&\int(0){\0}{}&\int(0){\0}{}&\int(0){\0}{}&\ln{1}&\cnt{4}&\bounds{1}{1}{1}{false}&&\cr
\l{3}&\int(2){\2}{}&\int(1){\2}{}&\int(0){\0}{}&\int(0){\0}{}&\ln{1}&\cnt{12}&\bounds{2}{2}{1}{false}&&\cr
\l{5}&\int(2){\2}{}&\int(0){\0}{}&\int(1){\2}{}&\int(0){\0}{}&\ln{1}&\cnt{12}&\bounds{2}{2}{1}{false}&&\cr
\l{7}&\int(2){\2}{}&\int(0){\0}{}&\int(0){\0}{}&\int(1){\2}{}&\ln{1}&\cnt{20}&\bounds{4}{4}{1}{false}&&\cr
\l{9}&\rat(11/7){4}{}&\rat(6/7){6}{}&\int(0){\0}{}&\rat(4/7){2}{}&\ln{1}&\cnt{40}&\bounds{8}{8}{1}{false}&&\cr
\l{11}&\rat(11/7){4}{}&\int(0){\0}{}&\rat(5/7){2}{}&\rat(5/7){6}{}&\ln{1}&\cnt{48}&\bounds{8}{8}{1}{false}&&\cr
\l{13}&\rat(10/7){3}{}&\rat(5/7){5}{}&\rat(6/7){1}{}&\int(0){\0}{}&\ln{1}&\cnt{36}&\bounds{6}{6}{1}{false}&&\cr
\l{15}&\rat(12/7){5}{}&\rat(4/7){4}{}&\int(0){\0}{}&\rat(5/7){6}{}&\ln{1}&\cnt{30}&\bounds{6}{6}{1}{false}&&\cr
\l{17}&\rat(12/7){5}{}&\rat(6/7){6}{}&\rat(3/7){4}{}&\int(0){\0}{}&\ln{1}&\cnt{20}&\bounds{4}{4}{1}{false}&&\cr
\l{19}&\rat(12/7){5}{}&\rat(1/7){1}{}&\rat(6/7){1}{}&\rat(2/7){1}{}&\ln{1}&\cnt{30}&\bounds{5}{5}{1}{false}&&\cr
\l{21}&\rat(12/7){5}{}&\int(0){\0}{}&\rat(1/7){6}{}&\rat(8/7){4}{}&\ln{1}&\cnt{30}&\bounds{5}{5}{1}{false}&&\cr
\endorbit
\orbit{4}{14}{100}
\h{14}&\rat(6/7){6}{}&\rat(12/7){3}{}&\rat(10/7){2}{}&\int(2){\2}{}&&\cnt{516}&\bounds{100}{100}{1}{false}&&\cr
\l{1}&\int(1){\2}{}&\int(0){\0}{}&\int(0){\0}{}&\int(2){\2}{}&\ln{1}&\cnt{6}&\bounds{1}{1}{1}{false}&&\cr
\l{3}&\rat(4/7){4}{}&\int(0){\0}{}&\rat(10/7){2}{}&\int(1){6}{}&\ln{1}&\cnt{15}&\bounds{3}{3}{1}{false}&&\cr
\l{5}&\rat(5/7){5}{}&\rat(9/7){4}{}&\int(0){\0}{}&\int(1){6}{}&\ln{1}&\cnt{24}&\bounds{4}{4}{1}{false}&&\cr
\l{7}&\rat(5/7){5}{}&\rat(4/7){1}{}&\rat(5/7){1}{}&\int(1){1}{}&\ln{1}&\cnt{36}&\bounds{6}{6}{1}{false}&&\cr
\l{9}&\rat(6/7){6}{}&\rat(12/7){3}{}&\rat(3/7){2}{}&\int(0){\0}{}&\ln{1}&\cnt{10}&\bounds{2}{2}{1}{false}&&\cr
\l{11}&\rat(6/7){6}{}&\rat(5/7){3}{}&\rat(10/7){2}{}&\int(0){\0}{}&\ln{1}&\cnt{12}&\bounds{3}{3}{1}{false}&&\cr
\l{13}&\rat(6/7){6}{}&\rat(8/7){2}{}&\int(0){\0}{}&\int(1){3}{}&\ln{1}&\cnt{30}&\bounds{6}{6}{1}{false}&&\cr
\l{15}&\rat(6/7){6}{}&\rat(6/7){5}{}&\rat(2/7){6}{}&\int(1){1}{}&\ln{1}&\cnt{30}&\bounds{6}{6}{1}{false}&&\cr
\l{17}&\rat(6/7){6}{}&\rat(4/7){1}{}&\rat(4/7){5}{}&\int(1){6}{}&\ln{1}&\cnt{30}&\bounds{6}{6}{1}{false}&&\cr
\l{19}&\rat(6/7){6}{}&\rat(3/7){6}{}&\rat(5/7){1}{}&\int(1){5}{}&\ln{1}&\cnt{40}&\bounds{8}{8}{1}{false}&&\cr
\l{21}&\rat(6/7){6}{}&\int(0){\0}{}&\rat(8/7){3}{}&\int(1){2}{}&\ln{1}&\cnt{25}&\bounds{5}{5}{1}{false}&&\cr
\endorbit
\orbit{4}{15}{96}
\h{15}&\rat(10/7){5}{}&\rat(10/7){5}{}&\rat(10/7){2}{}&\rat(12/7){3}{}&&\cnt{468}&\bounds{96}{100}{fail}{false}&&\cr
\l{1}&\rat(5/7){6}{}&\rat(5/7){6}{}&\rat(5/7){1}{}&\rat(6/7){5}{}&\ln{1\sdd}&\cnt{48}&\bounds{12}{16}{2}{false}&&\cr
\l{2}&\rat(8/7){4}{}&\rat(5/7){6}{}&\int(0){\0}{}&\rat(8/7){2}{}&\ln{3}&\cnt{30}&\bounds{6}{6}{1}{false}&&\cr
\l{4}&\rat(8/7){4}{}&\int(0){\0}{}&\rat(10/7){2}{}&\rat(3/7){6}{}&\ln{3}&\cnt{20}&\bounds{4}{4}{1}{false}&&\cr
\l{6}&\rat(6/7){3}{}&\rat(10/7){5}{}&\rat(5/7){1}{}&\int(0){\0}{}&\ln{3}&\cnt{20}&\bounds{4}{4}{1}{false}&&\cr
\endorbit
\orbit{4}{16}{50}
\h{16}&\int(6){\6}{\spec{1}}&\int(0){\0}{}&\int(0){\0}{}&\int(0){\0}{}&&\cnt{276}&\bounds{50}{50}{1}{false}&&\cr
\l{1}&\int(3){\4}{}&\int(0){\0}{}&\int(0){\0}{}&\int(0){\0}{}&\ln{1\sdd}&\cnt{24}&\bounds{8}{8}{1}{false}&&\cr
\l{2}&\int(3){\2}{}&\int(0){\2}{}&\int(0){\0}{}&\int(0){\0}{}&\ln{3\sdd}&\cnt{84}&\bounds{14}{14}{1}{false}&&\cr
\endorbit
\endorbits

\newcs{check-4A5_D4-1}{}
\newcs{check-4A5_D4-3}{}
\newcs{check-4A5_D4-44}{\noexpand\checkdone}
\newcs{check-4A5_D4-47}{\noexpand\checkdone}
\newcs{max-4A5_D4}{156}%
\newcs{count-4A5_D4}{93}%
\newcs{counts-4A5_D4}{\DO{1}{156}\DO{2}{120}\DO{3}{126}\DO{4}{108}\DO{5}{112}\DO{6}{106}\DO{7}{118}\DO{8}{114}\DO{9}{108}\DO{10}{92}\DO{11}{88}\DO{12}{112}\DO{13}{116}\DO{14}{108}\DO{15}{92}\DO{16}{112}\DO{17}{108}\DO{18}{108}\DO{19}{100}\DO{20}{108}\DO{21}{112}\DO{22}{104}\DO{23}{116}\DO{24}{118}\DO{25}{104}\DO{26}{96}\DO{27}{110}\DO{28}{102}\DO{29}{116}\DO{30}{104}\DO{31}{88}\DO{32}{100}\DO{33}{112}\DO{34}{110}\DO{35}{112}\DO{36}{114}\DO{37}{88}\DO{38}{82}\DO{39}{118}\DO{40}{118}\DO{41}{112}\DO{42}{98}\DO{43}{82}\DO{44}{122}\DO{45}{104}\DO{46}{92}\DO{47}{122}\DO{48}{108}\DO{49}{104}\DO{50}{92}\DO{51}{112}\DO{52}{118}\DO{53}{98}\DO{54}{90}\DO{55}{114}\DO{56}{108}\DO{57}{116}\DO{58}{108}\DO{59}{96}\DO{60}{92}\DO{61}{108}\DO{62}{114}\DO{63}{112}\DO{64}{100}\DO{65}{96}\DO{66}{92}\DO{67}{88}\DO{68}{84}\DO{69}{112}\DO{70}{104}\DO{71}{100}\DO{72}{92}\DO{73}{88}\DO{74}{118}\DO{75}{108}\DO{76}{104}\DO{77}{96}\DO{78}{80}\DO{79}{116}\DO{80}{106}\DO{81}{104}\DO{82}{92}\DO{83}{84}\DO{84}{108}\DO{85}{96}\DO{86}{96}\DO{87}{92}\DO{88}{80}\DO{89}{106}\DO{90}{80}\DO{91}{76}\DO{92}{52}\DO{93}{48}}%
\orbits
\orbit{5}{1}{156}
\h{1}&\rat(3/2){3}{}&\rat(3/2){3}{}&\rat(3/2){3}{}&\rat(3/2){3}{}&\root\int(0){\0}{}&&\cnt{456}&\bounds{156}{156}{1}{false}&&\cr
\l{1}&\rat(3/2){3}{}&\rat(3/2){3}{}&\int(0){\0}{}&\int(0){\0}{}&\int(0){1}{}&\ln{6\sd}&\cnt{4}&\bounds{2}{2}{1}{true}&&\cr
\l{2}&\rat(3/2){3}{}&\rat(1/2){5}{}&\rat(1/2){5}{}&\rat(1/2){5}{}&\int(0){\0}{}&\ln{8}&\cnt{27}&\bounds{9}{9}{1}{true}&&\cr
\endorbit
\orbit{5}{2}{120}
\h{2}&\root\rat(3/2){3}{\spec{1}}&\rat(3/2){3}{}&\rat(3/2){3}{}&\rat(3/2){3}{}&\int(0){\0}{}&&\cnt{372}&\bounds{120}{120}{1}{false}&&\cr
\l{1}&\rat(3/2){3}{}&\rat(3/2){3}{}&\int(0){\0}{}&\int(0){\0}{}&\int(0){1}{}&\ln{3}&\cnt{8}&\bounds{2}{2}{1}{true}&&\cr
\l{3}&\rat(3/2){3}{}&\rat(1/2){5}{}&\rat(1/2){5}{}&\rat(1/2){5}{}&\int(0){\0}{}&\ln{2}&\cnt{27}&\bounds{9}{9}{1}{true}&&\cr
\l{5}&\int(1){4}{}&\int(1){4}{}&\int(1){2}{}&\int(0){\0}{}&\int(0){\0}{}&\ln{3}&\cnt{27}&\bounds{9}{9}{1}{true}&&\cr
\l{7}&\int(1){2}{}&\int(1){4}{}&\int(0){\0}{}&\int(1){2}{}&\int(0){\0}{}&\ln{3}&\cnt{9}&\bounds{3}{3}{1}{true}&&\cr
\endorbit
\orbit{5}{3}{126}
\h{3}&\root\int(6){\6}{\spec{1}}&\int(0){\0}{}&\int(0){\0}{}&\int(0){\0}{}&\int(0){\0}{}&&\cnt{912}&\bounds{126}{168}{fail}{false}&&\cr
\l{1}&\int(3){3}{}&\int(0){3}{}&\int(0){\0}{}&\int(0){\0}{}&\int(0){1}{}&\ln{3\sdd}&\cnt{160}&\bounds{22}{32}{5}{true}&&\cr
\l{2}&\int(3){3}{}&\int(0){5}{}&\int(0){5}{}&\int(0){5}{}&\int(0){\0}{}&\ln{2\sd}&\cnt{216}&\bounds{30}{36}{3}{false}&&\cr
\endorbit
\orbit{5}{4}{108}
\h{4}&\int(6){\6}{}&\int(0){\0}{}&\int(0){\0}{}&\int(0){\0}{}&\root\int(0){\0}{}&&\cnt{672}&\bounds{108}{120}{fail}{false}&&\cr
\l{1}&\int(3){3}{}&\int(0){3}{}&\int(0){\0}{}&\int(0){\0}{}&\int(0){1}{}&\ln{3\sdd}&\cnt{80}&\bounds{12}{16}{5}{true}&&\cr
\l{2}&\int(3){3}{}&\int(0){5}{}&\int(0){5}{}&\int(0){5}{}&\int(0){\0}{}&\ln{2\sd}&\cnt{216}&\bounds{36}{36}{1}{false}&&\cr
\endorbit
\orbit{5}{5}{112}
\h{5}&\int(6){\6}{}&\root\int(0){\0}{}&\int(0){\0}{}&\int(0){\0}{}&\int(0){\0}{}&&\cnt{672}&\bounds{112}{128}{fail}{false}&&\cr
\l{1}&\int(3){3}{}&\int(0){3}{\spec{1}}&\int(0){\0}{}&\int(0){\0}{}&\int(0){1}{}&\ln{2\sd}&\cnt{32}&\bounds{8}{8}{1}{false}&&\cr
\l{2}&\int(3){3}{}&\int(0){5}{}&\int(0){5}{}&\int(0){5}{}&\int(0){\0}{}&\ln{2\sd}&\cnt{144}&\bounds{24}{24}{1}{false}&&\cr
\l{3}&\int(3){3}{}&\int(0){\0}{}&\int(0){3}{}&\int(0){\0}{}&\int(0){2}{}&\ln{2\sdd}&\cnt{160}&\bounds{24}{32}{2}{false}&&\cr
\endorbit
\orbit{5}{6}{106}
\h{6}&\rat(10/3){4}{\spec{2}}&\rat(4/3){4}{}&\rat(4/3){2}{}&\root\int(0){\0}{}&\int(0){\0}{}&&\cnt{528}&\bounds{106}{142}{fail}{false}&&\cr
\l{1}&\rat(5/3){2}{}&\rat(4/3){4}{}&\int(0){\0}{}&\int(0){2}{\spec{1}}&\int(0){\0}{}&\ln{2\sd}&\cnt{2}&\bounds{1}{1}{1}{true}&&\cr
\l{2}&\rat(5/3){2}{}&\rat(4/3){4}{}&\int(0){\0}{}&\int(0){2}{\spec{2}}&\int(0){\0}{}&\ln{2\sd}&\cnt{12}&\bounds{3}{4}{5}{true}&&\cr
\l{3}&\rat(5/3){2}{}&\rat(2/3){2}{}&\rat(2/3){4}{}&\int(0){\0}{}&\int(0){\0}{}&\ln{1\sdd}&\cnt{72}&\bounds{12}{24}{5}{true}&&\cr
\l{4}&\rat(5/3){2}{}&\rat(2/3){5}{}&\rat(2/3){1}{}&\int(0){\0}{}&\int(0){3}{}&\ln{1\sdd}&\cnt{64}&\bounds{10}{16}{5}{true}&&\cr
\l{5}&\int(3){\4}{}&\int(0){\0}{}&\int(0){\0}{}&\int(0){\0}{}&\int(0){\0}{}&\ln{1}&\cnt{6}&\bounds{2}{2}{1}{true}&&\cr
\l{7}&\int(2){\2}{}&\int(1){\2}{}&\int(0){\0}{}&\int(0){\0}{}&\int(0){\0}{}&\ln{2}&\cnt{24}&\bounds{6}{6}{1}{true}&&\cr
\l{9}&\int(2){3}{}&\int(1){3}{}&\int(0){\0}{}&\int(0){\0}{}&\int(0){1}{}&\ln{2}&\cnt{32}&\bounds{5}{8}{5}{true}&&\cr
\l{11}&\int(2){3}{}&\rat(2/3){5}{}&\rat(1/3){5}{}&\int(0){5}{}&\int(0){\0}{}&\ln{2}&\cnt{32}&\bounds{7}{8}{5}{true}&&\cr
\endorbit
\orbit{5}{7}{118}
\h{7}&\root\rat(10/3){4}{\spec{2}}&\rat(4/3){4}{}&\rat(4/3){2}{}&\int(0){\0}{}&\int(0){\0}{}&&\cnt{552}&\bounds{118}{126}{fail}{false}&&\cr
\l{1}&\rat(5/3){2}{}&\rat(4/3){4}{}&\int(0){\0}{}&\int(0){2}{}&\int(0){\0}{}&\ln{2\sd}&\cnt{30}&\bounds{6}{6}{1}{false}&&\cr
\l{2}&\rat(5/3){2}{}&\rat(2/3){2}{}&\rat(2/3){4}{}&\int(0){\0}{}&\int(0){\0}{}&\ln{1\sdd}&\cnt{72}&\bounds{16}{24}{2}{false}&&\cr
\l{3}&\rat(5/3){2}{}&\rat(2/3){5}{}&\rat(2/3){1}{}&\int(0){\0}{}&\int(0){3}{}&\ln{1\sdd}&\cnt{64}&\bounds{16}{16}{1}{false}&&\cr
\l{4}&\int(3){\4}{}&\int(0){\0}{}&\int(0){\0}{}&\int(0){\0}{}&\int(0){\0}{}&\ln{1}&\cnt{2}&\bounds{1}{1}{1}{false}&&\cr
\l{6}&\int(2){\2}{}&\int(1){\2}{}&\int(0){\0}{}&\int(0){\0}{}&\int(0){\0}{}&\ln{2}&\cnt{8}&\bounds{2}{2}{1}{false}&&\cr
\l{8}&\int(2){3}{}&\int(1){3}{}&\int(0){\0}{}&\int(0){\0}{}&\int(0){1}{}&\ln{2}&\cnt{32}&\bounds{8}{8}{1}{false}&&\cr
\l{10}&\int(2){3}{}&\rat(2/3){5}{}&\rat(1/3){5}{}&\int(0){5}{}&\int(0){\0}{}&\ln{2}&\cnt{48}&\bounds{8}{8}{1}{false}&&\cr
\endorbit
\orbit{5}{8}{114}
\h{8}&\rat(10/3){4}{\spec{2}}&\rat(4/3){4}{}&\rat(4/3){2}{}&\int(0){\0}{}&\root\int(0){\0}{}&&\cnt{528}&\bounds{114}{128}{fail}{false}&&\cr
\l{1}&\rat(5/3){2}{}&\rat(4/3){4}{}&\int(0){\0}{}&\int(0){2}{}&\int(0){\0}{}&\ln{2\sd}&\cnt{30}&\bounds{6}{6}{1}{false}&&\cr
\l{2}&\rat(5/3){2}{}&\rat(2/3){2}{}&\rat(2/3){4}{}&\int(0){\0}{}&\int(0){\0}{}&\ln{1\sdd}&\cnt{72}&\bounds{16}{24}{2}{false}&&\cr
\l{3}&\rat(5/3){2}{}&\rat(2/3){5}{}&\rat(2/3){1}{}&\int(0){\0}{}&\int(0){3}{}&\ln{1\sdd}&\cnt{32}&\bounds{10}{16}{2}{false}&&\cr
\l{4}&\int(3){\4}{}&\int(0){\0}{}&\int(0){\0}{}&\int(0){\0}{}&\int(0){\0}{}&\ln{1}&\cnt{6}&\bounds{2}{2}{1}{false}&&\cr
\l{6}&\int(2){\2}{}&\int(1){\2}{}&\int(0){\0}{}&\int(0){\0}{}&\int(0){\0}{}&\ln{2}&\cnt{24}&\bounds{6}{6}{1}{false}&&\cr
\l{8}&\int(2){3}{}&\int(1){3}{}&\int(0){\0}{}&\int(0){\0}{}&\int(0){1}{}&\ln{2}&\cnt{16}&\bounds{4}{4}{1}{false}&&\cr
\l{10}&\int(2){3}{}&\rat(2/3){5}{}&\rat(1/3){5}{}&\int(0){5}{}&\int(0){\0}{}&\ln{2}&\cnt{48}&\bounds{8}{8}{1}{false}&&\cr
\endorbit
\orbit{5}{9}{108}
\h{9}&\rat(10/3){4}{\spec{2}}&\root\rat(4/3){4}{\spec{1}}&\rat(4/3){2}{}&\int(0){\0}{}&\int(0){\0}{}&&\cnt{472}&\bounds{108}{108}{1}{false}&&\cr
\l{1}&\rat(5/3){2}{}&\rat(2/3){5}{}&\rat(2/3){1}{}&\int(0){\0}{}&\int(0){3}{}&\ln{1\sdd}&\cnt{64}&\bounds{16}{16}{1}{false}&&\cr
\l{2}&\int(3){\4}{}&\int(0){\0}{}&\int(0){\0}{}&\int(0){\0}{}&\int(0){\0}{}&\ln{1}&\cnt{6}&\bounds{2}{2}{1}{false}&&\cr
\l{4}&\int(2){\2}{}&\int(1){\2}{}&\int(0){\0}{}&\int(0){\0}{}&\int(0){\0}{}&\ln{1}&\cnt{12}&\bounds{4}{4}{1}{false}&&\cr
\l{6}&\int(2){\2}{}&\int(0){\0}{}&\int(1){\2}{}&\int(0){\0}{}&\int(0){\0}{}&\ln{1}&\cnt{24}&\bounds{6}{6}{1}{false}&&\cr
\l{8}&\int(2){3}{}&\int(1){3}{}&\int(0){\0}{}&\int(0){\0}{}&\int(0){1}{}&\ln{1}&\cnt{16}&\bounds{4}{4}{1}{false}&&\cr
\l{10}&\int(2){3}{}&\rat(2/3){5}{}&\rat(1/3){5}{}&\int(0){5}{}&\int(0){\0}{}&\ln{1}&\cnt{48}&\bounds{8}{8}{1}{false}&&\cr
\l{12}&\int(2){3}{}&\rat(1/3){1}{}&\rat(2/3){1}{}&\int(0){1}{}&\int(0){\0}{}&\ln{1}&\cnt{24}&\bounds{4}{4}{1}{false}&&\cr
\l{14}&\int(2){3}{}&\int(0){\0}{}&\int(1){3}{}&\int(0){\0}{}&\int(0){2}{}&\ln{1}&\cnt{32}&\bounds{8}{8}{1}{false}&&\cr
\l{16}&\rat(5/3){2}{}&\rat(4/3){4}{}&\int(0){\0}{}&\int(0){2}{}&\int(0){\0}{}&\ln{1}&\cnt{30}&\bounds{6}{6}{1}{false}&&\cr
\l{18}&\rat(5/3){2}{}&\rat(2/3){2}{\spec{1}}&\rat(2/3){4}{}&\int(0){\0}{}&\int(0){\0}{}&\ln{1}&\cnt{12}&\bounds{4}{4}{1}{false}&&\cr
\endorbit
\orbit{5}{10}{92}
\h{10}&\rat(10/3){4}{\spec{2}}&\root\rat(4/3){4}{\spec{2}}&\rat(4/3){2}{}&\int(0){\0}{}&\int(0){\0}{}&&\cnt{416}&\bounds{92}{100}{fail}{false}&&\cr
\l{1}&\rat(5/3){2}{}&\rat(2/3){2}{}&\rat(2/3){4}{}&\int(0){\0}{}&\int(0){\0}{}&\ln{1\sdd}&\cnt{72}&\bounds{16}{24}{2}{false}&&\cr
\l{2}&\int(3){\4}{}&\int(0){\0}{}&\int(0){\0}{}&\int(0){\0}{}&\int(0){\0}{}&\ln{1}&\cnt{6}&\bounds{2}{2}{1}{false}&&\cr
\l{4}&\int(2){\2}{}&\int(0){\0}{}&\int(1){\2}{}&\int(0){\0}{}&\int(0){\0}{}&\ln{1}&\cnt{24}&\bounds{6}{6}{1}{false}&&\cr
\l{6}&\int(2){3}{}&\int(1){3}{}&\int(0){\0}{}&\int(0){\0}{}&\int(0){1}{}&\ln{1}&\cnt{32}&\bounds{8}{8}{1}{false}&&\cr
\l{8}&\int(2){3}{}&\rat(1/3){1}{}&\rat(2/3){1}{}&\int(0){1}{}&\int(0){\0}{}&\ln{1}&\cnt{48}&\bounds{8}{8}{1}{false}&&\cr
\l{10}&\int(2){3}{}&\int(0){\0}{}&\int(1){3}{}&\int(0){\0}{}&\int(0){2}{}&\ln{1}&\cnt{32}&\bounds{8}{8}{1}{false}&&\cr
\l{12}&\rat(5/3){2}{}&\rat(4/3){4}{}&\int(0){\0}{}&\int(0){2}{}&\int(0){\0}{}&\ln{1}&\cnt{30}&\bounds{6}{6}{1}{false}&&\cr
\endorbit
\orbit{5}{11}{88}
\h{11}&\root\rat(10/3){4}{\spec{1}}&\rat(4/3){4}{}&\rat(4/3){2}{}&\int(0){\0}{}&\int(0){\0}{}&&\cnt{416}&\bounds{88}{88}{1}{false}&&\cr
\l{1}&\int(2){\2}{}&\int(1){\2}{}&\int(0){\0}{}&\int(0){\0}{}&\int(0){\0}{}&\ln{2}&\cnt{24}&\bounds{6}{6}{1}{false}&&\cr
\l{3}&\int(2){3}{}&\int(1){3}{}&\int(0){\0}{}&\int(0){\0}{}&\int(0){1}{}&\ln{2}&\cnt{32}&\bounds{8}{8}{1}{false}&&\cr
\l{5}&\int(2){3}{}&\rat(2/3){5}{}&\rat(1/3){5}{}&\int(0){5}{}&\int(0){\0}{}&\ln{2}&\cnt{48}&\bounds{8}{8}{1}{false}&&\cr
\endorbit
\orbit{5}{12}{112}
\h{12}&\rat(4/3){4}{}&\rat(4/3){4}{}&\rat(4/3){2}{}&\int(2){\2}{}&\root\int(0){\0}{}&&\cnt{480}&\bounds{112}{136}{fail}{false}&&\cr
\l{1}&\int(1){3}{}&\rat(2/3){5}{}&\rat(1/3){5}{}&\int(1){5}{}&\int(0){\0}{}&\ln{6\sd}&\cnt{32}&\bounds{7}{8}{5}{true}&&\cr
\l{2}&\rat(4/3){4}{}&\rat(2/3){2}{}&\int(0){\0}{}&\int(1){4}{}&\int(0){\0}{}&\ln{6\sd}&\cnt{24}&\bounds{5}{6}{5}{true}&&\cr
\l{3}&\rat(4/3){4}{}&\rat(2/3){5}{}&\int(0){\0}{}&\int(1){1}{}&\int(0){2}{}&\ln{6\sd}&\cnt{8}&\bounds{3}{4}{5}{true}&&\cr
\l{4}&\rat(2/3){5}{}&\rat(2/3){5}{}&\rat(2/3){1}{}&\int(1){3}{}&\int(0){\0}{}&\ln{1\sdd}&\cnt{48}&\bounds{10}{16}{5}{true}&&\cr
\l{5}&\int(1){\2}{}&\int(0){\0}{}&\int(0){\0}{}&\int(2){\2}{}&\int(0){\0}{}&\ln{3}&\cnt{8}&\bounds{2}{2}{1}{true}&&\cr
\endorbit
\orbit{5}{13}{116}
\h{13}&\rat(4/3){4}{}&\rat(4/3){4}{}&\rat(4/3){2}{}&\root\int(2){\2}{}&\int(0){\0}{}&&\cnt{424}&\bounds{116}{116}{1}{false}&&\cr
\l{1}&\int(1){3}{}&\rat(2/3){5}{}&\rat(1/3){5}{}&\int(1){5}{}&\int(0){\0}{}&\ln{6\sd}&\cnt{32}&\bounds{8}{8}{1}{false}&&\cr
\l{2}&\rat(4/3){4}{}&\rat(2/3){2}{}&\int(0){\0}{}&\int(1){4}{}&\int(0){\0}{}&\ln{6\sd}&\cnt{12}&\bounds{4}{4}{1}{false}&&\cr
\l{3}&\rat(4/3){4}{}&\rat(2/3){5}{}&\int(0){\0}{}&\int(1){1}{}&\int(0){2}{}&\ln{6\sd}&\cnt{16}&\bounds{4}{4}{1}{false}&&\cr
\l{4}&\rat(2/3){5}{}&\rat(2/3){5}{}&\rat(2/3){1}{}&\int(1){3}{\spec{1}}&\int(0){\0}{}&\ln{2\sd}&\cnt{8}&\bounds{4}{4}{1}{false}&&\cr
\l{5}&\int(1){\2}{}&\int(0){\0}{}&\int(0){\0}{}&\int(2){\2}{}&\int(0){\0}{}&\ln{3}&\cnt{8}&\bounds{2}{2}{1}{false}&&\cr
\endorbit
\orbit{5}{14}{108}
\h{14}&\root\rat(4/3){4}{\spec{1}}&\rat(4/3){4}{}&\rat(4/3){2}{}&\int(2){\2}{}&\int(0){\0}{}&&\cnt{424}&\bounds{108}{112}{fail}{false}&&\cr
\l{1}&\rat(2/3){5}{}&\int(1){3}{}&\rat(1/3){5}{}&\int(1){1}{}&\int(0){\0}{}&\ln{2\sd}&\cnt{32}&\bounds{8}{8}{1}{false}&&\cr
\l{2}&\rat(2/3){5}{}&\rat(4/3){4}{}&\int(0){\0}{}&\int(1){5}{}&\int(0){3}{}&\ln{2\sd}&\cnt{16}&\bounds{4}{4}{1}{false}&&\cr
\l{3}&\rat(2/3){5}{}&\rat(2/3){5}{}&\rat(2/3){1}{}&\int(1){3}{}&\int(0){\0}{}&\ln{1\sdd}&\cnt{48}&\bounds{12}{16}{2}{false}&&\cr
\l{4}&\int(1){\2}{}&\int(0){\0}{}&\int(0){\0}{}&\int(2){\2}{}&\int(0){\0}{}&\ln{1}&\cnt{4}&\bounds{2}{2}{1}{false}&&\cr
\l{6}&\int(1){3}{}&\rat(2/3){5}{}&\rat(1/3){5}{}&\int(1){5}{}&\int(0){\0}{}&\ln{2}&\cnt{16}&\bounds{4}{4}{1}{false}&&\cr
\l{8}&\rat(4/3){4}{}&\rat(4/3){4}{}&\rat(1/3){2}{}&\int(0){\0}{}&\int(0){\0}{}&\ln{2}&\cnt{8}&\bounds{2}{2}{1}{false}&&\cr
\l{10}&\rat(4/3){4}{}&\rat(2/3){2}{}&\int(0){\0}{}&\int(1){4}{}&\int(0){\0}{}&\ln{2}&\cnt{24}&\bounds{6}{6}{1}{false}&&\cr
\l{12}&\rat(4/3){4}{}&\rat(2/3){5}{}&\int(0){\0}{}&\int(1){1}{}&\int(0){2}{}&\ln{2}&\cnt{16}&\bounds{4}{4}{1}{false}&&\cr
\l{14}&\rat(2/3){2}{\spec{1}}&\rat(4/3){4}{}&\int(0){\0}{}&\int(1){2}{}&\int(0){\0}{}&\ln{2}&\cnt{4}&\bounds{1}{1}{1}{false}&&\cr
\endorbit
\orbit{5}{15}{92}
\h{15}&\root\rat(4/3){4}{\spec{2}}&\rat(4/3){4}{}&\rat(4/3){2}{}&\int(2){\2}{}&\int(0){\0}{}&&\cnt{368}&\bounds{92}{92}{1}{false}&&\cr
\l{1}&\rat(2/3){2}{}&\rat(4/3){4}{}&\int(0){\0}{}&\int(1){2}{}&\int(0){\0}{}&\ln{2\sd}&\cnt{24}&\bounds{6}{6}{1}{false}&&\cr
\l{2}&\int(1){3}{}&\rat(2/3){5}{}&\rat(1/3){5}{}&\int(1){5}{}&\int(0){\0}{}&\ln{2}&\cnt{32}&\bounds{8}{8}{1}{false}&&\cr
\l{4}&\rat(4/3){4}{}&\rat(4/3){4}{}&\rat(1/3){2}{}&\int(0){\0}{}&\int(0){\0}{}&\ln{2}&\cnt{8}&\bounds{2}{2}{1}{false}&&\cr
\l{6}&\rat(4/3){4}{}&\rat(2/3){2}{}&\int(0){\0}{}&\int(1){4}{}&\int(0){\0}{}&\ln{2}&\cnt{24}&\bounds{6}{6}{1}{false}&&\cr
\l{8}&\rat(4/3){4}{}&\rat(2/3){5}{}&\int(0){\0}{}&\int(1){1}{}&\int(0){2}{}&\ln{2}&\cnt{16}&\bounds{4}{4}{1}{false}&&\cr
\endorbit
\orbit{5}{16}{112}
\h{16}&\int(2){\2}{}&\int(0){\0}{}&\int(0){\0}{}&\int(0){\0}{}&\root\int(4){\4}{3}&&\cnt{784}&\bounds{112}{136}{fail}{false}&&\cr
\l{1}&\int(1){\2}{\spec{1}}&\int(0){\0}{}&\int(0){\0}{}&\int(0){\0}{}&\int(2){\2}{\spec{1}}&\ln{4\sd}&\cnt{4}&\bounds{1}{1}{1}{true}&&\cr
\l{2}&\int(1){3}{}&\int(0){\0}{}&\int(0){\0}{}&\int(0){3}{}&\int(2){3}{}&\ln{1\sdd}&\cnt{120}&\bounds{18}{24}{5}{true}&&\cr
\l{3}&\int(1){4}{}&\int(0){1}{}&\int(0){5}{}&\int(0){\0}{}&\int(2){3}{}&\ln{2\sd}&\cnt{144}&\bounds{21}{24}{3}{false}&&\cr
\l{4}&\int(1){5}{}&\int(0){4}{}&\int(0){\0}{}&\int(0){5}{}&\int(2){3}{}&\ln{4\sd}&\cnt{90}&\bounds{12}{15}{3}{false}&&\cr
\endorbit
\orbit{5}{17}{108}
\h{17}&\int(2){\2}{}&\int(0){\0}{}&\int(0){\0}{}&\root\int(0){\0}{}&\int(4){\4}{3}&&\cnt{624}&\bounds{108}{120}{fail}{false}&&\cr
\l{1}&\int(1){\2}{\spec{1}}&\int(0){\0}{}&\int(0){\0}{}&\int(0){\0}{}&\int(2){\2}{}&\ln{2\sd}&\cnt{24}&\bounds{5}{6}{5}{true}&&\cr
\l{2}&\int(1){3}{}&\int(0){\0}{}&\int(0){\0}{}&\int(0){3}{\spec{1}}&\int(2){3}{}&\ln{2\sd}&\cnt{24}&\bounds{5}{6}{5}{true}&&\cr
\l{3}&\int(1){4}{}&\int(0){1}{}&\int(0){5}{}&\int(0){\0}{}&\int(2){3}{}&\ln{2\sd}&\cnt{144}&\bounds{24}{24}{1}{false}&&\cr
\l{4}&\int(1){5}{}&\int(0){4}{}&\int(0){\0}{}&\int(0){5}{}&\int(2){3}{}&\ln{4\sd}&\cnt{60}&\bounds{10}{12}{5}{true}&&\cr
\endorbit
\orbit{5}{18}{108}
\h{18}&\int(2){\2}{}&\root\int(0){\0}{}&\int(0){\0}{}&\int(0){\0}{}&\int(4){\4}{3}&&\cnt{624}&\bounds{108}{112}{fail}{false}&&\cr
\l{1}&\int(1){\2}{\spec{1}}&\int(0){\0}{}&\int(0){\0}{}&\int(0){\0}{}&\int(2){\2}{}&\ln{2\sd}&\cnt{24}&\bounds{6}{6}{1}{false}&&\cr
\l{2}&\int(1){3}{}&\int(0){\0}{}&\int(0){\0}{}&\int(0){3}{}&\int(2){3}{}&\ln{1\sdd}&\cnt{120}&\bounds{20}{24}{2}{false}&&\cr
\l{3}&\int(1){4}{}&\int(0){1}{}&\int(0){5}{}&\int(0){\0}{}&\int(2){3}{}&\ln{2\sd}&\cnt{96}&\bounds{16}{16}{1}{false}&&\cr
\l{4}&\int(1){5}{}&\int(0){4}{\spec{1}}&\int(0){\0}{}&\int(0){5}{}&\int(2){3}{}&\ln{2\sd}&\cnt{36}&\bounds{6}{6}{1}{false}&&\cr
\l{5}&\int(1){5}{}&\int(0){4}{\spec{2}}&\int(0){\0}{}&\int(0){5}{}&\int(2){3}{}&\ln{2\sd}&\cnt{6}&\bounds{1}{1}{1}{false}&&\cr
\l{6}&\int(1){5}{}&\int(0){\0}{}&\int(0){2}{}&\int(0){1}{}&\int(2){3}{}&\ln{2\sd}&\cnt{90}&\bounds{15}{15}{1}{false}&&\cr
\endorbit
\orbit{5}{19}{100}
\h{19}&\root\int(2){\2}{}&\int(0){\0}{}&\int(0){\0}{}&\int(0){\0}{}&\int(4){\4}{3}&&\cnt{568}&\bounds{100}{100}{1}{false}&&\cr
\l{1}&\int(1){\2}{\spec{1}}&\int(0){\0}{}&\int(0){\0}{}&\int(0){\0}{}&\int(2){\2}{}&\ln{2\sd}&\cnt{12}&\bounds{4}{4}{1}{false}&&\cr
\l{2}&\int(1){3}{\spec{1}}&\int(0){\0}{}&\int(0){\0}{}&\int(0){3}{}&\int(2){3}{}&\ln{2\sd}&\cnt{20}&\bounds{4}{4}{1}{false}&&\cr
\l{3}&\int(1){4}{}&\int(0){1}{}&\int(0){5}{}&\int(0){\0}{}&\int(2){3}{}&\ln{2\sd}&\cnt{72}&\bounds{12}{12}{1}{false}&&\cr
\l{4}&\int(1){5}{}&\int(0){4}{}&\int(0){\0}{}&\int(0){5}{}&\int(2){3}{}&\ln{4\sd}&\cnt{90}&\bounds{15}{15}{1}{false}&&\cr
\endorbit
\orbit{5}{20}{108}
\h{20}&\int(2){\2}{}&\int(2){\2}{}&\int(2){\2}{}&\root\int(0){\0}{}&\int(0){\0}{}&&\cnt{528}&\bounds{108}{132}{fail}{false}&&\cr
\l{1}&\int(2){\2}{}&\int(1){\2}{\spec{1}}&\int(0){\0}{}&\int(0){\0}{}&\int(0){\0}{}&\ln{12\sd}&\cnt{4}&\bounds{1}{1}{1}{true}&&\cr
\l{2}&\int(1){3}{}&\int(1){5}{}&\int(1){5}{}&\int(0){5}{}&\int(0){\0}{}&\ln{6\sd}&\cnt{24}&\bounds{5}{6}{5}{true}&&\cr
\l{3}&\int(1){4}{}&\int(1){4}{}&\int(1){2}{}&\int(0){\0}{}&\int(0){\0}{}&\ln{2\sd}&\cnt{64}&\bounds{16}{16}{1}{false}&&\cr
\l{4}&\int(1){4}{}&\int(1){1}{}&\int(1){5}{}&\int(0){\0}{}&\int(0){3}{}&\ln{6\sd}&\cnt{32}&\bounds{5}{8}{5}{true}&&\cr
\l{5}&\int(1){5}{}&\int(1){5}{}&\int(1){1}{}&\int(0){3}{\spec{1}}&\int(0){\0}{}&\ln{4\sd}&\cnt{4}&\bounds{1}{1}{1}{true}&&\cr
\endorbit
\orbit{5}{21}{112}
\h{21}&\int(2){\2}{}&\int(2){\2}{}&\int(2){\2}{}&\int(0){\0}{}&\root\int(0){\0}{}&&\cnt{528}&\bounds{112}{112}{1}{false}&&\cr
\l{1}&\int(2){\2}{}&\int(1){\2}{\spec{1}}&\int(0){\0}{}&\int(0){\0}{}&\int(0){\0}{}&\ln{12\sd}&\cnt{4}&\bounds{1}{1}{1}{false}&&\cr
\l{2}&\int(1){3}{}&\int(1){5}{}&\int(1){5}{}&\int(0){5}{}&\int(0){\0}{}&\ln{6\sd}&\cnt{36}&\bounds{6}{6}{1}{false}&&\cr
\l{3}&\int(1){4}{}&\int(1){4}{}&\int(1){2}{}&\int(0){\0}{}&\int(0){\0}{}&\ln{2\sd}&\cnt{64}&\bounds{16}{16}{1}{false}&&\cr
\l{4}&\int(1){4}{}&\int(1){1}{}&\int(1){5}{}&\int(0){\0}{}&\int(0){3}{}&\ln{6\sd}&\cnt{16}&\bounds{4}{4}{1}{false}&&\cr
\l{5}&\int(1){5}{}&\int(1){5}{}&\int(1){1}{}&\int(0){3}{}&\int(0){\0}{}&\ln{2\sd}&\cnt{20}&\bounds{4}{4}{1}{false}&&\cr
\endorbit
\orbit{5}{22}{104}
\h{22}&\root\int(2){\2}{}&\int(2){\2}{}&\int(2){\2}{}&\int(0){\0}{}&\int(0){\0}{}&&\cnt{472}&\bounds{104}{104}{1}{false}&&\cr
\l{1}&\int(1){\2}{\spec{1}}&\int(2){\2}{}&\int(0){\0}{}&\int(0){\0}{}&\int(0){\0}{}&\ln{4\sd}&\cnt{2}&\bounds{1}{1}{1}{false}&&\cr
\l{2}&\int(1){3}{\spec{1}}&\int(1){5}{}&\int(1){5}{}&\int(0){5}{}&\int(0){\0}{}&\ln{4\sd}&\cnt{6}&\bounds{1}{1}{1}{false}&&\cr
\l{3}&\int(1){4}{}&\int(1){4}{}&\int(1){2}{}&\int(0){\0}{}&\int(0){\0}{}&\ln{2\sd}&\cnt{32}&\bounds{8}{8}{1}{false}&&\cr
\l{4}&\int(1){4}{}&\int(1){1}{}&\int(1){5}{}&\int(0){\0}{}&\int(0){3}{}&\ln{2\sd}&\cnt{16}&\bounds{4}{4}{1}{false}&&\cr
\l{5}&\int(1){5}{}&\int(1){3}{}&\int(1){5}{}&\int(0){1}{}&\int(0){\0}{}&\ln{4\sd}&\cnt{36}&\bounds{6}{6}{1}{false}&&\cr
\l{6}&\int(1){5}{}&\int(1){2}{}&\int(1){1}{}&\int(0){\0}{}&\int(0){2}{}&\ln{4\sd}&\cnt{32}&\bounds{8}{8}{1}{false}&&\cr
\l{7}&\int(1){5}{}&\int(1){5}{}&\int(1){1}{}&\int(0){3}{}&\int(0){\0}{}&\ln{2\sd}&\cnt{20}&\bounds{4}{4}{1}{false}&&\cr
\l{8}&\int(2){\2}{}&\int(1){\2}{\spec{1}}&\int(0){\0}{}&\int(0){\0}{}&\int(0){\0}{}&\ln{4}&\cnt{4}&\bounds{1}{1}{1}{false}&&\cr
\endorbit
\orbit{5}{23}{116}
\h{23}&\rat(5/6){5}{}&\rat(5/6){5}{}&\root\rat(10/3){4}{\spec{1}}&\int(0){\0}{}&\int(1){1}{}&&\cnt{656}&\bounds{116}{128}{fail}{false}&&\cr
\l{1}&\rat(2/3){4}{}&\rat(1/6){1}{}&\rat(5/3){5}{}&\int(0){\0}{}&\rat(1/2){3}{}&\ln{2\sd}&\cnt{100}&\bounds{20}{20}{1}{false}&&\cr
\l{2}&\rat(5/6){5}{}&\int(0){\0}{}&\rat(5/3){5}{}&\int(0){4}{}&\rat(1/2){2}{}&\ln{2\sd}&\cnt{60}&\bounds{10}{12}{5}{true}&&\cr
\l{3}&\int(1){\2}{}&\int(0){\0}{}&\int(2){\2}{}&\int(0){\0}{}&\int(0){\0}{}&\ln{2}&\cnt{15}&\bounds{3}{3}{1}{true}&&\cr
\l{5}&\rat(1/2){3}{}&\rat(5/6){5}{}&\rat(5/3){5}{}&\int(0){5}{}&\int(0){\0}{}&\ln{2}&\cnt{60}&\bounds{9}{10}{5}{true}&&\cr
\l{7}&\rat(5/6){5}{}&\rat(5/6){5}{}&\rat(4/3){4}{}&\int(0){\0}{}&\int(0){1}{}&\ln{1}&\cnt{18}&\bounds{4}{6}{5}{true}&&\cr
\endorbit
\orbit{5}{24}{118}
\h{24}&\rat(5/6){5}{}&\rat(5/6){5}{}&\rat(10/3){4}{\spec{1}}&\root\int(0){\0}{}&\int(1){1}{}&&\cnt{576}&\bounds{118}{118}{1}{false}&&\cr
\l{1}&\rat(2/3){4}{}&\rat(1/6){1}{}&\rat(5/3){5}{}&\int(0){\0}{}&\rat(1/2){3}{}&\ln{2\sd}&\cnt{100}&\bounds{20}{20}{1}{false}&&\cr
\l{2}&\rat(5/6){5}{}&\int(0){\0}{}&\rat(5/3){5}{}&\int(0){4}{\spec{1}}&\rat(1/2){2}{}&\ln{2\sd}&\cnt{24}&\bounds{6}{6}{1}{false}&&\cr
\l{3}&\rat(5/6){5}{}&\int(0){\0}{}&\rat(5/3){5}{}&\int(0){4}{\spec{2}}&\rat(1/2){2}{}&\ln{2\sd}&\cnt{4}&\bounds{1}{1}{1}{false}&&\cr
\l{4}&\int(1){\2}{}&\int(0){\0}{}&\int(2){\2}{}&\int(0){\0}{}&\int(0){\0}{}&\ln{2}&\cnt{25}&\bounds{5}{5}{1}{false}&&\cr
\l{6}&\rat(1/2){3}{}&\rat(5/6){5}{}&\rat(5/3){5}{}&\int(0){5}{}&\int(0){\0}{}&\ln{2}&\cnt{40}&\bounds{8}{8}{1}{false}&&\cr
\l{8}&\rat(5/6){5}{}&\rat(5/6){5}{}&\rat(4/3){4}{}&\int(0){\0}{}&\int(0){1}{}&\ln{1}&\cnt{30}&\bounds{6}{6}{1}{false}&&\cr
\endorbit
\orbit{5}{25}{104}
\h{25}&\root\rat(5/6){5}{}&\rat(5/6){5}{}&\rat(10/3){4}{\spec{1}}&\int(0){\0}{}&\int(1){1}{}&&\cnt{548}&\bounds{104}{104}{1}{false}&&\cr
\l{1}&\int(1){\2}{}&\int(0){\0}{}&\int(2){\2}{}&\int(0){\0}{}&\int(0){\0}{}&\ln{1}&\cnt{15}&\bounds{3}{3}{1}{false}&&\cr
\l{3}&\rat(1/2){3}{\spec{1}}&\rat(5/6){5}{}&\rat(5/3){5}{}&\int(0){5}{}&\int(0){\0}{}&\ln{1}&\cnt{18}&\bounds{3}{3}{1}{false}&&\cr
\l{5}&\rat(1/2){3}{\spec{2}}&\rat(5/6){5}{}&\rat(5/3){5}{}&\int(0){5}{}&\int(0){\0}{}&\ln{1}&\cnt{6}&\bounds{1}{1}{1}{false}&&\cr
\l{7}&\rat(2/3){4}{}&\rat(1/6){1}{}&\rat(5/3){5}{}&\int(0){\0}{}&\rat(1/2){3}{}&\ln{1}&\cnt{60}&\bounds{12}{12}{1}{false}&&\cr
\l{9}&\rat(5/6){5}{}&\rat(1/2){3}{}&\rat(5/3){5}{}&\int(0){1}{}&\int(0){\0}{}&\ln{1}&\cnt{60}&\bounds{10}{10}{1}{false}&&\cr
\l{11}&\rat(5/6){5}{}&\rat(5/6){5}{}&\rat(4/3){4}{}&\int(0){\0}{}&\int(0){1}{}&\ln{1}&\cnt{30}&\bounds{6}{6}{1}{false}&&\cr
\l{13}&\rat(5/6){5}{}&\rat(-1/6){5}{}&\rat(4/3){4}{}&\int(0){\0}{}&\int(1){1}{}&\ln{1}&\cnt{25}&\bounds{5}{5}{1}{false}&&\cr
\l{15}&\rat(5/6){5}{}&\int(0){\0}{}&\rat(5/3){5}{}&\int(0){4}{}&\rat(1/2){2}{}&\ln{1}&\cnt{60}&\bounds{12}{12}{1}{false}&&\cr
\endorbit
\orbit{5}{26}{96}
\h{26}&\rat(5/6){5}{}&\rat(5/6){5}{}&\rat(10/3){4}{\spec{1}}&\int(0){\0}{}&\root\int(1){1}{}&&\cnt{520}&\bounds{96}{96}{1}{false}&&\cr
\l{1}&\rat(2/3){4}{}&\rat(1/6){1}{}&\rat(5/3){5}{}&\int(0){\0}{}&\rat(1/2){3}{}&\ln{2\sd}&\cnt{50}&\bounds{10}{10}{1}{false}&&\cr
\l{2}&\rat(5/6){5}{}&\int(0){\0}{}&\rat(5/3){5}{}&\int(0){4}{}&\rat(1/2){2}{}&\ln{2\sd}&\cnt{30}&\bounds{6}{6}{1}{false}&&\cr
\l{3}&\int(1){\2}{}&\int(0){\0}{}&\int(2){\2}{}&\int(0){\0}{}&\int(0){\0}{}&\ln{2}&\cnt{25}&\bounds{5}{5}{1}{false}&&\cr
\l{5}&\rat(1/2){3}{}&\rat(5/6){5}{}&\rat(5/3){5}{}&\int(0){5}{}&\int(0){\0}{}&\ln{2}&\cnt{60}&\bounds{10}{10}{1}{false}&&\cr
\l{7}&\rat(5/6){5}{}&\rat(5/6){5}{}&\rat(4/3){4}{}&\int(0){\0}{}&\int(0){1}{\spec{1}}&\ln{2}&\cnt{5}&\bounds{1}{1}{1}{false}&&\cr
\endorbit
\orbit{5}{27}{110}
\h{27}&\root\int(4){\4}{\spec{1}}&\int(2){\2}{}&\int(0){\0}{}&\int(0){\0}{}&\int(0){\0}{}&&\cnt{680}&\bounds{110}{132}{fail}{false}&&\cr
\l{1}&\int(2){3}{}&\int(1){3}{}&\int(0){\0}{}&\int(0){\0}{}&\int(0){1}{}&\ln{1\sdd}&\cnt{96}&\bounds{14}{24}{5}{true}&&\cr
\l{2}&\int(2){3}{}&\int(1){5}{}&\int(0){5}{}&\int(0){5}{}&\int(0){\0}{}&\ln{2\sd}&\cnt{72}&\bounds{12}{12}{1}{false}&&\cr
\l{3}&\int(3){\4}{}&\int(0){\0}{}&\int(0){\0}{}&\int(0){\0}{}&\int(0){\0}{}&\ln{1}&\cnt{4}&\bounds{2}{2}{1}{true}&&\cr
\l{5}&\int(2){4}{}&\int(1){4}{}&\int(0){2}{}&\int(0){\0}{}&\int(0){\0}{}&\ln{2}&\cnt{60}&\bounds{10}{12}{5}{true}&&\cr
\l{7}&\int(2){4}{}&\int(1){5}{}&\int(0){\0}{}&\int(0){1}{}&\int(0){2}{}&\ln{2}&\cnt{48}&\bounds{7}{8}{5}{true}&&\cr
\endorbit
\orbit{5}{28}{102}
\h{28}&\int(4){\4}{}&\int(2){\2}{}&\root\int(0){\0}{}&\int(0){\0}{}&\int(0){\0}{}&&\cnt{576}&\bounds{102}{126}{fail}{false}&&\cr
\l{1}&\int(2){\2}{}&\int(1){\2}{\spec{1}}&\int(0){\0}{}&\int(0){\0}{}&\int(0){\0}{}&\ln{2\sd}&\cnt{16}&\bounds{4}{4}{1}{true}&&\cr
\l{2}&\int(2){3}{}&\int(1){3}{}&\int(0){\0}{}&\int(0){\0}{}&\int(0){1}{}&\ln{1\sdd}&\cnt{96}&\bounds{14}{24}{5}{true}&&\cr
\l{3}&\int(2){3}{}&\int(1){5}{}&\int(0){5}{}&\int(0){5}{}&\int(0){\0}{}&\ln{2\sd}&\cnt{48}&\bounds{8}{8}{1}{true}&&\cr
\l{4}&\int(2){4}{}&\int(1){4}{}&\int(0){2}{\spec{1}}&\int(0){\0}{}&\int(0){\0}{}&\ln{2\sd}&\cnt{4}&\bounds{1}{1}{1}{true}&&\cr
\l{5}&\int(2){4}{}&\int(1){4}{}&\int(0){2}{\spec{2}}&\int(0){\0}{}&\int(0){\0}{}&\ln{2\sd}&\cnt{24}&\bounds{5}{6}{5}{true}&&\cr
\l{6}&\int(2){4}{}&\int(1){2}{}&\int(0){\0}{}&\int(0){4}{}&\int(0){\0}{}&\ln{2\sd}&\cnt{60}&\bounds{10}{12}{5}{true}&&\cr
\l{7}&\int(2){4}{}&\int(1){5}{}&\int(0){\0}{}&\int(0){1}{}&\int(0){2}{}&\ln{2\sd}&\cnt{48}&\bounds{7}{8}{5}{true}&&\cr
\l{8}&\int(2){4}{}&\int(1){1}{}&\int(0){5}{}&\int(0){\0}{}&\int(0){3}{}&\ln{2\sd}&\cnt{32}&\bounds{5}{8}{5}{true}&&\cr
\l{9}&\int(3){\4}{\spec{1}}&\int(0){\0}{}&\int(0){\0}{}&\int(0){\0}{}&\int(0){\0}{}&\ln{2}&\cnt{4}&\bounds{2}{2}{1}{true}&&\cr
\endorbit
\orbit{5}{29}{116}
\h{29}&\int(4){\4}{}&\int(2){\2}{}&\int(0){\0}{}&\int(0){\0}{}&\root\int(0){\0}{}&&\cnt{576}&\bounds{116}{120}{fail}{false}&&\cr
\l{1}&\int(2){\2}{}&\int(1){\2}{\spec{1}}&\int(0){\0}{}&\int(0){\0}{}&\int(0){\0}{}&\ln{2\sd}&\cnt{16}&\bounds{4}{4}{1}{false}&&\cr
\l{2}&\int(2){3}{}&\int(1){3}{}&\int(0){\0}{}&\int(0){\0}{}&\int(0){1}{}&\ln{1\sdd}&\cnt{48}&\bounds{12}{16}{2}{false}&&\cr
\l{3}&\int(2){3}{}&\int(1){5}{}&\int(0){5}{}&\int(0){5}{}&\int(0){\0}{}&\ln{2\sd}&\cnt{72}&\bounds{12}{12}{1}{false}&&\cr
\l{4}&\int(2){4}{}&\int(1){4}{}&\int(0){2}{}&\int(0){\0}{}&\int(0){\0}{}&\ln{4\sd}&\cnt{60}&\bounds{12}{12}{1}{false}&&\cr
\l{5}&\int(2){4}{}&\int(1){5}{}&\int(0){\0}{}&\int(0){1}{}&\int(0){2}{}&\ln{4\sd}&\cnt{24}&\bounds{4}{4}{1}{false}&&\cr
\l{6}&\int(3){\4}{\spec{1}}&\int(0){\0}{}&\int(0){\0}{}&\int(0){\0}{}&\int(0){\0}{}&\ln{2}&\cnt{4}&\bounds{2}{2}{1}{false}&&\cr
\endorbit
\orbit{5}{30}{104}
\h{30}&\int(4){\4}{}&\root\int(2){\2}{}&\int(0){\0}{}&\int(0){\0}{}&\int(0){\0}{}&&\cnt{520}&\bounds{104}{104}{1}{false}&&\cr
\l{1}&\int(2){\2}{}&\int(1){\2}{\spec{1}}&\int(0){\0}{}&\int(0){\0}{}&\int(0){\0}{}&\ln{2\sd}&\cnt{8}&\bounds{4}{4}{1}{false}&&\cr
\l{2}&\int(2){3}{}&\int(1){3}{\spec{1}}&\int(0){\0}{}&\int(0){\0}{}&\int(0){1}{}&\ln{2\sd}&\cnt{16}&\bounds{4}{4}{1}{false}&&\cr
\l{3}&\int(2){3}{}&\int(1){5}{}&\int(0){5}{}&\int(0){5}{}&\int(0){\0}{}&\ln{2\sd}&\cnt{72}&\bounds{12}{12}{1}{false}&&\cr
\l{4}&\int(2){4}{}&\int(1){4}{}&\int(0){2}{}&\int(0){\0}{}&\int(0){\0}{}&\ln{4\sd}&\cnt{30}&\bounds{6}{6}{1}{false}&&\cr
\l{5}&\int(2){4}{}&\int(1){5}{}&\int(0){\0}{}&\int(0){1}{}&\int(0){2}{}&\ln{4\sd}&\cnt{48}&\bounds{8}{8}{1}{false}&&\cr
\l{6}&\int(3){\4}{\spec{1}}&\int(0){\0}{}&\int(0){\0}{}&\int(0){\0}{}&\int(0){\0}{}&\ln{2}&\cnt{4}&\bounds{2}{2}{1}{false}&&\cr
\endorbit
\orbit{5}{31}{88}
\h{31}&\root\int(4){\4}{\spec{2}}&\int(2){\2}{}&\int(0){\0}{}&\int(0){\0}{}&\int(0){\0}{}&&\cnt{464}&\bounds{88}{88}{1}{false}&&\cr
\l{1}&\int(2){\2}{}&\int(1){\2}{\spec{1}}&\int(0){\0}{}&\int(0){\0}{}&\int(0){\0}{}&\ln{2\sd}&\cnt{16}&\bounds{4}{4}{1}{false}&&\cr
\l{2}&\int(2){4}{}&\int(1){4}{}&\int(0){2}{}&\int(0){\0}{}&\int(0){\0}{}&\ln{4\sd}&\cnt{60}&\bounds{12}{12}{1}{false}&&\cr
\l{3}&\int(2){4}{}&\int(1){5}{}&\int(0){\0}{}&\int(0){1}{}&\int(0){2}{}&\ln{4\sd}&\cnt{48}&\bounds{8}{8}{1}{false}&&\cr
\endorbit
\orbit{5}{32}{100}
\h{32}&\rat(5/6){5}{}&\rat(5/6){5}{}&\rat(5/6){1}{}&\root\rat(7/2){3}{\spec{1}}&\int(0){\0}{}&&\cnt{604}&\bounds{100}{124}{fail}{false}&&\cr
\l{1}&\int(1){\2}{}&\int(0){\0}{}&\int(0){\0}{}&\int(2){\2}{}&\int(0){\0}{}&\ln{3}&\cnt{10}&\bounds{2}{2}{1}{true}&&\cr
\l{3}&\rat(1/2){3}{}&\rat(5/6){5}{}&\rat(1/6){5}{}&\rat(3/2){5}{}&\int(0){\0}{}&\ln{3}&\cnt{50}&\bounds{8}{10}{5}{true}&&\cr
\l{5}&\rat(2/3){4}{}&\int(0){\0}{}&\rat(5/6){1}{}&\rat(3/2){5}{}&\int(0){1}{}&\ln{3}&\cnt{40}&\bounds{6}{8}{5}{true}&&\cr
\l{7}&\rat(5/6){5}{}&\rat(5/6){5}{}&\rat(5/6){1}{}&\rat(1/2){3}{\spec{1}}&\int(0){\0}{}&\ln{1}&\cnt{1}&\bounds{1}{1}{1}{true}&&\cr
\l{9}&\rat(5/6){5}{}&\rat(5/6){5}{}&\rat(5/6){1}{}&\rat(1/2){3}{\spec{2}}&\int(0){\0}{}&\ln{1}&\cnt{1}&\bounds{1}{1}{1}{true}&&\cr
\endorbit
\orbit{5}{33}{112}
\h{33}&\rat(5/6){5}{}&\rat(5/6){5}{}&\rat(5/6){1}{}&\rat(7/2){3}{\spec{1}}&\root\int(0){\0}{}&&\cnt{552}&\bounds{112}{112}{1}{false}&&\cr
\l{1}&\int(1){\2}{}&\int(0){\0}{}&\int(0){\0}{}&\int(2){\2}{}&\int(0){\0}{}&\ln{3}&\cnt{20}&\bounds{4}{4}{1}{false}&&\cr
\l{3}&\rat(1/2){3}{}&\rat(5/6){5}{}&\rat(1/6){5}{}&\rat(3/2){5}{}&\int(0){\0}{}&\ln{3}&\cnt{50}&\bounds{10}{10}{1}{false}&&\cr
\l{5}&\rat(2/3){4}{}&\int(0){\0}{}&\rat(5/6){1}{}&\rat(3/2){5}{}&\int(0){1}{}&\ln{3}&\cnt{20}&\bounds{4}{4}{1}{false}&&\cr
\l{7}&\rat(5/6){5}{}&\rat(5/6){5}{}&\rat(5/6){1}{}&\rat(1/2){3}{}&\int(0){\0}{}&\ln{1}&\cnt{6}&\bounds{2}{2}{1}{false}&&\cr
\endorbit
\orbit{5}{34}{110}
\h{34}&\root\rat(5/6){5}{}&\rat(5/6){5}{}&\rat(5/6){1}{}&\rat(7/2){3}{\spec{1}}&\int(0){\0}{}&&\cnt{524}&\bounds{110}{110}{1}{false}&&\cr
\l{1}&\int(1){\2}{}&\int(0){\0}{}&\int(0){\0}{}&\int(2){\2}{}&\int(0){\0}{}&\ln{1}&\cnt{12}&\bounds{3}{3}{1}{false}&&\cr
\l{3}&\rat(1/2){3}{\spec{1}}&\rat(5/6){5}{}&\rat(1/6){5}{}&\rat(3/2){5}{}&\int(0){\0}{}&\ln{1}&\cnt{15}&\bounds{3}{3}{1}{false}&&\cr
\l{5}&\rat(1/2){3}{\spec{2}}&\rat(5/6){5}{}&\rat(1/6){5}{}&\rat(3/2){5}{}&\int(0){\0}{}&\ln{1}&\cnt{5}&\bounds{1}{1}{1}{false}&&\cr
\l{7}&\rat(2/3){4}{}&\rat(1/3){2}{}&\int(0){\0}{}&\int(2){4}{}&\int(0){\0}{}&\ln{1}&\cnt{30}&\bounds{6}{6}{1}{false}&&\cr
\l{9}&\rat(2/3){4}{}&\int(0){\0}{}&\rat(5/6){1}{}&\rat(3/2){5}{}&\int(0){1}{}&\ln{1}&\cnt{24}&\bounds{6}{6}{1}{false}&&\cr
\l{11}&\rat(5/6){5}{}&\rat(2/3){4}{}&\int(0){\0}{}&\rat(3/2){5}{}&\int(0){3}{}&\ln{1}&\cnt{40}&\bounds{8}{8}{1}{false}&&\cr
\l{13}&\rat(5/6){5}{}&\rat(5/6){5}{}&\rat(5/6){1}{}&\rat(1/2){3}{}&\int(0){\0}{}&\ln{1}&\cnt{6}&\bounds{2}{2}{1}{false}&&\cr
\l{15}&\rat(5/6){5}{}&\rat(5/6){5}{}&\rat(-1/6){1}{}&\rat(3/2){3}{}&\int(0){\0}{}&\ln{1}&\cnt{20}&\bounds{4}{4}{1}{false}&&\cr
\l{17}&\rat(5/6){5}{}&\rat(1/6){1}{}&\rat(1/2){3}{}&\rat(3/2){5}{}&\int(0){\0}{}&\ln{1}&\cnt{50}&\bounds{10}{10}{1}{false}&&\cr
\l{19}&\rat(5/6){5}{}&\rat(-1/6){5}{}&\rat(5/6){1}{}&\rat(3/2){3}{}&\int(0){\0}{}&\ln{1}&\cnt{20}&\bounds{4}{4}{1}{false}&&\cr
\l{21}&\rat(5/6){5}{}&\int(0){\0}{}&\rat(1/6){5}{}&\int(2){4}{}&\int(0){2}{}&\ln{1}&\cnt{40}&\bounds{8}{8}{1}{false}&&\cr
\endorbit
\orbit{5}{35}{112}
\h{35}&\root\int(4){\4}{\spec{1}}&\int(0){\0}{}&\int(0){\0}{}&\int(0){\0}{}&\int(2){\2}{}&&\cnt{680}&\bounds{112}{124}{fail}{false}&&\cr
\l{1}&\int(2){3}{}&\int(0){3}{}&\int(0){\0}{}&\int(0){\0}{}&\int(1){1}{}&\ln{3\sdd}&\cnt{80}&\bounds{12}{16}{5}{true}&&\cr
\l{2}&\int(3){\4}{}&\int(0){\0}{}&\int(0){\0}{}&\int(0){\0}{}&\int(0){\0}{}&\ln{1}&\cnt{4}&\bounds{2}{2}{1}{true}&&\cr
\l{4}&\int(2){4}{}&\int(0){5}{}&\int(0){\0}{}&\int(0){1}{}&\int(1){2}{}&\ln{3}&\cnt{72}&\bounds{12}{12}{1}{false}&&\cr
\endorbit
\orbit{5}{36}{114}
\h{36}&\int(4){\4}{}&\root\int(0){\0}{}&\int(0){\0}{}&\int(0){\0}{}&\int(2){\2}{}&&\cnt{576}&\bounds{114}{120}{fail}{false}&&\cr
\l{1}&\int(2){\2}{}&\int(0){\0}{}&\int(0){\0}{}&\int(0){\0}{}&\int(1){\2}{}&\ln{1\sdd}&\cnt{32}&\bounds{10}{16}{2}{false}&&\cr
\l{2}&\int(2){3}{}&\int(0){3}{\spec{1}}&\int(0){\0}{}&\int(0){\0}{}&\int(1){1}{}&\ln{2\sd}&\cnt{16}&\bounds{4}{4}{1}{false}&&\cr
\l{3}&\int(2){3}{}&\int(0){\0}{}&\int(0){3}{}&\int(0){\0}{}&\int(1){2}{}&\ln{2\sdd}&\cnt{80}&\bounds{16}{16}{1}{false}&&\cr
\l{4}&\int(2){4}{}&\int(0){5}{}&\int(0){\0}{}&\int(0){1}{}&\int(1){2}{}&\ln{4\sd}&\cnt{48}&\bounds{8}{8}{1}{false}&&\cr
\l{5}&\int(2){4}{}&\int(0){\0}{}&\int(0){1}{}&\int(0){5}{}&\int(1){1}{}&\ln{2\sd}&\cnt{72}&\bounds{12}{12}{1}{false}&&\cr
\l{6}&\int(3){\4}{\spec{1}}&\int(0){\0}{}&\int(0){\0}{}&\int(0){\0}{}&\int(0){\0}{}&\ln{2}&\cnt{4}&\bounds{2}{2}{1}{false}&&\cr
\endorbit
\orbit{5}{37}{88}
\h{37}&\int(4){\4}{}&\int(0){\0}{}&\int(0){\0}{}&\int(0){\0}{}&\root\int(2){\2}{\spec{1}}&&\cnt{464}&\bounds{88}{88}{1}{false}&&\cr
\l{1}&\int(2){3}{}&\int(0){3}{}&\int(0){\0}{}&\int(0){\0}{}&\int(1){1}{}&\ln{2\sdd}&\cnt{80}&\bounds{16}{16}{1}{false}&&\cr
\l{2}&\int(2){4}{}&\int(0){5}{}&\int(0){\0}{}&\int(0){1}{}&\int(1){2}{}&\ln{4\sd}&\cnt{72}&\bounds{12}{12}{1}{false}&&\cr
\l{3}&\int(3){\4}{\spec{1}}&\int(0){\0}{}&\int(0){\0}{}&\int(0){\0}{}&\int(0){\0}{}&\ln{2}&\cnt{4}&\bounds{2}{2}{1}{false}&&\cr
\endorbit
\orbit{5}{38}{82}
\h{38}&\root\int(4){\4}{\spec{2}}&\int(0){\0}{}&\int(0){\0}{}&\int(0){\0}{}&\int(2){\2}{}&&\cnt{464}&\bounds{82}{88}{fail}{false}&&\cr
\l{1}&\int(2){\2}{}&\int(0){\0}{}&\int(0){\0}{}&\int(0){\0}{}&\int(1){\2}{}&\ln{1\sdd}&\cnt{32}&\bounds{10}{16}{2}{false}&&\cr
\l{2}&\int(2){4}{}&\int(0){5}{}&\int(0){\0}{}&\int(0){1}{}&\int(1){2}{}&\ln{6\sd}&\cnt{72}&\bounds{12}{12}{1}{false}&&\cr
\endorbit
\orbit{5}{39}{118}
\h{39}&\root\rat(10/3){4}{\spec{1}}&\rat(4/3){4}{}&\rat(4/3){2}{}&\int(0){\0}{}&\int(0){\0}{}&&\cnt{656}&\bounds{118}{136}{fail}{false}&&\cr
\l{1}&\rat(5/3){5}{}&\int(1){3}{}&\rat(1/3){5}{}&\int(0){1}{}&\int(0){\0}{}&\ln{2\sd}&\cnt{96}&\bounds{16}{16}{1}{false}&&\cr
\l{2}&\rat(5/3){5}{}&\rat(4/3){4}{}&\int(0){\0}{}&\int(0){5}{}&\int(0){3}{}&\ln{2\sd}&\cnt{48}&\bounds{7}{8}{5}{true}&&\cr
\l{3}&\rat(5/3){5}{}&\rat(2/3){2}{}&\rat(2/3){1}{}&\int(0){\0}{}&\int(0){2}{}&\ln{2\sdd}&\cnt{96}&\bounds{18}{24}{2}{false}&&\cr
\l{4}&\rat(5/3){5}{}&\rat(2/3){5}{}&\rat(2/3){1}{}&\int(0){3}{}&\int(0){\0}{}&\ln{1\sdd}&\cnt{80}&\bounds{12}{16}{5}{true}&&\cr
\l{5}&\int(2){\2}{}&\int(1){\2}{}&\int(0){\0}{}&\int(0){\0}{}&\int(0){\0}{}&\ln{2}&\cnt{24}&\bounds{6}{6}{1}{true}&&\cr
\endorbit
\orbit{5}{40}{118}
\h{40}&\rat(10/3){4}{\spec{1}}&\rat(4/3){4}{}&\rat(4/3){2}{}&\root\int(0){\0}{}&\int(0){\0}{}&&\cnt{576}&\bounds{118}{136}{fail}{false}&&\cr
\l{1}&\rat(5/3){5}{}&\int(1){3}{}&\rat(1/3){5}{}&\int(0){1}{}&\int(0){\0}{}&\ln{2\sd}&\cnt{64}&\bounds{16}{16}{1}{false}&&\cr
\l{2}&\rat(5/3){5}{}&\rat(4/3){4}{}&\int(0){\0}{}&\int(0){5}{}&\int(0){3}{}&\ln{2\sd}&\cnt{32}&\bounds{5}{8}{5}{true}&&\cr
\l{3}&\rat(5/3){5}{}&\rat(2/3){2}{}&\rat(2/3){1}{}&\int(0){\0}{}&\int(0){2}{}&\ln{2\sdd}&\cnt{96}&\bounds{18}{24}{2}{false}&&\cr
\l{4}&\rat(5/3){5}{}&\rat(2/3){5}{}&\rat(2/3){1}{}&\int(0){3}{\spec{1}}&\int(0){\0}{}&\ln{2\sd}&\cnt{16}&\bounds{4}{4}{1}{true}&&\cr
\l{5}&\int(2){\2}{}&\int(1){\2}{}&\int(0){\0}{}&\int(0){\0}{}&\int(0){\0}{}&\ln{2}&\cnt{40}&\bounds{8}{8}{1}{true}&&\cr
\endorbit
\orbit{5}{41}{112}
\h{41}&\rat(10/3){4}{\spec{1}}&\rat(4/3){4}{}&\rat(4/3){2}{}&\int(0){\0}{}&\root\int(0){\0}{}&&\cnt{576}&\bounds{112}{120}{fail}{false}&&\cr
\l{1}&\rat(5/3){5}{}&\int(1){3}{}&\rat(1/3){5}{}&\int(0){1}{}&\int(0){\0}{}&\ln{2\sd}&\cnt{96}&\bounds{16}{16}{1}{false}&&\cr
\l{2}&\rat(5/3){5}{}&\rat(4/3){4}{}&\int(0){\0}{}&\int(0){5}{}&\int(0){3}{}&\ln{2\sd}&\cnt{24}&\bounds{4}{4}{1}{false}&&\cr
\l{3}&\rat(5/3){5}{}&\rat(2/3){2}{}&\rat(2/3){1}{}&\int(0){\0}{}&\int(0){2}{}&\ln{2\sdd}&\cnt{48}&\bounds{12}{16}{2}{false}&&\cr
\l{4}&\rat(5/3){5}{}&\rat(2/3){5}{}&\rat(2/3){1}{}&\int(0){3}{}&\int(0){\0}{}&\ln{1\sdd}&\cnt{80}&\bounds{16}{16}{1}{false}&&\cr
\l{5}&\int(2){\2}{}&\int(1){\2}{}&\int(0){\0}{}&\int(0){\0}{}&\int(0){\0}{}&\ln{2}&\cnt{40}&\bounds{8}{8}{1}{false}&&\cr
\endorbit
\orbit{5}{42}{98}
\h{42}&\rat(10/3){4}{\spec{1}}&\root\rat(4/3){4}{\spec{1}}&\rat(4/3){2}{}&\int(0){\0}{}&\int(0){\0}{}&&\cnt{520}&\bounds{98}{104}{fail}{false}&&\cr
\l{1}&\rat(5/3){5}{}&\rat(2/3){5}{}&\rat(2/3){4}{}&\int(0){\0}{}&\int(0){1}{}&\ln{1\sdd}&\cnt{96}&\bounds{18}{24}{2}{false}&&\cr
\l{2}&\rat(5/3){5}{}&\rat(2/3){5}{}&\rat(2/3){1}{}&\int(0){3}{}&\int(0){\0}{}&\ln{1\sdd}&\cnt{80}&\bounds{16}{16}{1}{false}&&\cr
\l{3}&\int(2){\2}{}&\int(1){\2}{}&\int(0){\0}{}&\int(0){\0}{}&\int(0){\0}{}&\ln{1}&\cnt{20}&\bounds{4}{4}{1}{false}&&\cr
\l{5}&\int(2){\2}{}&\int(0){\0}{}&\int(1){\2}{}&\int(0){\0}{}&\int(0){\0}{}&\ln{1}&\cnt{40}&\bounds{8}{8}{1}{false}&&\cr
\l{7}&\rat(5/3){5}{}&\int(1){3}{}&\rat(1/3){5}{}&\int(0){1}{}&\int(0){\0}{}&\ln{1}&\cnt{48}&\bounds{8}{8}{1}{false}&&\cr
\l{9}&\rat(5/3){5}{}&\rat(4/3){4}{}&\int(0){\0}{}&\int(0){5}{}&\int(0){3}{}&\ln{1}&\cnt{48}&\bounds{8}{8}{1}{false}&&\cr
\l{11}&\rat(5/3){5}{}&\rat(2/3){2}{\spec{1}}&\rat(2/3){1}{}&\int(0){\0}{}&\int(0){2}{}&\ln{1}&\cnt{16}&\bounds{4}{4}{1}{false}&&\cr
\endorbit
\orbit{5}{43}{82}
\h{43}&\rat(10/3){4}{\spec{1}}&\root\rat(4/3){4}{\spec{2}}&\rat(4/3){2}{}&\int(0){\0}{}&\int(0){\0}{}&&\cnt{464}&\bounds{82}{88}{fail}{false}&&\cr
\l{1}&\rat(5/3){5}{}&\rat(2/3){2}{}&\rat(2/3){1}{}&\int(0){\0}{}&\int(0){2}{}&\ln{1\sdd}&\cnt{96}&\bounds{18}{24}{2}{false}&&\cr
\l{2}&\int(2){\2}{}&\int(0){\0}{}&\int(1){\2}{}&\int(0){\0}{}&\int(0){\0}{}&\ln{1}&\cnt{40}&\bounds{8}{8}{1}{false}&&\cr
\l{4}&\rat(5/3){5}{}&\int(1){3}{}&\rat(1/3){5}{}&\int(0){1}{}&\int(0){\0}{}&\ln{1}&\cnt{96}&\bounds{16}{16}{1}{false}&&\cr
\l{6}&\rat(5/3){5}{}&\rat(4/3){4}{}&\int(0){\0}{}&\int(0){5}{}&\int(0){3}{}&\ln{1}&\cnt{48}&\bounds{8}{8}{1}{false}&&\cr
\endorbit
\orbit{5}{44}{122}
\h{44}&\rat(3/2){3}{}&\rat(3/2){3}{}&\root\int(0){\0}{}&\int(0){\0}{}&\int(3){1}{}&&\cnt{528}&\bounds{122}{124}{fail}{false}&&\cr
\l{1}&\rat(3/2){3}{}&\int(0){\0}{}&\int(0){3}{\spec{1}}&\int(0){\0}{}&\rat(3/2){2}{}&\ln{4\sd}&\cnt{4}&\bounds{1}{1}{1}{true}&&\cr
\l{2}&\rat(3/2){3}{}&\int(0){\0}{}&\int(0){\0}{}&\int(0){3}{}&\rat(3/2){3}{}&\ln{2\sd}&\cnt{20}&\bounds{3}{4}{5}{true}&&\cr
\l{3}&\int(1){4}{}&\rat(1/2){5}{}&\int(0){\0}{}&\int(0){1}{}&\rat(3/2){2}{}&\ln{4\sd}&\cnt{54}&\bounds{9}{9}{1}{true}&&\cr
\l{4}&\int(1){4}{}&\rat(1/2){1}{}&\int(0){5}{}&\int(0){\0}{}&\rat(3/2){3}{}&\ln{4\sd}&\cnt{36}&\bounds{9}{9}{1}{true}&&\cr
\l{5}&\int(1){\2}{}&\int(0){\0}{}&\int(0){\0}{}&\int(0){\0}{}&\int(2){\2}{}&\ln{2}&\cnt{27}&\bounds{9}{9}{1}{true}&&\cr
\l{7}&\rat(3/2){3}{}&\rat(3/2){3}{}&\int(0){\0}{}&\int(0){\0}{}&\int(0){1}{\spec{1}}&\ln{2}&\cnt{1}&\bounds{1}{1}{1}{true}&&\cr
\endorbit
\orbit{5}{45}{104}
\h{45}&\rat(3/2){3}{}&\rat(3/2){3}{}&\int(0){\0}{}&\int(0){\0}{}&\root\int(3){1}{}&&\cnt{552}&\bounds{104}{104}{1}{false}&&\cr
\l{1}&\rat(3/2){3}{}&\int(0){\0}{}&\int(0){3}{}&\int(0){\0}{}&\rat(3/2){2}{}&\ln{4\sd}&\cnt{20}&\bounds{4}{4}{1}{false}&&\cr
\l{2}&\int(1){4}{}&\rat(1/2){5}{}&\int(0){\0}{}&\int(0){1}{}&\rat(3/2){2}{}&\ln{8\sd}&\cnt{54}&\bounds{9}{9}{1}{false}&&\cr
\l{3}&\int(1){\2}{}&\int(0){\0}{}&\int(0){\0}{}&\int(0){\0}{}&\int(2){\2}{}&\ln{2}&\cnt{9}&\bounds{3}{3}{1}{false}&&\cr
\l{5}&\rat(3/2){3}{}&\rat(3/2){3}{}&\int(0){\0}{}&\int(0){\0}{}&\int(0){1}{\spec{1}}&\ln{2}&\cnt{1}&\bounds{1}{1}{1}{false}&&\cr
\endorbit
\orbit{5}{46}{92}
\h{46}&\root\rat(3/2){3}{\spec{1}}&\rat(3/2){3}{}&\int(0){\0}{}&\int(0){\0}{}&\int(3){1}{}&&\cnt{444}&\bounds{92}{92}{1}{false}&&\cr
\l{1}&\int(1){\2}{}&\int(0){\0}{}&\int(0){\0}{}&\int(0){\0}{}&\int(2){\2}{}&\ln{1}&\cnt{9}&\bounds{3}{3}{1}{false}&&\cr
\l{3}&\rat(3/2){3}{}&\rat(3/2){3}{}&\int(0){\0}{}&\int(0){\0}{}&\int(0){1}{\spec{1}}&\ln{2}&\cnt{1}&\bounds{1}{1}{1}{false}&&\cr
\l{5}&\rat(3/2){3}{}&\rat(1/2){3}{}&\int(0){\0}{}&\int(0){\0}{}&\int(1){1}{}&\ln{1}&\cnt{27}&\bounds{9}{9}{1}{false}&&\cr
\l{7}&\rat(3/2){3}{}&\int(0){\0}{}&\int(0){3}{}&\int(0){\0}{}&\rat(3/2){2}{}&\ln{2}&\cnt{20}&\bounds{4}{4}{1}{false}&&\cr
\l{9}&\int(1){4}{}&\rat(1/2){5}{}&\int(0){\0}{}&\int(0){1}{}&\rat(3/2){2}{}&\ln{2}&\cnt{54}&\bounds{9}{9}{1}{false}&&\cr
\l{11}&\int(1){2}{}&\rat(1/2){5}{}&\int(0){1}{}&\int(0){\0}{}&\rat(3/2){3}{}&\ln{2}&\cnt{18}&\bounds{3}{3}{1}{false}&&\cr
\endorbit
\orbit{5}{47}{122}
\h{47}&\rat(3/2){3}{}&\rat(3/2){3}{}&\int(2){\2}{}&\root\int(0){\0}{}&\int(1){1}{}&&\cnt{480}&\bounds{122}{124}{fail}{false}&&\cr
\l{1}&\rat(3/2){3}{}&\int(0){\0}{}&\int(1){3}{}&\int(0){\0}{}&\rat(1/2){2}{}&\ln{2\sd}&\cnt{24}&\bounds{5}{6}{5}{true}&&\cr
\l{2}&\int(1){4}{}&\rat(1/2){1}{}&\int(1){5}{}&\int(0){\0}{}&\rat(1/2){3}{}&\ln{4\sd}&\cnt{36}&\bounds{9}{9}{1}{true}&&\cr
\l{3}&\int(1){\2}{}&\int(0){\0}{}&\int(2){\2}{}&\int(0){\0}{}&\int(0){\0}{}&\ln{2}&\cnt{9}&\bounds{3}{3}{1}{true}&&\cr
\l{5}&\rat(3/2){3}{}&\rat(3/2){3}{}&\int(0){\0}{}&\int(0){\0}{}&\int(0){1}{}&\ln{1}&\cnt{6}&\bounds{2}{2}{1}{true}&&\cr
\l{7}&\rat(3/2){3}{}&\rat(1/2){5}{}&\int(1){5}{}&\int(0){5}{}&\int(0){\0}{}&\ln{4}&\cnt{12}&\bounds{3}{3}{1}{true}&&\cr
\l{9}&\int(1){4}{}&\int(1){4}{}&\int(1){2}{}&\int(0){\0}{}&\int(0){\0}{}&\ln{2}&\cnt{36}&\bounds{9}{9}{1}{true}&&\cr
\endorbit
\orbit{5}{48}{108}
\h{48}&\rat(3/2){3}{}&\rat(3/2){3}{}&\int(2){\2}{}&\int(0){\0}{}&\root\int(1){1}{}&&\cnt{424}&\bounds{108}{108}{1}{false}&&\cr
\l{1}&\rat(3/2){3}{}&\int(0){\0}{}&\int(1){3}{}&\int(0){\0}{}&\rat(1/2){2}{}&\ln{2\sd}&\cnt{12}&\bounds{4}{4}{1}{false}&&\cr
\l{2}&\int(1){4}{}&\rat(1/2){1}{}&\int(1){5}{}&\int(0){\0}{}&\rat(1/2){3}{}&\ln{4\sd}&\cnt{18}&\bounds{6}{6}{1}{false}&&\cr
\l{3}&\int(1){\2}{}&\int(0){\0}{}&\int(2){\2}{}&\int(0){\0}{}&\int(0){\0}{}&\ln{2}&\cnt{9}&\bounds{3}{3}{1}{false}&&\cr
\l{5}&\rat(3/2){3}{}&\rat(3/2){3}{}&\int(0){\0}{}&\int(0){\0}{}&\int(0){1}{\spec{1}}&\ln{2}&\cnt{1}&\bounds{1}{1}{1}{false}&&\cr
\l{7}&\rat(3/2){3}{}&\rat(1/2){5}{}&\int(1){5}{}&\int(0){5}{}&\int(0){\0}{}&\ln{4}&\cnt{18}&\bounds{3}{3}{1}{false}&&\cr
\l{9}&\int(1){4}{}&\int(1){4}{}&\int(1){2}{}&\int(0){\0}{}&\int(0){\0}{}&\ln{2}&\cnt{36}&\bounds{9}{9}{1}{false}&&\cr
\endorbit
\orbit{5}{49}{104}
\h{49}&\rat(3/2){3}{}&\rat(3/2){3}{}&\root\int(2){\2}{}&\int(0){\0}{}&\int(1){1}{}&&\cnt{424}&\bounds{104}{104}{1}{false}&&\cr
\l{1}&\rat(3/2){3}{}&\int(0){\0}{}&\int(1){3}{\spec{1}}&\int(0){\0}{}&\rat(1/2){2}{}&\ln{4\sd}&\cnt{4}&\bounds{1}{1}{1}{false}&&\cr
\l{2}&\int(1){4}{}&\rat(1/2){1}{}&\int(1){5}{}&\int(0){\0}{}&\rat(1/2){3}{}&\ln{4\sd}&\cnt{36}&\bounds{9}{9}{1}{false}&&\cr
\l{3}&\int(1){\2}{}&\int(0){\0}{}&\int(2){\2}{}&\int(0){\0}{}&\int(0){\0}{}&\ln{2}&\cnt{9}&\bounds{3}{3}{1}{false}&&\cr
\l{5}&\rat(3/2){3}{}&\rat(3/2){3}{}&\int(0){\0}{}&\int(0){\0}{}&\int(0){1}{}&\ln{1}&\cnt{6}&\bounds{2}{2}{1}{false}&&\cr
\l{7}&\rat(3/2){3}{}&\rat(1/2){5}{}&\int(1){5}{}&\int(0){5}{}&\int(0){\0}{}&\ln{4}&\cnt{18}&\bounds{3}{3}{1}{false}&&\cr
\l{9}&\int(1){4}{}&\int(1){4}{}&\int(1){2}{}&\int(0){\0}{}&\int(0){\0}{}&\ln{2}&\cnt{18}&\bounds{6}{6}{1}{false}&&\cr
\endorbit
\orbit{5}{50}{92}
\h{50}&\root\rat(3/2){3}{\spec{1}}&\rat(3/2){3}{}&\int(2){\2}{}&\int(0){\0}{}&\int(1){1}{}&&\cnt{396}&\bounds{92}{92}{1}{false}&&\cr
\l{1}&\int(1){\2}{}&\int(0){\0}{}&\int(2){\2}{}&\int(0){\0}{}&\int(0){\0}{}&\ln{1}&\cnt{3}&\bounds{1}{1}{1}{false}&&\cr
\l{3}&\rat(3/2){3}{}&\rat(3/2){3}{}&\int(0){\0}{}&\int(0){\0}{}&\int(0){1}{}&\ln{1}&\cnt{6}&\bounds{2}{2}{1}{false}&&\cr
\l{5}&\rat(3/2){3}{}&\rat(1/2){3}{}&\int(0){\0}{}&\int(0){\0}{}&\int(1){1}{}&\ln{1}&\cnt{9}&\bounds{3}{3}{1}{false}&&\cr
\l{7}&\rat(3/2){3}{}&\rat(1/2){5}{}&\int(1){5}{}&\int(0){5}{}&\int(0){\0}{}&\ln{1}&\cnt{18}&\bounds{3}{3}{1}{false}&&\cr
\l{9}&\rat(3/2){3}{}&\rat(1/2){1}{}&\int(1){1}{}&\int(0){1}{}&\int(0){\0}{}&\ln{1}&\cnt{18}&\bounds{3}{3}{1}{false}&&\cr
\l{11}&\rat(3/2){3}{}&\int(0){\0}{}&\int(1){3}{}&\int(0){\0}{}&\rat(1/2){2}{}&\ln{1}&\cnt{24}&\bounds{6}{6}{1}{false}&&\cr
\l{13}&\int(1){4}{}&\int(1){4}{}&\int(1){2}{}&\int(0){\0}{}&\int(0){\0}{}&\ln{1}&\cnt{36}&\bounds{9}{9}{1}{false}&&\cr
\l{15}&\int(1){4}{}&\rat(1/2){1}{}&\int(1){5}{}&\int(0){\0}{}&\rat(1/2){3}{}&\ln{1}&\cnt{36}&\bounds{9}{9}{1}{false}&&\cr
\l{17}&\int(1){4}{}&\int(0){\0}{}&\int(1){1}{}&\int(0){5}{}&\int(1){1}{}&\ln{1}&\cnt{18}&\bounds{3}{3}{1}{false}&&\cr
\l{19}&\int(1){2}{}&\int(1){2}{}&\int(1){4}{}&\int(0){\0}{}&\int(0){\0}{}&\ln{1}&\cnt{12}&\bounds{3}{3}{1}{false}&&\cr
\l{21}&\int(1){2}{}&\rat(1/2){5}{}&\int(1){1}{}&\int(0){\0}{}&\rat(1/2){3}{}&\ln{1}&\cnt{12}&\bounds{3}{3}{1}{false}&&\cr
\l{23}&\int(1){2}{}&\int(0){\0}{}&\int(1){5}{}&\int(0){1}{}&\int(1){1}{}&\ln{1}&\cnt{6}&\bounds{1}{1}{1}{false}&&\cr
\endorbit
\orbit{5}{51}{112}
\h{51}&\rat(7/2){3}{\spec{1}}&\rat(3/2){3}{}&\root\int(0){\0}{}&\int(0){\0}{}&\int(1){1}{}&&\cnt{552}&\bounds{112}{122}{fail}{false}&&\cr
\l{1}&\int(3){\4}{}&\int(0){\0}{}&\int(0){\0}{}&\int(0){\0}{}&\int(0){\0}{}&\ln{1}&\cnt{6}&\bounds{2}{2}{1}{true}&&\cr
\l{3}&\int(2){\2}{}&\int(1){\2}{}&\int(0){\0}{}&\int(0){\0}{}&\int(0){\0}{}&\ln{1}&\cnt{36}&\bounds{9}{9}{1}{false}&&\cr
\l{5}&\int(2){\2}{}&\int(0){\0}{}&\int(0){\0}{}&\int(0){\0}{}&\int(1){\2}{}&\ln{1}&\cnt{24}&\bounds{5}{6}{5}{true}&&\cr
\l{7}&\int(2){4}{}&\int(1){4}{}&\int(0){2}{\spec{1}}&\int(0){\0}{}&\int(0){\0}{}&\ln{1}&\cnt{3}&\bounds{1}{1}{1}{true}&&\cr
\l{9}&\int(2){4}{}&\int(1){4}{}&\int(0){2}{\spec{2}}&\int(0){\0}{}&\int(0){\0}{}&\ln{1}&\cnt{18}&\bounds{4}{6}{5}{true}&&\cr
\l{11}&\int(2){4}{}&\int(1){2}{}&\int(0){\0}{}&\int(0){4}{}&\int(0){\0}{}&\ln{1}&\cnt{45}&\bounds{9}{9}{1}{true}&&\cr
\l{13}&\int(2){4}{}&\rat(1/2){5}{}&\int(0){\0}{}&\int(0){1}{}&\rat(1/2){2}{}&\ln{1}&\cnt{72}&\bounds{12}{12}{1}{false}&&\cr
\l{15}&\int(2){4}{}&\rat(1/2){1}{}&\int(0){5}{}&\int(0){\0}{}&\rat(1/2){3}{}&\ln{1}&\cnt{48}&\bounds{10}{12}{5}{true}&&\cr
\l{17}&\int(2){4}{}&\int(0){\0}{}&\int(0){1}{}&\int(0){5}{}&\int(1){1}{}&\ln{1}&\cnt{24}&\bounds{4}{4}{1}{true}&&\cr
\endorbit
\orbit{5}{52}{118}
\h{52}&\root\rat(7/2){3}{\spec{1}}&\rat(3/2){3}{}&\int(0){\0}{}&\int(0){\0}{}&\int(1){1}{}&&\cnt{604}&\bounds{118}{120}{fail}{false}&&\cr
\l{1}&\int(3){\4}{\spec{1}}&\int(0){\0}{}&\int(0){\0}{}&\int(0){\0}{}&\int(0){\0}{}&\ln{1}&\cnt{1}&\bounds{1}{1}{1}{true}&&\cr
\l{3}&\int(3){\4}{\spec{2}}&\int(0){\0}{}&\int(0){\0}{}&\int(0){\0}{}&\int(0){\0}{}&\ln{1}&\cnt{1}&\bounds{1}{1}{1}{true}&&\cr
\l{5}&\int(2){\2}{}&\int(1){\2}{}&\int(0){\0}{}&\int(0){\0}{}&\int(0){\0}{}&\ln{1}&\cnt{18}&\bounds{6}{6}{1}{true}&&\cr
\l{7}&\int(2){\2}{}&\int(0){\0}{}&\int(0){\0}{}&\int(0){\0}{}&\int(1){\2}{}&\ln{1}&\cnt{12}&\bounds{3}{4}{5}{true}&&\cr
\l{9}&\int(2){4}{}&\int(1){4}{}&\int(0){2}{}&\int(0){\0}{}&\int(0){\0}{}&\ln{2}&\cnt{45}&\bounds{9}{9}{1}{true}&&\cr
\l{11}&\int(2){4}{}&\rat(1/2){5}{}&\int(0){\0}{}&\int(0){1}{}&\rat(1/2){2}{}&\ln{2}&\cnt{72}&\bounds{12}{12}{1}{false}&&\cr
\l{13}&\int(2){4}{}&\int(0){\0}{}&\int(0){1}{}&\int(0){5}{}&\int(1){1}{}&\ln{1}&\cnt{36}&\bounds{6}{6}{1}{true}&&\cr
\endorbit
\orbit{5}{53}{98}
\h{53}&\rat(7/2){3}{\spec{1}}&\rat(3/2){3}{}&\int(0){\0}{}&\int(0){\0}{}&\root\int(1){1}{}&&\cnt{496}&\bounds{98}{98}{1}{false}&&\cr
\l{1}&\int(3){\4}{}&\int(0){\0}{}&\int(0){\0}{}&\int(0){\0}{}&\int(0){\0}{}&\ln{1}&\cnt{6}&\bounds{2}{2}{1}{false}&&\cr
\l{3}&\int(2){\2}{}&\int(1){\2}{}&\int(0){\0}{}&\int(0){\0}{}&\int(0){\0}{}&\ln{1}&\cnt{36}&\bounds{9}{9}{1}{false}&&\cr
\l{5}&\int(2){\2}{}&\int(0){\0}{}&\int(0){\0}{}&\int(0){\0}{}&\int(1){\2}{\spec{1}}&\ln{2}&\cnt{4}&\bounds{1}{1}{1}{false}&&\cr
\l{7}&\int(2){4}{}&\int(1){4}{}&\int(0){2}{}&\int(0){\0}{}&\int(0){\0}{}&\ln{2}&\cnt{45}&\bounds{9}{9}{1}{false}&&\cr
\l{9}&\int(2){4}{}&\rat(1/2){5}{}&\int(0){\0}{}&\int(0){1}{}&\rat(1/2){2}{}&\ln{2}&\cnt{36}&\bounds{6}{6}{1}{false}&&\cr
\l{11}&\int(2){4}{}&\int(0){\0}{}&\int(0){1}{}&\int(0){5}{}&\int(1){1}{}&\ln{1}&\cnt{36}&\bounds{6}{6}{1}{false}&&\cr
\endorbit
\orbit{5}{54}{90}
\h{54}&\rat(7/2){3}{\spec{1}}&\root\rat(3/2){3}{\spec{1}}&\int(0){\0}{}&\int(0){\0}{}&\int(1){1}{}&&\cnt{468}&\bounds{90}{90}{1}{false}&&\cr
\l{1}&\int(3){\4}{}&\int(0){\0}{}&\int(0){\0}{}&\int(0){\0}{}&\int(0){\0}{}&\ln{1}&\cnt{6}&\bounds{2}{2}{1}{false}&&\cr
\l{3}&\int(2){\2}{}&\int(1){\2}{}&\int(0){\0}{}&\int(0){\0}{}&\int(0){\0}{}&\ln{1}&\cnt{12}&\bounds{3}{3}{1}{false}&&\cr
\l{5}&\int(2){\2}{}&\int(0){\0}{}&\int(0){\0}{}&\int(0){\0}{}&\int(1){\2}{}&\ln{1}&\cnt{24}&\bounds{6}{6}{1}{false}&&\cr
\l{7}&\int(2){4}{}&\int(1){4}{}&\int(0){2}{}&\int(0){\0}{}&\int(0){\0}{}&\ln{1}&\cnt{45}&\bounds{9}{9}{1}{false}&&\cr
\l{9}&\int(2){4}{}&\int(1){2}{}&\int(0){\0}{}&\int(0){4}{}&\int(0){\0}{}&\ln{1}&\cnt{15}&\bounds{3}{3}{1}{false}&&\cr
\l{11}&\int(2){4}{}&\rat(1/2){5}{}&\int(0){\0}{}&\int(0){1}{}&\rat(1/2){2}{}&\ln{1}&\cnt{72}&\bounds{12}{12}{1}{false}&&\cr
\l{13}&\int(2){4}{}&\rat(1/2){1}{}&\int(0){5}{}&\int(0){\0}{}&\rat(1/2){3}{}&\ln{1}&\cnt{24}&\bounds{4}{4}{1}{false}&&\cr
\l{15}&\int(2){4}{}&\int(0){\0}{}&\int(0){1}{}&\int(0){5}{}&\int(1){1}{}&\ln{1}&\cnt{36}&\bounds{6}{6}{1}{false}&&\cr
\endorbit
\orbit{5}{55}{114}
\h{55}&\int(2){\2}{}&\int(2){\2}{}&\int(0){\0}{}&\int(0){\0}{}&\int(2){\2}{}&&\cnt{624}&\bounds{114}{128}{fail}{false}&&\cr
\l{1}&\int(2){\2}{}&\int(0){\0}{}&\int(0){\0}{}&\int(0){\0}{}&\int(1){\2}{}&\ln{2\sd}&\cnt{8}&\bounds{3}{4}{5}{true}&&\cr
\l{2}&\int(1){3}{}&\int(1){3}{}&\int(0){\0}{}&\int(0){\0}{}&\int(1){1}{}&\ln{1\sdd}&\cnt{72}&\bounds{12}{24}{5}{true}&&\cr
\l{3}&\int(1){4}{}&\int(1){5}{}&\int(0){\0}{}&\int(0){1}{}&\int(1){2}{}&\ln{8\sd}&\cnt{48}&\bounds{8}{8}{1}{true}&&\cr
\l{4}&\int(1){5}{}&\int(1){5}{}&\int(0){4}{}&\int(0){\0}{}&\int(1){1}{}&\ln{4\sd}&\cnt{30}&\bounds{6}{6}{1}{true}&&\cr
\l{5}&\int(2){\2}{}&\int(1){\2}{\spec{1}}&\int(0){\0}{}&\int(0){\0}{}&\int(0){\0}{}&\ln{4}&\cnt{4}&\bounds{1}{1}{1}{true}&&\cr
\endorbit
\orbit{5}{56}{108}
\h{56}&\rat(17/6){5}{\spec{2}}&\root\rat(5/6){5}{}&\rat(5/6){1}{}&\rat(3/2){3}{}&\int(0){\0}{}&&\cnt{476}&\bounds{108}{120}{fail}{false}&&\cr
\l{1}&\int(3){\4}{}&\int(0){\0}{}&\int(0){\0}{}&\int(0){\0}{}&\int(0){\0}{}&\ln{1}&\cnt{3}&\bounds{1}{1}{1}{true}&&\cr
\l{3}&\int(2){\2}{}&\int(1){\2}{}&\int(0){\0}{}&\int(0){\0}{}&\int(0){\0}{}&\ln{1}&\cnt{6}&\bounds{2}{2}{1}{true}&&\cr
\l{5}&\int(2){\2}{}&\int(0){\0}{}&\int(1){\2}{}&\int(0){\0}{}&\int(0){\0}{}&\ln{1}&\cnt{10}&\bounds{2}{2}{1}{true}&&\cr
\l{7}&\int(2){\2}{}&\int(0){\0}{}&\int(0){\0}{}&\int(1){\2}{}&\int(0){\0}{}&\ln{1}&\cnt{18}&\bounds{6}{6}{1}{true}&&\cr
\l{9}&\rat(3/2){3}{}&\rat(5/6){5}{}&\rat(1/6){5}{}&\rat(1/2){5}{}&\int(0){\0}{}&\ln{1}&\cnt{45}&\bounds{9}{9}{1}{true}&&\cr
\l{11}&\rat(3/2){3}{}&\rat(1/6){1}{}&\rat(5/6){1}{}&\rat(1/2){1}{}&\int(0){\0}{}&\ln{1}&\cnt{27}&\bounds{9}{9}{1}{true}&&\cr
\l{13}&\rat(3/2){3}{}&\int(0){\0}{}&\int(0){\0}{}&\rat(3/2){3}{}&\int(0){3}{}&\ln{1}&\cnt{24}&\bounds{4}{6}{5}{true}&&\cr
\l{15}&\rat(5/3){4}{}&\rat(2/3){4}{}&\rat(2/3){2}{}&\int(0){\0}{}&\int(0){\0}{}&\ln{1}&\cnt{15}&\bounds{3}{3}{1}{true}&&\cr
\l{17}&\rat(5/3){4}{}&\rat(1/3){2}{\spec{1}}&\int(0){\0}{}&\int(1){4}{}&\int(0){\0}{}&\ln{1}&\cnt{3}&\bounds{1}{1}{1}{true}&&\cr
\l{19}&\rat(5/3){4}{}&\rat(1/3){2}{\spec{2}}&\int(0){\0}{}&\int(1){4}{}&\int(0){\0}{}&\ln{1}&\cnt{9}&\bounds{3}{3}{1}{true}&&\cr
\l{21}&\rat(5/3){4}{}&\rat(5/6){5}{}&\int(0){\0}{}&\rat(1/2){1}{}&\int(0){2}{}&\ln{1}&\cnt{24}&\bounds{4}{6}{5}{true}&&\cr
\l{23}&\rat(5/3){4}{}&\int(0){\0}{}&\rat(1/3){4}{}&\int(1){2}{}&\int(0){\0}{}&\ln{1}&\cnt{30}&\bounds{6}{6}{1}{true}&&\cr
\l{25}&\rat(5/3){4}{}&\int(0){\0}{}&\rat(5/6){1}{}&\rat(1/2){5}{}&\int(0){1}{}&\ln{1}&\cnt{24}&\bounds{4}{6}{5}{true}&&\cr
\endorbit
\orbit{5}{57}{116}
\h{57}&\rat(17/6){5}{\spec{2}}&\rat(5/6){5}{}&\rat(5/6){1}{}&\rat(3/2){3}{}&\root\int(0){\0}{}&&\cnt{504}&\bounds{116}{116}{1}{false}&&\cr
\l{1}&\int(3){\4}{}&\int(0){\0}{}&\int(0){\0}{}&\int(0){\0}{}&\int(0){\0}{}&\ln{1}&\cnt{3}&\bounds{1}{1}{1}{false}&&\cr
\l{3}&\int(2){\2}{}&\int(1){\2}{}&\int(0){\0}{}&\int(0){\0}{}&\int(0){\0}{}&\ln{2}&\cnt{10}&\bounds{2}{2}{1}{false}&&\cr
\l{5}&\int(2){\2}{}&\int(0){\0}{}&\int(0){\0}{}&\int(1){\2}{}&\int(0){\0}{}&\ln{1}&\cnt{18}&\bounds{6}{6}{1}{false}&&\cr
\l{7}&\rat(3/2){3}{}&\rat(5/6){5}{}&\rat(1/6){5}{}&\rat(1/2){5}{}&\int(0){\0}{}&\ln{2}&\cnt{45}&\bounds{9}{9}{1}{false}&&\cr
\l{9}&\rat(3/2){3}{}&\int(0){\0}{}&\int(0){\0}{}&\rat(3/2){3}{}&\int(0){3}{}&\ln{1}&\cnt{12}&\bounds{4}{4}{1}{false}&&\cr
\l{11}&\rat(5/3){4}{}&\rat(2/3){4}{}&\rat(2/3){2}{}&\int(0){\0}{}&\int(0){\0}{}&\ln{1}&\cnt{25}&\bounds{5}{5}{1}{false}&&\cr
\l{13}&\rat(5/3){4}{}&\rat(1/3){2}{}&\int(0){\0}{}&\int(1){4}{}&\int(0){\0}{}&\ln{2}&\cnt{30}&\bounds{6}{6}{1}{false}&&\cr
\l{15}&\rat(5/3){4}{}&\rat(5/6){5}{}&\int(0){\0}{}&\rat(1/2){1}{}&\int(0){2}{}&\ln{2}&\cnt{12}&\bounds{4}{4}{1}{false}&&\cr
\endorbit
\orbit{5}{58}{108}
\h{58}&\root\rat(17/6){5}{\spec{2}}&\rat(5/6){5}{}&\rat(5/6){1}{}&\rat(3/2){3}{}&\int(0){\0}{}&&\cnt{500}&\bounds{108}{108}{1}{false}&&\cr
\l{1}&\int(3){\4}{}&\int(0){\0}{}&\int(0){\0}{}&\int(0){\0}{}&\int(0){\0}{}&\ln{1}&\cnt{3}&\bounds{1}{1}{1}{false}&&\cr
\l{3}&\rat(3/2){3}{}&\rat(5/6){5}{}&\rat(1/6){5}{}&\rat(1/2){5}{}&\int(0){\0}{}&\ln{2}&\cnt{45}&\bounds{9}{9}{1}{false}&&\cr
\l{5}&\rat(3/2){3}{}&\int(0){\0}{}&\int(0){\0}{}&\rat(3/2){3}{}&\int(0){3}{}&\ln{1}&\cnt{24}&\bounds{6}{6}{1}{false}&&\cr
\l{7}&\rat(5/3){4}{}&\rat(2/3){4}{}&\rat(2/3){2}{}&\int(0){\0}{}&\int(0){\0}{}&\ln{1}&\cnt{25}&\bounds{5}{5}{1}{false}&&\cr
\l{9}&\rat(5/3){4}{}&\rat(1/3){2}{}&\int(0){\0}{}&\int(1){4}{}&\int(0){\0}{}&\ln{2}&\cnt{30}&\bounds{6}{6}{1}{false}&&\cr
\l{11}&\rat(5/3){4}{}&\rat(5/6){5}{}&\int(0){\0}{}&\rat(1/2){1}{}&\int(0){2}{}&\ln{2}&\cnt{24}&\bounds{6}{6}{1}{false}&&\cr
\endorbit
\orbit{5}{59}{96}
\h{59}&\root\rat(17/6){5}{\spec{1}}&\rat(5/6){5}{}&\rat(5/6){1}{}&\rat(3/2){3}{}&\int(0){\0}{}&&\cnt{420}&\bounds{96}{96}{1}{false}&&\cr
\l{1}&\int(3){\4}{}&\int(0){\0}{}&\int(0){\0}{}&\int(0){\0}{}&\int(0){\0}{}&\ln{1}&\cnt{1}&\bounds{1}{1}{1}{false}&&\cr
\l{3}&\int(2){\2}{}&\int(1){\2}{}&\int(0){\0}{}&\int(0){\0}{}&\int(0){\0}{}&\ln{2}&\cnt{10}&\bounds{2}{2}{1}{false}&&\cr
\l{5}&\int(2){\2}{}&\int(0){\0}{}&\int(0){\0}{}&\int(1){\2}{}&\int(0){\0}{}&\ln{1}&\cnt{18}&\bounds{6}{6}{1}{false}&&\cr
\l{7}&\rat(3/2){3}{}&\rat(5/6){5}{}&\rat(1/6){5}{}&\rat(1/2){5}{}&\int(0){\0}{}&\ln{2}&\cnt{15}&\bounds{3}{3}{1}{false}&&\cr
\l{9}&\rat(3/2){3}{}&\int(0){\0}{}&\int(0){\0}{}&\rat(3/2){3}{}&\int(0){3}{}&\ln{1}&\cnt{8}&\bounds{2}{2}{1}{false}&&\cr
\l{11}&\rat(5/3){4}{}&\rat(2/3){4}{}&\rat(2/3){2}{}&\int(0){\0}{}&\int(0){\0}{}&\ln{1}&\cnt{25}&\bounds{5}{5}{1}{false}&&\cr
\l{13}&\rat(5/3){4}{}&\rat(1/3){2}{}&\int(0){\0}{}&\int(1){4}{}&\int(0){\0}{}&\ln{2}&\cnt{30}&\bounds{6}{6}{1}{false}&&\cr
\l{15}&\rat(5/3){4}{}&\rat(5/6){5}{}&\int(0){\0}{}&\rat(1/2){1}{}&\int(0){2}{}&\ln{2}&\cnt{24}&\bounds{6}{6}{1}{false}&&\cr
\endorbit
\orbit{5}{60}{92}
\h{60}&\rat(17/6){5}{\spec{2}}&\rat(5/6){5}{}&\rat(5/6){1}{}&\root\rat(3/2){3}{\spec{1}}&\int(0){\0}{}&&\cnt{420}&\bounds{92}{92}{1}{false}&&\cr
\l{1}&\int(3){\4}{}&\int(0){\0}{}&\int(0){\0}{}&\int(0){\0}{}&\int(0){\0}{}&\ln{1}&\cnt{3}&\bounds{1}{1}{1}{false}&&\cr
\l{3}&\int(2){\2}{}&\int(1){\2}{}&\int(0){\0}{}&\int(0){\0}{}&\int(0){\0}{}&\ln{1}&\cnt{10}&\bounds{2}{2}{1}{false}&&\cr
\l{5}&\int(2){\2}{}&\int(0){\0}{}&\int(1){\2}{}&\int(0){\0}{}&\int(0){\0}{}&\ln{1}&\cnt{10}&\bounds{2}{2}{1}{false}&&\cr
\l{7}&\int(2){\2}{}&\int(0){\0}{}&\int(0){\0}{}&\int(1){\2}{}&\int(0){\0}{}&\ln{1}&\cnt{6}&\bounds{2}{2}{1}{false}&&\cr
\l{9}&\rat(3/2){3}{}&\rat(5/6){5}{}&\rat(1/6){5}{}&\rat(1/2){5}{}&\int(0){\0}{}&\ln{1}&\cnt{45}&\bounds{9}{9}{1}{false}&&\cr
\l{11}&\rat(3/2){3}{}&\rat(1/6){1}{}&\rat(5/6){1}{}&\rat(1/2){1}{}&\int(0){\0}{}&\ln{1}&\cnt{15}&\bounds{3}{3}{1}{false}&&\cr
\l{13}&\rat(3/2){3}{}&\int(0){\0}{}&\int(0){\0}{}&\rat(3/2){3}{}&\int(0){3}{}&\ln{1}&\cnt{24}&\bounds{6}{6}{1}{false}&&\cr
\l{15}&\rat(5/3){4}{}&\rat(2/3){4}{}&\rat(2/3){2}{}&\int(0){\0}{}&\int(0){\0}{}&\ln{1}&\cnt{25}&\bounds{5}{5}{1}{false}&&\cr
\l{17}&\rat(5/3){4}{}&\rat(1/3){2}{}&\int(0){\0}{}&\int(1){4}{}&\int(0){\0}{}&\ln{1}&\cnt{30}&\bounds{6}{6}{1}{false}&&\cr
\l{19}&\rat(5/3){4}{}&\rat(5/6){5}{}&\int(0){\0}{}&\rat(1/2){1}{}&\int(0){2}{}&\ln{1}&\cnt{8}&\bounds{2}{2}{1}{false}&&\cr
\l{21}&\rat(5/3){4}{}&\int(0){\0}{}&\rat(1/3){4}{}&\int(1){2}{}&\int(0){\0}{}&\ln{1}&\cnt{10}&\bounds{2}{2}{1}{false}&&\cr
\l{23}&\rat(5/3){4}{}&\int(0){\0}{}&\rat(5/6){1}{}&\rat(1/2){5}{}&\int(0){1}{}&\ln{1}&\cnt{24}&\bounds{6}{6}{1}{false}&&\cr
\endorbit
\orbit{5}{61}{108}
\h{61}&\rat(4/3){4}{}&\rat(4/3){4}{}&\rat(4/3){2}{}&\int(0){\0}{}&\int(2){\2}{}&&\cnt{528}&\bounds{108}{132}{fail}{false}&&\cr
\l{1}&\rat(4/3){4}{}&\rat(2/3){5}{}&\int(0){\0}{}&\int(0){1}{}&\int(1){2}{}&\ln{6\sd}&\cnt{24}&\bounds{4}{4}{1}{true}&&\cr
\l{2}&\rat(2/3){2}{}&\rat(2/3){5}{}&\rat(2/3){1}{}&\int(0){\0}{}&\int(1){3}{}&\ln{3\sdd}&\cnt{48}&\bounds{10}{16}{5}{true}&&\cr
\l{3}&\int(1){\2}{}&\int(0){\0}{}&\int(0){\0}{}&\int(0){\0}{}&\int(2){\2}{}&\ln{3}&\cnt{8}&\bounds{2}{2}{1}{true}&&\cr
\l{5}&\int(1){3}{}&\int(1){3}{}&\int(0){\0}{}&\int(0){\0}{}&\int(1){1}{}&\ln{3}&\cnt{32}&\bounds{7}{8}{5}{true}&&\cr
\endorbit
\orbit{5}{62}{114}
\h{62}&\rat(17/6){5}{\spec{2}}&\rat(5/6){5}{}&\rat(4/3){4}{}&\root\int(0){\0}{}&\int(1){1}{}&&\cnt{504}&\bounds{114}{120}{fail}{false}&&\cr
\l{1}&\int(3){\4}{}&\int(0){\0}{}&\int(0){\0}{}&\int(0){\0}{}&\int(0){\0}{}&\ln{1}&\cnt{3}&\bounds{1}{1}{1}{true}&&\cr
\l{3}&\int(2){\2}{}&\int(1){\2}{}&\int(0){\0}{}&\int(0){\0}{}&\int(0){\0}{}&\ln{1}&\cnt{10}&\bounds{2}{2}{1}{true}&&\cr
\l{5}&\int(2){\2}{}&\int(0){\0}{}&\int(1){\2}{}&\int(0){\0}{}&\int(0){\0}{}&\ln{1}&\cnt{16}&\bounds{4}{4}{1}{true}&&\cr
\l{7}&\int(2){\2}{}&\int(0){\0}{}&\int(0){\0}{}&\int(0){\0}{}&\int(1){\2}{}&\ln{1}&\cnt{12}&\bounds{3}{4}{5}{true}&&\cr
\l{9}&\rat(3/2){3}{}&\rat(1/2){3}{}&\int(0){\0}{}&\int(0){\0}{}&\int(1){1}{}&\ln{1}&\cnt{30}&\bounds{6}{6}{1}{true}&&\cr
\l{11}&\rat(3/2){3}{}&\rat(5/6){5}{}&\rat(2/3){5}{}&\int(0){5}{}&\int(0){\0}{}&\ln{1}&\cnt{24}&\bounds{6}{6}{1}{true}&&\cr
\l{13}&\rat(3/2){3}{}&\int(0){\0}{}&\int(1){3}{}&\int(0){\0}{}&\rat(1/2){2}{}&\ln{1}&\cnt{48}&\bounds{10}{12}{5}{true}&&\cr
\l{15}&\rat(5/3){4}{}&\rat(2/3){4}{}&\rat(2/3){2}{}&\int(0){\0}{}&\int(0){\0}{}&\ln{1}&\cnt{30}&\bounds{6}{6}{1}{true}&&\cr
\l{17}&\rat(5/3){4}{}&\rat(5/6){5}{}&\int(0){\0}{}&\int(0){1}{}&\rat(1/2){2}{}&\ln{1}&\cnt{16}&\bounds{4}{4}{1}{true}&&\cr
\l{19}&\rat(5/3){4}{}&\rat(1/6){1}{}&\rat(2/3){5}{}&\int(0){\0}{}&\rat(1/2){3}{}&\ln{1}&\cnt{40}&\bounds{8}{8}{1}{true}&&\cr
\l{21}&\rat(5/3){4}{}&\int(0){\0}{}&\rat(4/3){4}{}&\int(0){2}{\spec{1}}&\int(0){\0}{}&\ln{1}&\cnt{1}&\bounds{1}{1}{1}{true}&&\cr
\l{23}&\rat(5/3){4}{}&\int(0){\0}{}&\rat(4/3){4}{}&\int(0){2}{\spec{2}}&\int(0){\0}{}&\ln{1}&\cnt{6}&\bounds{2}{2}{1}{true}&&\cr
\l{25}&\rat(5/3){4}{}&\int(0){\0}{}&\rat(1/3){1}{}&\int(0){5}{}&\int(1){1}{}&\ln{1}&\cnt{16}&\bounds{4}{4}{1}{true}&&\cr
\endorbit
\orbit{5}{63}{112}
\h{63}&\rat(17/6){5}{\spec{2}}&\root\rat(5/6){5}{}&\rat(4/3){4}{}&\int(0){\0}{}&\int(1){1}{}&&\cnt{476}&\bounds{112}{112}{1}{false}&&\cr
\l{1}&\int(3){\4}{}&\int(0){\0}{}&\int(0){\0}{}&\int(0){\0}{}&\int(0){\0}{}&\ln{1}&\cnt{3}&\bounds{1}{1}{1}{false}&&\cr
\l{3}&\int(2){\2}{}&\int(1){\2}{}&\int(0){\0}{}&\int(0){\0}{}&\int(0){\0}{}&\ln{1}&\cnt{6}&\bounds{2}{2}{1}{false}&&\cr
\l{5}&\int(2){\2}{}&\int(0){\0}{}&\int(1){\2}{}&\int(0){\0}{}&\int(0){\0}{}&\ln{1}&\cnt{16}&\bounds{4}{4}{1}{false}&&\cr
\l{7}&\int(2){\2}{}&\int(0){\0}{}&\int(0){\0}{}&\int(0){\0}{}&\int(1){\2}{}&\ln{1}&\cnt{12}&\bounds{4}{4}{1}{false}&&\cr
\l{9}&\rat(3/2){3}{}&\rat(1/2){3}{\spec{1}}&\int(0){\0}{}&\int(0){\0}{}&\int(1){1}{}&\ln{1}&\cnt{9}&\bounds{3}{3}{1}{false}&&\cr
\l{11}&\rat(3/2){3}{}&\rat(1/2){3}{\spec{2}}&\int(0){\0}{}&\int(0){\0}{}&\int(1){1}{}&\ln{1}&\cnt{3}&\bounds{1}{1}{1}{false}&&\cr
\l{13}&\rat(3/2){3}{}&\rat(5/6){5}{}&\rat(2/3){5}{}&\int(0){5}{}&\int(0){\0}{}&\ln{1}&\cnt{36}&\bounds{6}{6}{1}{false}&&\cr
\l{15}&\rat(3/2){3}{}&\int(0){\0}{}&\int(1){3}{}&\int(0){\0}{}&\rat(1/2){2}{}&\ln{1}&\cnt{48}&\bounds{12}{12}{1}{false}&&\cr
\l{17}&\rat(5/3){4}{}&\rat(2/3){4}{}&\rat(2/3){2}{}&\int(0){\0}{}&\int(0){\0}{}&\ln{1}&\cnt{18}&\bounds{6}{6}{1}{false}&&\cr
\l{19}&\rat(5/3){4}{}&\rat(5/6){5}{}&\int(0){\0}{}&\int(0){1}{}&\rat(1/2){2}{}&\ln{1}&\cnt{24}&\bounds{4}{4}{1}{false}&&\cr
\l{21}&\rat(5/3){4}{}&\rat(1/6){1}{}&\rat(2/3){5}{}&\int(0){\0}{}&\rat(1/2){3}{}&\ln{1}&\cnt{24}&\bounds{6}{6}{1}{false}&&\cr
\l{23}&\rat(5/3){4}{}&\int(0){\0}{}&\rat(4/3){4}{}&\int(0){2}{}&\int(0){\0}{}&\ln{1}&\cnt{15}&\bounds{3}{3}{1}{false}&&\cr
\l{25}&\rat(5/3){4}{}&\int(0){\0}{}&\rat(1/3){1}{}&\int(0){5}{}&\int(1){1}{}&\ln{1}&\cnt{24}&\bounds{4}{4}{1}{false}&&\cr
\endorbit
\orbit{5}{64}{100}
\h{64}&\root\rat(17/6){5}{\spec{2}}&\rat(5/6){5}{}&\rat(4/3){4}{}&\int(0){\0}{}&\int(1){1}{}&&\cnt{500}&\bounds{100}{100}{1}{false}&&\cr
\l{1}&\int(3){\4}{}&\int(0){\0}{}&\int(0){\0}{}&\int(0){\0}{}&\int(0){\0}{}&\ln{1}&\cnt{3}&\bounds{1}{1}{1}{false}&&\cr
\l{3}&\rat(3/2){3}{}&\rat(1/2){3}{}&\int(0){\0}{}&\int(0){\0}{}&\int(1){1}{}&\ln{1}&\cnt{30}&\bounds{6}{6}{1}{false}&&\cr
\l{5}&\rat(3/2){3}{}&\rat(5/6){5}{}&\rat(2/3){5}{}&\int(0){5}{}&\int(0){\0}{}&\ln{1}&\cnt{36}&\bounds{6}{6}{1}{false}&&\cr
\l{7}&\rat(3/2){3}{}&\int(0){\0}{}&\int(1){3}{}&\int(0){\0}{}&\rat(1/2){2}{}&\ln{1}&\cnt{48}&\bounds{12}{12}{1}{false}&&\cr
\l{9}&\rat(5/3){4}{}&\rat(2/3){4}{}&\rat(2/3){2}{}&\int(0){\0}{}&\int(0){\0}{}&\ln{1}&\cnt{30}&\bounds{6}{6}{1}{false}&&\cr
\l{11}&\rat(5/3){4}{}&\rat(5/6){5}{}&\int(0){\0}{}&\int(0){1}{}&\rat(1/2){2}{}&\ln{1}&\cnt{24}&\bounds{4}{4}{1}{false}&&\cr
\l{13}&\rat(5/3){4}{}&\rat(1/6){1}{}&\rat(2/3){5}{}&\int(0){\0}{}&\rat(1/2){3}{}&\ln{1}&\cnt{40}&\bounds{8}{8}{1}{false}&&\cr
\l{15}&\rat(5/3){4}{}&\int(0){\0}{}&\rat(4/3){4}{}&\int(0){2}{}&\int(0){\0}{}&\ln{1}&\cnt{15}&\bounds{3}{3}{1}{false}&&\cr
\l{17}&\rat(5/3){4}{}&\int(0){\0}{}&\rat(1/3){1}{}&\int(0){5}{}&\int(1){1}{}&\ln{1}&\cnt{24}&\bounds{4}{4}{1}{false}&&\cr
\endorbit
\orbit{5}{65}{96}
\h{65}&\rat(17/6){5}{\spec{2}}&\rat(5/6){5}{}&\root\rat(4/3){4}{\spec{1}}&\int(0){\0}{}&\int(1){1}{}&&\cnt{448}&\bounds{96}{96}{1}{false}&&\cr
\l{1}&\int(3){\4}{}&\int(0){\0}{}&\int(0){\0}{}&\int(0){\0}{}&\int(0){\0}{}&\ln{1}&\cnt{3}&\bounds{1}{1}{1}{false}&&\cr
\l{3}&\int(2){\2}{}&\int(1){\2}{}&\int(0){\0}{}&\int(0){\0}{}&\int(0){\0}{}&\ln{1}&\cnt{10}&\bounds{2}{2}{1}{false}&&\cr
\l{5}&\int(2){\2}{}&\int(0){\0}{}&\int(1){\2}{}&\int(0){\0}{}&\int(0){\0}{}&\ln{1}&\cnt{8}&\bounds{4}{4}{1}{false}&&\cr
\l{7}&\int(2){\2}{}&\int(0){\0}{}&\int(0){\0}{}&\int(0){\0}{}&\int(1){\2}{}&\ln{1}&\cnt{12}&\bounds{4}{4}{1}{false}&&\cr
\l{9}&\rat(3/2){3}{}&\rat(1/2){3}{}&\int(0){\0}{}&\int(0){\0}{}&\int(1){1}{}&\ln{1}&\cnt{30}&\bounds{6}{6}{1}{false}&&\cr
\l{11}&\rat(3/2){3}{}&\rat(5/6){5}{}&\rat(2/3){5}{}&\int(0){5}{}&\int(0){\0}{}&\ln{1}&\cnt{36}&\bounds{6}{6}{1}{false}&&\cr
\l{13}&\rat(3/2){3}{}&\int(0){\0}{}&\int(1){3}{}&\int(0){\0}{}&\rat(1/2){2}{}&\ln{1}&\cnt{24}&\bounds{6}{6}{1}{false}&&\cr
\l{15}&\rat(5/3){4}{}&\rat(2/3){4}{}&\rat(2/3){2}{\spec{1}}&\int(0){\0}{}&\int(0){\0}{}&\ln{1}&\cnt{5}&\bounds{1}{1}{1}{false}&&\cr
\l{17}&\rat(5/3){4}{}&\rat(2/3){4}{}&\rat(2/3){2}{\spec{2}}&\int(0){\0}{}&\int(0){\0}{}&\ln{1}&\cnt{5}&\bounds{1}{1}{1}{false}&&\cr
\l{19}&\rat(5/3){4}{}&\rat(5/6){5}{}&\int(0){\0}{}&\int(0){1}{}&\rat(1/2){2}{}&\ln{1}&\cnt{24}&\bounds{4}{4}{1}{false}&&\cr
\l{21}&\rat(5/3){4}{}&\rat(1/6){1}{}&\rat(2/3){5}{}&\int(0){\0}{}&\rat(1/2){3}{}&\ln{1}&\cnt{40}&\bounds{8}{8}{1}{false}&&\cr
\l{23}&\rat(5/3){4}{}&\int(0){\0}{}&\rat(4/3){4}{}&\int(0){2}{}&\int(0){\0}{}&\ln{1}&\cnt{15}&\bounds{3}{3}{1}{false}&&\cr
\l{25}&\rat(5/3){4}{}&\int(0){\0}{}&\rat(1/3){1}{}&\int(0){5}{}&\int(1){1}{}&\ln{1}&\cnt{12}&\bounds{2}{2}{1}{false}&&\cr
\endorbit
\orbit{5}{66}{92}
\h{66}&\rat(17/6){5}{\spec{2}}&\rat(5/6){5}{}&\rat(4/3){4}{}&\int(0){\0}{}&\root\int(1){1}{}&&\cnt{448}&\bounds{92}{92}{1}{false}&&\cr
\l{1}&\int(3){\4}{}&\int(0){\0}{}&\int(0){\0}{}&\int(0){\0}{}&\int(0){\0}{}&\ln{1}&\cnt{3}&\bounds{1}{1}{1}{false}&&\cr
\l{3}&\int(2){\2}{}&\int(1){\2}{}&\int(0){\0}{}&\int(0){\0}{}&\int(0){\0}{}&\ln{1}&\cnt{10}&\bounds{2}{2}{1}{false}&&\cr
\l{5}&\int(2){\2}{}&\int(0){\0}{}&\int(1){\2}{}&\int(0){\0}{}&\int(0){\0}{}&\ln{1}&\cnt{16}&\bounds{4}{4}{1}{false}&&\cr
\l{7}&\int(2){\2}{}&\int(0){\0}{}&\int(0){\0}{}&\int(0){\0}{}&\int(1){\2}{\spec{1}}&\ln{1}&\cnt{2}&\bounds{1}{1}{1}{false}&&\cr
\l{9}&\int(2){\2}{}&\int(0){\0}{}&\int(0){\0}{}&\int(0){\0}{}&\int(1){\2}{\spec{2}}&\ln{1}&\cnt{2}&\bounds{1}{1}{1}{false}&&\cr
\l{11}&\rat(3/2){3}{}&\rat(1/2){3}{}&\int(0){\0}{}&\int(0){\0}{}&\int(1){1}{}&\ln{1}&\cnt{30}&\bounds{6}{6}{1}{false}&&\cr
\l{13}&\rat(3/2){3}{}&\rat(5/6){5}{}&\rat(2/3){5}{}&\int(0){5}{}&\int(0){\0}{}&\ln{1}&\cnt{36}&\bounds{6}{6}{1}{false}&&\cr
\l{15}&\rat(3/2){3}{}&\int(0){\0}{}&\int(1){3}{}&\int(0){\0}{}&\rat(1/2){2}{}&\ln{1}&\cnt{24}&\bounds{6}{6}{1}{false}&&\cr
\l{17}&\rat(5/3){4}{}&\rat(2/3){4}{}&\rat(2/3){2}{}&\int(0){\0}{}&\int(0){\0}{}&\ln{1}&\cnt{30}&\bounds{6}{6}{1}{false}&&\cr
\l{19}&\rat(5/3){4}{}&\rat(5/6){5}{}&\int(0){\0}{}&\int(0){1}{}&\rat(1/2){2}{}&\ln{1}&\cnt{12}&\bounds{2}{2}{1}{false}&&\cr
\l{21}&\rat(5/3){4}{}&\rat(1/6){1}{}&\rat(2/3){5}{}&\int(0){\0}{}&\rat(1/2){3}{}&\ln{1}&\cnt{20}&\bounds{4}{4}{1}{false}&&\cr
\l{23}&\rat(5/3){4}{}&\int(0){\0}{}&\rat(4/3){4}{}&\int(0){2}{}&\int(0){\0}{}&\ln{1}&\cnt{15}&\bounds{3}{3}{1}{false}&&\cr
\l{25}&\rat(5/3){4}{}&\int(0){\0}{}&\rat(1/3){1}{}&\int(0){5}{}&\int(1){1}{}&\ln{1}&\cnt{24}&\bounds{4}{4}{1}{false}&&\cr
\endorbit
\orbit{5}{67}{88}
\h{67}&\root\rat(17/6){5}{\spec{1}}&\rat(5/6){5}{}&\rat(4/3){4}{}&\int(0){\0}{}&\int(1){1}{}&&\cnt{420}&\bounds{88}{88}{1}{false}&&\cr
\l{1}&\int(3){\4}{}&\int(0){\0}{}&\int(0){\0}{}&\int(0){\0}{}&\int(0){\0}{}&\ln{1}&\cnt{1}&\bounds{1}{1}{1}{false}&&\cr
\l{3}&\int(2){\2}{}&\int(1){\2}{}&\int(0){\0}{}&\int(0){\0}{}&\int(0){\0}{}&\ln{1}&\cnt{10}&\bounds{2}{2}{1}{false}&&\cr
\l{5}&\int(2){\2}{}&\int(0){\0}{}&\int(1){\2}{}&\int(0){\0}{}&\int(0){\0}{}&\ln{1}&\cnt{16}&\bounds{4}{4}{1}{false}&&\cr
\l{7}&\int(2){\2}{}&\int(0){\0}{}&\int(0){\0}{}&\int(0){\0}{}&\int(1){\2}{}&\ln{1}&\cnt{12}&\bounds{4}{4}{1}{false}&&\cr
\l{9}&\rat(3/2){3}{}&\rat(1/2){3}{}&\int(0){\0}{}&\int(0){\0}{}&\int(1){1}{}&\ln{1}&\cnt{10}&\bounds{2}{2}{1}{false}&&\cr
\l{11}&\rat(3/2){3}{}&\rat(5/6){5}{}&\rat(2/3){5}{}&\int(0){5}{}&\int(0){\0}{}&\ln{1}&\cnt{12}&\bounds{2}{2}{1}{false}&&\cr
\l{13}&\rat(3/2){3}{}&\int(0){\0}{}&\int(1){3}{}&\int(0){\0}{}&\rat(1/2){2}{}&\ln{1}&\cnt{16}&\bounds{4}{4}{1}{false}&&\cr
\l{15}&\rat(5/3){4}{}&\rat(2/3){4}{}&\rat(2/3){2}{}&\int(0){\0}{}&\int(0){\0}{}&\ln{1}&\cnt{30}&\bounds{6}{6}{1}{false}&&\cr
\l{17}&\rat(5/3){4}{}&\rat(5/6){5}{}&\int(0){\0}{}&\int(0){1}{}&\rat(1/2){2}{}&\ln{1}&\cnt{24}&\bounds{4}{4}{1}{false}&&\cr
\l{19}&\rat(5/3){4}{}&\rat(1/6){1}{}&\rat(2/3){5}{}&\int(0){\0}{}&\rat(1/2){3}{}&\ln{1}&\cnt{40}&\bounds{8}{8}{1}{false}&&\cr
\l{21}&\rat(5/3){4}{}&\int(0){\0}{}&\rat(4/3){4}{}&\int(0){2}{}&\int(0){\0}{}&\ln{1}&\cnt{15}&\bounds{3}{3}{1}{false}&&\cr
\l{23}&\rat(5/3){4}{}&\int(0){\0}{}&\rat(1/3){1}{}&\int(0){5}{}&\int(1){1}{}&\ln{1}&\cnt{24}&\bounds{4}{4}{1}{false}&&\cr
\endorbit
\orbit{5}{68}{84}
\h{68}&\rat(17/6){5}{\spec{2}}&\rat(5/6){5}{}&\root\rat(4/3){4}{\spec{2}}&\int(0){\0}{}&\int(1){1}{}&&\cnt{392}&\bounds{84}{84}{1}{false}&&\cr
\l{1}&\int(3){\4}{}&\int(0){\0}{}&\int(0){\0}{}&\int(0){\0}{}&\int(0){\0}{}&\ln{1}&\cnt{3}&\bounds{1}{1}{1}{false}&&\cr
\l{3}&\int(2){\2}{}&\int(1){\2}{}&\int(0){\0}{}&\int(0){\0}{}&\int(0){\0}{}&\ln{1}&\cnt{10}&\bounds{2}{2}{1}{false}&&\cr
\l{5}&\int(2){\2}{}&\int(0){\0}{}&\int(0){\0}{}&\int(0){\0}{}&\int(1){\2}{}&\ln{1}&\cnt{12}&\bounds{4}{4}{1}{false}&&\cr
\l{7}&\rat(3/2){3}{}&\rat(1/2){3}{}&\int(0){\0}{}&\int(0){\0}{}&\int(1){1}{}&\ln{1}&\cnt{30}&\bounds{6}{6}{1}{false}&&\cr
\l{9}&\rat(3/2){3}{}&\int(0){\0}{}&\int(1){3}{}&\int(0){\0}{}&\rat(1/2){2}{}&\ln{1}&\cnt{48}&\bounds{12}{12}{1}{false}&&\cr
\l{11}&\rat(5/3){4}{}&\rat(2/3){4}{}&\rat(2/3){2}{}&\int(0){\0}{}&\int(0){\0}{}&\ln{1}&\cnt{30}&\bounds{6}{6}{1}{false}&&\cr
\l{13}&\rat(5/3){4}{}&\rat(5/6){5}{}&\int(0){\0}{}&\int(0){1}{}&\rat(1/2){2}{}&\ln{1}&\cnt{24}&\bounds{4}{4}{1}{false}&&\cr
\l{15}&\rat(5/3){4}{}&\int(0){\0}{}&\rat(4/3){4}{}&\int(0){2}{}&\int(0){\0}{}&\ln{1}&\cnt{15}&\bounds{3}{3}{1}{false}&&\cr
\l{17}&\rat(5/3){4}{}&\int(0){\0}{}&\rat(1/3){1}{}&\int(0){5}{}&\int(1){1}{}&\ln{1}&\cnt{24}&\bounds{4}{4}{1}{false}&&\cr
\endorbit
\orbit{5}{69}{112}
\h{69}&\root\rat(5/6){5}{}&\rat(4/3){4}{}&\rat(3/2){3}{}&\rat(4/3){2}{}&\int(1){2}{}&&\cnt{428}&\bounds{112}{120}{fail}{false}&&\cr
\l{1}&\rat(1/2){3}{\spec{1}}&\int(0){\0}{}&\rat(3/2){3}{}&\int(0){\0}{}&\int(1){2}{}&\ln{1}&\cnt{3}&\bounds{1}{1}{1}{true}&&\cr
\l{3}&\rat(1/2){3}{\spec{2}}&\int(0){\0}{}&\rat(3/2){3}{}&\int(0){\0}{}&\int(1){2}{}&\ln{1}&\cnt{1}&\bounds{1}{1}{1}{true}&&\cr
\l{5}&\rat(2/3){4}{}&\rat(4/3){4}{}&\int(1){2}{}&\int(0){\0}{}&\int(0){\0}{}&\ln{2}&\cnt{9}&\bounds{3}{3}{1}{true}&&\cr
\l{7}&\rat(2/3){4}{}&\rat(2/3){5}{}&\int(0){\0}{}&\rat(2/3){1}{}&\int(1){2}{}&\ln{1}&\cnt{12}&\bounds{4}{4}{1}{true}&&\cr
\l{9}&\rat(5/6){5}{}&\int(1){3}{}&\rat(1/2){5}{}&\rat(2/3){1}{}&\int(0){\0}{}&\ln{2}&\cnt{24}&\bounds{6}{6}{1}{true}&&\cr
\l{11}&\rat(5/6){5}{}&\rat(4/3){4}{}&\int(0){\0}{}&\rat(1/3){5}{}&\rat(1/2){3}{}&\ln{2}&\cnt{16}&\bounds{4}{4}{1}{true}&&\cr
\l{13}&\rat(5/6){5}{}&\rat(2/3){2}{}&\rat(1/2){1}{}&\int(0){\0}{}&\int(1){2}{}&\ln{2}&\cnt{18}&\bounds{4}{6}{5}{true}&&\cr
\l{15}&\rat(5/6){5}{}&\rat(2/3){5}{}&\int(1){4}{}&\int(0){\0}{}&\rat(1/2){1}{}&\ln{2}&\cnt{24}&\bounds{6}{6}{1}{true}&&\cr
\l{17}&\rat(5/6){5}{}&\rat(1/3){1}{}&\rat(3/2){3}{}&\rat(1/3){5}{}&\int(0){\0}{}&\ln{1}&\cnt{16}&\bounds{4}{4}{1}{true}&&\cr
\endorbit
\orbit{5}{70}{104}
\h{70}&\rat(5/6){5}{}&\rat(4/3){4}{}&\rat(3/2){3}{}&\rat(4/3){2}{}&\root\int(1){2}{}&&\cnt{400}&\bounds{104}{104}{1}{false}&&\cr
\l{1}&\rat(1/2){3}{}&\int(0){\0}{}&\rat(3/2){3}{}&\int(0){\0}{}&\int(1){2}{}&\ln{1}&\cnt{10}&\bounds{2}{2}{1}{false}&&\cr
\l{3}&\rat(2/3){4}{}&\rat(4/3){4}{}&\int(1){2}{}&\int(0){\0}{}&\int(0){\0}{}&\ln{2}&\cnt{15}&\bounds{3}{3}{1}{false}&&\cr
\l{5}&\rat(2/3){4}{}&\rat(2/3){5}{}&\int(0){\0}{}&\rat(2/3){1}{}&\int(1){2}{}&\ln{1}&\cnt{20}&\bounds{4}{4}{1}{false}&&\cr
\l{7}&\rat(5/6){5}{}&\int(1){3}{}&\rat(1/2){5}{}&\rat(2/3){1}{}&\int(0){\0}{}&\ln{2}&\cnt{24}&\bounds{6}{6}{1}{false}&&\cr
\l{9}&\rat(5/6){5}{}&\rat(4/3){4}{}&\int(0){\0}{}&\rat(1/3){5}{}&\rat(1/2){3}{}&\ln{2}&\cnt{8}&\bounds{2}{2}{1}{false}&&\cr
\l{11}&\rat(5/6){5}{}&\rat(2/3){2}{}&\rat(1/2){1}{}&\int(0){\0}{}&\int(1){2}{}&\ln{2}&\cnt{18}&\bounds{6}{6}{1}{false}&&\cr
\l{13}&\rat(5/6){5}{}&\rat(2/3){5}{}&\int(1){4}{}&\int(0){\0}{}&\rat(1/2){1}{}&\ln{2}&\cnt{12}&\bounds{4}{4}{1}{false}&&\cr
\l{15}&\rat(5/6){5}{}&\rat(1/3){1}{}&\rat(3/2){3}{}&\rat(1/3){5}{}&\int(0){\0}{}&\ln{1}&\cnt{16}&\bounds{4}{4}{1}{false}&&\cr
\endorbit
\orbit{5}{71}{100}
\h{71}&\rat(5/6){5}{}&\root\rat(4/3){4}{\spec{1}}&\rat(3/2){3}{}&\rat(4/3){2}{}&\int(1){2}{}&&\cnt{400}&\bounds{100}{100}{1}{false}&&\cr
\l{1}&\rat(1/2){3}{}&\int(0){\0}{}&\rat(3/2){3}{}&\int(0){\0}{}&\int(1){2}{}&\ln{1}&\cnt{10}&\bounds{2}{2}{1}{false}&&\cr
\l{3}&\rat(2/3){4}{}&\rat(4/3){4}{}&\int(1){2}{}&\int(0){\0}{}&\int(0){\0}{}&\ln{1}&\cnt{15}&\bounds{3}{3}{1}{false}&&\cr
\l{5}&\rat(2/3){4}{}&\rat(2/3){5}{}&\int(0){\0}{}&\rat(2/3){1}{}&\int(1){2}{}&\ln{1}&\cnt{20}&\bounds{4}{4}{1}{false}&&\cr
\l{7}&\rat(2/3){4}{}&\int(0){\0}{}&\int(1){4}{}&\rat(4/3){2}{}&\int(0){\0}{}&\ln{1}&\cnt{15}&\bounds{3}{3}{1}{false}&&\cr
\l{9}&\rat(5/6){5}{}&\int(1){3}{}&\rat(1/2){5}{}&\rat(2/3){1}{}&\int(0){\0}{}&\ln{1}&\cnt{12}&\bounds{4}{4}{1}{false}&&\cr
\l{11}&\rat(5/6){5}{}&\rat(4/3){4}{}&\int(0){\0}{}&\rat(1/3){5}{}&\rat(1/2){3}{}&\ln{1}&\cnt{16}&\bounds{4}{4}{1}{false}&&\cr
\l{13}&\rat(5/6){5}{}&\rat(2/3){2}{\spec{1}}&\rat(1/2){1}{}&\int(0){\0}{}&\int(1){2}{}&\ln{1}&\cnt{3}&\bounds{1}{1}{1}{false}&&\cr
\l{15}&\rat(5/6){5}{}&\rat(2/3){2}{\spec{2}}&\rat(1/2){1}{}&\int(0){\0}{}&\int(1){2}{}&\ln{1}&\cnt{3}&\bounds{1}{1}{1}{false}&&\cr
\l{17}&\rat(5/6){5}{}&\rat(2/3){5}{}&\int(1){4}{}&\int(0){\0}{}&\rat(1/2){1}{}&\ln{1}&\cnt{24}&\bounds{6}{6}{1}{false}&&\cr
\l{19}&\rat(5/6){5}{}&\rat(2/3){5}{}&\rat(1/2){1}{}&\int(1){3}{}&\int(0){\0}{}&\ln{1}&\cnt{24}&\bounds{6}{6}{1}{false}&&\cr
\l{21}&\rat(5/6){5}{}&\rat(1/3){1}{}&\rat(3/2){3}{}&\rat(1/3){5}{}&\int(0){\0}{}&\ln{1}&\cnt{8}&\bounds{2}{2}{1}{false}&&\cr
\l{23}&\rat(5/6){5}{}&\rat(1/3){1}{}&\int(0){\0}{}&\rat(4/3){2}{}&\rat(1/2){1}{}&\ln{1}&\cnt{8}&\bounds{2}{2}{1}{false}&&\cr
\l{25}&\rat(5/6){5}{}&\int(0){\0}{}&\int(1){2}{}&\rat(2/3){1}{}&\rat(1/2){3}{}&\ln{1}&\cnt{24}&\bounds{6}{6}{1}{false}&&\cr
\l{27}&\rat(5/6){5}{}&\int(0){\0}{}&\rat(1/2){5}{}&\rat(2/3){4}{}&\int(1){2}{}&\ln{1}&\cnt{18}&\bounds{6}{6}{1}{false}&&\cr
\endorbit
\orbit{5}{72}{92}
\h{72}&\rat(5/6){5}{}&\rat(4/3){4}{}&\root\rat(3/2){3}{\spec{1}}&\rat(4/3){2}{}&\int(1){2}{}&&\cnt{372}&\bounds{92}{92}{1}{false}&&\cr
\l{1}&\rat(1/2){3}{}&\int(0){\0}{}&\rat(3/2){3}{}&\int(0){\0}{}&\int(1){2}{}&\ln{1}&\cnt{10}&\bounds{2}{2}{1}{false}&&\cr
\l{3}&\rat(2/3){4}{}&\rat(4/3){4}{}&\int(1){2}{}&\int(0){\0}{}&\int(0){\0}{}&\ln{1}&\cnt{5}&\bounds{1}{1}{1}{false}&&\cr
\l{5}&\rat(2/3){4}{}&\rat(2/3){5}{}&\int(0){\0}{}&\rat(2/3){1}{}&\int(1){2}{}&\ln{1}&\cnt{20}&\bounds{4}{4}{1}{false}&&\cr
\l{7}&\rat(2/3){4}{}&\int(0){\0}{}&\int(1){4}{}&\rat(4/3){2}{}&\int(0){\0}{}&\ln{1}&\cnt{15}&\bounds{3}{3}{1}{false}&&\cr
\l{9}&\rat(5/6){5}{}&\int(1){3}{}&\rat(1/2){5}{}&\rat(2/3){1}{}&\int(0){\0}{}&\ln{1}&\cnt{24}&\bounds{6}{6}{1}{false}&&\cr
\l{11}&\rat(5/6){5}{}&\rat(4/3){4}{}&\int(0){\0}{}&\rat(1/3){5}{}&\rat(1/2){3}{}&\ln{1}&\cnt{16}&\bounds{4}{4}{1}{false}&&\cr
\l{13}&\rat(5/6){5}{}&\rat(2/3){2}{}&\rat(1/2){1}{}&\int(0){\0}{}&\int(1){2}{}&\ln{1}&\cnt{6}&\bounds{2}{2}{1}{false}&&\cr
\l{15}&\rat(5/6){5}{}&\rat(2/3){5}{}&\int(1){4}{}&\int(0){\0}{}&\rat(1/2){1}{}&\ln{1}&\cnt{24}&\bounds{6}{6}{1}{false}&&\cr
\l{17}&\rat(5/6){5}{}&\rat(2/3){5}{}&\rat(1/2){1}{}&\int(1){3}{}&\int(0){\0}{}&\ln{1}&\cnt{8}&\bounds{2}{2}{1}{false}&&\cr
\l{19}&\rat(5/6){5}{}&\rat(1/3){1}{}&\rat(3/2){3}{}&\rat(1/3){5}{}&\int(0){\0}{}&\ln{1}&\cnt{16}&\bounds{4}{4}{1}{false}&&\cr
\l{21}&\rat(5/6){5}{}&\rat(1/3){1}{}&\int(0){\0}{}&\rat(4/3){2}{}&\rat(1/2){1}{}&\ln{1}&\cnt{16}&\bounds{4}{4}{1}{false}&&\cr
\l{23}&\rat(5/6){5}{}&\int(0){\0}{}&\int(1){2}{}&\rat(2/3){1}{}&\rat(1/2){3}{}&\ln{1}&\cnt{8}&\bounds{2}{2}{1}{false}&&\cr
\l{25}&\rat(5/6){5}{}&\int(0){\0}{}&\rat(1/2){5}{}&\rat(2/3){4}{}&\int(1){2}{}&\ln{1}&\cnt{18}&\bounds{6}{6}{1}{false}&&\cr
\endorbit
\orbit{5}{73}{88}
\h{73}&\rat(5/6){5}{}&\root\rat(4/3){4}{\spec{2}}&\rat(3/2){3}{}&\rat(4/3){2}{}&\int(1){2}{}&&\cnt{344}&\bounds{88}{88}{1}{false}&&\cr
\l{1}&\rat(1/2){3}{}&\int(0){\0}{}&\rat(3/2){3}{}&\int(0){\0}{}&\int(1){2}{}&\ln{1}&\cnt{10}&\bounds{2}{2}{1}{false}&&\cr
\l{3}&\rat(2/3){4}{}&\rat(4/3){4}{}&\int(1){2}{}&\int(0){\0}{}&\int(0){\0}{}&\ln{1}&\cnt{15}&\bounds{3}{3}{1}{false}&&\cr
\l{5}&\rat(2/3){4}{}&\int(0){\0}{}&\int(1){4}{}&\rat(4/3){2}{}&\int(0){\0}{}&\ln{1}&\cnt{15}&\bounds{3}{3}{1}{false}&&\cr
\l{7}&\rat(5/6){5}{}&\int(1){3}{}&\rat(1/2){5}{}&\rat(2/3){1}{}&\int(0){\0}{}&\ln{1}&\cnt{24}&\bounds{6}{6}{1}{false}&&\cr
\l{9}&\rat(5/6){5}{}&\rat(4/3){4}{}&\int(0){\0}{}&\rat(1/3){5}{}&\rat(1/2){3}{}&\ln{1}&\cnt{16}&\bounds{4}{4}{1}{false}&&\cr
\l{11}&\rat(5/6){5}{}&\rat(2/3){2}{}&\rat(1/2){1}{}&\int(0){\0}{}&\int(1){2}{}&\ln{1}&\cnt{18}&\bounds{6}{6}{1}{false}&&\cr
\l{13}&\rat(5/6){5}{}&\rat(1/3){1}{}&\rat(3/2){3}{}&\rat(1/3){5}{}&\int(0){\0}{}&\ln{1}&\cnt{16}&\bounds{4}{4}{1}{false}&&\cr
\l{15}&\rat(5/6){5}{}&\rat(1/3){1}{}&\int(0){\0}{}&\rat(4/3){2}{}&\rat(1/2){1}{}&\ln{1}&\cnt{16}&\bounds{4}{4}{1}{false}&&\cr
\l{17}&\rat(5/6){5}{}&\int(0){\0}{}&\int(1){2}{}&\rat(2/3){1}{}&\rat(1/2){3}{}&\ln{1}&\cnt{24}&\bounds{6}{6}{1}{false}&&\cr
\l{19}&\rat(5/6){5}{}&\int(0){\0}{}&\rat(1/2){5}{}&\rat(2/3){4}{}&\int(1){2}{}&\ln{1}&\cnt{18}&\bounds{6}{6}{1}{false}&&\cr
\endorbit
\orbit{5}{74}{118}
\h{74}&\rat(5/6){5}{}&\rat(5/6){5}{}&\rat(4/3){4}{}&\root\int(0){\0}{}&\int(3){1}{}&&\cnt{528}&\bounds{118}{118}{1}{false}&&\cr
\l{1}&\rat(2/3){4}{}&\rat(1/6){1}{}&\rat(2/3){5}{}&\int(0){\0}{}&\rat(3/2){3}{}&\ln{2\sd}&\cnt{50}&\bounds{10}{10}{1}{false}&&\cr
\l{2}&\rat(5/6){5}{}&\int(0){\0}{}&\rat(2/3){2}{}&\int(0){1}{}&\rat(3/2){3}{}&\ln{2\sd}&\cnt{24}&\bounds{6}{6}{1}{false}&&\cr
\l{3}&\rat(5/6){5}{}&\int(0){\0}{}&\rat(2/3){5}{}&\int(0){4}{\spec{1}}&\rat(3/2){2}{}&\ln{2\sd}&\cnt{12}&\bounds{4}{4}{1}{false}&&\cr
\l{4}&\rat(5/6){5}{}&\int(0){\0}{}&\rat(2/3){5}{}&\int(0){4}{\spec{2}}&\rat(3/2){2}{}&\ln{2\sd}&\cnt{2}&\bounds{1}{1}{1}{false}&&\cr
\l{5}&\int(1){\2}{}&\int(0){\0}{}&\int(0){\0}{}&\int(0){\0}{}&\int(2){\2}{}&\ln{2}&\cnt{15}&\bounds{3}{3}{1}{false}&&\cr
\l{7}&\rat(1/2){3}{}&\int(0){\0}{}&\int(1){3}{}&\int(0){\0}{}&\rat(3/2){2}{}&\ln{2}&\cnt{40}&\bounds{8}{8}{1}{false}&&\cr
\l{9}&\rat(2/3){4}{}&\rat(5/6){5}{}&\int(0){\0}{}&\int(0){1}{}&\rat(3/2){2}{}&\ln{2}&\cnt{20}&\bounds{4}{4}{1}{false}&&\cr
\l{11}&\rat(5/6){5}{}&\rat(5/6){5}{}&\rat(4/3){4}{}&\int(0){\0}{}&\int(0){1}{\spec{1}}&\ln{2}&\cnt{1}&\bounds{1}{1}{1}{false}&&\cr
\l{13}&\rat(5/6){5}{}&\rat(5/6){5}{}&\rat(1/3){4}{}&\int(0){\0}{}&\int(1){1}{}&\ln{1}&\cnt{24}&\bounds{6}{6}{1}{false}&&\cr
\endorbit
\orbit{5}{75}{108}
\h{75}&\rat(5/6){5}{}&\rat(5/6){5}{}&\rat(4/3){4}{}&\int(0){\0}{}&\root\int(3){1}{}&&\cnt{552}&\bounds{108}{108}{1}{false}&&\cr
\l{1}&\rat(2/3){4}{}&\rat(1/6){1}{}&\rat(2/3){5}{}&\int(0){\0}{}&\rat(3/2){3}{}&\ln{2\sd}&\cnt{50}&\bounds{10}{10}{1}{false}&&\cr
\l{2}&\rat(5/6){5}{}&\int(0){\0}{}&\rat(2/3){2}{}&\int(0){1}{}&\rat(3/2){3}{}&\ln{2\sd}&\cnt{36}&\bounds{6}{6}{1}{false}&&\cr
\l{3}&\rat(5/6){5}{}&\int(0){\0}{}&\rat(2/3){5}{}&\int(0){4}{}&\rat(3/2){2}{}&\ln{2\sd}&\cnt{30}&\bounds{6}{6}{1}{false}&&\cr
\l{4}&\int(1){\2}{}&\int(0){\0}{}&\int(0){\0}{}&\int(0){\0}{}&\int(2){\2}{}&\ln{2}&\cnt{5}&\bounds{1}{1}{1}{false}&&\cr
\l{6}&\rat(1/2){3}{}&\int(0){\0}{}&\int(1){3}{}&\int(0){\0}{}&\rat(3/2){2}{}&\ln{2}&\cnt{40}&\bounds{8}{8}{1}{false}&&\cr
\l{8}&\rat(2/3){4}{}&\rat(5/6){5}{}&\int(0){\0}{}&\int(0){1}{}&\rat(3/2){2}{}&\ln{2}&\cnt{30}&\bounds{5}{5}{1}{false}&&\cr
\l{10}&\rat(5/6){5}{}&\rat(5/6){5}{}&\rat(4/3){4}{}&\int(0){\0}{}&\int(0){1}{\spec{1}}&\ln{2}&\cnt{1}&\bounds{1}{1}{1}{false}&&\cr
\l{12}&\rat(5/6){5}{}&\rat(5/6){5}{}&\rat(1/3){4}{}&\int(0){\0}{}&\int(1){1}{}&\ln{1}&\cnt{8}&\bounds{2}{2}{1}{false}&&\cr
\endorbit
\orbit{5}{76}{104}
\h{76}&\root\rat(5/6){5}{}&\rat(5/6){5}{}&\rat(4/3){4}{}&\int(0){\0}{}&\int(3){1}{}&&\cnt{500}&\bounds{104}{104}{1}{false}&&\cr
\l{1}&\int(1){\2}{}&\int(0){\0}{}&\int(0){\0}{}&\int(0){\0}{}&\int(2){\2}{}&\ln{1}&\cnt{9}&\bounds{3}{3}{1}{false}&&\cr
\l{3}&\rat(1/2){3}{\spec{1}}&\int(0){\0}{}&\int(1){3}{}&\int(0){\0}{}&\rat(3/2){2}{}&\ln{1}&\cnt{12}&\bounds{3}{3}{1}{false}&&\cr
\l{5}&\rat(1/2){3}{\spec{2}}&\int(0){\0}{}&\int(1){3}{}&\int(0){\0}{}&\rat(3/2){2}{}&\ln{1}&\cnt{4}&\bounds{1}{1}{1}{false}&&\cr
\l{7}&\rat(2/3){4}{}&\rat(5/6){5}{}&\int(0){\0}{}&\int(0){1}{}&\rat(3/2){2}{}&\ln{1}&\cnt{18}&\bounds{3}{3}{1}{false}&&\cr
\l{9}&\rat(2/3){4}{}&\rat(1/6){1}{}&\rat(2/3){5}{}&\int(0){\0}{}&\rat(3/2){3}{}&\ln{1}&\cnt{30}&\bounds{6}{6}{1}{false}&&\cr
\l{11}&\rat(5/6){5}{}&\rat(2/3){4}{}&\int(0){\0}{}&\int(0){5}{}&\rat(3/2){3}{}&\ln{1}&\cnt{30}&\bounds{5}{5}{1}{false}&&\cr
\l{13}&\rat(5/6){5}{}&\rat(1/3){2}{}&\rat(1/3){1}{}&\int(0){\0}{}&\rat(3/2){2}{}&\ln{1}&\cnt{40}&\bounds{8}{8}{1}{false}&&\cr
\l{15}&\rat(5/6){5}{}&\rat(5/6){5}{}&\rat(4/3){4}{}&\int(0){\0}{}&\int(0){1}{\spec{1}}&\ln{1}&\cnt{1}&\bounds{1}{1}{1}{false}&&\cr
\l{17}&\rat(5/6){5}{}&\rat(5/6){5}{}&\rat(4/3){4}{}&\int(0){\0}{}&\int(0){1}{\spec{2}}&\ln{1}&\cnt{1}&\bounds{1}{1}{1}{false}&&\cr
\l{19}&\rat(5/6){5}{}&\rat(5/6){5}{}&\rat(1/3){4}{}&\int(0){\0}{}&\int(1){1}{}&\ln{1}&\cnt{24}&\bounds{6}{6}{1}{false}&&\cr
\l{21}&\rat(5/6){5}{}&\rat(-1/6){5}{}&\rat(4/3){4}{}&\int(0){\0}{}&\int(1){1}{}&\ln{1}&\cnt{15}&\bounds{3}{3}{1}{false}&&\cr
\l{23}&\rat(5/6){5}{}&\int(0){\0}{}&\rat(2/3){2}{}&\int(0){1}{}&\rat(3/2){3}{}&\ln{1}&\cnt{36}&\bounds{6}{6}{1}{false}&&\cr
\l{25}&\rat(5/6){5}{}&\int(0){\0}{}&\rat(2/3){5}{}&\int(0){4}{}&\rat(3/2){2}{}&\ln{1}&\cnt{30}&\bounds{6}{6}{1}{false}&&\cr
\endorbit
\orbit{5}{77}{96}
\h{77}&\rat(5/6){5}{}&\rat(5/6){5}{}&\root\rat(4/3){4}{\spec{1}}&\int(0){\0}{}&\int(3){1}{}&&\cnt{472}&\bounds{96}{96}{1}{false}&&\cr
\l{1}&\rat(2/3){4}{}&\rat(1/6){1}{}&\rat(2/3){5}{}&\int(0){\0}{}&\rat(3/2){3}{}&\ln{2\sd}&\cnt{50}&\bounds{10}{10}{1}{false}&&\cr
\l{2}&\rat(5/6){5}{}&\int(0){\0}{}&\rat(2/3){5}{}&\int(0){4}{}&\rat(3/2){2}{}&\ln{2\sd}&\cnt{30}&\bounds{6}{6}{1}{false}&&\cr
\l{3}&\int(1){\2}{}&\int(0){\0}{}&\int(0){\0}{}&\int(0){\0}{}&\int(2){\2}{}&\ln{2}&\cnt{15}&\bounds{3}{3}{1}{false}&&\cr
\l{5}&\rat(1/2){3}{}&\int(0){\0}{}&\int(1){3}{}&\int(0){\0}{}&\rat(3/2){2}{}&\ln{2}&\cnt{20}&\bounds{4}{4}{1}{false}&&\cr
\l{7}&\rat(2/3){4}{}&\rat(5/6){5}{}&\int(0){\0}{}&\int(0){1}{}&\rat(3/2){2}{}&\ln{2}&\cnt{30}&\bounds{5}{5}{1}{false}&&\cr
\l{9}&\rat(5/6){5}{}&\rat(5/6){5}{}&\rat(4/3){4}{}&\int(0){\0}{}&\int(0){1}{\spec{1}}&\ln{2}&\cnt{1}&\bounds{1}{1}{1}{false}&&\cr
\l{11}&\rat(5/6){5}{}&\rat(5/6){5}{}&\rat(1/3){4}{}&\int(0){\0}{}&\int(1){1}{}&\ln{1}&\cnt{12}&\bounds{4}{4}{1}{false}&&\cr
\l{13}&\rat(5/6){5}{}&\int(0){\0}{}&\rat(2/3){2}{\spec{1}}&\int(0){1}{}&\rat(3/2){3}{}&\ln{2}&\cnt{6}&\bounds{1}{1}{1}{false}&&\cr
\endorbit
\orbit{5}{78}{80}
\h{78}&\rat(5/6){5}{}&\rat(5/6){5}{}&\root\rat(4/3){4}{\spec{2}}&\int(0){\0}{}&\int(3){1}{}&&\cnt{416}&\bounds{80}{80}{1}{false}&&\cr
\l{1}&\rat(5/6){5}{}&\int(0){\0}{}&\rat(2/3){2}{}&\int(0){1}{}&\rat(3/2){3}{}&\ln{2\sd}&\cnt{36}&\bounds{6}{6}{1}{false}&&\cr
\l{2}&\int(1){\2}{}&\int(0){\0}{}&\int(0){\0}{}&\int(0){\0}{}&\int(2){\2}{}&\ln{2}&\cnt{15}&\bounds{3}{3}{1}{false}&&\cr
\l{4}&\rat(1/2){3}{}&\int(0){\0}{}&\int(1){3}{}&\int(0){\0}{}&\rat(3/2){2}{}&\ln{2}&\cnt{40}&\bounds{8}{8}{1}{false}&&\cr
\l{6}&\rat(2/3){4}{}&\rat(5/6){5}{}&\int(0){\0}{}&\int(0){1}{}&\rat(3/2){2}{}&\ln{2}&\cnt{30}&\bounds{5}{5}{1}{false}&&\cr
\l{8}&\rat(5/6){5}{}&\rat(5/6){5}{}&\rat(4/3){4}{}&\int(0){\0}{}&\int(0){1}{\spec{1}}&\ln{2}&\cnt{1}&\bounds{1}{1}{1}{false}&&\cr
\endorbit
\orbit{5}{79}{116}
\h{79}&\rat(5/6){5}{}&\rat(5/6){5}{}&\rat(10/3){4}{\spec{2}}&\root\int(0){\0}{}&\int(1){1}{}&&\cnt{528}&\bounds{116}{116}{1}{false}&&\cr
\l{1}&\rat(5/6){5}{}&\int(0){\0}{}&\rat(5/3){2}{}&\int(0){1}{}&\rat(1/2){3}{}&\ln{2\sd}&\cnt{32}&\bounds{8}{8}{1}{false}&&\cr
\l{2}&\int(1){\2}{}&\int(0){\0}{}&\int(2){\2}{}&\int(0){\0}{}&\int(0){\0}{}&\ln{2}&\cnt{15}&\bounds{3}{3}{1}{false}&&\cr
\l{4}&\rat(1/2){3}{}&\int(0){\0}{}&\int(2){3}{}&\int(0){\0}{}&\rat(1/2){2}{}&\ln{2}&\cnt{40}&\bounds{8}{8}{1}{false}&&\cr
\l{6}&\rat(2/3){4}{}&\rat(2/3){4}{}&\rat(5/3){2}{}&\int(0){\0}{}&\int(0){\0}{}&\ln{1}&\cnt{50}&\bounds{10}{10}{1}{false}&&\cr
\l{8}&\rat(2/3){4}{}&\int(0){\0}{}&\rat(4/3){1}{}&\int(0){5}{}&\int(1){1}{}&\ln{2}&\cnt{20}&\bounds{4}{4}{1}{false}&&\cr
\l{10}&\rat(5/6){5}{}&\rat(5/6){5}{}&\rat(4/3){4}{}&\int(0){\0}{}&\int(0){1}{}&\ln{1}&\cnt{18}&\bounds{6}{6}{1}{false}&&\cr
\l{12}&\rat(5/6){5}{}&\rat(5/6){5}{}&\rat(1/3){4}{}&\int(0){\0}{}&\int(1){1}{}&\ln{1}&\cnt{6}&\bounds{2}{2}{1}{false}&&\cr
\l{14}&\rat(5/6){5}{}&\rat(5/6){5}{}&\rat(4/3){1}{}&\int(0){3}{\spec{1}}&\int(0){\0}{}&\ln{2}&\cnt{4}&\bounds{1}{1}{1}{false}&&\cr
\endorbit
\orbit{5}{80}{106}
\h{80}&\rat(5/6){5}{}&\rat(5/6){5}{}&\root\rat(10/3){4}{\spec{2}}&\int(0){\0}{}&\int(1){1}{}&&\cnt{552}&\bounds{106}{106}{1}{false}&&\cr
\l{1}&\rat(5/6){5}{}&\int(0){\0}{}&\rat(5/3){2}{}&\int(0){1}{}&\rat(1/2){3}{}&\ln{2\sd}&\cnt{48}&\bounds{8}{8}{1}{false}&&\cr
\l{2}&\int(1){\2}{}&\int(0){\0}{}&\int(2){\2}{}&\int(0){\0}{}&\int(0){\0}{}&\ln{2}&\cnt{5}&\bounds{1}{1}{1}{false}&&\cr
\l{4}&\rat(1/2){3}{}&\int(0){\0}{}&\int(2){3}{}&\int(0){\0}{}&\rat(1/2){2}{}&\ln{2}&\cnt{40}&\bounds{8}{8}{1}{false}&&\cr
\l{6}&\rat(2/3){4}{}&\rat(2/3){4}{}&\rat(5/3){2}{}&\int(0){\0}{}&\int(0){\0}{}&\ln{1}&\cnt{50}&\bounds{10}{10}{1}{false}&&\cr
\l{8}&\rat(2/3){4}{}&\int(0){\0}{}&\rat(4/3){1}{}&\int(0){5}{}&\int(1){1}{}&\ln{2}&\cnt{30}&\bounds{5}{5}{1}{false}&&\cr
\l{10}&\rat(5/6){5}{}&\rat(5/6){5}{}&\rat(4/3){4}{}&\int(0){\0}{}&\int(0){1}{}&\ln{1}&\cnt{6}&\bounds{2}{2}{1}{false}&&\cr
\l{12}&\rat(5/6){5}{}&\rat(5/6){5}{}&\rat(1/3){4}{}&\int(0){\0}{}&\int(1){1}{}&\ln{1}&\cnt{2}&\bounds{1}{1}{1}{false}&&\cr
\l{14}&\rat(5/6){5}{}&\rat(5/6){5}{}&\rat(4/3){1}{}&\int(0){3}{}&\int(0){\0}{}&\ln{1}&\cnt{20}&\bounds{4}{4}{1}{false}&&\cr
\endorbit
\orbit{5}{81}{104}
\h{81}&\root\rat(5/6){5}{}&\rat(5/6){5}{}&\rat(10/3){4}{\spec{2}}&\int(0){\0}{}&\int(1){1}{}&&\cnt{500}&\bounds{104}{104}{1}{false}&&\cr
\l{1}&\int(1){\2}{}&\int(0){\0}{}&\int(2){\2}{}&\int(0){\0}{}&\int(0){\0}{}&\ln{1}&\cnt{9}&\bounds{3}{3}{1}{false}&&\cr
\l{3}&\rat(1/2){3}{\spec{1}}&\int(0){\0}{}&\int(2){3}{}&\int(0){\0}{}&\rat(1/2){2}{}&\ln{1}&\cnt{12}&\bounds{3}{3}{1}{false}&&\cr
\l{5}&\rat(1/2){3}{\spec{2}}&\int(0){\0}{}&\int(2){3}{}&\int(0){\0}{}&\rat(1/2){2}{}&\ln{1}&\cnt{4}&\bounds{1}{1}{1}{false}&&\cr
\l{7}&\rat(2/3){4}{}&\rat(2/3){4}{}&\rat(5/3){2}{}&\int(0){\0}{}&\int(0){\0}{}&\ln{1}&\cnt{30}&\bounds{6}{6}{1}{false}&&\cr
\l{9}&\rat(2/3){4}{}&\int(0){\0}{}&\rat(4/3){1}{}&\int(0){5}{}&\int(1){1}{}&\ln{1}&\cnt{18}&\bounds{3}{3}{1}{false}&&\cr
\l{11}&\rat(5/6){5}{}&\rat(1/3){2}{}&\rat(4/3){1}{}&\int(0){\0}{}&\rat(1/2){2}{}&\ln{1}&\cnt{40}&\bounds{8}{8}{1}{false}&&\cr
\l{13}&\rat(5/6){5}{}&\rat(5/6){5}{}&\rat(4/3){4}{}&\int(0){\0}{}&\int(0){1}{}&\ln{1}&\cnt{18}&\bounds{6}{6}{1}{false}&&\cr
\l{15}&\rat(5/6){5}{}&\rat(5/6){5}{}&\rat(1/3){4}{}&\int(0){\0}{}&\int(1){1}{}&\ln{1}&\cnt{6}&\bounds{2}{2}{1}{false}&&\cr
\l{17}&\rat(5/6){5}{}&\rat(5/6){5}{}&\rat(4/3){1}{}&\int(0){3}{}&\int(0){\0}{}&\ln{1}&\cnt{20}&\bounds{4}{4}{1}{false}&&\cr
\l{19}&\rat(5/6){5}{}&\rat(1/6){1}{}&\int(2){3}{}&\int(0){5}{}&\int(0){\0}{}&\ln{1}&\cnt{30}&\bounds{5}{5}{1}{false}&&\cr
\l{21}&\rat(5/6){5}{}&\rat(-1/6){5}{}&\rat(4/3){4}{}&\int(0){\0}{}&\int(1){1}{}&\ln{1}&\cnt{15}&\bounds{3}{3}{1}{false}&&\cr
\l{23}&\rat(5/6){5}{}&\int(0){\0}{}&\rat(5/3){2}{}&\int(0){1}{}&\rat(1/2){3}{}&\ln{1}&\cnt{48}&\bounds{8}{8}{1}{false}&&\cr
\endorbit
\orbit{5}{82}{92}
\h{82}&\rat(5/6){5}{}&\rat(5/6){5}{}&\rat(10/3){4}{\spec{2}}&\int(0){\0}{}&\root\int(1){1}{}&&\cnt{472}&\bounds{92}{92}{1}{false}&&\cr
\l{1}&\rat(5/6){5}{}&\int(0){\0}{}&\rat(5/3){2}{}&\int(0){1}{}&\rat(1/2){3}{}&\ln{2\sd}&\cnt{24}&\bounds{4}{4}{1}{false}&&\cr
\l{2}&\int(1){\2}{}&\int(0){\0}{}&\int(2){\2}{}&\int(0){\0}{}&\int(0){\0}{}&\ln{2}&\cnt{15}&\bounds{3}{3}{1}{false}&&\cr
\l{4}&\rat(1/2){3}{}&\int(0){\0}{}&\int(2){3}{}&\int(0){\0}{}&\rat(1/2){2}{}&\ln{2}&\cnt{20}&\bounds{4}{4}{1}{false}&&\cr
\l{6}&\rat(2/3){4}{}&\rat(2/3){4}{}&\rat(5/3){2}{}&\int(0){\0}{}&\int(0){\0}{}&\ln{1}&\cnt{50}&\bounds{10}{10}{1}{false}&&\cr
\l{8}&\rat(2/3){4}{}&\int(0){\0}{}&\rat(4/3){1}{}&\int(0){5}{}&\int(1){1}{}&\ln{2}&\cnt{30}&\bounds{5}{5}{1}{false}&&\cr
\l{10}&\rat(5/6){5}{}&\rat(5/6){5}{}&\rat(4/3){4}{}&\int(0){\0}{}&\int(0){1}{\spec{1}}&\ln{2}&\cnt{3}&\bounds{1}{1}{1}{false}&&\cr
\l{12}&\rat(5/6){5}{}&\rat(5/6){5}{}&\rat(1/3){4}{}&\int(0){\0}{}&\int(1){1}{}&\ln{1}&\cnt{6}&\bounds{2}{2}{1}{false}&&\cr
\l{14}&\rat(5/6){5}{}&\rat(5/6){5}{}&\rat(4/3){1}{}&\int(0){3}{}&\int(0){\0}{}&\ln{1}&\cnt{20}&\bounds{4}{4}{1}{false}&&\cr
\endorbit
\orbit{5}{83}{84}
\h{83}&\rat(5/6){5}{}&\rat(5/6){5}{}&\root\rat(10/3){4}{\spec{1}}&\int(0){\0}{}&\int(1){1}{}&&\cnt{416}&\bounds{84}{84}{1}{false}&&\cr
\l{1}&\int(1){\2}{}&\int(0){\0}{}&\int(2){\2}{}&\int(0){\0}{}&\int(0){\0}{}&\ln{2}&\cnt{15}&\bounds{3}{3}{1}{false}&&\cr
\l{3}&\rat(1/2){3}{}&\int(0){\0}{}&\int(2){3}{}&\int(0){\0}{}&\rat(1/2){2}{}&\ln{2}&\cnt{40}&\bounds{8}{8}{1}{false}&&\cr
\l{5}&\rat(2/3){4}{}&\int(0){\0}{}&\rat(4/3){1}{}&\int(0){5}{}&\int(1){1}{}&\ln{2}&\cnt{30}&\bounds{5}{5}{1}{false}&&\cr
\l{7}&\rat(5/6){5}{}&\rat(5/6){5}{}&\rat(4/3){4}{}&\int(0){\0}{}&\int(0){1}{}&\ln{1}&\cnt{18}&\bounds{6}{6}{1}{false}&&\cr
\l{9}&\rat(5/6){5}{}&\rat(5/6){5}{}&\rat(4/3){1}{}&\int(0){3}{}&\int(0){\0}{}&\ln{1}&\cnt{20}&\bounds{4}{4}{1}{false}&&\cr
\endorbit
\orbit{5}{84}{108}
\h{84}&\root\rat(5/6){5}{}&\rat(5/6){5}{}&\rat(4/3){4}{}&\int(2){\2}{}&\int(1){1}{}&&\cnt{452}&\bounds{108}{108}{1}{false}&&\cr
\l{1}&\int(1){\2}{}&\int(0){\0}{}&\int(0){\0}{}&\int(2){\2}{}&\int(0){\0}{}&\ln{1}&\cnt{3}&\bounds{1}{1}{1}{false}&&\cr
\l{3}&\rat(1/2){3}{\spec{1}}&\rat(5/6){5}{}&\rat(2/3){5}{}&\int(1){5}{}&\int(0){\0}{}&\ln{1}&\cnt{6}&\bounds{2}{2}{1}{false}&&\cr
\l{5}&\rat(1/2){3}{\spec{2}}&\rat(5/6){5}{}&\rat(2/3){5}{}&\int(1){5}{}&\int(0){\0}{}&\ln{1}&\cnt{2}&\bounds{1}{1}{1}{false}&&\cr
\l{7}&\rat(2/3){4}{}&\rat(5/6){5}{}&\int(0){\0}{}&\int(1){1}{}&\rat(1/2){2}{}&\ln{1}&\cnt{12}&\bounds{3}{3}{1}{false}&&\cr
\l{9}&\rat(2/3){4}{}&\int(0){\0}{}&\rat(4/3){4}{}&\int(1){2}{}&\int(0){\0}{}&\ln{1}&\cnt{12}&\bounds{3}{3}{1}{false}&&\cr
\l{11}&\rat(2/3){4}{}&\int(0){\0}{}&\rat(1/3){1}{}&\int(1){5}{}&\int(1){1}{}&\ln{1}&\cnt{12}&\bounds{3}{3}{1}{false}&&\cr
\l{13}&\rat(5/6){5}{}&\rat(1/2){3}{}&\rat(2/3){5}{}&\int(1){1}{}&\int(0){\0}{}&\ln{1}&\cnt{20}&\bounds{4}{4}{1}{false}&&\cr
\l{15}&\rat(5/6){5}{}&\rat(2/3){4}{}&\int(0){\0}{}&\int(1){5}{}&\rat(1/2){3}{}&\ln{1}&\cnt{20}&\bounds{4}{4}{1}{false}&&\cr
\l{17}&\rat(5/6){5}{}&\rat(5/6){5}{}&\rat(4/3){4}{}&\int(0){\0}{}&\int(0){1}{}&\ln{1}&\cnt{6}&\bounds{2}{2}{1}{false}&&\cr
\l{19}&\rat(5/6){5}{}&\rat(5/6){5}{}&\rat(1/3){4}{}&\int(0){\0}{}&\int(1){1}{}&\ln{1}&\cnt{8}&\bounds{2}{2}{1}{false}&&\cr
\l{21}&\rat(5/6){5}{}&\rat(5/6){5}{}&\rat(1/3){1}{}&\int(1){3}{}&\int(0){\0}{}&\ln{1}&\cnt{24}&\bounds{6}{6}{1}{false}&&\cr
\l{23}&\rat(5/6){5}{}&\rat(1/6){1}{}&\int(1){3}{}&\int(1){5}{}&\int(0){\0}{}&\ln{1}&\cnt{20}&\bounds{4}{4}{1}{false}&&\cr
\l{25}&\rat(5/6){5}{}&\rat(1/6){1}{}&\int(0){\0}{}&\int(1){2}{}&\int(1){1}{}&\ln{1}&\cnt{20}&\bounds{4}{4}{1}{false}&&\cr
\l{27}&\rat(5/6){5}{}&\rat(-1/6){5}{}&\rat(4/3){4}{}&\int(0){\0}{}&\int(1){1}{}&\ln{1}&\cnt{5}&\bounds{1}{1}{1}{false}&&\cr
\l{29}&\rat(5/6){5}{}&\int(0){\0}{}&\rat(2/3){2}{}&\int(1){1}{}&\rat(1/2){3}{}&\ln{1}&\cnt{24}&\bounds{6}{6}{1}{false}&&\cr
\l{31}&\rat(5/6){5}{}&\int(0){\0}{}&\rat(2/3){5}{}&\int(1){4}{}&\rat(1/2){2}{}&\ln{1}&\cnt{32}&\bounds{8}{8}{1}{false}&&\cr
\endorbit
\orbit{5}{85}{96}
\h{85}&\rat(5/6){5}{}&\rat(5/6){5}{}&\root\rat(4/3){4}{\spec{1}}&\int(2){\2}{}&\int(1){1}{}&&\cnt{424}&\bounds{96}{96}{1}{false}&&\cr
\l{1}&\rat(5/6){5}{}&\int(0){\0}{}&\rat(2/3){5}{}&\int(1){4}{}&\rat(1/2){2}{}&\ln{2\sd}&\cnt{32}&\bounds{8}{8}{1}{false}&&\cr
\l{2}&\int(1){\2}{}&\int(0){\0}{}&\int(0){\0}{}&\int(2){\2}{}&\int(0){\0}{}&\ln{2}&\cnt{5}&\bounds{1}{1}{1}{false}&&\cr
\l{4}&\rat(1/2){3}{}&\rat(5/6){5}{}&\rat(2/3){5}{}&\int(1){5}{}&\int(0){\0}{}&\ln{2}&\cnt{20}&\bounds{4}{4}{1}{false}&&\cr
\l{6}&\rat(2/3){4}{}&\rat(5/6){5}{}&\int(0){\0}{}&\int(1){1}{}&\rat(1/2){2}{}&\ln{2}&\cnt{20}&\bounds{4}{4}{1}{false}&&\cr
\l{8}&\rat(2/3){4}{}&\int(0){\0}{}&\rat(4/3){4}{}&\int(1){2}{}&\int(0){\0}{}&\ln{2}&\cnt{20}&\bounds{4}{4}{1}{false}&&\cr
\l{10}&\rat(2/3){4}{}&\int(0){\0}{}&\rat(1/3){1}{}&\int(1){5}{}&\int(1){1}{}&\ln{2}&\cnt{10}&\bounds{2}{2}{1}{false}&&\cr
\l{12}&\rat(5/6){5}{}&\rat(5/6){5}{}&\rat(4/3){4}{}&\int(0){\0}{}&\int(0){1}{}&\ln{1}&\cnt{6}&\bounds{2}{2}{1}{false}&&\cr
\l{14}&\rat(5/6){5}{}&\rat(5/6){5}{}&\rat(1/3){4}{}&\int(0){\0}{}&\int(1){1}{}&\ln{1}&\cnt{4}&\bounds{2}{2}{1}{false}&&\cr
\l{16}&\rat(5/6){5}{}&\rat(5/6){5}{}&\rat(1/3){1}{}&\int(1){3}{}&\int(0){\0}{}&\ln{1}&\cnt{12}&\bounds{4}{4}{1}{false}&&\cr
\l{18}&\rat(5/6){5}{}&\int(0){\0}{}&\rat(2/3){2}{\spec{1}}&\int(1){1}{}&\rat(1/2){3}{}&\ln{2}&\cnt{4}&\bounds{1}{1}{1}{false}&&\cr
\endorbit
\orbit{5}{86}{96}
\h{86}&\rat(5/6){5}{}&\rat(5/6){5}{}&\rat(4/3){4}{}&\int(2){\2}{}&\root\int(1){1}{}&&\cnt{424}&\bounds{96}{96}{1}{false}&&\cr
\l{1}&\rat(5/6){5}{}&\int(0){\0}{}&\rat(2/3){2}{}&\int(1){1}{}&\rat(1/2){3}{}&\ln{2\sd}&\cnt{12}&\bounds{4}{4}{1}{false}&&\cr
\l{2}&\rat(5/6){5}{}&\int(0){\0}{}&\rat(2/3){5}{}&\int(1){4}{}&\rat(1/2){2}{}&\ln{2\sd}&\cnt{16}&\bounds{4}{4}{1}{false}&&\cr
\l{3}&\int(1){\2}{}&\int(0){\0}{}&\int(0){\0}{}&\int(2){\2}{}&\int(0){\0}{}&\ln{2}&\cnt{5}&\bounds{1}{1}{1}{false}&&\cr
\l{5}&\rat(1/2){3}{}&\rat(5/6){5}{}&\rat(2/3){5}{}&\int(1){5}{}&\int(0){\0}{}&\ln{2}&\cnt{20}&\bounds{4}{4}{1}{false}&&\cr
\l{7}&\rat(2/3){4}{}&\rat(5/6){5}{}&\int(0){\0}{}&\int(1){1}{}&\rat(1/2){2}{}&\ln{2}&\cnt{10}&\bounds{2}{2}{1}{false}&&\cr
\l{9}&\rat(2/3){4}{}&\int(0){\0}{}&\rat(4/3){4}{}&\int(1){2}{}&\int(0){\0}{}&\ln{2}&\cnt{20}&\bounds{4}{4}{1}{false}&&\cr
\l{11}&\rat(2/3){4}{}&\int(0){\0}{}&\rat(1/3){1}{}&\int(1){5}{}&\int(1){1}{}&\ln{2}&\cnt{20}&\bounds{4}{4}{1}{false}&&\cr
\l{13}&\rat(5/6){5}{}&\rat(5/6){5}{}&\rat(4/3){4}{}&\int(0){\0}{}&\int(0){1}{\spec{1}}&\ln{2}&\cnt{1}&\bounds{1}{1}{1}{false}&&\cr
\l{15}&\rat(5/6){5}{}&\rat(5/6){5}{}&\rat(1/3){4}{}&\int(0){\0}{}&\int(1){1}{}&\ln{1}&\cnt{8}&\bounds{2}{2}{1}{false}&&\cr
\l{17}&\rat(5/6){5}{}&\rat(5/6){5}{}&\rat(1/3){1}{}&\int(1){3}{}&\int(0){\0}{}&\ln{1}&\cnt{24}&\bounds{6}{6}{1}{false}&&\cr
\endorbit
\orbit{5}{87}{92}
\h{87}&\rat(5/6){5}{}&\rat(5/6){5}{}&\rat(4/3){4}{}&\root\int(2){\2}{}&\int(1){1}{}&&\cnt{424}&\bounds{92}{92}{1}{false}&&\cr
\l{1}&\rat(5/6){5}{}&\int(0){\0}{}&\rat(2/3){2}{}&\int(1){1}{}&\rat(1/2){3}{}&\ln{2\sd}&\cnt{24}&\bounds{6}{6}{1}{false}&&\cr
\l{2}&\rat(5/6){5}{}&\int(0){\0}{}&\rat(2/3){5}{}&\int(1){4}{}&\rat(1/2){2}{}&\ln{2\sd}&\cnt{16}&\bounds{4}{4}{1}{false}&&\cr
\l{3}&\int(1){\2}{}&\int(0){\0}{}&\int(0){\0}{}&\int(2){\2}{}&\int(0){\0}{}&\ln{2}&\cnt{5}&\bounds{1}{1}{1}{false}&&\cr
\l{5}&\rat(1/2){3}{}&\rat(5/6){5}{}&\rat(2/3){5}{}&\int(1){5}{}&\int(0){\0}{}&\ln{2}&\cnt{20}&\bounds{4}{4}{1}{false}&&\cr
\l{7}&\rat(2/3){4}{}&\rat(5/6){5}{}&\int(0){\0}{}&\int(1){1}{}&\rat(1/2){2}{}&\ln{2}&\cnt{20}&\bounds{4}{4}{1}{false}&&\cr
\l{9}&\rat(2/3){4}{}&\int(0){\0}{}&\rat(4/3){4}{}&\int(1){2}{}&\int(0){\0}{}&\ln{2}&\cnt{10}&\bounds{2}{2}{1}{false}&&\cr
\l{11}&\rat(2/3){4}{}&\int(0){\0}{}&\rat(1/3){1}{}&\int(1){5}{}&\int(1){1}{}&\ln{2}&\cnt{20}&\bounds{4}{4}{1}{false}&&\cr
\l{13}&\rat(5/6){5}{}&\rat(5/6){5}{}&\rat(4/3){4}{}&\int(0){\0}{}&\int(0){1}{}&\ln{1}&\cnt{6}&\bounds{2}{2}{1}{false}&&\cr
\l{15}&\rat(5/6){5}{}&\rat(5/6){5}{}&\rat(1/3){4}{}&\int(0){\0}{}&\int(1){1}{}&\ln{1}&\cnt{8}&\bounds{2}{2}{1}{false}&&\cr
\l{17}&\rat(5/6){5}{}&\rat(5/6){5}{}&\rat(1/3){1}{}&\int(1){3}{\spec{1}}&\int(0){\0}{}&\ln{2}&\cnt{4}&\bounds{1}{1}{1}{false}&&\cr
\endorbit
\orbit{5}{88}{80}
\h{88}&\rat(5/6){5}{}&\rat(5/6){5}{}&\root\rat(4/3){4}{\spec{2}}&\int(2){\2}{}&\int(1){1}{}&&\cnt{368}&\bounds{80}{80}{1}{false}&&\cr
\l{1}&\rat(5/6){5}{}&\int(0){\0}{}&\rat(2/3){2}{}&\int(1){1}{}&\rat(1/2){3}{}&\ln{2\sd}&\cnt{24}&\bounds{6}{6}{1}{false}&&\cr
\l{2}&\int(1){\2}{}&\int(0){\0}{}&\int(0){\0}{}&\int(2){\2}{}&\int(0){\0}{}&\ln{2}&\cnt{5}&\bounds{1}{1}{1}{false}&&\cr
\l{4}&\rat(2/3){4}{}&\rat(5/6){5}{}&\int(0){\0}{}&\int(1){1}{}&\rat(1/2){2}{}&\ln{2}&\cnt{20}&\bounds{4}{4}{1}{false}&&\cr
\l{6}&\rat(2/3){4}{}&\int(0){\0}{}&\rat(4/3){4}{}&\int(1){2}{}&\int(0){\0}{}&\ln{2}&\cnt{20}&\bounds{4}{4}{1}{false}&&\cr
\l{8}&\rat(2/3){4}{}&\int(0){\0}{}&\rat(1/3){1}{}&\int(1){5}{}&\int(1){1}{}&\ln{2}&\cnt{20}&\bounds{4}{4}{1}{false}&&\cr
\l{10}&\rat(5/6){5}{}&\rat(5/6){5}{}&\rat(4/3){4}{}&\int(0){\0}{}&\int(0){1}{}&\ln{1}&\cnt{6}&\bounds{2}{2}{1}{false}&&\cr
\l{12}&\rat(5/6){5}{}&\rat(5/6){5}{}&\rat(1/3){1}{}&\int(1){3}{}&\int(0){\0}{}&\ln{1}&\cnt{24}&\bounds{6}{6}{1}{false}&&\cr
\endorbit
\orbit{5}{89}{106}
\h{89}&\root\rat(5/6){5}{}&\rat(5/6){5}{}&\rat(5/6){1}{}&\rat(3/2){3}{}&\int(2){\2}{}&&\cnt{452}&\bounds{106}{106}{1}{false}&&\cr
\l{1}&\int(1){\2}{}&\int(0){\0}{}&\int(0){\0}{}&\int(0){\0}{}&\int(2){\2}{}&\ln{1}&\cnt{3}&\bounds{1}{1}{1}{false}&&\cr
\l{3}&\rat(1/2){3}{\spec{1}}&\int(0){\0}{}&\int(0){\0}{}&\rat(3/2){3}{}&\int(1){3}{}&\ln{1}&\cnt{6}&\bounds{2}{2}{1}{false}&&\cr
\l{5}&\rat(1/2){3}{\spec{2}}&\int(0){\0}{}&\int(0){\0}{}&\rat(3/2){3}{}&\int(1){3}{}&\ln{1}&\cnt{2}&\bounds{1}{1}{1}{false}&&\cr
\l{7}&\rat(2/3){4}{}&\rat(5/6){5}{}&\int(0){\0}{}&\rat(1/2){1}{}&\int(1){2}{}&\ln{2}&\cnt{18}&\bounds{6}{6}{1}{false}&&\cr
\l{9}&\rat(5/6){5}{}&\rat(2/3){4}{}&\int(0){\0}{}&\rat(1/2){5}{}&\int(1){3}{}&\ln{2}&\cnt{30}&\bounds{6}{6}{1}{false}&&\cr
\l{11}&\rat(5/6){5}{}&\rat(1/3){2}{}&\rat(5/6){1}{}&\int(0){\0}{}&\int(1){2}{}&\ln{2}&\cnt{20}&\bounds{4}{4}{1}{false}&&\cr
\l{13}&\rat(5/6){5}{}&\rat(5/6){5}{}&\rat(5/6){1}{}&\rat(1/2){3}{}&\int(0){\0}{}&\ln{1}&\cnt{9}&\bounds{3}{3}{1}{false}&&\cr
\l{15}&\rat(5/6){5}{}&\rat(5/6){5}{}&\rat(-1/6){1}{}&\rat(3/2){3}{}&\int(0){\0}{}&\ln{2}&\cnt{5}&\bounds{1}{1}{1}{false}&&\cr
\l{17}&\rat(5/6){5}{}&\rat(1/6){1}{}&\int(0){\0}{}&\int(1){2}{}&\int(1){1}{}&\ln{2}&\cnt{30}&\bounds{6}{6}{1}{false}&&\cr
\endorbit
\orbit{5}{90}{80}
\h{90}&\rat(5/6){5}{}&\rat(5/6){5}{}&\rat(5/6){1}{}&\root\rat(3/2){3}{\spec{1}}&\int(2){\2}{}&&\cnt{396}&\bounds{80}{80}{1}{false}&&\cr
\l{1}&\int(1){\2}{}&\int(0){\0}{}&\int(0){\0}{}&\int(0){\0}{}&\int(2){\2}{}&\ln{3}&\cnt{5}&\bounds{1}{1}{1}{false}&&\cr
\l{3}&\rat(1/2){3}{}&\int(0){\0}{}&\int(0){\0}{}&\rat(3/2){3}{}&\int(1){3}{}&\ln{3}&\cnt{20}&\bounds{4}{4}{1}{false}&&\cr
\l{5}&\rat(2/3){4}{}&\rat(5/6){5}{}&\int(0){\0}{}&\rat(1/2){1}{}&\int(1){2}{}&\ln{3}&\cnt{10}&\bounds{2}{2}{1}{false}&&\cr
\l{7}&\rat(2/3){4}{}&\int(0){\0}{}&\rat(5/6){1}{}&\rat(1/2){5}{}&\int(1){1}{}&\ln{3}&\cnt{30}&\bounds{6}{6}{1}{false}&&\cr
\l{9}&\rat(5/6){5}{}&\rat(5/6){5}{}&\rat(5/6){1}{}&\rat(1/2){3}{}&\int(0){\0}{}&\ln{1}&\cnt{3}&\bounds{1}{1}{1}{false}&&\cr
\endorbit
\orbit{5}{91}{76}
\h{91}&\rat(5/6){5}{}&\rat(5/6){5}{}&\rat(5/6){1}{}&\rat(3/2){3}{}&\root\int(2){\2}{\spec{1}}&&\cnt{368}&\bounds{76}{76}{1}{false}&&\cr
\l{1}&\int(1){\2}{}&\int(0){\0}{}&\int(0){\0}{}&\int(0){\0}{}&\int(2){\2}{}&\ln{1}&\cnt{5}&\bounds{1}{1}{1}{false}&&\cr
\l{3}&\rat(2/3){4}{}&\rat(5/6){5}{}&\int(0){\0}{}&\rat(1/2){1}{}&\int(1){2}{}&\ln{2}&\cnt{30}&\bounds{6}{6}{1}{false}&&\cr
\l{5}&\rat(5/6){5}{}&\rat(1/3){2}{}&\rat(5/6){1}{}&\int(0){\0}{}&\int(1){2}{}&\ln{2}&\cnt{20}&\bounds{4}{4}{1}{false}&&\cr
\l{7}&\rat(5/6){5}{}&\rat(5/6){5}{}&\rat(5/6){1}{}&\rat(1/2){3}{}&\int(0){\0}{}&\ln{1}&\cnt{9}&\bounds{3}{3}{1}{false}&&\cr
\l{9}&\rat(5/6){5}{}&\rat(5/6){5}{}&\rat(-1/6){1}{}&\rat(3/2){3}{}&\int(0){\0}{}&\ln{2}&\cnt{5}&\bounds{1}{1}{1}{false}&&\cr
\l{11}&\rat(5/6){5}{}&\rat(1/6){1}{}&\int(0){\0}{}&\int(1){2}{}&\int(1){1}{}&\ln{2}&\cnt{30}&\bounds{6}{6}{1}{false}&&\cr
\endorbit
\orbit{5}{92}{52}
\h{92}&\int(6){\6}{\spec{1}}&\int(0){\0}{}&\int(0){\0}{}&\int(0){\0}{}&\int(0){\0}{}&&\cnt{240}&\bounds{52}{52}{1}{false}&&\cr
\l{1}&\int(3){\4}{}&\int(0){\0}{}&\int(0){\0}{}&\int(0){\0}{}&\int(0){\0}{}&\ln{1\sdd}&\cnt{12}&\bounds{6}{6}{1}{false}&&\cr
\l{2}&\int(3){\2}{}&\int(0){\2}{}&\int(0){\0}{}&\int(0){\0}{}&\int(0){\0}{}&\ln{3\sdd}&\cnt{60}&\bounds{12}{12}{1}{false}&&\cr
\l{3}&\int(3){\2}{}&\int(0){\0}{}&\int(0){\0}{}&\int(0){\0}{}&\int(0){\2}{}&\ln{1\sdd}&\cnt{48}&\bounds{10}{10}{1}{false}&&\cr
\endorbit
\orbit{5}{93}{48}
\h{93}&\int(0){\0}{}&\int(0){\0}{}&\int(0){\0}{}&\int(0){\0}{}&\int(6){\6}{\spec{1}}&&\cnt{240}&\bounds{48}{48}{1}{false}&&\cr
\l{1}&\int(0){\2}{}&\int(0){\0}{}&\int(0){\0}{}&\int(0){\0}{}&\int(3){\2}{}&\ln{4\sdd}&\cnt{60}&\bounds{12}{12}{1}{false}&&\cr
\endorbit
\endorbits

\subsubsection{Configuration \hlink{4A5_D4}{1}}\label{s.4A5_D4-1}
There are two sets $\fL\in\Bnd_{17}(\oorb2)$, defined by the patterns~$\pat$
such that
$\pat|_{\oorb2}=\const=7$ or $9$.
The former has
rank $19$, and its only nontrivial
geometric $\oorb2$-proper (see \autoref{ss.ext.max})
extension has $\pat(\orb)=\pat(\orb^*)=2$
for a pair of dual orbits $\orb,\orb^*\subset\oorb1$ and $\pat(\orb')=0$ for
all other orbits $\orb'\subset\oorb1$.
The other set, which is denoted $\Lmax1$, is maximal, see
\autoref{ss.maximal}.
This set of size $144$ is determined by the pattern
(see \autoref{rem.pattern.unique})
\[
\pat=\pattern{0,9,*}\quad\mbox{(not $\Q$-complete)}.
\maxlabel1
\]


\subsubsection{Configuration \hlink{4A5_D4}{3}}\label{s.4A5_D4-3}
It is not practical to compute the admissible sets for all orbits;
thus, we argue as in \autoref{ss.brute.force.blocks} and
only compute admissible
subsets $\fL\subset\orb\subset\oorb1$ of size at least~$18$.
This suffices to show that there is a unique set $\fL\in\Bnd_4(\oorb1)$,
with the pattern $\pat$ taking values $(22,22,18)$ on $\oorb1$ and
identical~$0$ on~$\oorb2$.
This set is maximal (see \autoref{ss.maximal}).


\newcs{max-6D4}{120}%
\newcs{count-6D4}{13}%
\newcs{counts-6D4}{\DO{1}{108}\DO{2}{112}\DO{3}{114}\DO{4}{108}\DO{5}{90}\DO{6}{118}\DO{7}{112}\DO{8}{96}\DO{9}{102}\DO{10}{88}\DO{11}{80}\DO{12}{120}\DO{13}{50}}%
\orbits
\orbit{6}{1}{108}
\h{1}&\int(2){\2}{}&\int(2){\2}{}&\int(2){\2}{}&\root\int(0){\0}{}&\int(0){\0}{}&\int(0){\0}{}&&\cnt{528}&\bounds{108}{168}{fail}{false}&&\cr
\l{1}&\int(2){\2}{}&\int(1){\2}{}&\int(0){\0}{}&\int(0){\0}{}&\int(0){\0}{}&\int(0){\0}{}&\ln{6\sd}&\cnt{8}&\bounds{3}{4}{5}{true}&&\cr
\l{2}&\int(1){3}{}&\int(1){3}{}&\int(1){3}{}&\int(0){3}{}&\int(0){\0}{}&\int(0){\0}{}&\ln{3\sdd}&\cnt{32}&\bounds{10}{16}{2}{true}&&\cr
\l{3}&\int(1){3}{}&\int(1){2}{}&\int(1){1}{}&\int(0){\0}{}&\int(0){1}{}&\int(0){\0}{}&\ln{6\sdd}&\cnt{64}&\bounds{10}{16}{5}{true}&&\cr
\endorbit
\orbit{6}{2}{112}
\h{2}&\root\int(2){\2}{\spec{1}}&\int(2){\2}{}&\int(2){\2}{}&\int(0){\0}{}&\int(0){\0}{}&\int(0){\0}{}&&\cnt{416}&\bounds{112}{112}{1}{false}&&\cr
\l{1}&\int(1){2}{}&\int(1){3}{}&\int(1){1}{}&\int(0){\0}{}&\int(0){\0}{}&\int(0){1}{}&\ln{6\sdd}&\cnt{64}&\bounds{16}{16}{1}{false}&&\cr
\l{2}&\int(2){\2}{}&\int(1){\2}{}&\int(0){\0}{}&\int(0){\0}{}&\int(0){\0}{}&\int(0){\0}{}&\ln{2}&\cnt{8}&\bounds{4}{4}{1}{false}&&\cr
\endorbit
\orbit{6}{3}{114}
\h{3}&\int(3){3}{}&\int(1){3}{}&\int(1){3}{}&\int(1){3}{}&\root\int(0){\0}{}&\int(0){\0}{}&&\cnt{528}&\bounds{114}{144}{fail}{false}&&\cr
\l{1}&\rat(3/2){2}{}&\int(1){3}{}&\rat(1/2){1}{}&\int(0){\0}{}&\int(0){\0}{}&\int(0){1}{}&\ln{6\sd}&\cnt{32}&\bounds{5}{8}{5}{true}&&\cr
\l{2}&\rat(3/2){2}{}&\int(1){3}{}&\int(0){\0}{}&\rat(1/2){1}{}&\int(0){1}{}&\int(0){\0}{}&\ln{6\sd}&\cnt{16}&\bounds{4}{4}{1}{true}&&\cr
\l{3}&\rat(3/2){2}{}&\rat(1/2){2}{}&\rat(1/2){2}{}&\rat(1/2){2}{}&\int(0){\0}{}&\int(0){\0}{}&\ln{2\sd}&\cnt{64}&\bounds{16}{16}{1}{false}&&\cr
\l{4}&\int(3){\4}{2}&\int(0){\0}{}&\int(0){\0}{}&\int(0){\0}{}&\int(0){\0}{}&\int(0){\0}{}&\ln{2}&\cnt{1}&\bounds{1}{1}{1}{true}&&\cr
\l{6}&\int(2){\2}{}&\int(1){\2}{}&\int(0){\0}{}&\int(0){\0}{}&\int(0){\0}{}&\int(0){\0}{}&\ln{3}&\cnt{18}&\bounds{4}{6}{5}{true}&&\cr
\endorbit
\orbit{6}{4}{108}
\h{4}&\root\int(3){3}{}&\int(1){3}{}&\int(1){3}{}&\int(1){3}{}&\int(0){\0}{}&\int(0){\0}{}&&\cnt{552}&\bounds{108}{144}{fail}{false}&&\cr
\l{1}&\rat(3/2){2}{}&\int(1){3}{}&\rat(1/2){1}{}&\int(0){\0}{}&\int(0){\0}{}&\int(0){1}{}&\ln{12\sd}&\cnt{32}&\bounds{5}{8}{5}{true}&&\cr
\l{2}&\rat(3/2){2}{}&\rat(1/2){2}{}&\rat(1/2){2}{}&\rat(1/2){2}{}&\int(0){\0}{}&\int(0){\0}{}&\ln{2\sd}&\cnt{64}&\bounds{16}{16}{1}{false}&&\cr
\l{3}&\int(3){\4}{2}&\int(0){\0}{}&\int(0){\0}{}&\int(0){\0}{}&\int(0){\0}{}&\int(0){\0}{}&\ln{2}&\cnt{1}&\bounds{1}{1}{1}{true}&&\cr
\l{5}&\int(2){\2}{}&\int(1){\2}{}&\int(0){\0}{}&\int(0){\0}{}&\int(0){\0}{}&\int(0){\0}{}&\ln{3}&\cnt{6}&\bounds{2}{2}{1}{true}&&\cr
\endorbit
\orbit{6}{5}{90}
\h{5}&\int(3){3}{}&\root\int(1){3}{}&\int(1){3}{}&\int(1){3}{}&\int(0){\0}{}&\int(0){\0}{}&&\cnt{472}&\bounds{90}{128}{fail}{false}&&\cr
\l{1}&\rat(3/2){2}{}&\rat(1/2){2}{}&\rat(1/2){2}{}&\rat(1/2){2}{}&\int(0){\0}{}&\int(0){\0}{}&\ln{2\sd}&\cnt{32}&\bounds{7}{8}{5}{true}&&\cr
\l{2}&\rat(3/2){2}{}&\rat(1/2){1}{}&\int(1){3}{}&\int(0){\0}{}&\int(0){3}{}&\int(0){\0}{}&\ln{4\sd}&\cnt{16}&\bounds{3}{4}{5}{true}&&\cr
\l{3}&\int(3){\4}{2}&\int(0){\0}{}&\int(0){\0}{}&\int(0){\0}{}&\int(0){\0}{}&\int(0){\0}{}&\ln{2}&\cnt{1}&\bounds{1}{1}{1}{true}&&\cr
\l{5}&\int(2){\2}{}&\int(1){\2}{\spec{1}}&\int(0){\0}{}&\int(0){\0}{}&\int(0){\0}{}&\int(0){\0}{}&\ln{2}&\cnt{3}&\bounds{1}{1}{1}{true}&&\cr
\l{7}&\int(2){\2}{}&\int(0){\0}{}&\int(1){\2}{}&\int(0){\0}{}&\int(0){\0}{}&\int(0){\0}{}&\ln{2}&\cnt{18}&\bounds{4}{6}{5}{true}&&\cr
\l{9}&\rat(3/2){2}{}&\int(1){3}{}&\rat(1/2){1}{}&\int(0){\0}{}&\int(0){\0}{}&\int(0){1}{}&\ln{4}&\cnt{32}&\bounds{5}{8}{5}{true}&&\cr
\endorbit
\orbit{6}{6}{118}
\h{6}&\int(1){3}{}&\int(1){3}{}&\int(1){3}{}&\int(1){3}{}&\int(2){\2}{}&\root\int(0){\0}{}&&\cnt{480}&\bounds{118}{136}{fail}{false}&&\cr
\l{1}&\int(1){3}{}&\int(1){3}{}&\int(0){\0}{}&\int(0){\0}{}&\int(1){3}{}&\int(0){3}{}&\ln{6\sd}&\cnt{8}&\bounds{3}{4}{5}{true}&&\cr
\l{2}&\int(1){3}{}&\rat(1/2){2}{}&\rat(1/2){1}{}&\int(0){\0}{}&\int(1){1}{}&\int(0){\0}{}&\ln{12\sd}&\cnt{32}&\bounds{7}{8}{5}{true}&&\cr
\l{3}&\int(1){\2}{}&\int(0){\0}{}&\int(0){\0}{}&\int(0){\0}{}&\int(2){\2}{}&\int(0){\0}{}&\ln{4}&\cnt{6}&\bounds{2}{2}{1}{true}&&\cr
\endorbit
\orbit{6}{7}{112}
\h{7}&\root\int(1){3}{}&\int(1){3}{}&\int(1){3}{}&\int(1){3}{}&\int(2){\2}{}&\int(0){\0}{}&&\cnt{424}&\bounds{112}{112}{1}{false}&&\cr
\l{1}&\rat(1/2){2}{}&\int(1){3}{}&\int(0){\0}{}&\rat(1/2){1}{}&\int(1){1}{}&\int(0){\0}{}&\ln{6\sd}&\cnt{16}&\bounds{4}{4}{1}{false}&&\cr
\l{2}&\int(1){\2}{\spec{1}}&\int(0){\0}{}&\int(0){\0}{}&\int(0){\0}{}&\int(2){\2}{}&\int(0){\0}{}&\ln{2}&\cnt{1}&\bounds{1}{1}{1}{false}&&\cr
\l{4}&\int(1){3}{}&\int(1){3}{}&\int(1){3}{}&\int(0){3}{}&\int(0){\0}{}&\int(0){\0}{}&\ln{3}&\cnt{6}&\bounds{2}{2}{1}{false}&&\cr
\l{6}&\int(1){3}{}&\int(1){3}{}&\int(0){\0}{}&\int(0){\0}{}&\int(1){3}{}&\int(0){3}{}&\ln{3}&\cnt{16}&\bounds{4}{4}{1}{false}&&\cr
\l{8}&\int(1){3}{}&\rat(1/2){2}{}&\rat(1/2){1}{}&\int(0){\0}{}&\int(1){1}{}&\int(0){\0}{}&\ln{3}&\cnt{32}&\bounds{8}{8}{1}{false}&&\cr
\endorbit
\orbit{6}{8}{96}
\h{8}&\int(1){3}{}&\int(1){3}{}&\int(1){3}{}&\int(1){3}{}&\root\int(2){\2}{\spec{1}}&\int(0){\0}{}&&\cnt{368}&\bounds{96}{96}{1}{false}&&\cr
\l{1}&\int(1){3}{}&\rat(1/2){2}{}&\rat(1/2){1}{}&\int(0){\0}{}&\int(1){1}{}&\int(0){\0}{}&\ln{8\sd}&\cnt{32}&\bounds{8}{8}{1}{false}&&\cr
\l{2}&\int(1){3}{}&\int(0){\0}{}&\int(1){3}{}&\int(0){\0}{}&\int(1){2}{}&\int(0){1}{}&\ln{4\sd}&\cnt{16}&\bounds{4}{4}{1}{false}&&\cr
\l{3}&\int(1){\2}{}&\int(0){\0}{}&\int(0){\0}{}&\int(0){\0}{}&\int(2){\2}{}&\int(0){\0}{}&\ln{4}&\cnt{6}&\bounds{2}{2}{1}{false}&&\cr
\endorbit
\orbit{6}{9}{102}
\h{9}&\root\int(4){\4}{3}&\int(2){\2}{}&\int(0){\0}{}&\int(0){\0}{}&\int(0){\0}{}&\int(0){\0}{}&&\cnt{784}&\bounds{102}{200}{fail}{false}&&\cr
\l{1}&\int(2){\2}{\spec{1}}&\int(1){\2}{}&\int(0){\0}{}&\int(0){\0}{}&\int(0){\0}{}&\int(0){\0}{}&\ln{2\sd}&\cnt{8}&\bounds{3}{4}{5}{true}&&\cr
\l{2}&\int(2){3}{}&\int(1){3}{}&\int(0){3}{}&\int(0){3}{}&\int(0){\0}{}&\int(0){\0}{}&\ln{6\sdd}&\cnt{128}&\bounds{16}{32}{5}{true}&&\cr
\endorbit
\orbit{6}{10}{88}
\h{10}&\int(4){\4}{3}&\int(2){\2}{}&\root\int(0){\0}{}&\int(0){\0}{}&\int(0){\0}{}&\int(0){\0}{}&&\cnt{624}&\bounds{88}{160}{fail}{false}&&\cr
\l{1}&\int(2){\2}{}&\int(1){\2}{}&\int(0){\0}{}&\int(0){\0}{}&\int(0){\0}{}&\int(0){\0}{}&\ln{1\sdd}&\cnt{48}&\bounds{10}{16}{5}{true}&&\cr
\l{2}&\int(2){3}{}&\int(1){3}{}&\int(0){3}{}&\int(0){3}{}&\int(0){\0}{}&\int(0){\0}{}&\ln{3\sdd}&\cnt{64}&\bounds{10}{16}{5}{true}&&\cr
\l{3}&\int(2){3}{}&\int(1){3}{}&\int(0){\0}{}&\int(0){\0}{}&\int(0){3}{}&\int(0){3}{}&\ln{3\sdd}&\cnt{128}&\bounds{16}{32}{5}{true}&&\cr
\endorbit
\orbit{6}{11}{80}
\h{11}&\int(4){\4}{3}&\root\int(2){\2}{\spec{1}}&\int(0){\0}{}&\int(0){\0}{}&\int(0){\0}{}&\int(0){\0}{}&&\cnt{512}&\bounds{80}{128}{fail}{false}&&\cr
\l{1}&\int(2){3}{}&\int(1){2}{}&\int(0){1}{}&\int(0){\0}{}&\int(0){1}{}&\int(0){\0}{}&\ln{4\sdd}&\cnt{128}&\bounds{20}{32}{2}{false}&&\cr
\endorbit
\orbit{6}{12}{120}
\h{12}&\int(1){3}{}&\int(1){3}{}&\int(1){2}{}&\int(1){2}{}&\int(1){1}{}&\int(1){1}{}&&\cnt{480}&\bounds{120}{120}{1}{false}&&\cr
\l{1}&\int(1){3}{}&\int(1){3}{}&\rat(1/2){3}{}&\rat(1/2){3}{}&\int(0){\0}{}&\int(0){\0}{}&\ln{30\sd}&\cnt{16}&\bounds{4}{4}{1}{true}&&\cr
\endorbit
\orbit{6}{13}{50}
\h{13}&\int(6){\6}{\spec{1}}&\int(0){\0}{}&\int(0){\0}{}&\int(0){\0}{}&\int(0){\0}{}&\int(0){\0}{}&&\cnt{240}&\bounds{50}{50}{1}{false}&&\cr
\l{1}&\int(3){\2}{}&\int(0){\2}{}&\int(0){\0}{}&\int(0){\0}{}&\int(0){\0}{}&\int(0){\0}{}&\ln{5\sdd}&\cnt{48}&\bounds{10}{10}{1}{false}&&\cr
\endorbit
\endorbits

\newcs{check-6A4-1}{}
\newcs{check-6A4-7}{\noexpand\checkdone}
\newcs{check-6A4-21}{\noexpand\checkdone}
\newcs{check-6A4-27}{\noexpand\checkdone}
\newcs{max-6A4}{150}%
\newcs{count-6A4}{39}%
\newcs{counts-6A4}{\DO{1}{150}\DO{2}{118}\DO{3}{112}\DO{4}{118}\DO{5}{112}\DO{6}{108}\DO{7}{132}\DO{8}{116}\DO{9}{108}\DO{10}{106}\DO{11}{104}\DO{12}{118}\DO{13}{114}\DO{14}{112}\DO{15}{104}\DO{16}{96}\DO{17}{116}\DO{18}{110}\DO{19}{114}\DO{20}{96}\DO{21}{126}\DO{22}{114}\DO{23}{112}\DO{24}{100}\DO{25}{94}\DO{26}{92}\DO{27}{122}\DO{28}{110}\DO{29}{102}\DO{30}{120}\DO{31}{120}\DO{32}{110}\DO{33}{100}\DO{34}{96}\DO{35}{110}\DO{36}{110}\DO{37}{120}\DO{38}{100}\DO{39}{51}}%
\orbits
\orbit{6}{1}{150}
\h{1}&\rat(6/5){3}{}&\rat(6/5){3}{}&\rat(6/5){2}{}&\rat(6/5){2}{}&\rat(6/5){3}{}&\root\int(0){\0}{}&&\cnt{452}&\bounds{150}{156}{fail}{false}&&\cr
\l{1}&\rat(6/5){3}{}&\rat(6/5){3}{}&\rat(3/5){1}{}&\int(0){\0}{}&\int(0){\0}{}&\int(0){1}{}&\ln{10\sd}&\cnt{6}&\bounds{2}{2}{1}{true}&&\cr
\l{2}&\rat(6/5){3}{}&\rat(4/5){2}{}&\int(0){\0}{}&\rat(2/5){4}{}&\rat(3/5){4}{}&\int(0){\0}{}&\ln{20\sd}&\cnt{18}&\bounds{6}{6}{1}{true}&&\cr
\l{3}&\rat(3/5){4}{}&\rat(3/5){4}{}&\rat(3/5){1}{}&\rat(3/5){1}{}&\rat(3/5){4}{}&\int(0){\0}{}&\ln{1\sdd}&\cnt{32}&\bounds{10}{16}{2}{true}&&\cr
\endorbit
\orbit{6}{2}{118}
\h{2}&\root\rat(6/5){3}{\spec{1}}&\rat(6/5){3}{}&\rat(6/5){2}{}&\rat(6/5){2}{}&\rat(6/5){3}{}&\int(0){\0}{}&&\cnt{396}&\bounds{118}{124}{fail}{false}&&\cr
\l{1}&\rat(3/5){4}{}&\rat(6/5){3}{}&\rat(2/5){4}{}&\rat(4/5){3}{}&\int(0){\0}{}&\int(0){\0}{}&\ln{4\sd}&\cnt{18}&\bounds{6}{6}{1}{false}&&\cr
\l{2}&\rat(3/5){4}{}&\rat(6/5){3}{}&\int(0){\0}{}&\int(0){\0}{}&\rat(6/5){3}{}&\int(0){4}{}&\ln{2\sd}&\cnt{10}&\bounds{2}{2}{1}{false}&&\cr
\l{3}&\rat(3/5){4}{}&\rat(3/5){4}{}&\rat(3/5){1}{}&\rat(3/5){1}{}&\rat(3/5){4}{}&\int(0){\0}{}&\ln{1\sdd}&\cnt{32}&\bounds{10}{16}{2}{false}&&\cr
\l{4}&\rat(6/5){3}{}&\rat(6/5){3}{}&\rat(3/5){1}{}&\int(0){\0}{}&\int(0){\0}{}&\int(0){1}{}&\ln{4}&\cnt{10}&\bounds{2}{2}{1}{false}&&\cr
\l{6}&\rat(6/5){3}{}&\rat(4/5){2}{}&\int(0){\0}{}&\rat(2/5){4}{}&\rat(3/5){4}{}&\int(0){\0}{}&\ln{4}&\cnt{18}&\bounds{6}{6}{1}{false}&&\cr
\l{8}&\rat(4/5){2}{}&\rat(6/5){3}{}&\int(0){\0}{}&\rat(3/5){1}{}&\rat(2/5){1}{}&\int(0){\0}{}&\ln{4}&\cnt{6}&\bounds{2}{2}{1}{false}&&\cr
\endorbit
\orbit{6}{3}{112}
\h{3}&\root\rat(6/5){3}{\spec{2}}&\rat(6/5){3}{}&\rat(6/5){2}{}&\rat(6/5){2}{}&\rat(6/5){3}{}&\int(0){\0}{}&&\cnt{368}&\bounds{112}{112}{1}{false}&&\cr
\l{1}&\rat(6/5){3}{}&\rat(6/5){3}{}&\rat(3/5){1}{}&\int(0){\0}{}&\int(0){\0}{}&\int(0){1}{}&\ln{4}&\cnt{10}&\bounds{2}{2}{1}{false}&&\cr
\l{3}&\rat(6/5){3}{}&\rat(4/5){2}{}&\int(0){\0}{}&\rat(2/5){4}{}&\rat(3/5){4}{}&\int(0){\0}{}&\ln{4}&\cnt{18}&\bounds{6}{6}{1}{false}&&\cr
\l{5}&\rat(4/5){2}{}&\rat(6/5){3}{}&\int(0){\0}{}&\rat(3/5){1}{}&\rat(2/5){1}{}&\int(0){\0}{}&\ln{4}&\cnt{18}&\bounds{6}{6}{1}{false}&&\cr
\endorbit
\orbit{6}{4}{118}
\h{4}&\root\int(4){\4}{\spec{1}}&\int(2){\2}{}&\int(0){\0}{}&\int(0){\0}{}&\int(0){\0}{}&\int(0){\0}{}&&\cnt{704}&\bounds{118}{142}{fail}{false}&&\cr
\l{1}&\int(3){\4}{}&\int(0){\0}{}&\int(0){\0}{}&\int(0){\0}{}&\int(0){\0}{}&\int(0){\0}{}&\ln{1}&\cnt{2}&\bounds{1}{1}{1}{true}&&\cr
\l{3}&\int(2){3}{}&\int(1){3}{}&\int(0){1}{}&\int(0){\0}{}&\int(0){\0}{}&\int(0){1}{}&\ln{2}&\cnt{75}&\bounds{13}{15}{3}{false}&&\cr
\l{5}&\int(2){3}{}&\int(1){4}{}&\int(0){4}{}&\int(0){\0}{}&\int(0){2}{}&\int(0){\0}{}&\ln{4}&\cnt{50}&\bounds{8}{10}{5}{true}&&\cr
\endorbit
\orbit{6}{5}{112}
\h{5}&\int(4){\4}{}&\int(2){\2}{}&\root\int(0){\0}{}&\int(0){\0}{}&\int(0){\0}{}&\int(0){\0}{}&&\cnt{572}&\bounds{112}{120}{fail}{false}&&\cr
\l{1}&\int(2){\2}{}&\int(1){\2}{\spec{1}}&\int(0){\0}{}&\int(0){\0}{}&\int(0){\0}{}&\int(0){\0}{}&\ln{2\sd}&\cnt{12}&\bounds{4}{4}{1}{true}&&\cr
\l{2}&\int(2){3}{}&\int(1){3}{}&\int(0){1}{}&\int(0){\0}{}&\int(0){\0}{}&\int(0){1}{}&\ln{2\sd}&\cnt{45}&\bounds{9}{9}{1}{true}&&\cr
\l{3}&\int(2){3}{}&\int(1){2}{}&\int(0){\0}{}&\int(0){4}{}&\int(0){4}{}&\int(0){\0}{}&\ln{2\sd}&\cnt{75}&\bounds{15}{15}{1}{false}&&\cr
\l{4}&\int(2){3}{}&\int(1){4}{}&\int(0){4}{}&\int(0){\0}{}&\int(0){2}{}&\int(0){\0}{}&\ln{2\sd}&\cnt{30}&\bounds{6}{6}{1}{true}&&\cr
\l{5}&\int(2){3}{}&\int(1){4}{}&\int(0){\0}{}&\int(0){2}{}&\int(0){\0}{}&\int(0){4}{}&\ln{2\sd}&\cnt{50}&\bounds{8}{10}{5}{true}&&\cr
\l{6}&\int(2){3}{}&\int(1){1}{}&\int(0){3}{\spec{1}}&\int(0){1}{}&\int(0){\0}{}&\int(0){\0}{}&\ln{2\sd}&\cnt{15}&\bounds{3}{3}{1}{true}&&\cr
\l{7}&\int(2){3}{}&\int(1){1}{}&\int(0){3}{\spec{2}}&\int(0){1}{}&\int(0){\0}{}&\int(0){\0}{}&\ln{2\sd}&\cnt{5}&\bounds{1}{1}{1}{true}&&\cr
\l{8}&\int(2){3}{}&\int(1){1}{}&\int(0){\0}{}&\int(0){\0}{}&\int(0){1}{}&\int(0){3}{}&\ln{2\sd}&\cnt{50}&\bounds{8}{10}{5}{true}&&\cr
\l{9}&\int(3){\4}{\spec{1}}&\int(0){\0}{}&\int(0){\0}{}&\int(0){\0}{}&\int(0){\0}{}&\int(0){\0}{}&\ln{2}&\cnt{2}&\bounds{1}{1}{1}{true}&&\cr
\endorbit
\orbit{6}{6}{108}
\h{6}&\int(4){\4}{}&\root\int(2){\2}{}&\int(0){\0}{}&\int(0){\0}{}&\int(0){\0}{}&\int(0){\0}{}&&\cnt{516}&\bounds{108}{108}{1}{false}&&\cr
\l{1}&\int(2){\2}{}&\int(1){\2}{\spec{1}}&\int(0){\0}{}&\int(0){\0}{}&\int(0){\0}{}&\int(0){\0}{}&\ln{2\sd}&\cnt{4}&\bounds{2}{2}{1}{false}&&\cr
\l{2}&\int(2){3}{}&\int(1){3}{}&\int(0){1}{}&\int(0){\0}{}&\int(0){\0}{}&\int(0){1}{}&\ln{4\sd}&\cnt{25}&\bounds{5}{5}{1}{false}&&\cr
\l{3}&\int(2){3}{}&\int(1){4}{}&\int(0){4}{}&\int(0){\0}{}&\int(0){2}{}&\int(0){\0}{}&\ln{8\sd}&\cnt{50}&\bounds{10}{10}{1}{false}&&\cr
\l{4}&\int(3){\4}{\spec{1}}&\int(0){\0}{}&\int(0){\0}{}&\int(0){\0}{}&\int(0){\0}{}&\int(0){\0}{}&\ln{2}&\cnt{2}&\bounds{1}{1}{1}{false}&&\cr
\endorbit
\orbit{6}{7}{132}
\h{7}&\rat(4/5){4}{}&\rat(4/5){4}{}&\rat(6/5){3}{}&\int(2){\2}{}&\root\int(0){\0}{}&\rat(6/5){3}{}&&\cnt{476}&\bounds{132}{140}{fail}{false}&&\cr
\l{1}&\rat(3/5){3}{}&\rat(1/5){1}{}&\rat(6/5){3}{}&\int(1){1}{}&\int(0){\0}{}&\int(0){\0}{}&\ln{2\sd}&\cnt{16}&\bounds{4}{4}{1}{true}&&\cr
\l{2}&\rat(4/5){4}{}&\int(0){\0}{}&\rat(6/5){3}{}&\int(1){4}{}&\int(0){3}{\spec{1}}&\int(0){\0}{}&\ln{2\sd}&\cnt{3}&\bounds{1}{1}{1}{true}&&\cr
\l{3}&\rat(4/5){4}{}&\int(0){\0}{}&\rat(6/5){3}{}&\int(1){4}{}&\int(0){3}{\spec{2}}&\int(0){\0}{}&\ln{2\sd}&\cnt{1}&\bounds{1}{1}{1}{true}&&\cr
\l{4}&\rat(4/5){4}{}&\int(0){\0}{}&\rat(4/5){2}{}&\int(1){2}{}&\int(0){\0}{}&\rat(2/5){1}{}&\ln{2\sd}&\cnt{27}&\bounds{9}{9}{1}{true}&&\cr
\l{5}&\rat(4/5){4}{}&\int(0){\0}{}&\rat(3/5){4}{}&\int(1){1}{}&\int(0){1}{}&\rat(3/5){4}{}&\ln{2\sd}&\cnt{12}&\bounds{4}{4}{1}{true}&&\cr
\l{6}&\rat(4/5){4}{}&\int(0){\0}{}&\int(0){\0}{}&\int(1){3}{}&\int(0){4}{}&\rat(6/5){3}{}&\ln{2\sd}&\cnt{9}&\bounds{3}{3}{1}{true}&&\cr
\l{7}&\int(1){\2}{}&\int(0){\0}{}&\int(0){\0}{}&\int(2){\2}{}&\int(0){\0}{}&\int(0){\0}{}&\ln{2}&\cnt{4}&\bounds{1}{1}{1}{true}&&\cr
\l{9}&\rat(3/5){3}{}&\rat(4/5){4}{}&\int(0){\0}{}&\int(1){2}{}&\int(0){\0}{}&\rat(3/5){4}{}&\ln{2}&\cnt{24}&\bounds{6}{6}{1}{true}&&\cr
\l{11}&\rat(3/5){3}{}&\int(0){\0}{}&\rat(3/5){4}{}&\int(1){4}{}&\int(0){\0}{}&\rat(4/5){2}{}&\ln{2}&\cnt{24}&\bounds{6}{6}{1}{true}&&\cr
\l{13}&\rat(2/5){2}{}&\rat(4/5){4}{}&\rat(4/5){2}{}&\int(1){4}{}&\int(0){\0}{}&\int(0){\0}{}&\ln{2}&\cnt{18}&\bounds{4}{6}{5}{true}&&\cr
\l{15}&\rat(4/5){4}{}&\rat(4/5){4}{}&\rat(6/5){3}{}&\int(0){\0}{}&\int(0){\0}{}&\rat(1/5){3}{}&\ln{2}&\cnt{6}&\bounds{2}{2}{1}{true}&&\cr
\l{17}&\rat(4/5){4}{}&\rat(4/5){4}{}&\rat(2/5){1}{}&\int(1){1}{}&\int(0){4}{}&\int(0){\0}{}&\ln{2}&\cnt{9}&\bounds{3}{3}{1}{true}&&\cr
\endorbit
\orbit{6}{8}{116}
\h{8}&\root\rat(4/5){4}{}&\rat(4/5){4}{}&\rat(6/5){3}{}&\int(2){\2}{}&\int(0){\0}{}&\rat(6/5){3}{}&&\cnt{448}&\bounds{116}{120}{fail}{false}&&\cr
\l{1}&\int(1){\2}{}&\int(0){\0}{}&\int(0){\0}{}&\int(2){\2}{}&\int(0){\0}{}&\int(0){\0}{}&\ln{1}&\cnt{2}&\bounds{1}{1}{1}{true}&&\cr
\l{3}&\rat(3/5){3}{}&\rat(4/5){4}{}&\int(0){\0}{}&\int(1){2}{}&\int(0){\0}{}&\rat(3/5){4}{}&\ln{1}&\cnt{12}&\bounds{4}{4}{1}{true}&&\cr
\l{5}&\rat(3/5){3}{}&\rat(1/5){1}{}&\rat(6/5){3}{}&\int(1){1}{}&\int(0){\0}{}&\int(0){\0}{}&\ln{1}&\cnt{8}&\bounds{2}{2}{1}{true}&&\cr
\l{7}&\rat(3/5){3}{}&\int(0){\0}{}&\rat(3/5){4}{}&\int(1){4}{}&\int(0){\0}{}&\rat(4/5){2}{}&\ln{1}&\cnt{12}&\bounds{4}{4}{1}{true}&&\cr
\l{9}&\rat(2/5){2}{\spec{1}}&\rat(4/5){4}{}&\rat(4/5){2}{}&\int(1){4}{}&\int(0){\0}{}&\int(0){\0}{}&\ln{1}&\cnt{3}&\bounds{1}{1}{1}{true}&&\cr
\l{11}&\rat(2/5){2}{\spec{1}}&\int(0){\0}{}&\rat(2/5){1}{}&\int(1){1}{}&\int(0){\0}{}&\rat(6/5){3}{}&\ln{1}&\cnt{3}&\bounds{1}{1}{1}{true}&&\cr
\l{13}&\rat(4/5){4}{}&\rat(3/5){3}{}&\rat(3/5){4}{}&\int(1){3}{}&\int(0){\0}{}&\int(0){\0}{}&\ln{1}&\cnt{24}&\bounds{6}{6}{1}{true}&&\cr
\l{15}&\rat(4/5){4}{}&\rat(2/5){2}{}&\int(0){\0}{}&\int(1){1}{}&\int(0){\0}{}&\rat(4/5){2}{}&\ln{1}&\cnt{18}&\bounds{4}{6}{5}{true}&&\cr
\l{17}&\rat(4/5){4}{}&\rat(4/5){4}{}&\rat(6/5){3}{}&\int(0){\0}{}&\int(0){\0}{}&\rat(1/5){3}{}&\ln{1}&\cnt{6}&\bounds{2}{2}{1}{true}&&\cr
\l{19}&\rat(4/5){4}{}&\rat(4/5){4}{}&\rat(1/5){3}{}&\int(0){\0}{}&\int(0){\0}{}&\rat(6/5){3}{}&\ln{1}&\cnt{6}&\bounds{2}{2}{1}{true}&&\cr
\l{21}&\rat(4/5){4}{}&\rat(4/5){4}{}&\rat(2/5){1}{}&\int(1){1}{}&\int(0){4}{}&\int(0){\0}{}&\ln{1}&\cnt{15}&\bounds{3}{3}{1}{true}&&\cr
\l{23}&\rat(4/5){4}{}&\rat(4/5){4}{}&\int(0){\0}{}&\int(1){4}{}&\int(0){1}{}&\rat(2/5){1}{}&\ln{1}&\cnt{15}&\bounds{3}{3}{1}{true}&&\cr
\l{25}&\rat(4/5){4}{}&\rat(1/5){1}{}&\rat(2/5){1}{}&\int(1){4}{}&\int(0){\0}{}&\rat(3/5){4}{}&\ln{1}&\cnt{24}&\bounds{6}{6}{1}{true}&&\cr
\l{27}&\rat(4/5){4}{}&\rat(-1/5){4}{}&\rat(6/5){3}{}&\int(0){\0}{}&\int(0){\0}{}&\rat(6/5){3}{}&\ln{1}&\cnt{4}&\bounds{1}{1}{1}{true}&&\cr
\l{29}&\rat(4/5){4}{}&\int(0){\0}{}&\rat(6/5){3}{}&\int(1){4}{}&\int(0){3}{}&\int(0){\0}{}&\ln{1}&\cnt{10}&\bounds{2}{2}{1}{true}&&\cr
\l{31}&\rat(4/5){4}{}&\int(0){\0}{}&\rat(4/5){2}{}&\int(1){2}{}&\int(0){\0}{}&\rat(2/5){1}{}&\ln{1}&\cnt{27}&\bounds{9}{9}{1}{true}&&\cr
\l{33}&\rat(4/5){4}{}&\int(0){\0}{}&\rat(3/5){4}{}&\int(1){1}{}&\int(0){1}{}&\rat(3/5){4}{}&\ln{1}&\cnt{20}&\bounds{4}{4}{1}{true}&&\cr
\l{35}&\rat(4/5){4}{}&\int(0){\0}{}&\int(0){\0}{}&\int(1){3}{}&\int(0){4}{}&\rat(6/5){3}{}&\ln{1}&\cnt{15}&\bounds{3}{3}{1}{true}&&\cr
\endorbit
\orbit{6}{9}{108}
\h{9}&\rat(4/5){4}{}&\rat(4/5){4}{}&\rat(6/5){3}{}&\root\int(2){\2}{}&\int(0){\0}{}&\rat(6/5){3}{}&&\cnt{420}&\bounds{108}{108}{1}{false}&&\cr
\l{1}&\rat(3/5){3}{}&\rat(1/5){1}{}&\rat(6/5){3}{}&\int(1){1}{}&\int(0){\0}{}&\int(0){\0}{}&\ln{2\sd}&\cnt{16}&\bounds{4}{4}{1}{false}&&\cr
\l{2}&\rat(4/5){4}{}&\int(0){\0}{}&\rat(6/5){3}{}&\int(1){4}{}&\int(0){3}{}&\int(0){\0}{}&\ln{2\sd}&\cnt{10}&\bounds{2}{2}{1}{false}&&\cr
\l{3}&\rat(4/5){4}{}&\int(0){\0}{}&\rat(4/5){2}{}&\int(1){2}{}&\int(0){\0}{}&\rat(2/5){1}{}&\ln{2\sd}&\cnt{9}&\bounds{3}{3}{1}{false}&&\cr
\l{4}&\rat(4/5){4}{}&\int(0){\0}{}&\rat(3/5){4}{}&\int(1){1}{}&\int(0){1}{}&\rat(3/5){4}{}&\ln{2\sd}&\cnt{20}&\bounds{4}{4}{1}{false}&&\cr
\l{5}&\rat(4/5){4}{}&\int(0){\0}{}&\int(0){\0}{}&\int(1){3}{}&\int(0){4}{}&\rat(6/5){3}{}&\ln{2\sd}&\cnt{5}&\bounds{1}{1}{1}{false}&&\cr
\l{6}&\int(1){\2}{}&\int(0){\0}{}&\int(0){\0}{}&\int(2){\2}{}&\int(0){\0}{}&\int(0){\0}{}&\ln{2}&\cnt{4}&\bounds{1}{1}{1}{false}&&\cr
\l{8}&\rat(3/5){3}{}&\rat(4/5){4}{}&\int(0){\0}{}&\int(1){2}{}&\int(0){\0}{}&\rat(3/5){4}{}&\ln{2}&\cnt{8}&\bounds{2}{2}{1}{false}&&\cr
\l{10}&\rat(3/5){3}{}&\int(0){\0}{}&\rat(3/5){4}{}&\int(1){4}{}&\int(0){\0}{}&\rat(4/5){2}{}&\ln{2}&\cnt{24}&\bounds{6}{6}{1}{false}&&\cr
\l{12}&\rat(2/5){2}{}&\rat(4/5){4}{}&\rat(4/5){2}{}&\int(1){4}{}&\int(0){\0}{}&\int(0){\0}{}&\ln{2}&\cnt{18}&\bounds{6}{6}{1}{false}&&\cr
\l{14}&\rat(4/5){4}{}&\rat(4/5){4}{}&\rat(6/5){3}{}&\int(0){\0}{}&\int(0){\0}{}&\rat(1/5){3}{}&\ln{2}&\cnt{6}&\bounds{2}{2}{1}{false}&&\cr
\l{16}&\rat(4/5){4}{}&\rat(4/5){4}{}&\rat(2/5){1}{}&\int(1){1}{}&\int(0){4}{}&\int(0){\0}{}&\ln{2}&\cnt{15}&\bounds{3}{3}{1}{false}&&\cr
\endorbit
\orbit{6}{10}{106}
\h{10}&\rat(4/5){4}{}&\rat(4/5){4}{}&\root\rat(6/5){3}{\spec{1}}&\int(2){\2}{}&\int(0){\0}{}&\rat(6/5){3}{}&&\cnt{420}&\bounds{106}{106}{1}{false}&&\cr
\l{1}&\int(1){\2}{}&\int(0){\0}{}&\int(0){\0}{}&\int(2){\2}{}&\int(0){\0}{}&\int(0){\0}{}&\ln{1}&\cnt{4}&\bounds{1}{1}{1}{false}&&\cr
\l{3}&\rat(3/5){3}{}&\rat(4/5){4}{}&\int(0){\0}{}&\int(1){2}{}&\int(0){\0}{}&\rat(3/5){4}{}&\ln{1}&\cnt{24}&\bounds{6}{6}{1}{false}&&\cr
\l{5}&\rat(3/5){3}{}&\rat(1/5){1}{}&\rat(6/5){3}{}&\int(1){1}{}&\int(0){\0}{}&\int(0){\0}{}&\ln{1}&\cnt{16}&\bounds{4}{4}{1}{false}&&\cr
\l{7}&\rat(3/5){3}{}&\int(0){\0}{}&\rat(3/5){4}{}&\int(1){4}{}&\int(0){\0}{}&\rat(4/5){2}{}&\ln{1}&\cnt{24}&\bounds{6}{6}{1}{false}&&\cr
\l{9}&\rat(2/5){2}{}&\rat(4/5){4}{}&\rat(4/5){2}{}&\int(1){4}{}&\int(0){\0}{}&\int(0){\0}{}&\ln{1}&\cnt{6}&\bounds{2}{2}{1}{false}&&\cr
\l{11}&\rat(4/5){4}{}&\rat(3/5){3}{}&\rat(3/5){4}{}&\int(1){3}{}&\int(0){\0}{}&\int(0){\0}{}&\ln{1}&\cnt{24}&\bounds{6}{6}{1}{false}&&\cr
\l{13}&\rat(4/5){4}{}&\rat(2/5){2}{}&\int(0){\0}{}&\int(1){1}{}&\int(0){\0}{}&\rat(4/5){2}{}&\ln{1}&\cnt{18}&\bounds{6}{6}{1}{false}&&\cr
\l{15}&\rat(4/5){4}{}&\rat(4/5){4}{}&\rat(6/5){3}{}&\int(0){\0}{}&\int(0){\0}{}&\rat(1/5){3}{}&\ln{1}&\cnt{6}&\bounds{2}{2}{1}{false}&&\cr
\l{17}&\rat(4/5){4}{}&\rat(4/5){4}{}&\rat(1/5){3}{}&\int(0){\0}{}&\int(0){\0}{}&\rat(6/5){3}{}&\ln{1}&\cnt{2}&\bounds{1}{1}{1}{false}&&\cr
\l{19}&\rat(4/5){4}{}&\rat(4/5){4}{}&\rat(2/5){1}{}&\int(1){1}{}&\int(0){4}{}&\int(0){\0}{}&\ln{1}&\cnt{5}&\bounds{1}{1}{1}{false}&&\cr
\l{21}&\rat(4/5){4}{}&\rat(4/5){4}{}&\int(0){\0}{}&\int(1){4}{}&\int(0){1}{}&\rat(2/5){1}{}&\ln{1}&\cnt{15}&\bounds{3}{3}{1}{false}&&\cr
\l{23}&\rat(4/5){4}{}&\rat(1/5){1}{}&\rat(2/5){1}{}&\int(1){4}{}&\int(0){\0}{}&\rat(3/5){4}{}&\ln{1}&\cnt{8}&\bounds{2}{2}{1}{false}&&\cr
\l{25}&\rat(4/5){4}{}&\rat(-1/5){4}{}&\rat(6/5){3}{}&\int(0){\0}{}&\int(0){\0}{}&\rat(6/5){3}{}&\ln{1}&\cnt{4}&\bounds{1}{1}{1}{false}&&\cr
\l{27}&\rat(4/5){4}{}&\int(0){\0}{}&\rat(6/5){3}{}&\int(1){4}{}&\int(0){3}{}&\int(0){\0}{}&\ln{1}&\cnt{10}&\bounds{2}{2}{1}{false}&&\cr
\l{29}&\rat(4/5){4}{}&\int(0){\0}{}&\rat(4/5){2}{}&\int(1){2}{}&\int(0){\0}{}&\rat(2/5){1}{}&\ln{1}&\cnt{9}&\bounds{3}{3}{1}{false}&&\cr
\l{31}&\rat(4/5){4}{}&\int(0){\0}{}&\rat(3/5){4}{}&\int(1){1}{}&\int(0){1}{}&\rat(3/5){4}{}&\ln{1}&\cnt{20}&\bounds{4}{4}{1}{false}&&\cr
\l{33}&\rat(4/5){4}{}&\int(0){\0}{}&\int(0){\0}{}&\int(1){3}{}&\int(0){4}{}&\rat(6/5){3}{}&\ln{1}&\cnt{15}&\bounds{3}{3}{1}{false}&&\cr
\endorbit
\orbit{6}{11}{104}
\h{11}&\rat(4/5){4}{}&\rat(4/5){4}{}&\root\rat(6/5){3}{\spec{2}}&\int(2){\2}{}&\int(0){\0}{}&\rat(6/5){3}{}&&\cnt{392}&\bounds{104}{104}{1}{false}&&\cr
\l{1}&\int(1){\2}{}&\int(0){\0}{}&\int(0){\0}{}&\int(2){\2}{}&\int(0){\0}{}&\int(0){\0}{}&\ln{1}&\cnt{4}&\bounds{1}{1}{1}{false}&&\cr
\l{3}&\rat(3/5){3}{}&\rat(4/5){4}{}&\int(0){\0}{}&\int(1){2}{}&\int(0){\0}{}&\rat(3/5){4}{}&\ln{1}&\cnt{24}&\bounds{6}{6}{1}{false}&&\cr
\l{5}&\rat(3/5){3}{}&\rat(1/5){1}{}&\rat(6/5){3}{}&\int(1){1}{}&\int(0){\0}{}&\int(0){\0}{}&\ln{1}&\cnt{16}&\bounds{4}{4}{1}{false}&&\cr
\l{7}&\rat(2/5){2}{}&\rat(4/5){4}{}&\rat(4/5){2}{}&\int(1){4}{}&\int(0){\0}{}&\int(0){\0}{}&\ln{1}&\cnt{18}&\bounds{6}{6}{1}{false}&&\cr
\l{9}&\rat(4/5){4}{}&\rat(2/5){2}{}&\int(0){\0}{}&\int(1){1}{}&\int(0){\0}{}&\rat(4/5){2}{}&\ln{1}&\cnt{18}&\bounds{6}{6}{1}{false}&&\cr
\l{11}&\rat(4/5){4}{}&\rat(4/5){4}{}&\rat(6/5){3}{}&\int(0){\0}{}&\int(0){\0}{}&\rat(1/5){3}{}&\ln{1}&\cnt{6}&\bounds{2}{2}{1}{false}&&\cr
\l{13}&\rat(4/5){4}{}&\rat(4/5){4}{}&\rat(2/5){1}{}&\int(1){1}{}&\int(0){4}{}&\int(0){\0}{}&\ln{1}&\cnt{15}&\bounds{3}{3}{1}{false}&&\cr
\l{15}&\rat(4/5){4}{}&\rat(4/5){4}{}&\int(0){\0}{}&\int(1){4}{}&\int(0){1}{}&\rat(2/5){1}{}&\ln{1}&\cnt{15}&\bounds{3}{3}{1}{false}&&\cr
\l{17}&\rat(4/5){4}{}&\rat(1/5){1}{}&\rat(2/5){1}{}&\int(1){4}{}&\int(0){\0}{}&\rat(3/5){4}{}&\ln{1}&\cnt{24}&\bounds{6}{6}{1}{false}&&\cr
\l{19}&\rat(4/5){4}{}&\rat(-1/5){4}{}&\rat(6/5){3}{}&\int(0){\0}{}&\int(0){\0}{}&\rat(6/5){3}{}&\ln{1}&\cnt{4}&\bounds{1}{1}{1}{false}&&\cr
\l{21}&\rat(4/5){4}{}&\int(0){\0}{}&\rat(6/5){3}{}&\int(1){4}{}&\int(0){3}{}&\int(0){\0}{}&\ln{1}&\cnt{10}&\bounds{2}{2}{1}{false}&&\cr
\l{23}&\rat(4/5){4}{}&\int(0){\0}{}&\rat(4/5){2}{}&\int(1){2}{}&\int(0){\0}{}&\rat(2/5){1}{}&\ln{1}&\cnt{27}&\bounds{9}{9}{1}{false}&&\cr
\l{25}&\rat(4/5){4}{}&\int(0){\0}{}&\int(0){\0}{}&\int(1){3}{}&\int(0){4}{}&\rat(6/5){3}{}&\ln{1}&\cnt{15}&\bounds{3}{3}{1}{false}&&\cr
\endorbit
\orbit{6}{12}{118}
\h{12}&\rat(4/5){4}{}&\rat(4/5){4}{}&\root\rat(16/5){3}{\spec{1}}&\int(0){\0}{}&\int(0){\0}{}&\rat(6/5){3}{}&&\cnt{628}&\bounds{118}{140}{fail}{false}&&\cr
\l{1}&\rat(2/5){2}{}&\rat(2/5){2}{}&\rat(8/5){4}{}&\int(0){\0}{}&\int(0){\0}{}&\rat(3/5){4}{}&\ln{1\sdd}&\cnt{72}&\bounds{12}{24}{5}{true}&&\cr
\l{2}&\rat(4/5){4}{}&\int(0){\0}{}&\rat(8/5){4}{}&\int(0){1}{}&\int(0){1}{}&\rat(3/5){4}{}&\ln{2\sd}&\cnt{50}&\bounds{9}{10}{5}{true}&&\cr
\l{3}&\int(1){\2}{}&\int(0){\0}{}&\int(2){\2}{}&\int(0){\0}{}&\int(0){\0}{}&\int(0){\0}{}&\ln{2}&\cnt{8}&\bounds{2}{2}{1}{true}&&\cr
\l{5}&\rat(3/5){3}{}&\rat(4/5){4}{}&\rat(8/5){4}{}&\int(0){\0}{}&\int(0){2}{}&\int(0){\0}{}&\ln{2}&\cnt{40}&\bounds{7}{8}{5}{true}&&\cr
\l{7}&\rat(3/5){3}{}&\int(0){\0}{}&\rat(8/5){4}{}&\int(0){4}{}&\int(0){\0}{}&\rat(4/5){2}{}&\ln{2}&\cnt{60}&\bounds{11}{12}{5}{true}&&\cr
\l{9}&\rat(4/5){4}{}&\rat(4/5){4}{}&\rat(6/5){3}{}&\int(0){\0}{}&\int(0){\0}{}&\rat(1/5){3}{}&\ln{1}&\cnt{12}&\bounds{4}{4}{1}{true}&&\cr
\endorbit
\orbit{6}{13}{114}
\h{13}&\rat(4/5){4}{}&\rat(4/5){4}{}&\rat(16/5){3}{\spec{1}}&\root\int(0){\0}{}&\int(0){\0}{}&\rat(6/5){3}{}&&\cnt{548}&\bounds{114}{130}{fail}{false}&&\cr
\l{1}&\rat(2/5){2}{}&\rat(2/5){2}{}&\rat(8/5){4}{}&\int(0){\0}{}&\int(0){\0}{}&\rat(3/5){4}{}&\ln{1\sdd}&\cnt{72}&\bounds{12}{24}{5}{true}&&\cr
\l{2}&\int(1){\2}{}&\int(0){\0}{}&\int(2){\2}{}&\int(0){\0}{}&\int(0){\0}{}&\int(0){\0}{}&\ln{1}&\cnt{16}&\bounds{4}{4}{1}{true}&&\cr
\l{4}&\rat(3/5){3}{}&\rat(4/5){4}{}&\rat(8/5){4}{}&\int(0){\0}{}&\int(0){2}{}&\int(0){\0}{}&\ln{1}&\cnt{40}&\bounds{7}{8}{5}{true}&&\cr
\l{6}&\rat(3/5){3}{}&\int(0){\0}{}&\rat(8/5){4}{}&\int(0){4}{}&\int(0){\0}{}&\rat(4/5){2}{}&\ln{1}&\cnt{36}&\bounds{9}{9}{1}{true}&&\cr
\l{8}&\rat(4/5){4}{}&\rat(3/5){3}{}&\rat(8/5){4}{}&\int(0){3}{\spec{1}}&\int(0){\0}{}&\int(0){\0}{}&\ln{1}&\cnt{12}&\bounds{3}{3}{1}{true}&&\cr
\l{10}&\rat(4/5){4}{}&\rat(3/5){3}{}&\rat(8/5){4}{}&\int(0){3}{\spec{2}}&\int(0){\0}{}&\int(0){\0}{}&\ln{1}&\cnt{4}&\bounds{1}{1}{1}{true}&&\cr
\l{12}&\rat(4/5){4}{}&\rat(4/5){4}{}&\rat(6/5){3}{}&\int(0){\0}{}&\int(0){\0}{}&\rat(1/5){3}{}&\ln{1}&\cnt{24}&\bounds{6}{6}{1}{true}&&\cr
\l{14}&\rat(4/5){4}{}&\rat(1/5){1}{}&\rat(8/5){4}{}&\int(0){\0}{}&\int(0){4}{}&\rat(2/5){1}{}&\ln{1}&\cnt{60}&\bounds{11}{12}{5}{true}&&\cr
\l{16}&\rat(4/5){4}{}&\rat(-1/5){4}{}&\rat(6/5){3}{}&\int(0){\0}{}&\int(0){\0}{}&\rat(6/5){3}{}&\ln{1}&\cnt{16}&\bounds{4}{4}{1}{true}&&\cr
\l{18}&\rat(4/5){4}{}&\int(0){\0}{}&\rat(8/5){4}{}&\int(0){1}{}&\int(0){1}{}&\rat(3/5){4}{}&\ln{1}&\cnt{30}&\bounds{6}{6}{1}{true}&&\cr
\endorbit
\orbit{6}{14}{112}
\h{14}&\root\rat(4/5){4}{}&\rat(4/5){4}{}&\rat(16/5){3}{\spec{1}}&\int(0){\0}{}&\int(0){\0}{}&\rat(6/5){3}{}&&\cnt{520}&\bounds{112}{112}{1}{false}&&\cr
\l{1}&\int(1){\2}{}&\int(0){\0}{}&\int(2){\2}{}&\int(0){\0}{}&\int(0){\0}{}&\int(0){\0}{}&\ln{1}&\cnt{8}&\bounds{2}{2}{1}{false}&&\cr
\l{3}&\rat(3/5){3}{}&\rat(4/5){4}{}&\rat(8/5){4}{}&\int(0){\0}{}&\int(0){2}{}&\int(0){\0}{}&\ln{1}&\cnt{20}&\bounds{4}{4}{1}{false}&&\cr
\l{5}&\rat(3/5){3}{}&\int(0){\0}{}&\rat(8/5){4}{}&\int(0){4}{}&\int(0){\0}{}&\rat(4/5){2}{}&\ln{1}&\cnt{30}&\bounds{6}{6}{1}{false}&&\cr
\l{7}&\rat(2/5){2}{\spec{1}}&\rat(2/5){2}{}&\rat(8/5){4}{}&\int(0){\0}{}&\int(0){\0}{}&\rat(3/5){4}{}&\ln{1}&\cnt{12}&\bounds{4}{4}{1}{false}&&\cr
\l{9}&\rat(4/5){4}{}&\rat(3/5){3}{}&\rat(8/5){4}{}&\int(0){3}{}&\int(0){\0}{}&\int(0){\0}{}&\ln{1}&\cnt{40}&\bounds{8}{8}{1}{false}&&\cr
\l{11}&\rat(4/5){4}{}&\rat(4/5){4}{}&\rat(6/5){3}{}&\int(0){\0}{}&\int(0){\0}{}&\rat(1/5){3}{}&\ln{1}&\cnt{24}&\bounds{6}{6}{1}{false}&&\cr
\l{13}&\rat(4/5){4}{}&\rat(1/5){1}{}&\rat(8/5){4}{}&\int(0){\0}{}&\int(0){4}{}&\rat(2/5){1}{}&\ln{1}&\cnt{60}&\bounds{12}{12}{1}{false}&&\cr
\l{15}&\rat(4/5){4}{}&\rat(-1/5){4}{}&\rat(6/5){3}{}&\int(0){\0}{}&\int(0){\0}{}&\rat(6/5){3}{}&\ln{1}&\cnt{16}&\bounds{4}{4}{1}{false}&&\cr
\l{17}&\rat(4/5){4}{}&\int(0){\0}{}&\rat(8/5){4}{}&\int(0){1}{}&\int(0){1}{}&\rat(3/5){4}{}&\ln{1}&\cnt{50}&\bounds{10}{10}{1}{false}&&\cr
\endorbit
\orbit{6}{15}{104}
\h{15}&\rat(4/5){4}{}&\rat(4/5){4}{}&\rat(16/5){3}{\spec{1}}&\int(0){\0}{}&\int(0){\0}{}&\root\rat(6/5){3}{\spec{1}}&&\cnt{492}&\bounds{104}{112}{fail}{false}&&\cr
\l{1}&\rat(2/5){2}{}&\rat(2/5){2}{}&\rat(8/5){4}{}&\int(0){\0}{}&\int(0){\0}{}&\rat(3/5){4}{}&\ln{1\sdd}&\cnt{72}&\bounds{16}{24}{2}{false}&&\cr
\l{2}&\rat(4/5){4}{}&\int(0){\0}{}&\rat(8/5){4}{}&\int(0){1}{}&\int(0){1}{}&\rat(3/5){4}{}&\ln{2\sd}&\cnt{50}&\bounds{10}{10}{1}{false}&&\cr
\l{3}&\int(1){\2}{}&\int(0){\0}{}&\int(2){\2}{}&\int(0){\0}{}&\int(0){\0}{}&\int(0){\0}{}&\ln{2}&\cnt{16}&\bounds{4}{4}{1}{false}&&\cr
\l{5}&\rat(3/5){3}{}&\rat(4/5){4}{}&\rat(8/5){4}{}&\int(0){\0}{}&\int(0){2}{}&\int(0){\0}{}&\ln{2}&\cnt{40}&\bounds{8}{8}{1}{false}&&\cr
\l{7}&\rat(3/5){3}{}&\int(0){\0}{}&\rat(8/5){4}{}&\int(0){4}{}&\int(0){\0}{}&\rat(4/5){2}{}&\ln{2}&\cnt{20}&\bounds{4}{4}{1}{false}&&\cr
\l{9}&\rat(4/5){4}{}&\rat(4/5){4}{}&\rat(6/5){3}{}&\int(0){\0}{}&\int(0){\0}{}&\rat(1/5){3}{}&\ln{1}&\cnt{8}&\bounds{2}{2}{1}{false}&&\cr
\endorbit
\orbit{6}{16}{96}
\h{16}&\rat(4/5){4}{}&\rat(4/5){4}{}&\rat(16/5){3}{\spec{1}}&\int(0){\0}{}&\int(0){\0}{}&\root\rat(6/5){3}{\spec{2}}&&\cnt{464}&\bounds{96}{96}{1}{false}&&\cr
\l{1}&\int(1){\2}{}&\int(0){\0}{}&\int(2){\2}{}&\int(0){\0}{}&\int(0){\0}{}&\int(0){\0}{}&\ln{2}&\cnt{16}&\bounds{4}{4}{1}{false}&&\cr
\l{3}&\rat(3/5){3}{}&\rat(4/5){4}{}&\rat(8/5){4}{}&\int(0){\0}{}&\int(0){2}{}&\int(0){\0}{}&\ln{2}&\cnt{40}&\bounds{8}{8}{1}{false}&&\cr
\l{5}&\rat(3/5){3}{}&\int(0){\0}{}&\rat(8/5){4}{}&\int(0){4}{}&\int(0){\0}{}&\rat(4/5){2}{}&\ln{2}&\cnt{60}&\bounds{12}{12}{1}{false}&&\cr
\endorbit
\orbit{6}{17}{116}
\h{17}&\rat(14/5){4}{\spec{2}}&\rat(4/5){4}{}&\rat(4/5){1}{}&\rat(4/5){1}{}&\rat(4/5){4}{}&\root\int(0){\0}{}&&\cnt{500}&\bounds{116}{128}{fail}{false}&&\cr
\l{1}&\rat(7/5){2}{}&\rat(3/5){3}{}&\int(0){\0}{}&\rat(4/5){1}{}&\rat(1/5){1}{}&\int(0){\0}{}&\ln{4\sd}&\cnt{32}&\bounds{7}{8}{5}{true}&&\cr
\l{2}&\rat(7/5){2}{}&\rat(4/5){4}{}&\int(0){\0}{}&\int(0){\0}{}&\rat(4/5){4}{}&\int(0){2}{\spec{1}}&\ln{2\sd}&\cnt{2}&\bounds{1}{1}{1}{true}&&\cr
\l{3}&\rat(7/5){2}{}&\rat(4/5){4}{}&\int(0){\0}{}&\int(0){\0}{}&\rat(4/5){4}{}&\int(0){2}{\spec{2}}&\ln{2\sd}&\cnt{6}&\bounds{2}{2}{1}{true}&&\cr
\l{4}&\int(3){\4}{}&\int(0){\0}{}&\int(0){\0}{}&\int(0){\0}{}&\int(0){\0}{}&\int(0){\0}{}&\ln{1}&\cnt{2}&\bounds{1}{1}{1}{true}&&\cr
\l{6}&\int(2){\2}{}&\int(1){\2}{}&\int(0){\0}{}&\int(0){\0}{}&\int(0){\0}{}&\int(0){\0}{}&\ln{4}&\cnt{8}&\bounds{2}{2}{1}{true}&&\cr
\l{8}&\rat(8/5){3}{}&\rat(3/5){3}{}&\rat(4/5){1}{}&\int(0){\0}{}&\int(0){\0}{}&\int(0){1}{}&\ln{4}&\cnt{12}&\bounds{3}{3}{1}{true}&&\cr
\l{10}&\rat(8/5){3}{}&\rat(2/5){2}{}&\int(0){\0}{}&\rat(1/5){4}{}&\rat(4/5){4}{}&\int(0){\0}{}&\ln{4}&\cnt{24}&\bounds{5}{6}{5}{true}&&\cr
\endorbit
\orbit{6}{18}{110}
\h{18}&\root\rat(14/5){4}{\spec{2}}&\rat(4/5){4}{}&\rat(4/5){1}{}&\rat(4/5){1}{}&\rat(4/5){4}{}&\int(0){\0}{}&&\cnt{524}&\bounds{110}{122}{fail}{false}&&\cr
\l{1}&\rat(7/5){2}{}&\rat(3/5){3}{}&\int(0){\0}{}&\rat(4/5){1}{}&\rat(1/5){1}{}&\int(0){\0}{}&\ln{4\sd}&\cnt{32}&\bounds{7}{8}{5}{true}&&\cr
\l{2}&\rat(7/5){2}{}&\rat(4/5){4}{}&\int(0){\0}{}&\int(0){\0}{}&\rat(4/5){4}{}&\int(0){2}{}&\ln{2\sd}&\cnt{20}&\bounds{4}{4}{1}{true}&&\cr
\l{3}&\int(3){\4}{}&\int(0){\0}{}&\int(0){\0}{}&\int(0){\0}{}&\int(0){\0}{}&\int(0){\0}{}&\ln{1}&\cnt{2}&\bounds{1}{1}{1}{true}&&\cr
\l{5}&\rat(8/5){3}{}&\rat(3/5){3}{}&\rat(4/5){1}{}&\int(0){\0}{}&\int(0){\0}{}&\int(0){1}{}&\ln{4}&\cnt{20}&\bounds{4}{4}{1}{true}&&\cr
\l{7}&\rat(8/5){3}{}&\rat(2/5){2}{}&\int(0){\0}{}&\rat(1/5){4}{}&\rat(4/5){4}{}&\int(0){\0}{}&\ln{4}&\cnt{24}&\bounds{5}{6}{5}{true}&&\cr
\endorbit
\orbit{6}{19}{114}
\h{19}&\rat(14/5){4}{\spec{2}}&\root\rat(4/5){4}{}&\rat(4/5){1}{}&\rat(4/5){1}{}&\rat(4/5){4}{}&\int(0){\0}{}&&\cnt{472}&\bounds{114}{114}{1}{false}&&\cr
\l{1}&\int(3){\4}{}&\int(0){\0}{}&\int(0){\0}{}&\int(0){\0}{}&\int(0){\0}{}&\int(0){\0}{}&\ln{1}&\cnt{2}&\bounds{1}{1}{1}{false}&&\cr
\l{3}&\int(2){\2}{}&\int(1){\2}{}&\int(0){\0}{}&\int(0){\0}{}&\int(0){\0}{}&\int(0){\0}{}&\ln{1}&\cnt{4}&\bounds{2}{2}{1}{false}&&\cr
\l{5}&\int(2){\2}{}&\int(0){\0}{}&\int(1){\2}{}&\int(0){\0}{}&\int(0){\0}{}&\int(0){\0}{}&\ln{1}&\cnt{8}&\bounds{2}{2}{1}{false}&&\cr
\l{7}&\int(2){\2}{}&\int(0){\0}{}&\int(0){\0}{}&\int(1){\2}{}&\int(0){\0}{}&\int(0){\0}{}&\ln{1}&\cnt{8}&\bounds{2}{2}{1}{false}&&\cr
\l{9}&\int(2){\2}{}&\int(0){\0}{}&\int(0){\0}{}&\int(0){\0}{}&\int(1){\2}{}&\int(0){\0}{}&\ln{1}&\cnt{8}&\bounds{2}{2}{1}{false}&&\cr
\l{11}&\rat(8/5){3}{}&\rat(3/5){3}{}&\rat(4/5){1}{}&\int(0){\0}{}&\int(0){\0}{}&\int(0){1}{}&\ln{1}&\cnt{10}&\bounds{2}{2}{1}{false}&&\cr
\l{13}&\rat(8/5){3}{}&\rat(2/5){2}{\spec{1}}&\int(0){\0}{}&\rat(1/5){4}{}&\rat(4/5){4}{}&\int(0){\0}{}&\ln{1}&\cnt{4}&\bounds{1}{1}{1}{false}&&\cr
\l{15}&\rat(8/5){3}{}&\rat(2/5){2}{\spec{2}}&\int(0){\0}{}&\rat(1/5){4}{}&\rat(4/5){4}{}&\int(0){\0}{}&\ln{1}&\cnt{4}&\bounds{1}{1}{1}{false}&&\cr
\l{17}&\rat(8/5){3}{}&\rat(4/5){4}{}&\rat(1/5){4}{}&\int(0){\0}{}&\rat(2/5){2}{}&\int(0){\0}{}&\ln{1}&\cnt{24}&\bounds{6}{6}{1}{false}&&\cr
\l{19}&\rat(8/5){3}{}&\rat(4/5){4}{}&\int(0){\0}{}&\rat(3/5){2}{}&\int(0){\0}{}&\int(0){4}{}&\ln{1}&\cnt{20}&\bounds{4}{4}{1}{false}&&\cr
\l{21}&\rat(8/5){3}{}&\rat(1/5){1}{}&\rat(2/5){3}{}&\rat(4/5){1}{}&\int(0){\0}{}&\int(0){\0}{}&\ln{1}&\cnt{12}&\bounds{4}{4}{1}{false}&&\cr
\l{23}&\rat(8/5){3}{}&\int(0){\0}{}&\rat(3/5){2}{}&\int(0){\0}{}&\rat(4/5){4}{}&\int(0){4}{}&\ln{1}&\cnt{20}&\bounds{4}{4}{1}{false}&&\cr
\l{25}&\rat(8/5){3}{}&\int(0){\0}{}&\rat(4/5){1}{}&\rat(2/5){3}{}&\rat(1/5){1}{}&\int(0){\0}{}&\ln{1}&\cnt{24}&\bounds{6}{6}{1}{false}&&\cr
\l{27}&\rat(8/5){3}{}&\int(0){\0}{}&\int(0){\0}{}&\rat(4/5){1}{}&\rat(3/5){3}{}&\int(0){1}{}&\ln{1}&\cnt{20}&\bounds{4}{4}{1}{false}&&\cr
\l{29}&\rat(7/5){2}{}&\rat(3/5){3}{}&\int(0){\0}{}&\rat(4/5){1}{}&\rat(1/5){1}{}&\int(0){\0}{}&\ln{1}&\cnt{16}&\bounds{4}{4}{1}{false}&&\cr
\l{31}&\rat(7/5){2}{}&\rat(4/5){4}{}&\rat(3/5){2}{}&\rat(1/5){4}{}&\int(0){\0}{}&\int(0){\0}{}&\ln{1}&\cnt{32}&\bounds{8}{8}{1}{false}&&\cr
\l{33}&\rat(7/5){2}{}&\rat(4/5){4}{}&\int(0){\0}{}&\int(0){\0}{}&\rat(4/5){4}{}&\int(0){2}{}&\ln{1}&\cnt{20}&\bounds{4}{4}{1}{false}&&\cr
\endorbit
\orbit{6}{20}{96}
\h{20}&\root\rat(14/5){4}{\spec{1}}&\rat(4/5){4}{}&\rat(4/5){1}{}&\rat(4/5){1}{}&\rat(4/5){4}{}&\int(0){\0}{}&&\cnt{416}&\bounds{96}{96}{1}{false}&&\cr
\l{1}&\int(2){\2}{}&\int(1){\2}{}&\int(0){\0}{}&\int(0){\0}{}&\int(0){\0}{}&\int(0){\0}{}&\ln{4}&\cnt{8}&\bounds{2}{2}{1}{false}&&\cr
\l{3}&\rat(8/5){3}{}&\rat(3/5){3}{}&\rat(4/5){1}{}&\int(0){\0}{}&\int(0){\0}{}&\int(0){1}{}&\ln{4}&\cnt{20}&\bounds{4}{4}{1}{false}&&\cr
\l{5}&\rat(8/5){3}{}&\rat(2/5){2}{}&\int(0){\0}{}&\rat(1/5){4}{}&\rat(4/5){4}{}&\int(0){\0}{}&\ln{4}&\cnt{24}&\bounds{6}{6}{1}{false}&&\cr
\endorbit
\orbit{6}{21}{126}
\h{21}&\rat(14/5){4}{\spec{2}}&\rat(4/5){4}{}&\rat(6/5){3}{}&\root\int(0){\0}{}&\int(0){\0}{}&\rat(6/5){3}{}&&\cnt{500}&\bounds{126}{132}{fail}{false}&&\cr
\l{1}&\rat(7/5){2}{}&\rat(2/5){2}{}&\rat(3/5){4}{}&\int(0){\0}{}&\int(0){\0}{}&\rat(3/5){4}{}&\ln{1\sdd}&\cnt{48}&\bounds{10}{16}{5}{true}&&\cr
\l{2}&\int(3){\4}{}&\int(0){\0}{}&\int(0){\0}{}&\int(0){\0}{}&\int(0){\0}{}&\int(0){\0}{}&\ln{1}&\cnt{2}&\bounds{1}{1}{1}{true}&&\cr
\l{4}&\int(2){\2}{}&\int(1){\2}{}&\int(0){\0}{}&\int(0){\0}{}&\int(0){\0}{}&\int(0){\0}{}&\ln{1}&\cnt{8}&\bounds{2}{2}{1}{true}&&\cr
\l{6}&\int(2){\2}{}&\int(0){\0}{}&\int(1){\2}{}&\int(0){\0}{}&\int(0){\0}{}&\int(0){\0}{}&\ln{1}&\cnt{12}&\bounds{4}{4}{1}{true}&&\cr
\l{8}&\int(2){\2}{}&\int(0){\0}{}&\int(0){\0}{}&\int(0){\0}{}&\int(0){\0}{}&\int(1){\2}{}&\ln{1}&\cnt{12}&\bounds{4}{4}{1}{true}&&\cr
\l{10}&\rat(8/5){3}{}&\rat(3/5){3}{}&\rat(2/5){1}{}&\int(0){\0}{}&\int(0){\0}{}&\rat(2/5){1}{}&\ln{1}&\cnt{36}&\bounds{9}{9}{1}{true}&&\cr
\l{12}&\rat(8/5){3}{}&\rat(4/5){4}{}&\rat(3/5){4}{}&\int(0){\0}{}&\int(0){2}{}&\int(0){\0}{}&\ln{1}&\cnt{20}&\bounds{4}{4}{1}{true}&&\cr
\l{14}&\rat(8/5){3}{}&\rat(4/5){4}{}&\int(0){\0}{}&\int(0){2}{\spec{1}}&\int(0){\0}{}&\rat(3/5){4}{}&\ln{1}&\cnt{2}&\bounds{1}{1}{1}{true}&&\cr
\l{16}&\rat(8/5){3}{}&\rat(4/5){4}{}&\int(0){\0}{}&\int(0){2}{\spec{2}}&\int(0){\0}{}&\rat(3/5){4}{}&\ln{1}&\cnt{6}&\bounds{2}{2}{1}{true}&&\cr
\l{18}&\rat(8/5){3}{}&\rat(1/5){1}{}&\rat(6/5){3}{}&\int(0){1}{}&\int(0){\0}{}&\int(0){\0}{}&\ln{1}&\cnt{12}&\bounds{3}{3}{1}{true}&&\cr
\l{20}&\rat(8/5){3}{}&\rat(1/5){1}{}&\int(0){\0}{}&\int(0){\0}{}&\int(0){1}{}&\rat(6/5){3}{}&\ln{1}&\cnt{20}&\bounds{4}{4}{1}{true}&&\cr
\l{22}&\rat(8/5){3}{}&\int(0){\0}{}&\rat(4/5){2}{}&\int(0){\0}{}&\int(0){4}{}&\rat(3/5){4}{}&\ln{1}&\cnt{30}&\bounds{6}{6}{1}{true}&&\cr
\l{24}&\rat(8/5){3}{}&\int(0){\0}{}&\rat(3/5){4}{}&\int(0){4}{}&\int(0){\0}{}&\rat(4/5){2}{}&\ln{1}&\cnt{18}&\bounds{6}{6}{1}{true}&&\cr
\l{26}&\rat(7/5){2}{}&\rat(4/5){4}{}&\rat(4/5){2}{}&\int(0){4}{}&\int(0){\0}{}&\int(0){\0}{}&\ln{1}&\cnt{18}&\bounds{6}{6}{1}{true}&&\cr
\l{28}&\rat(7/5){2}{}&\rat(4/5){4}{}&\int(0){\0}{}&\int(0){\0}{}&\int(0){4}{}&\rat(4/5){2}{}&\ln{1}&\cnt{30}&\bounds{6}{6}{1}{true}&&\cr
\endorbit
\orbit{6}{22}{114}
\h{22}&\rat(14/5){4}{\spec{2}}&\root\rat(4/5){4}{}&\rat(6/5){3}{}&\int(0){\0}{}&\int(0){\0}{}&\rat(6/5){3}{}&&\cnt{472}&\bounds{114}{114}{1}{false}&&\cr
\l{1}&\int(3){\4}{}&\int(0){\0}{}&\int(0){\0}{}&\int(0){\0}{}&\int(0){\0}{}&\int(0){\0}{}&\ln{1}&\cnt{2}&\bounds{1}{1}{1}{false}&&\cr
\l{3}&\int(2){\2}{}&\int(1){\2}{}&\int(0){\0}{}&\int(0){\0}{}&\int(0){\0}{}&\int(0){\0}{}&\ln{1}&\cnt{4}&\bounds{2}{2}{1}{false}&&\cr
\l{5}&\int(2){\2}{}&\int(0){\0}{}&\int(1){\2}{}&\int(0){\0}{}&\int(0){\0}{}&\int(0){\0}{}&\ln{2}&\cnt{12}&\bounds{4}{4}{1}{false}&&\cr
\l{7}&\rat(8/5){3}{}&\rat(3/5){3}{}&\rat(2/5){1}{}&\int(0){\0}{}&\int(0){\0}{}&\rat(2/5){1}{}&\ln{1}&\cnt{18}&\bounds{6}{6}{1}{false}&&\cr
\l{9}&\rat(8/5){3}{}&\rat(4/5){4}{}&\rat(3/5){4}{}&\int(0){\0}{}&\int(0){2}{}&\int(0){\0}{}&\ln{2}&\cnt{20}&\bounds{4}{4}{1}{false}&&\cr
\l{11}&\rat(8/5){3}{}&\rat(1/5){1}{}&\rat(6/5){3}{}&\int(0){1}{}&\int(0){\0}{}&\int(0){\0}{}&\ln{2}&\cnt{10}&\bounds{2}{2}{1}{false}&&\cr
\l{13}&\rat(8/5){3}{}&\int(0){\0}{}&\rat(4/5){2}{}&\int(0){\0}{}&\int(0){4}{}&\rat(3/5){4}{}&\ln{2}&\cnt{30}&\bounds{6}{6}{1}{false}&&\cr
\l{15}&\rat(7/5){2}{}&\rat(2/5){2}{\spec{1}}&\rat(3/5){4}{}&\int(0){\0}{}&\int(0){\0}{}&\rat(3/5){4}{}&\ln{1}&\cnt{8}&\bounds{4}{4}{1}{false}&&\cr
\l{17}&\rat(7/5){2}{}&\rat(4/5){4}{}&\rat(4/5){2}{}&\int(0){4}{}&\int(0){\0}{}&\int(0){\0}{}&\ln{2}&\cnt{30}&\bounds{6}{6}{1}{false}&&\cr
\endorbit
\orbit{6}{23}{112}
\h{23}&\root\rat(14/5){4}{\spec{2}}&\rat(4/5){4}{}&\rat(6/5){3}{}&\int(0){\0}{}&\int(0){\0}{}&\rat(6/5){3}{}&&\cnt{524}&\bounds{112}{116}{fail}{false}&&\cr
\l{1}&\rat(7/5){2}{}&\rat(2/5){2}{}&\rat(3/5){4}{}&\int(0){\0}{}&\int(0){\0}{}&\rat(3/5){4}{}&\ln{1\sdd}&\cnt{48}&\bounds{12}{16}{2}{false}&&\cr
\l{2}&\int(3){\4}{}&\int(0){\0}{}&\int(0){\0}{}&\int(0){\0}{}&\int(0){\0}{}&\int(0){\0}{}&\ln{1}&\cnt{2}&\bounds{1}{1}{1}{false}&&\cr
\l{4}&\rat(8/5){3}{}&\rat(3/5){3}{}&\rat(2/5){1}{}&\int(0){\0}{}&\int(0){\0}{}&\rat(2/5){1}{}&\ln{1}&\cnt{36}&\bounds{9}{9}{1}{false}&&\cr
\l{6}&\rat(8/5){3}{}&\rat(4/5){4}{}&\rat(3/5){4}{}&\int(0){\0}{}&\int(0){2}{}&\int(0){\0}{}&\ln{2}&\cnt{20}&\bounds{4}{4}{1}{false}&&\cr
\l{8}&\rat(8/5){3}{}&\rat(1/5){1}{}&\rat(6/5){3}{}&\int(0){1}{}&\int(0){\0}{}&\int(0){\0}{}&\ln{2}&\cnt{20}&\bounds{4}{4}{1}{false}&&\cr
\l{10}&\rat(8/5){3}{}&\int(0){\0}{}&\rat(4/5){2}{}&\int(0){\0}{}&\int(0){4}{}&\rat(3/5){4}{}&\ln{2}&\cnt{30}&\bounds{6}{6}{1}{false}&&\cr
\l{12}&\rat(7/5){2}{}&\rat(4/5){4}{}&\rat(4/5){2}{}&\int(0){4}{}&\int(0){\0}{}&\int(0){\0}{}&\ln{2}&\cnt{30}&\bounds{6}{6}{1}{false}&&\cr
\endorbit
\orbit{6}{24}{100}
\h{24}&\rat(14/5){4}{\spec{2}}&\rat(4/5){4}{}&\root\rat(6/5){3}{\spec{1}}&\int(0){\0}{}&\int(0){\0}{}&\rat(6/5){3}{}&&\cnt{444}&\bounds{100}{104}{fail}{false}&&\cr
\l{1}&\rat(7/5){2}{}&\rat(2/5){2}{}&\rat(3/5){4}{}&\int(0){\0}{}&\int(0){\0}{}&\rat(3/5){4}{}&\ln{1\sdd}&\cnt{48}&\bounds{12}{16}{2}{false}&&\cr
\l{2}&\int(3){\4}{}&\int(0){\0}{}&\int(0){\0}{}&\int(0){\0}{}&\int(0){\0}{}&\int(0){\0}{}&\ln{1}&\cnt{2}&\bounds{1}{1}{1}{false}&&\cr
\l{4}&\int(2){\2}{}&\int(1){\2}{}&\int(0){\0}{}&\int(0){\0}{}&\int(0){\0}{}&\int(0){\0}{}&\ln{1}&\cnt{8}&\bounds{2}{2}{1}{false}&&\cr
\l{6}&\int(2){\2}{}&\int(0){\0}{}&\int(1){\2}{}&\int(0){\0}{}&\int(0){\0}{}&\int(0){\0}{}&\ln{1}&\cnt{4}&\bounds{2}{2}{1}{false}&&\cr
\l{8}&\int(2){\2}{}&\int(0){\0}{}&\int(0){\0}{}&\int(0){\0}{}&\int(0){\0}{}&\int(1){\2}{}&\ln{1}&\cnt{12}&\bounds{4}{4}{1}{false}&&\cr
\l{10}&\rat(8/5){3}{}&\rat(3/5){3}{}&\rat(2/5){1}{}&\int(0){\0}{}&\int(0){\0}{}&\rat(2/5){1}{}&\ln{1}&\cnt{12}&\bounds{3}{3}{1}{false}&&\cr
\l{12}&\rat(8/5){3}{}&\rat(4/5){4}{}&\rat(3/5){4}{}&\int(0){\0}{}&\int(0){2}{}&\int(0){\0}{}&\ln{1}&\cnt{20}&\bounds{4}{4}{1}{false}&&\cr
\l{14}&\rat(8/5){3}{}&\rat(4/5){4}{}&\int(0){\0}{}&\int(0){2}{}&\int(0){\0}{}&\rat(3/5){4}{}&\ln{1}&\cnt{20}&\bounds{4}{4}{1}{false}&&\cr
\l{16}&\rat(8/5){3}{}&\rat(1/5){1}{}&\rat(6/5){3}{}&\int(0){1}{}&\int(0){\0}{}&\int(0){\0}{}&\ln{1}&\cnt{20}&\bounds{4}{4}{1}{false}&&\cr
\l{18}&\rat(8/5){3}{}&\rat(1/5){1}{}&\int(0){\0}{}&\int(0){\0}{}&\int(0){1}{}&\rat(6/5){3}{}&\ln{1}&\cnt{20}&\bounds{4}{4}{1}{false}&&\cr
\l{20}&\rat(8/5){3}{}&\int(0){\0}{}&\rat(4/5){2}{}&\int(0){\0}{}&\int(0){4}{}&\rat(3/5){4}{}&\ln{1}&\cnt{10}&\bounds{2}{2}{1}{false}&&\cr
\l{22}&\rat(8/5){3}{}&\int(0){\0}{}&\rat(3/5){4}{}&\int(0){4}{}&\int(0){\0}{}&\rat(4/5){2}{}&\ln{1}&\cnt{30}&\bounds{6}{6}{1}{false}&&\cr
\l{24}&\rat(7/5){2}{}&\rat(4/5){4}{}&\rat(4/5){2}{}&\int(0){4}{}&\int(0){\0}{}&\int(0){\0}{}&\ln{1}&\cnt{10}&\bounds{2}{2}{1}{false}&&\cr
\l{26}&\rat(7/5){2}{}&\rat(4/5){4}{}&\int(0){\0}{}&\int(0){\0}{}&\int(0){4}{}&\rat(4/5){2}{}&\ln{1}&\cnt{30}&\bounds{6}{6}{1}{false}&&\cr
\endorbit
\orbit{6}{25}{94}
\h{25}&\root\rat(14/5){4}{\spec{1}}&\rat(4/5){4}{}&\rat(6/5){3}{}&\int(0){\0}{}&\int(0){\0}{}&\rat(6/5){3}{}&&\cnt{416}&\bounds{94}{94}{1}{false}&&\cr
\l{1}&\int(2){\2}{}&\int(1){\2}{}&\int(0){\0}{}&\int(0){\0}{}&\int(0){\0}{}&\int(0){\0}{}&\ln{1}&\cnt{8}&\bounds{2}{2}{1}{false}&&\cr
\l{3}&\int(2){\2}{}&\int(0){\0}{}&\int(1){\2}{}&\int(0){\0}{}&\int(0){\0}{}&\int(0){\0}{}&\ln{2}&\cnt{12}&\bounds{4}{4}{1}{false}&&\cr
\l{5}&\rat(8/5){3}{}&\rat(3/5){3}{}&\rat(2/5){1}{}&\int(0){\0}{}&\int(0){\0}{}&\rat(2/5){1}{}&\ln{1}&\cnt{36}&\bounds{9}{9}{1}{false}&&\cr
\l{7}&\rat(8/5){3}{}&\rat(4/5){4}{}&\rat(3/5){4}{}&\int(0){\0}{}&\int(0){2}{}&\int(0){\0}{}&\ln{2}&\cnt{20}&\bounds{4}{4}{1}{false}&&\cr
\l{9}&\rat(8/5){3}{}&\rat(1/5){1}{}&\rat(6/5){3}{}&\int(0){1}{}&\int(0){\0}{}&\int(0){\0}{}&\ln{2}&\cnt{20}&\bounds{4}{4}{1}{false}&&\cr
\l{11}&\rat(8/5){3}{}&\int(0){\0}{}&\rat(4/5){2}{}&\int(0){\0}{}&\int(0){4}{}&\rat(3/5){4}{}&\ln{2}&\cnt{30}&\bounds{6}{6}{1}{false}&&\cr
\endorbit
\orbit{6}{26}{92}
\h{26}&\rat(14/5){4}{\spec{2}}&\rat(4/5){4}{}&\root\rat(6/5){3}{\spec{2}}&\int(0){\0}{}&\int(0){\0}{}&\rat(6/5){3}{}&&\cnt{416}&\bounds{92}{92}{1}{false}&&\cr
\l{1}&\int(3){\4}{}&\int(0){\0}{}&\int(0){\0}{}&\int(0){\0}{}&\int(0){\0}{}&\int(0){\0}{}&\ln{1}&\cnt{2}&\bounds{1}{1}{1}{false}&&\cr
\l{3}&\int(2){\2}{}&\int(1){\2}{}&\int(0){\0}{}&\int(0){\0}{}&\int(0){\0}{}&\int(0){\0}{}&\ln{1}&\cnt{8}&\bounds{2}{2}{1}{false}&&\cr
\l{5}&\int(2){\2}{}&\int(0){\0}{}&\int(0){\0}{}&\int(0){\0}{}&\int(0){\0}{}&\int(1){\2}{}&\ln{1}&\cnt{12}&\bounds{4}{4}{1}{false}&&\cr
\l{7}&\rat(8/5){3}{}&\rat(3/5){3}{}&\rat(2/5){1}{}&\int(0){\0}{}&\int(0){\0}{}&\rat(2/5){1}{}&\ln{1}&\cnt{36}&\bounds{9}{9}{1}{false}&&\cr
\l{9}&\rat(8/5){3}{}&\rat(4/5){4}{}&\int(0){\0}{}&\int(0){2}{}&\int(0){\0}{}&\rat(3/5){4}{}&\ln{1}&\cnt{20}&\bounds{4}{4}{1}{false}&&\cr
\l{11}&\rat(8/5){3}{}&\rat(1/5){1}{}&\rat(6/5){3}{}&\int(0){1}{}&\int(0){\0}{}&\int(0){\0}{}&\ln{1}&\cnt{20}&\bounds{4}{4}{1}{false}&&\cr
\l{13}&\rat(8/5){3}{}&\rat(1/5){1}{}&\int(0){\0}{}&\int(0){\0}{}&\int(0){1}{}&\rat(6/5){3}{}&\ln{1}&\cnt{20}&\bounds{4}{4}{1}{false}&&\cr
\l{15}&\rat(8/5){3}{}&\int(0){\0}{}&\rat(4/5){2}{}&\int(0){\0}{}&\int(0){4}{}&\rat(3/5){4}{}&\ln{1}&\cnt{30}&\bounds{6}{6}{1}{false}&&\cr
\l{17}&\rat(7/5){2}{}&\rat(4/5){4}{}&\rat(4/5){2}{}&\int(0){4}{}&\int(0){\0}{}&\int(0){\0}{}&\ln{1}&\cnt{30}&\bounds{6}{6}{1}{false}&&\cr
\l{19}&\rat(7/5){2}{}&\rat(4/5){4}{}&\int(0){\0}{}&\int(0){\0}{}&\int(0){4}{}&\rat(4/5){2}{}&\ln{1}&\cnt{30}&\bounds{6}{6}{1}{false}&&\cr
\endorbit
\orbit{6}{27}{122}
\h{27}&\root\rat(4/5){4}{}&\rat(4/5){4}{}&\rat(4/5){4}{}&\rat(6/5){2}{}&\rat(6/5){3}{}&\rat(6/5){2}{}&&\cnt{424}&\bounds{122}{128}{fail}{false}&&\cr
\l{1}&\rat(3/5){3}{}&\rat(4/5){4}{}&\rat(4/5){4}{}&\int(0){\0}{}&\rat(4/5){2}{}&\int(0){\0}{}&\ln{1}&\cnt{6}&\bounds{2}{2}{1}{true}&&\cr
\l{3}&\rat(3/5){3}{}&\rat(4/5){4}{}&\int(0){\0}{}&\rat(6/5){2}{}&\int(0){\0}{}&\rat(2/5){4}{}&\ln{2}&\cnt{6}&\bounds{2}{2}{1}{true}&&\cr
\l{5}&\rat(3/5){3}{}&\int(0){\0}{}&\int(0){\0}{}&\rat(3/5){1}{}&\rat(6/5){3}{}&\rat(3/5){1}{}&\ln{1}&\cnt{8}&\bounds{3}{4}{5}{true}&&\cr
\l{7}&\rat(2/5){2}{\spec{1}}&\rat(4/5){4}{}&\int(0){\0}{}&\int(0){\0}{}&\rat(3/5){4}{}&\rat(6/5){2}{}&\ln{2}&\cnt{2}&\bounds{1}{1}{1}{true}&&\cr
\l{9}&\rat(4/5){4}{}&\rat(3/5){3}{}&\rat(4/5){4}{}&\rat(4/5){3}{}&\int(0){\0}{}&\int(0){\0}{}&\ln{2}&\cnt{12}&\bounds{3}{3}{1}{true}&&\cr
\l{11}&\rat(4/5){4}{}&\rat(3/5){3}{}&\int(0){\0}{}&\int(0){\0}{}&\rat(6/5){3}{}&\rat(2/5){4}{}&\ln{2}&\cnt{12}&\bounds{3}{3}{1}{true}&&\cr
\l{13}&\rat(4/5){4}{}&\rat(2/5){2}{}&\int(0){\0}{}&\rat(3/5){1}{}&\int(0){\0}{}&\rat(6/5){2}{}&\ln{2}&\cnt{12}&\bounds{3}{4}{5}{true}&&\cr
\l{15}&\rat(4/5){4}{}&\rat(4/5){4}{}&\rat(1/5){1}{}&\rat(3/5){1}{}&\rat(3/5){4}{}&\int(0){\0}{}&\ln{2}&\cnt{16}&\bounds{4}{4}{1}{true}&&\cr
\l{17}&\rat(4/5){4}{}&\rat(4/5){4}{}&\int(0){\0}{}&\rat(2/5){4}{}&\rat(2/5){1}{}&\rat(3/5){1}{}&\ln{2}&\cnt{18}&\bounds{6}{6}{1}{true}&&\cr
\l{19}&\rat(4/5){4}{}&\rat(1/5){1}{}&\int(0){\0}{}&\rat(6/5){2}{}&\rat(4/5){2}{}&\int(0){\0}{}&\ln{2}&\cnt{12}&\bounds{3}{3}{1}{true}&&\cr
\l{21}&\rat(4/5){4}{}&\int(0){\0}{}&\int(0){\0}{}&\rat(4/5){3}{}&\rat(3/5){4}{}&\rat(4/5){3}{}&\ln{1}&\cnt{18}&\bounds{6}{6}{1}{true}&&\cr
\endorbit
\orbit{6}{28}{110}
\h{28}&\rat(4/5){4}{}&\rat(4/5){4}{}&\rat(4/5){4}{}&\root\rat(6/5){2}{\spec{2}}&\rat(6/5){3}{}&\rat(6/5){2}{}&&\cnt{396}&\bounds{110}{110}{1}{false}&&\cr
\l{1}&\rat(4/5){4}{}&\rat(2/5){2}{}&\int(0){\0}{}&\rat(3/5){1}{}&\int(0){\0}{}&\rat(6/5){2}{}&\ln{2\sd}&\cnt{12}&\bounds{4}{4}{1}{false}&&\cr
\l{2}&\rat(3/5){3}{}&\rat(4/5){4}{}&\rat(4/5){4}{}&\int(0){\0}{}&\rat(4/5){2}{}&\int(0){\0}{}&\ln{2}&\cnt{12}&\bounds{3}{3}{1}{false}&&\cr
\l{4}&\rat(3/5){3}{}&\rat(4/5){4}{}&\int(0){\0}{}&\rat(6/5){2}{}&\int(0){\0}{}&\rat(2/5){4}{}&\ln{2}&\cnt{12}&\bounds{3}{3}{1}{false}&&\cr
\l{6}&\rat(3/5){3}{}&\int(0){\0}{}&\rat(4/5){4}{}&\rat(2/5){4}{}&\int(0){\0}{}&\rat(6/5){2}{}&\ln{2}&\cnt{4}&\bounds{1}{1}{1}{false}&&\cr
\l{8}&\rat(3/5){3}{}&\int(0){\0}{}&\int(0){\0}{}&\rat(3/5){1}{}&\rat(6/5){3}{}&\rat(3/5){1}{}&\ln{2}&\cnt{16}&\bounds{4}{4}{1}{false}&&\cr
\l{10}&\rat(2/5){2}{}&\rat(4/5){4}{}&\int(0){\0}{}&\int(0){\0}{}&\rat(3/5){4}{}&\rat(6/5){2}{}&\ln{2}&\cnt{12}&\bounds{4}{4}{1}{false}&&\cr
\l{12}&\rat(4/5){4}{}&\rat(3/5){3}{}&\rat(4/5){4}{}&\rat(4/5){3}{}&\int(0){\0}{}&\int(0){\0}{}&\ln{1}&\cnt{4}&\bounds{1}{1}{1}{false}&&\cr
\l{14}&\rat(4/5){4}{}&\rat(3/5){3}{}&\int(0){\0}{}&\int(0){\0}{}&\rat(6/5){3}{}&\rat(2/5){4}{}&\ln{2}&\cnt{12}&\bounds{3}{3}{1}{false}&&\cr
\l{16}&\rat(4/5){4}{}&\rat(4/5){4}{}&\int(0){\0}{}&\rat(2/5){4}{}&\rat(2/5){1}{}&\rat(3/5){1}{}&\ln{2}&\cnt{6}&\bounds{2}{2}{1}{false}&&\cr
\l{18}&\rat(4/5){4}{}&\rat(1/5){1}{}&\rat(4/5){4}{}&\int(0){\0}{}&\rat(3/5){4}{}&\rat(3/5){1}{}&\ln{1}&\cnt{16}&\bounds{4}{4}{1}{false}&&\cr
\l{20}&\rat(4/5){4}{}&\int(0){\0}{}&\rat(4/5){4}{}&\rat(3/5){1}{}&\rat(2/5){1}{}&\rat(2/5){4}{}&\ln{1}&\cnt{18}&\bounds{6}{6}{1}{false}&&\cr
\endorbit
\orbit{6}{29}{102}
\h{29}&\rat(4/5){4}{}&\rat(4/5){4}{}&\rat(4/5){4}{}&\root\rat(6/5){2}{\spec{1}}&\rat(6/5){3}{}&\rat(6/5){2}{}&&\cnt{368}&\bounds{102}{102}{1}{false}&&\cr
\l{1}&\rat(3/5){3}{}&\rat(4/5){4}{}&\rat(4/5){4}{}&\int(0){\0}{}&\rat(4/5){2}{}&\int(0){\0}{}&\ln{2}&\cnt{12}&\bounds{3}{3}{1}{false}&&\cr
\l{3}&\rat(3/5){3}{}&\rat(4/5){4}{}&\int(0){\0}{}&\rat(6/5){2}{}&\int(0){\0}{}&\rat(2/5){4}{}&\ln{2}&\cnt{12}&\bounds{3}{3}{1}{false}&&\cr
\l{5}&\rat(3/5){3}{}&\int(0){\0}{}&\rat(4/5){4}{}&\rat(2/5){4}{}&\int(0){\0}{}&\rat(6/5){2}{}&\ln{2}&\cnt{12}&\bounds{3}{3}{1}{false}&&\cr
\l{7}&\rat(2/5){2}{}&\rat(4/5){4}{}&\int(0){\0}{}&\int(0){\0}{}&\rat(3/5){4}{}&\rat(6/5){2}{}&\ln{2}&\cnt{12}&\bounds{4}{4}{1}{false}&&\cr
\l{9}&\rat(4/5){4}{}&\rat(3/5){3}{}&\rat(4/5){4}{}&\rat(4/5){3}{}&\int(0){\0}{}&\int(0){\0}{}&\ln{1}&\cnt{12}&\bounds{3}{3}{1}{false}&&\cr
\l{11}&\rat(4/5){4}{}&\rat(3/5){3}{}&\int(0){\0}{}&\int(0){\0}{}&\rat(6/5){3}{}&\rat(2/5){4}{}&\ln{2}&\cnt{12}&\bounds{3}{3}{1}{false}&&\cr
\l{13}&\rat(4/5){4}{}&\rat(4/5){4}{}&\int(0){\0}{}&\rat(2/5){4}{}&\rat(2/5){1}{}&\rat(3/5){1}{}&\ln{2}&\cnt{18}&\bounds{6}{6}{1}{false}&&\cr
\l{15}&\rat(4/5){4}{}&\rat(1/5){1}{}&\rat(4/5){4}{}&\int(0){\0}{}&\rat(3/5){4}{}&\rat(3/5){1}{}&\ln{1}&\cnt{16}&\bounds{4}{4}{1}{false}&&\cr
\endorbit
\orbit{6}{30}{120}
\h{30}&\rat(4/5){4}{}&\rat(4/5){4}{}&\rat(16/5){3}{\spec{2}}&\root\int(0){\0}{}&\int(0){\0}{}&\rat(6/5){3}{}&&\cnt{524}&\bounds{120}{128}{fail}{false}&&\cr
\l{1}&\int(1){\2}{}&\int(0){\0}{}&\int(2){\2}{}&\int(0){\0}{}&\int(0){\0}{}&\int(0){\0}{}&\ln{1}&\cnt{12}&\bounds{3}{3}{1}{true}&&\cr
\l{3}&\rat(3/5){3}{}&\rat(3/5){3}{}&\rat(7/5){1}{}&\int(0){\0}{}&\int(0){\0}{}&\rat(2/5){1}{}&\ln{1}&\cnt{48}&\bounds{10}{12}{5}{true}&&\cr
\l{5}&\rat(3/5){3}{}&\int(0){\0}{}&\rat(9/5){2}{}&\int(0){\0}{}&\int(0){4}{}&\rat(3/5){4}{}&\ln{1}&\cnt{40}&\bounds{8}{8}{1}{true}&&\cr
\l{7}&\rat(2/5){2}{}&\rat(4/5){4}{}&\rat(9/5){2}{}&\int(0){4}{}&\int(0){\0}{}&\int(0){\0}{}&\ln{1}&\cnt{18}&\bounds{4}{6}{5}{true}&&\cr
\l{9}&\rat(4/5){4}{}&\rat(2/5){2}{}&\rat(9/5){2}{}&\int(0){\0}{}&\int(0){1}{}&\int(0){\0}{}&\ln{1}&\cnt{30}&\bounds{6}{6}{1}{true}&&\cr
\l{11}&\rat(4/5){4}{}&\rat(4/5){4}{}&\rat(6/5){3}{}&\int(0){\0}{}&\int(0){\0}{}&\rat(1/5){3}{}&\ln{1}&\cnt{18}&\bounds{6}{6}{1}{true}&&\cr
\l{13}&\rat(4/5){4}{}&\rat(4/5){4}{}&\rat(1/5){3}{}&\int(0){\0}{}&\int(0){\0}{}&\rat(6/5){3}{}&\ln{1}&\cnt{3}&\bounds{1}{1}{1}{true}&&\cr
\l{15}&\rat(4/5){4}{}&\rat(4/5){4}{}&\rat(7/5){1}{}&\int(0){1}{}&\int(0){4}{}&\int(0){\0}{}&\ln{1}&\cnt{15}&\bounds{3}{3}{1}{true}&&\cr
\l{17}&\rat(4/5){4}{}&\rat(1/5){1}{}&\rat(7/5){1}{}&\int(0){4}{}&\int(0){\0}{}&\rat(3/5){4}{}&\ln{1}&\cnt{24}&\bounds{6}{6}{1}{true}&&\cr
\l{19}&\rat(4/5){4}{}&\rat(-1/5){4}{}&\rat(6/5){3}{}&\int(0){\0}{}&\int(0){\0}{}&\rat(6/5){3}{}&\ln{1}&\cnt{12}&\bounds{3}{3}{1}{true}&&\cr
\l{21}&\rat(4/5){4}{}&\int(0){\0}{}&\rat(9/5){2}{}&\int(0){2}{\spec{1}}&\int(0){\0}{}&\rat(2/5){1}{}&\ln{1}&\cnt{3}&\bounds{1}{1}{1}{true}&&\cr
\l{23}&\rat(4/5){4}{}&\int(0){\0}{}&\rat(9/5){2}{}&\int(0){2}{\spec{2}}&\int(0){\0}{}&\rat(2/5){1}{}&\ln{1}&\cnt{9}&\bounds{3}{3}{1}{true}&&\cr
\l{25}&\rat(4/5){4}{}&\int(0){\0}{}&\rat(7/5){1}{}&\int(0){\0}{}&\int(0){2}{}&\rat(4/5){2}{}&\ln{1}&\cnt{30}&\bounds{6}{6}{1}{true}&&\cr
\endorbit
\orbit{6}{31}{120}
\h{31}&\rat(4/5){4}{}&\rat(4/5){4}{}&\root\rat(16/5){3}{\spec{2}}&\int(0){\0}{}&\int(0){\0}{}&\rat(6/5){3}{}&&\cnt{576}&\bounds{120}{124}{fail}{false}&&\cr
\l{1}&\int(1){\2}{}&\int(0){\0}{}&\int(2){\2}{}&\int(0){\0}{}&\int(0){\0}{}&\int(0){\0}{}&\ln{2}&\cnt{4}&\bounds{1}{1}{1}{true}&&\cr
\l{3}&\rat(3/5){3}{}&\rat(3/5){3}{}&\rat(7/5){1}{}&\int(0){\0}{}&\int(0){\0}{}&\rat(2/5){1}{}&\ln{1}&\cnt{48}&\bounds{10}{12}{5}{true}&&\cr
\l{5}&\rat(3/5){3}{}&\int(0){\0}{}&\rat(9/5){2}{}&\int(0){\0}{}&\int(0){4}{}&\rat(3/5){4}{}&\ln{2}&\cnt{40}&\bounds{8}{8}{1}{true}&&\cr
\l{7}&\rat(2/5){2}{}&\rat(4/5){4}{}&\rat(9/5){2}{}&\int(0){4}{}&\int(0){\0}{}&\int(0){\0}{}&\ln{2}&\cnt{30}&\bounds{6}{6}{1}{true}&&\cr
\l{9}&\rat(4/5){4}{}&\rat(4/5){4}{}&\rat(6/5){3}{}&\int(0){\0}{}&\int(0){\0}{}&\rat(1/5){3}{}&\ln{1}&\cnt{6}&\bounds{2}{2}{1}{true}&&\cr
\l{11}&\rat(4/5){4}{}&\rat(4/5){4}{}&\rat(1/5){3}{}&\int(0){\0}{}&\int(0){\0}{}&\rat(6/5){3}{}&\ln{1}&\cnt{1}&\bounds{1}{1}{1}{true}&&\cr
\l{13}&\rat(4/5){4}{}&\rat(4/5){4}{}&\rat(7/5){1}{}&\int(0){1}{}&\int(0){4}{}&\int(0){\0}{}&\ln{1}&\cnt{25}&\bounds{5}{5}{1}{true}&&\cr
\l{15}&\rat(4/5){4}{}&\int(0){\0}{}&\rat(9/5){2}{}&\int(0){2}{}&\int(0){\0}{}&\rat(2/5){1}{}&\ln{2}&\cnt{30}&\bounds{6}{6}{1}{true}&&\cr
\endorbit
\orbit{6}{32}{110}
\h{32}&\root\rat(4/5){4}{}&\rat(4/5){4}{}&\rat(16/5){3}{\spec{2}}&\int(0){\0}{}&\int(0){\0}{}&\rat(6/5){3}{}&&\cnt{496}&\bounds{110}{110}{1}{false}&&\cr
\l{1}&\int(1){\2}{}&\int(0){\0}{}&\int(2){\2}{}&\int(0){\0}{}&\int(0){\0}{}&\int(0){\0}{}&\ln{1}&\cnt{6}&\bounds{2}{2}{1}{false}&&\cr
\l{3}&\rat(3/5){3}{}&\rat(3/5){3}{}&\rat(7/5){1}{}&\int(0){\0}{}&\int(0){\0}{}&\rat(2/5){1}{}&\ln{1}&\cnt{24}&\bounds{6}{6}{1}{false}&&\cr
\l{5}&\rat(3/5){3}{}&\int(0){\0}{}&\rat(9/5){2}{}&\int(0){\0}{}&\int(0){4}{}&\rat(3/5){4}{}&\ln{1}&\cnt{20}&\bounds{4}{4}{1}{false}&&\cr
\l{7}&\rat(2/5){2}{\spec{1}}&\rat(4/5){4}{}&\rat(9/5){2}{}&\int(0){4}{}&\int(0){\0}{}&\int(0){\0}{}&\ln{1}&\cnt{5}&\bounds{1}{1}{1}{false}&&\cr
\l{9}&\rat(2/5){2}{\spec{1}}&\int(0){\0}{}&\rat(7/5){1}{}&\int(0){1}{}&\int(0){\0}{}&\rat(6/5){3}{}&\ln{1}&\cnt{5}&\bounds{1}{1}{1}{false}&&\cr
\l{11}&\rat(4/5){4}{}&\rat(2/5){2}{}&\rat(9/5){2}{}&\int(0){\0}{}&\int(0){1}{}&\int(0){\0}{}&\ln{1}&\cnt{30}&\bounds{6}{6}{1}{false}&&\cr
\l{13}&\rat(4/5){4}{}&\rat(4/5){4}{}&\rat(6/5){3}{}&\int(0){\0}{}&\int(0){\0}{}&\rat(1/5){3}{}&\ln{1}&\cnt{18}&\bounds{6}{6}{1}{false}&&\cr
\l{15}&\rat(4/5){4}{}&\rat(4/5){4}{}&\rat(1/5){3}{}&\int(0){\0}{}&\int(0){\0}{}&\rat(6/5){3}{}&\ln{1}&\cnt{3}&\bounds{1}{1}{1}{false}&&\cr
\l{17}&\rat(4/5){4}{}&\rat(4/5){4}{}&\rat(7/5){1}{}&\int(0){1}{}&\int(0){4}{}&\int(0){\0}{}&\ln{1}&\cnt{25}&\bounds{5}{5}{1}{false}&&\cr
\l{19}&\rat(4/5){4}{}&\rat(1/5){1}{}&\rat(7/5){1}{}&\int(0){4}{}&\int(0){\0}{}&\rat(3/5){4}{}&\ln{1}&\cnt{40}&\bounds{8}{8}{1}{false}&&\cr
\l{21}&\rat(4/5){4}{}&\rat(-1/5){4}{}&\rat(6/5){3}{}&\int(0){\0}{}&\int(0){\0}{}&\rat(6/5){3}{}&\ln{1}&\cnt{12}&\bounds{3}{3}{1}{false}&&\cr
\l{23}&\rat(4/5){4}{}&\int(0){\0}{}&\rat(9/5){2}{}&\int(0){2}{}&\int(0){\0}{}&\rat(2/5){1}{}&\ln{1}&\cnt{30}&\bounds{6}{6}{1}{false}&&\cr
\l{25}&\rat(4/5){4}{}&\int(0){\0}{}&\rat(7/5){1}{}&\int(0){\0}{}&\int(0){2}{}&\rat(4/5){2}{}&\ln{1}&\cnt{30}&\bounds{6}{6}{1}{false}&&\cr
\endorbit
\orbit{6}{33}{100}
\h{33}&\rat(4/5){4}{}&\rat(4/5){4}{}&\rat(16/5){3}{\spec{2}}&\int(0){\0}{}&\int(0){\0}{}&\root\rat(6/5){3}{\spec{1}}&&\cnt{468}&\bounds{100}{100}{1}{false}&&\cr
\l{1}&\int(1){\2}{}&\int(0){\0}{}&\int(2){\2}{}&\int(0){\0}{}&\int(0){\0}{}&\int(0){\0}{}&\ln{2}&\cnt{12}&\bounds{3}{3}{1}{false}&&\cr
\l{3}&\rat(3/5){3}{}&\rat(3/5){3}{}&\rat(7/5){1}{}&\int(0){\0}{}&\int(0){\0}{}&\rat(2/5){1}{}&\ln{1}&\cnt{16}&\bounds{4}{4}{1}{false}&&\cr
\l{5}&\rat(3/5){3}{}&\int(0){\0}{}&\rat(9/5){2}{}&\int(0){\0}{}&\int(0){4}{}&\rat(3/5){4}{}&\ln{2}&\cnt{40}&\bounds{8}{8}{1}{false}&&\cr
\l{7}&\rat(2/5){2}{}&\rat(4/5){4}{}&\rat(9/5){2}{}&\int(0){4}{}&\int(0){\0}{}&\int(0){\0}{}&\ln{2}&\cnt{30}&\bounds{6}{6}{1}{false}&&\cr
\l{9}&\rat(4/5){4}{}&\rat(4/5){4}{}&\rat(6/5){3}{}&\int(0){\0}{}&\int(0){\0}{}&\rat(1/5){3}{}&\ln{1}&\cnt{6}&\bounds{2}{2}{1}{false}&&\cr
\l{11}&\rat(4/5){4}{}&\rat(4/5){4}{}&\rat(1/5){3}{}&\int(0){\0}{}&\int(0){\0}{}&\rat(6/5){3}{}&\ln{1}&\cnt{3}&\bounds{1}{1}{1}{false}&&\cr
\l{13}&\rat(4/5){4}{}&\rat(4/5){4}{}&\rat(7/5){1}{}&\int(0){1}{}&\int(0){4}{}&\int(0){\0}{}&\ln{1}&\cnt{25}&\bounds{5}{5}{1}{false}&&\cr
\l{15}&\rat(4/5){4}{}&\int(0){\0}{}&\rat(9/5){2}{}&\int(0){2}{}&\int(0){\0}{}&\rat(2/5){1}{}&\ln{2}&\cnt{10}&\bounds{2}{2}{1}{false}&&\cr
\endorbit
\orbit{6}{34}{96}
\h{34}&\rat(4/5){4}{}&\rat(4/5){4}{}&\rat(16/5){3}{\spec{2}}&\int(0){\0}{}&\int(0){\0}{}&\root\rat(6/5){3}{\spec{2}}&&\cnt{440}&\bounds{96}{96}{1}{false}&&\cr
\l{1}&\int(1){\2}{}&\int(0){\0}{}&\int(2){\2}{}&\int(0){\0}{}&\int(0){\0}{}&\int(0){\0}{}&\ln{2}&\cnt{12}&\bounds{3}{3}{1}{false}&&\cr
\l{3}&\rat(3/5){3}{}&\rat(3/5){3}{}&\rat(7/5){1}{}&\int(0){\0}{}&\int(0){\0}{}&\rat(2/5){1}{}&\ln{1}&\cnt{48}&\bounds{12}{12}{1}{false}&&\cr
\l{5}&\rat(2/5){2}{}&\rat(4/5){4}{}&\rat(9/5){2}{}&\int(0){4}{}&\int(0){\0}{}&\int(0){\0}{}&\ln{2}&\cnt{30}&\bounds{6}{6}{1}{false}&&\cr
\l{7}&\rat(4/5){4}{}&\rat(4/5){4}{}&\rat(1/5){3}{}&\int(0){\0}{}&\int(0){\0}{}&\rat(6/5){3}{}&\ln{1}&\cnt{3}&\bounds{1}{1}{1}{false}&&\cr
\l{9}&\rat(4/5){4}{}&\rat(4/5){4}{}&\rat(7/5){1}{}&\int(0){1}{}&\int(0){4}{}&\int(0){\0}{}&\ln{1}&\cnt{25}&\bounds{5}{5}{1}{false}&&\cr
\l{11}&\rat(4/5){4}{}&\int(0){\0}{}&\rat(9/5){2}{}&\int(0){2}{}&\int(0){\0}{}&\rat(2/5){1}{}&\ln{2}&\cnt{30}&\bounds{6}{6}{1}{false}&&\cr
\endorbit
\orbit{6}{35}{110}
\h{35}&\root\rat(4/5){4}{}&\rat(4/5){4}{}&\rat(4/5){1}{}&\rat(4/5){1}{}&\rat(4/5){4}{}&\int(2){\2}{}&&\cnt{448}&\bounds{110}{126}{fail}{false}&&\cr
\l{1}&\int(1){\2}{}&\int(0){\0}{}&\int(0){\0}{}&\int(0){\0}{}&\int(0){\0}{}&\int(2){\2}{}&\ln{1}&\cnt{2}&\bounds{1}{1}{1}{true}&&\cr
\l{3}&\rat(3/5){3}{}&\rat(3/5){3}{}&\rat(4/5){1}{}&\int(0){\0}{}&\int(0){\0}{}&\int(1){1}{}&\ln{4}&\cnt{8}&\bounds{2}{2}{1}{true}&&\cr
\l{5}&\rat(2/5){2}{\spec{1}}&\rat(4/5){4}{}&\int(0){\0}{}&\int(0){\0}{}&\rat(4/5){4}{}&\int(1){2}{}&\ln{2}&\cnt{3}&\bounds{1}{1}{1}{true}&&\cr
\l{7}&\rat(4/5){4}{}&\rat(3/5){3}{}&\int(0){\0}{}&\int(0){\0}{}&\rat(3/5){3}{}&\int(1){4}{}&\ln{2}&\cnt{16}&\bounds{4}{4}{1}{true}&&\cr
\l{9}&\rat(4/5){4}{}&\rat(2/5){2}{}&\int(0){\0}{}&\rat(4/5){1}{}&\int(0){\0}{}&\int(1){2}{}&\ln{4}&\cnt{18}&\bounds{4}{6}{5}{true}&&\cr
\l{11}&\rat(4/5){4}{}&\rat(4/5){4}{}&\rat(4/5){1}{}&\rat(4/5){1}{}&\rat(-1/5){4}{}&\int(0){\0}{}&\ln{4}&\cnt{4}&\bounds{1}{1}{1}{true}&&\cr
\l{13}&\rat(4/5){4}{}&\rat(4/5){4}{}&\int(0){\0}{}&\rat(1/5){4}{}&\rat(1/5){1}{}&\int(1){1}{}&\ln{4}&\cnt{16}&\bounds{4}{4}{1}{true}&&\cr
\endorbit
\orbit{6}{36}{110}
\h{36}&\rat(4/5){4}{}&\rat(4/5){4}{}&\rat(4/5){1}{}&\rat(4/5){1}{}&\rat(4/5){4}{}&\root\int(2){\2}{}&&\cnt{420}&\bounds{110}{110}{1}{false}&&\cr
\l{1}&\rat(2/5){2}{}&\rat(4/5){4}{}&\int(0){\0}{}&\int(0){\0}{}&\rat(4/5){4}{}&\int(1){2}{}&\ln{10\sd}&\cnt{6}&\bounds{2}{2}{1}{false}&&\cr
\l{2}&\int(1){\2}{}&\int(0){\0}{}&\int(0){\0}{}&\int(0){\0}{}&\int(0){\0}{}&\int(2){\2}{}&\ln{5}&\cnt{4}&\bounds{1}{1}{1}{false}&&\cr
\l{4}&\rat(3/5){3}{}&\rat(3/5){3}{}&\rat(4/5){1}{}&\int(0){\0}{}&\int(0){\0}{}&\int(1){1}{}&\ln{10}&\cnt{16}&\bounds{4}{4}{1}{false}&&\cr
\endorbit
\orbit{6}{37}{120}
\h{37}&\int(2){\2}{}&\int(2){\2}{}&\int(2){\2}{}&\root\int(0){\0}{}&\int(0){\0}{}&\int(0){\0}{}&&\cnt{524}&\bounds{120}{120}{1}{false}&&\cr
\l{1}&\int(2){\2}{}&\int(1){\2}{\spec{1}}&\int(0){\0}{}&\int(0){\0}{}&\int(0){\0}{}&\int(0){\0}{}&\ln{4\sd}&\cnt{3}&\bounds{1}{1}{1}{true}&&\cr
\l{2}&\int(1){3}{}&\int(1){3}{}&\int(1){1}{}&\int(0){\0}{}&\int(0){\0}{}&\int(0){1}{}&\ln{4\sd}&\cnt{45}&\bounds{9}{9}{1}{true}&&\cr
\l{3}&\int(1){3}{}&\int(1){4}{}&\int(1){4}{}&\int(0){\0}{}&\int(0){2}{}&\int(0){\0}{}&\ln{4\sd}&\cnt{30}&\bounds{6}{6}{1}{true}&&\cr
\l{4}&\int(1){3}{}&\int(1){1}{}&\int(1){3}{}&\int(0){1}{}&\int(0){\0}{}&\int(0){\0}{}&\ln{2\sd}&\cnt{27}&\bounds{9}{9}{1}{true}&&\cr
\l{5}&\int(1){4}{}&\int(1){3}{}&\int(1){4}{}&\int(0){3}{\spec{1}}&\int(0){\0}{}&\int(0){\0}{}&\ln{2\sd}&\cnt{9}&\bounds{3}{3}{1}{true}&&\cr
\l{6}&\int(1){4}{}&\int(1){3}{}&\int(1){4}{}&\int(0){3}{\spec{2}}&\int(0){\0}{}&\int(0){\0}{}&\ln{2\sd}&\cnt{3}&\bounds{1}{1}{1}{true}&&\cr
\l{7}&\int(1){4}{}&\int(1){4}{}&\int(1){1}{}&\int(0){1}{}&\int(0){4}{}&\int(0){\0}{}&\ln{4\sd}&\cnt{15}&\bounds{3}{3}{1}{true}&&\cr
\l{8}&\int(1){4}{}&\int(1){1}{}&\int(1){4}{}&\int(0){\0}{}&\int(0){4}{}&\int(0){1}{}&\ln{2\sd}&\cnt{25}&\bounds{5}{5}{1}{true}&&\cr
\l{9}&\int(2){\2}{}&\int(0){\0}{}&\int(1){\2}{\spec{1}}&\int(0){\0}{}&\int(0){\0}{}&\int(0){\0}{}&\ln{4}&\cnt{3}&\bounds{1}{1}{1}{true}&&\cr
\endorbit
\orbit{6}{38}{100}
\h{38}&\root\int(2){\2}{}&\int(2){\2}{}&\int(2){\2}{}&\int(0){\0}{}&\int(0){\0}{}&\int(0){\0}{}&&\cnt{468}&\bounds{100}{100}{1}{false}&&\cr
\l{1}&\int(1){\2}{\spec{1}}&\int(2){\2}{}&\int(0){\0}{}&\int(0){\0}{}&\int(0){\0}{}&\int(0){\0}{}&\ln{4\sd}&\cnt{1}&\bounds{1}{1}{1}{false}&&\cr
\l{2}&\int(1){3}{}&\int(1){3}{}&\int(1){1}{}&\int(0){\0}{}&\int(0){\0}{}&\int(0){1}{}&\ln{4\sd}&\cnt{15}&\bounds{3}{3}{1}{false}&&\cr
\l{3}&\int(1){3}{}&\int(1){4}{}&\int(1){4}{}&\int(0){\0}{}&\int(0){2}{}&\int(0){\0}{}&\ln{2\sd}&\cnt{10}&\bounds{2}{2}{1}{false}&&\cr
\l{4}&\int(1){4}{}&\int(1){3}{}&\int(1){4}{}&\int(0){3}{}&\int(0){\0}{}&\int(0){\0}{}&\ln{4\sd}&\cnt{30}&\bounds{6}{6}{1}{false}&&\cr
\l{5}&\int(1){4}{}&\int(1){2}{}&\int(1){2}{}&\int(0){\0}{}&\int(0){1}{}&\int(0){\0}{}&\ln{2\sd}&\cnt{45}&\bounds{9}{9}{1}{false}&&\cr
\l{6}&\int(1){4}{}&\int(1){4}{}&\int(1){1}{}&\int(0){1}{}&\int(0){4}{}&\int(0){\0}{}&\ln{4\sd}&\cnt{25}&\bounds{5}{5}{1}{false}&&\cr
\l{7}&\int(1){4}{}&\int(1){1}{}&\int(1){1}{}&\int(0){4}{}&\int(0){\0}{}&\int(0){4}{}&\ln{2\sd}&\cnt{25}&\bounds{5}{5}{1}{false}&&\cr
\l{8}&\int(2){\2}{}&\int(1){\2}{\spec{1}}&\int(0){\0}{}&\int(0){\0}{}&\int(0){\0}{}&\int(0){\0}{}&\ln{4}&\cnt{3}&\bounds{1}{1}{1}{false}&&\cr
\endorbit
\orbit{6}{39}{51}
\h{39}&\int(6){\6}{\spec{1}}&\int(0){\0}{}&\int(0){\0}{}&\int(0){\0}{}&\int(0){\0}{}&\int(0){\0}{}&&\cnt{204}&\bounds{51}{51}{1}{false}&&\cr
\l{1}&\int(3){\4}{}&\int(0){\0}{}&\int(0){\0}{}&\int(0){\0}{}&\int(0){\0}{}&\int(0){\0}{}&\ln{1\sdd}&\cnt{4}&\bounds{1}{1}{1}{false}&&\cr
\l{2}&\int(3){\2}{}&\int(0){\2}{}&\int(0){\0}{}&\int(0){\0}{}&\int(0){\0}{}&\int(0){\0}{}&\ln{5\sdd}&\cnt{40}&\bounds{10}{10}{1}{false}&&\cr
\endorbit
\endorbits

\subsubsection{Configuration \hlink{6A4}{1}}\label{s.6A4-1}
There is a unique set $\Lmisc1\in\Bnd_{28}(\Orb)$;
it has $130$ lines and is described by the pattern
(see \autoref{rem.pattern.unique})
\[
\pat=\pattern{0,6,10}\quad\mbox{(not $\Q$-complete)}.
\misclabel1{130}
\]
This case completes the technical details of the proof of \autoref{th.<=6}.
\endproof


\def\ltitle{}

\section{The lattice $\N(8\bA\sb3)$}\label{S.8A3}
Starting from this section, we relax the goal to
\[*
\ls|\fL|\ge M:=132.
\]
The result of this section is the following theorem.

\chardef\mincount130
\newcs{check-8A3-1}{\noexpand\checkshow}
\newcs{check-8A3-4}{}
\newcs{check-8A3-5}{\noexpand\checkdone}
\newcs{check-8A3-7}{\noexpand\checkdone}
\newcs{check-8A3-10}{\noexpand\checkdone}
\newcs{check-8A3-14}{\noexpand\checkdone}
\newcs{check-8A3-15}{\noexpand\checkdone}
\newcs{check-8A3-18}{\noexpand\checkdone}
\newcs{check-8A3-19}{\noexpand\checkdone}
\newcs{check-8A3-21}{\noexpand\checkdone}
\newcs{check-8A3-27}{\noexpand\checkdone}

\theorem\label{th.8A3}
Let $\fL\subset\N(8\bA_3)$ be a geometric set. Then, unless
\roster*
\item
$\ls|\fL|=132$
and $\fL$ is conjugate to $\Lsub1$, see \eqref{eq.Lsub.1},
or $\Lsub2$, see \eqref{eq.Lsub.2},
\endroster
one has $\ls|\fL|\le130$.
\endtheorem

\proof
We proceed as in \autoref{S.few}, listing pairs $\hbar,\rbar\in\N(8\bA_3)$
and respective orbits (and following the notation of \autoref{S.few}).
\newcs{max-8A3}{160}%
\newcs{count-8A3}{28}%
\newcs{counts-8A3}{\DO{1}{132}\DO{2}{112}\DO{3}{96}\DO{4}{160}\DO{5}{136}\DO{6}{124}\DO{7}{136}\DO{8}{126}\DO{9}{106}\DO{10}{148}\DO{11}{128}\DO{12}{120}\DO{13}{118}\DO{14}{140}\DO{15}{138}\DO{16}{124}\DO{17}{116}\DO{18}{146}\DO{19}{140}\DO{20}{122}\DO{21}{140}\DO{22}{112}\DO{23}{112}\DO{24}{126}\DO{25}{122}\DO{26}{104}\DO{27}{140}\DO{28}{42}}%
\orbits
\orbit{8}{1}{132}
\h{1}&\root\int(4){\4}{\spec{1}}&\int(2){\2}{}&\int(0){\0}{}&\int(0){\0}{}&\int(0){\0}{}&\int(0){\0}{}&\int(0){\0}{}&\int(0){\0}{}&&\cnt{728}&\bounds{132}{200}{fail}{false}&&\cr
\l{1}&\int(2){2}{}&\int(1){2}{}&\int(0){2}{}&\int(0){\0}{}&\int(0){\0}{}&\int(0){\0}{}&\int(0){2}{}&\int(0){\0}{}&\ln{3\sdd}&\cnt{72}&\bounds{12}{24}{5}{true}&&\cr
\l{2}&\int(2){2}{}&\int(1){3}{}&\int(0){3}{}&\int(0){\0}{}&\int(0){3}{}&\int(0){\0}{}&\int(0){\0}{}&\int(0){1}{}&\ln{8\sd}&\cnt{64}&\bounds{12}{16}{3}{false}&&\cr
\endorbit
\orbit{8}{2}{112}
\h{2}&\int(4){\4}{}&\int(2){\2}{}&\root\int(0){\0}{}&\int(0){\0}{}&\int(0){\0}{}&\int(0){\0}{}&\int(0){\0}{}&\int(0){\0}{}&&\cnt{568}&\bounds{112}{160}{fail}{false}&&\cr
\l{1}&\int(2){\2}{}&\int(1){\2}{\spec{1}}&\int(0){\0}{}&\int(0){\0}{}&\int(0){\0}{}&\int(0){\0}{}&\int(0){\0}{}&\int(0){\0}{}&\ln{2\sd}&\cnt{8}&\bounds{3}{4}{5}{true}&&\cr
\l{2}&\int(2){2}{}&\int(1){2}{}&\int(0){2}{\spec{1}}&\int(0){\0}{}&\int(0){\0}{}&\int(0){\0}{}&\int(0){2}{}&\int(0){\0}{}&\ln{2\sd}&\cnt{12}&\bounds{3}{4}{5}{true}&&\cr
\l{3}&\int(2){2}{}&\int(1){2}{}&\int(0){\0}{}&\int(0){2}{}&\int(0){2}{}&\int(0){\0}{}&\int(0){\0}{}&\int(0){\0}{}&\ln{2\sdd}&\cnt{72}&\bounds{12}{24}{5}{true}&&\cr
\l{4}&\int(2){2}{}&\int(1){3}{}&\int(0){3}{}&\int(0){\0}{}&\int(0){3}{}&\int(0){\0}{}&\int(0){\0}{}&\int(0){1}{}&\ln{4\sd}&\cnt{32}&\bounds{7}{8}{5}{true}&&\cr
\l{5}&\int(2){2}{}&\int(1){3}{}&\int(0){\0}{}&\int(0){3}{}&\int(0){\0}{}&\int(0){\0}{}&\int(0){1}{}&\int(0){3}{}&\ln{4\sd}&\cnt{64}&\bounds{12}{16}{3}{false}&&\cr
\endorbit
\orbit{8}{3}{96}
\h{3}&\int(4){\4}{}&\root\int(2){\2}{}&\int(0){\0}{}&\int(0){\0}{}&\int(0){\0}{}&\int(0){\0}{}&\int(0){\0}{}&\int(0){\0}{}&&\cnt{512}&\bounds{96}{128}{fail}{false}&&\cr
\l{1}&\int(2){2}{}&\int(1){3}{}&\int(0){3}{}&\int(0){\0}{}&\int(0){3}{}&\int(0){\0}{}&\int(0){\0}{}&\int(0){1}{}&\ln{8\sd}&\cnt{64}&\bounds{12}{16}{3}{false}&&\cr
\endorbit
\orbit{8}{4}{160}
\h{4}&\rat(3/4){3}{}&\rat(3/4){3}{}&\rat(3/4){3}{}&\int(1){2}{}&\int(1){2}{}&\root\int(0){\0}{}&\rat(3/4){1}{}&\int(1){2}{}&&\cnt{448}&\bounds{160}{160}{1}{false}&&\cr
\l{1}&\rat(3/4){3}{}&\rat(3/4){3}{}&\int(0){\0}{}&\int(1){2}{}&\int(0){\0}{}&\int(0){3}{}&\int(0){\0}{}&\rat(1/2){1}{}&\ln{12\sd}&\cnt{4}&\bounds{2}{2}{1}{true}&&\cr
\l{2}&\rat(1/2){2}{}&\rat(1/2){2}{}&\int(0){\0}{}&\int(1){2}{}&\int(1){2}{}&\int(0){\0}{}&\int(0){\0}{}&\int(0){\0}{}&\ln{6}&\cnt{9}&\bounds{3}{3}{1}{true}&&\cr
\l{4}&\rat(1/2){2}{}&\rat(3/4){3}{}&\rat(3/4){3}{}&\int(0){\0}{}&\rat(1/2){3}{}&\int(0){\0}{}&\int(0){\0}{}&\rat(1/2){1}{}&\ln{12}&\cnt{12}&\bounds{4}{4}{1}{true}&&\cr
\l{6}&\rat(3/4){3}{}&\rat(3/4){3}{}&\rat(3/4){3}{}&\int(0){\0}{}&\int(0){\0}{}&\int(0){2}{\spec{1}}&\rat(3/4){1}{}&\int(0){\0}{}&\ln{2}&\cnt{1}&\bounds{1}{1}{1}{true}&&\cr
\endorbit
\orbit{8}{5}{136}
\h{5}&\root\rat(3/4){3}{}&\rat(3/4){3}{}&\rat(3/4){3}{}&\int(1){2}{}&\int(1){2}{}&\int(0){\0}{}&\rat(3/4){1}{}&\int(1){2}{}&&\cnt{420}&\bounds{136}{136}{1}{false}&&\cr
\l{1}&\rat(1/2){2}{}&\rat(1/2){2}{}&\int(0){\0}{}&\int(1){2}{}&\int(1){2}{}&\int(0){\0}{}&\int(0){\0}{}&\int(0){\0}{}&\ln{3}&\cnt{3}&\bounds{1}{1}{1}{true}&&\cr
\l{3}&\rat(1/2){2}{}&\rat(3/4){3}{}&\rat(3/4){3}{}&\int(0){\0}{}&\rat(1/2){3}{}&\int(0){\0}{}&\int(0){\0}{}&\rat(1/2){1}{}&\ln{3}&\cnt{4}&\bounds{2}{2}{1}{true}&&\cr
\l{5}&\rat(3/4){3}{}&\rat(1/2){2}{}&\rat(3/4){3}{}&\rat(1/2){1}{}&\int(0){\0}{}&\int(0){\0}{}&\int(0){\0}{}&\rat(1/2){3}{}&\ln{6}&\cnt{12}&\bounds{4}{4}{1}{true}&&\cr
\l{7}&\rat(3/4){3}{}&\rat(3/4){3}{}&\rat(3/4){3}{}&\int(0){\0}{}&\int(0){\0}{}&\int(0){2}{}&\rat(3/4){1}{}&\int(0){\0}{}&\ln{1}&\cnt{6}&\bounds{2}{2}{1}{true}&&\cr
\l{9}&\rat(3/4){3}{}&\rat(3/4){3}{}&\rat(1/4){1}{}&\int(0){\0}{}&\int(0){\0}{}&\int(0){\0}{}&\rat(1/4){3}{}&\int(1){2}{}&\ln{3}&\cnt{9}&\bounds{3}{3}{1}{true}&&\cr
\l{11}&\rat(3/4){3}{}&\rat(3/4){3}{}&\int(0){\0}{}&\int(1){2}{}&\int(0){\0}{}&\int(0){3}{}&\int(0){\0}{}&\rat(1/2){1}{}&\ln{6}&\cnt{8}&\bounds{2}{2}{1}{true}&&\cr
\l{13}&\rat(3/4){3}{}&\rat(1/4){1}{}&\int(0){\0}{}&\rat(1/2){3}{}&\rat(1/2){3}{}&\int(0){\0}{}&\int(0){\0}{}&\int(1){2}{}&\ln{3}&\cnt{12}&\bounds{4}{4}{1}{true}&&\cr
\endorbit
\orbit{8}{6}{124}
\h{6}&\rat(3/4){3}{}&\rat(3/4){3}{}&\rat(3/4){3}{}&\root\int(1){2}{\spec{1}}&\int(1){2}{}&\int(0){\0}{}&\rat(3/4){1}{}&\int(1){2}{}&&\cnt{392}&\bounds{124}{124}{1}{false}&&\cr
\l{1}&\rat(3/4){3}{}&\int(0){\0}{}&\int(0){\0}{}&\rat(1/2){3}{}&\int(1){2}{}&\int(0){3}{}&\rat(3/4){1}{}&\int(0){\0}{}&\ln{2\sd}&\cnt{8}&\bounds{2}{2}{1}{true}&&\cr
\l{2}&\rat(1/2){2}{}&\rat(1/2){2}{}&\int(0){\0}{}&\int(1){2}{}&\int(1){2}{}&\int(0){\0}{}&\int(0){\0}{}&\int(0){\0}{}&\ln{4}&\cnt{9}&\bounds{3}{3}{1}{true}&&\cr
\l{4}&\rat(1/2){2}{}&\rat(3/4){3}{}&\rat(3/4){3}{}&\int(0){\0}{}&\rat(1/2){3}{}&\int(0){\0}{}&\int(0){\0}{}&\rat(1/2){1}{}&\ln{4}&\cnt{12}&\bounds{4}{4}{1}{true}&&\cr
\l{6}&\rat(1/2){2}{}&\rat(3/4){3}{}&\int(0){\0}{}&\rat(1/2){3}{}&\int(0){\0}{}&\int(0){\0}{}&\rat(3/4){1}{}&\rat(1/2){3}{}&\ln{4}&\cnt{12}&\bounds{4}{4}{1}{true}&&\cr
\l{8}&\rat(1/2){2}{}&\int(0){\0}{}&\int(0){\0}{}&\int(0){\0}{}&\int(1){2}{}&\int(0){\0}{}&\rat(1/2){2}{}&\int(1){2}{}&\ln{2}&\cnt{9}&\bounds{3}{3}{1}{true}&&\cr
\l{10}&\rat(3/4){3}{}&\rat(3/4){3}{}&\rat(3/4){3}{}&\int(0){\0}{}&\int(0){\0}{}&\int(0){2}{}&\rat(3/4){1}{}&\int(0){\0}{}&\ln{1}&\cnt{6}&\bounds{2}{2}{1}{true}&&\cr
\l{12}&\rat(3/4){3}{}&\rat(3/4){3}{}&\int(0){\0}{}&\int(1){2}{}&\int(0){\0}{}&\int(0){3}{}&\int(0){\0}{}&\rat(1/2){1}{}&\ln{4}&\cnt{8}&\bounds{2}{2}{1}{true}&&\cr
\endorbit
\orbit{8}{7}{136}
\h{7}&\root\int(3){2}{\spec{1}}&\int(1){2}{}&\int(1){2}{}&\int(0){\0}{}&\int(0){\0}{}&\int(0){\0}{}&\int(1){2}{}&\int(0){\0}{}&&\cnt{600}&\bounds{136}{172}{fail}{false}&&\cr
\l{1}&\rat(3/2){3}{}&\int(1){2}{}&\rat(1/2){3}{}&\int(0){1}{}&\int(0){\0}{}&\int(0){\0}{}&\int(0){\0}{}&\int(0){3}{}&\ln{12\sd}&\cnt{32}&\bounds{7}{8}{5}{true}&&\cr
\l{2}&\rat(3/2){3}{}&\rat(1/2){3}{}&\rat(1/2){3}{}&\int(0){\0}{}&\int(0){\0}{}&\int(0){2}{}&\rat(1/2){1}{}&\int(0){\0}{}&\ln{4\sdd}&\cnt{48}&\bounds{10}{16}{5}{true}&&\cr
\l{3}&\int(2){\2}{}&\int(1){\2}{}&\int(0){\0}{}&\int(0){\0}{}&\int(0){\0}{}&\int(0){\0}{}&\int(0){\0}{}&\int(0){\0}{}&\ln{3}&\cnt{4}&\bounds{2}{2}{1}{true}&&\cr
\endorbit
\orbit{8}{8}{126}
\h{8}&\int(3){2}{\spec{1}}&\int(1){2}{}&\int(1){2}{}&\root\int(0){\0}{}&\int(0){\0}{}&\int(0){\0}{}&\int(1){2}{}&\int(0){\0}{}&&\cnt{520}&\bounds{126}{152}{fail}{false}&&\cr
\l{1}&\rat(3/2){3}{}&\int(1){2}{}&\rat(1/2){3}{}&\int(0){1}{}&\int(0){\0}{}&\int(0){\0}{}&\int(0){\0}{}&\int(0){3}{}&\ln{6\sd}&\cnt{16}&\bounds{4}{4}{1}{true}&&\cr
\l{2}&\rat(3/2){3}{}&\int(1){2}{}&\rat(1/2){1}{}&\int(0){\0}{}&\int(0){3}{}&\int(0){3}{}&\int(0){\0}{}&\int(0){\0}{}&\ln{6\sd}&\cnt{32}&\bounds{7}{8}{5}{true}&&\cr
\l{3}&\rat(3/2){3}{}&\rat(1/2){3}{}&\rat(1/2){3}{}&\int(0){\0}{}&\int(0){\0}{}&\int(0){2}{}&\rat(1/2){1}{}&\int(0){\0}{}&\ln{3\sdd}&\cnt{48}&\bounds{10}{16}{5}{true}&&\cr
\l{4}&\rat(3/2){3}{}&\rat(1/2){1}{}&\rat(1/2){1}{}&\int(0){2}{\spec{1}}&\int(0){\0}{}&\int(0){\0}{}&\rat(1/2){1}{}&\int(0){\0}{}&\ln{2\sd}&\cnt{8}&\bounds{3}{4}{5}{true}&&\cr
\l{5}&\int(2){\2}{}&\int(1){\2}{}&\int(0){\0}{}&\int(0){\0}{}&\int(0){\0}{}&\int(0){\0}{}&\int(0){\0}{}&\int(0){\0}{}&\ln{3}&\cnt{12}&\bounds{4}{4}{1}{true}&&\cr
\endorbit
\orbit{8}{9}{106}
\h{9}&\int(3){2}{\spec{1}}&\root\int(1){2}{\spec{1}}&\int(1){2}{}&\int(0){\0}{}&\int(0){\0}{}&\int(0){\0}{}&\int(1){2}{}&\int(0){\0}{}&&\cnt{464}&\bounds{106}{128}{fail}{false}&&\cr
\l{1}&\rat(3/2){3}{}&\rat(1/2){3}{}&\int(1){2}{}&\int(0){3}{}&\int(0){1}{}&\int(0){\0}{}&\int(0){\0}{}&\int(0){\0}{}&\ln{2\sd}&\cnt{32}&\bounds{7}{8}{5}{true}&&\cr
\l{2}&\rat(3/2){3}{}&\rat(1/2){3}{}&\rat(1/2){3}{}&\int(0){\0}{}&\int(0){\0}{}&\int(0){2}{}&\rat(1/2){1}{}&\int(0){\0}{}&\ln{2\sdd}&\cnt{48}&\bounds{10}{16}{3}{false}&&\cr
\l{3}&\int(2){\2}{}&\int(0){\0}{}&\int(1){\2}{}&\int(0){\0}{}&\int(0){\0}{}&\int(0){\0}{}&\int(0){\0}{}&\int(0){\0}{}&\ln{2}&\cnt{12}&\bounds{4}{4}{1}{true}&&\cr
\l{5}&\rat(3/2){3}{}&\int(1){2}{}&\rat(1/2){3}{}&\int(0){1}{}&\int(0){\0}{}&\int(0){\0}{}&\int(0){\0}{}&\int(0){3}{}&\ln{4}&\cnt{32}&\bounds{7}{8}{5}{true}&&\cr
\endorbit
\orbit{8}{10}{148}
\h{10}&\rat(3/4){3}{}&\rat(3/4){3}{}&\rat(3/4){3}{}&\int(2){\2}{}&\root\int(0){\0}{}&\int(1){2}{}&\rat(3/4){1}{}&\int(0){\0}{}&&\cnt{472}&\bounds{148}{152}{fail}{false}&&\cr
\l{1}&\rat(1/2){2}{}&\rat(3/4){3}{}&\rat(1/4){1}{}&\int(1){1}{}&\int(0){\0}{}&\rat(1/2){1}{}&\int(0){\0}{}&\int(0){\0}{}&\ln{4\sd}&\cnt{18}&\bounds{6}{6}{1}{true}&&\cr
\l{2}&\rat(1/2){2}{}&\rat(1/4){1}{}&\rat(3/4){3}{}&\int(1){3}{}&\int(0){\0}{}&\rat(1/2){3}{}&\int(0){\0}{}&\int(0){\0}{}&\ln{4\sd}&\cnt{18}&\bounds{6}{6}{1}{true}&&\cr
\l{3}&\rat(3/4){3}{}&\rat(3/4){3}{}&\int(0){\0}{}&\int(1){2}{}&\int(0){\0}{}&\rat(1/2){3}{}&\int(0){\0}{}&\int(0){1}{}&\ln{2\sd}&\cnt{16}&\bounds{4}{4}{1}{true}&&\cr
\l{4}&\rat(3/4){3}{}&\int(0){\0}{}&\rat(3/4){3}{}&\int(1){2}{}&\int(0){3}{}&\rat(1/2){1}{}&\int(0){\0}{}&\int(0){\0}{}&\ln{2\sd}&\cnt{8}&\bounds{3}{4}{5}{true}&&\cr
\l{5}&\rat(3/4){3}{}&\int(0){\0}{}&\int(0){\0}{}&\int(1){3}{}&\int(0){2}{\spec{1}}&\rat(1/2){3}{}&\rat(3/4){1}{}&\int(0){\0}{}&\ln{4\sd}&\cnt{2}&\bounds{1}{1}{1}{true}&&\cr
\l{6}&\rat(3/4){3}{}&\int(0){\0}{}&\int(0){\0}{}&\int(1){1}{}&\int(0){\0}{}&\rat(1/2){1}{}&\rat(3/4){1}{}&\int(0){2}{}&\ln{2\sd}&\cnt{12}&\bounds{3}{4}{5}{true}&&\cr
\l{7}&\int(1){\2}{}&\int(0){\0}{}&\int(0){\0}{}&\int(2){\2}{}&\int(0){\0}{}&\int(0){\0}{}&\int(0){\0}{}&\int(0){\0}{}&\ln{4}&\cnt{3}&\bounds{1}{1}{1}{true}&&\cr
\l{9}&\rat(1/2){2}{}&\rat(3/4){3}{}&\int(0){\0}{}&\int(1){3}{}&\int(0){\0}{}&\int(0){\0}{}&\rat(3/4){1}{}&\int(0){3}{}&\ln{4}&\cnt{12}&\bounds{3}{3}{1}{true}&&\cr
\l{11}&\rat(1/2){2}{}&\int(0){\0}{}&\rat(3/4){3}{}&\int(1){1}{}&\int(0){1}{}&\int(0){\0}{}&\rat(3/4){1}{}&\int(0){\0}{}&\ln{4}&\cnt{6}&\bounds{2}{2}{1}{true}&&\cr
\l{13}&\rat(1/2){2}{}&\int(0){\0}{}&\int(0){\0}{}&\int(1){2}{}&\int(0){\0}{}&\int(1){2}{}&\rat(1/2){2}{}&\int(0){\0}{}&\ln{2}&\cnt{18}&\bounds{6}{6}{1}{true}&&\cr
\l{15}&\rat(3/4){3}{}&\rat(3/4){3}{}&\rat(3/4){3}{}&\int(0){\0}{}&\int(0){\0}{}&\int(0){2}{}&\rat(3/4){1}{}&\int(0){\0}{}&\ln{1}&\cnt{4}&\bounds{2}{2}{1}{true}&&\cr
\endorbit
\orbit{8}{11}{128}
\h{11}&\root\rat(3/4){3}{}&\rat(3/4){3}{}&\rat(3/4){3}{}&\int(2){\2}{}&\int(0){\0}{}&\int(1){2}{}&\rat(3/4){1}{}&\int(0){\0}{}&&\cnt{444}&\bounds{128}{132}{fail}{false}&&\cr
\l{1}&\int(1){\2}{}&\int(0){\0}{}&\int(0){\0}{}&\int(2){\2}{}&\int(0){\0}{}&\int(0){\0}{}&\int(0){\0}{}&\int(0){\0}{}&\ln{1}&\cnt{1}&\bounds{1}{1}{1}{true}&&\cr
\l{3}&\rat(1/2){2}{}&\rat(3/4){3}{}&\rat(1/4){1}{}&\int(1){1}{}&\int(0){\0}{}&\rat(1/2){1}{}&\int(0){\0}{}&\int(0){\0}{}&\ln{2}&\cnt{6}&\bounds{2}{2}{1}{true}&&\cr
\l{5}&\rat(1/2){2}{}&\rat(3/4){3}{}&\int(0){\0}{}&\int(1){3}{}&\int(0){\0}{}&\int(0){\0}{}&\rat(3/4){1}{}&\int(0){3}{}&\ln{2}&\cnt{4}&\bounds{1}{1}{1}{true}&&\cr
\l{7}&\rat(1/2){2}{}&\int(0){\0}{}&\int(0){\0}{}&\int(1){2}{}&\int(0){\0}{}&\int(1){2}{}&\rat(1/2){2}{}&\int(0){\0}{}&\ln{1}&\cnt{6}&\bounds{2}{2}{1}{true}&&\cr
\l{9}&\rat(3/4){3}{}&\rat(1/2){2}{}&\rat(3/4){3}{}&\int(1){1}{}&\int(0){\0}{}&\int(0){\0}{}&\int(0){\0}{}&\int(0){3}{}&\ln{2}&\cnt{12}&\bounds{3}{3}{1}{true}&&\cr
\l{11}&\rat(3/4){3}{}&\rat(1/2){2}{}&\int(0){\0}{}&\int(1){3}{}&\int(0){\0}{}&\rat(1/2){1}{}&\rat(1/4){3}{}&\int(0){\0}{}&\ln{2}&\cnt{18}&\bounds{6}{6}{1}{true}&&\cr
\l{13}&\rat(3/4){3}{}&\rat(3/4){3}{}&\rat(3/4){3}{}&\int(0){\0}{}&\int(0){\0}{}&\int(1){2}{}&\rat(-1/4){1}{}&\int(0){\0}{}&\ln{1}&\cnt{3}&\bounds{1}{1}{1}{true}&&\cr
\l{15}&\rat(3/4){3}{}&\rat(3/4){3}{}&\rat(3/4){3}{}&\int(0){\0}{}&\int(0){\0}{}&\int(0){2}{}&\rat(3/4){1}{}&\int(0){\0}{}&\ln{1}&\cnt{4}&\bounds{2}{2}{1}{true}&&\cr
\l{17}&\rat(3/4){3}{}&\rat(3/4){3}{}&\rat(-1/4){3}{}&\int(0){\0}{}&\int(0){\0}{}&\int(1){2}{}&\rat(3/4){1}{}&\int(0){\0}{}&\ln{2}&\cnt{3}&\bounds{1}{1}{1}{true}&&\cr
\l{19}&\rat(3/4){3}{}&\rat(3/4){3}{}&\int(0){\0}{}&\int(1){2}{}&\int(0){\0}{}&\rat(1/2){3}{}&\int(0){\0}{}&\int(0){1}{}&\ln{2}&\cnt{16}&\bounds{4}{4}{1}{true}&&\cr
\l{21}&\rat(3/4){3}{}&\rat(3/4){3}{}&\int(0){\0}{}&\int(1){1}{}&\int(0){3}{}&\int(0){\0}{}&\rat(1/2){2}{}&\int(0){\0}{}&\ln{2}&\cnt{12}&\bounds{3}{3}{1}{true}&&\cr
\l{23}&\rat(3/4){3}{}&\rat(1/4){1}{}&\rat(1/4){1}{}&\int(1){2}{}&\int(0){\0}{}&\int(0){\0}{}&\rat(3/4){1}{}&\int(0){\0}{}&\ln{1}&\cnt{18}&\bounds{6}{6}{1}{true}&&\cr
\l{25}&\rat(3/4){3}{}&\rat(1/4){1}{}&\int(0){\0}{}&\int(1){1}{}&\int(0){1}{}&\int(1){2}{}&\int(0){\0}{}&\int(0){\0}{}&\ln{2}&\cnt{12}&\bounds{3}{3}{1}{true}&&\cr
\l{27}&\rat(3/4){3}{}&\int(0){\0}{}&\int(0){\0}{}&\int(1){3}{}&\int(0){2}{}&\rat(1/2){3}{}&\rat(3/4){1}{}&\int(0){\0}{}&\ln{2}&\cnt{12}&\bounds{3}{4}{5}{true}&&\cr
\endorbit
\orbit{8}{12}{120}
\h{12}&\rat(3/4){3}{}&\rat(3/4){3}{}&\rat(3/4){3}{}&\root\int(2){\2}{}&\int(0){\0}{}&\int(1){2}{}&\rat(3/4){1}{}&\int(0){\0}{}&&\cnt{416}&\bounds{120}{124}{fail}{false}&&\cr
\l{1}&\rat(1/2){2}{}&\rat(3/4){3}{}&\rat(1/4){1}{}&\int(1){1}{}&\int(0){\0}{}&\rat(1/2){1}{}&\int(0){\0}{}&\int(0){\0}{}&\ln{8\sd}&\cnt{18}&\bounds{6}{6}{1}{true}&&\cr
\l{2}&\rat(3/4){3}{}&\int(0){\0}{}&\int(0){\0}{}&\int(1){3}{}&\int(0){2}{}&\rat(1/2){3}{}&\rat(3/4){1}{}&\int(0){\0}{}&\ln{4\sd}&\cnt{12}&\bounds{3}{4}{5}{true}&&\cr
\l{3}&\int(1){\2}{}&\int(0){\0}{}&\int(0){\0}{}&\int(2){\2}{}&\int(0){\0}{}&\int(0){\0}{}&\int(0){\0}{}&\int(0){\0}{}&\ln{4}&\cnt{3}&\bounds{1}{1}{1}{true}&&\cr
\l{5}&\rat(1/2){2}{}&\rat(3/4){3}{}&\int(0){\0}{}&\int(1){3}{}&\int(0){\0}{}&\int(0){\0}{}&\rat(3/4){1}{}&\int(0){3}{}&\ln{8}&\cnt{12}&\bounds{3}{3}{1}{true}&&\cr
\l{7}&\rat(3/4){3}{}&\rat(3/4){3}{}&\rat(3/4){3}{}&\int(0){\0}{}&\int(0){\0}{}&\int(0){2}{}&\rat(3/4){1}{}&\int(0){\0}{}&\ln{1}&\cnt{4}&\bounds{2}{2}{1}{true}&&\cr
\endorbit
\orbit{8}{13}{118}
\h{13}&\rat(3/4){3}{}&\rat(3/4){3}{}&\rat(3/4){3}{}&\int(2){\2}{}&\int(0){\0}{}&\root\int(1){2}{\spec{1}}&\rat(3/4){1}{}&\int(0){\0}{}&&\cnt{416}&\bounds{118}{120}{fail}{false}&&\cr
\l{1}&\rat(1/2){2}{}&\rat(1/4){1}{}&\rat(3/4){3}{}&\int(1){3}{}&\int(0){\0}{}&\rat(1/2){3}{}&\int(0){\0}{}&\int(0){\0}{}&\ln{4\sd}&\cnt{18}&\bounds{6}{6}{1}{true}&&\cr
\l{2}&\rat(3/4){3}{}&\rat(3/4){3}{}&\int(0){\0}{}&\int(1){2}{}&\int(0){\0}{}&\rat(1/2){3}{}&\int(0){\0}{}&\int(0){1}{}&\ln{2\sd}&\cnt{16}&\bounds{4}{4}{1}{true}&&\cr
\l{3}&\rat(3/4){3}{}&\int(0){\0}{}&\int(0){\0}{}&\int(1){3}{}&\int(0){2}{}&\rat(1/2){3}{}&\rat(3/4){1}{}&\int(0){\0}{}&\ln{2\sd}&\cnt{12}&\bounds{3}{4}{5}{true}&&\cr
\l{4}&\int(1){\2}{}&\int(0){\0}{}&\int(0){\0}{}&\int(2){\2}{}&\int(0){\0}{}&\int(0){\0}{}&\int(0){\0}{}&\int(0){\0}{}&\ln{4}&\cnt{3}&\bounds{1}{1}{1}{true}&&\cr
\l{6}&\rat(1/2){2}{}&\rat(3/4){3}{}&\int(0){\0}{}&\int(1){3}{}&\int(0){\0}{}&\int(0){\0}{}&\rat(3/4){1}{}&\int(0){3}{}&\ln{4}&\cnt{12}&\bounds{3}{3}{1}{true}&&\cr
\l{8}&\rat(1/2){2}{}&\int(0){\0}{}&\rat(3/4){3}{}&\int(1){1}{}&\int(0){1}{}&\int(0){\0}{}&\rat(3/4){1}{}&\int(0){\0}{}&\ln{4}&\cnt{12}&\bounds{3}{3}{1}{true}&&\cr
\l{10}&\rat(1/2){2}{}&\int(0){\0}{}&\int(0){\0}{}&\int(1){2}{}&\int(0){\0}{}&\int(1){2}{}&\rat(1/2){2}{}&\int(0){\0}{}&\ln{2}&\cnt{18}&\bounds{6}{6}{1}{true}&&\cr
\endorbit
\orbit{8}{14}{140}
\h{14}&\root\rat(11/4){3}{\spec{2}}&\rat(3/4){3}{}&\rat(3/4){3}{}&\int(0){\0}{}&\int(0){\0}{}&\int(1){2}{}&\rat(3/4){1}{}&\int(0){\0}{}&&\cnt{548}&\bounds{140}{152}{fail}{false}&&\cr
\l{1}&\int(3){\4}{}&\int(0){\0}{}&\int(0){\0}{}&\int(0){\0}{}&\int(0){\0}{}&\int(0){\0}{}&\int(0){\0}{}&\int(0){\0}{}&\ln{1}&\cnt{1}&\bounds{1}{1}{1}{true}&&\cr
\l{3}&\rat(3/2){2}{}&\rat(1/2){2}{}&\rat(1/2){2}{}&\int(0){\0}{}&\int(0){\0}{}&\int(0){\0}{}&\rat(1/2){2}{}&\int(0){\0}{}&\ln{1}&\cnt{27}&\bounds{9}{9}{1}{true}&&\cr
\l{5}&\rat(3/2){2}{}&\rat(1/2){2}{}&\int(0){\0}{}&\int(0){\0}{}&\int(0){\0}{}&\int(1){2}{}&\int(0){\0}{}&\int(0){2}{}&\ln{3}&\cnt{18}&\bounds{4}{6}{5}{true}&&\cr
\l{7}&\rat(3/2){2}{}&\rat(3/4){3}{}&\rat(3/4){3}{}&\int(0){\0}{}&\int(0){3}{}&\int(0){\0}{}&\int(0){\0}{}&\int(0){1}{}&\ln{3}&\cnt{16}&\bounds{4}{4}{1}{true}&&\cr
\l{9}&\rat(3/2){2}{}&\rat(3/4){3}{}&\rat(1/4){1}{}&\int(0){1}{}&\int(0){\0}{}&\rat(1/2){1}{}&\int(0){\0}{}&\int(0){\0}{}&\ln{6}&\cnt{24}&\bounds{6}{6}{1}{true}&&\cr
\endorbit
\orbit{8}{15}{138}
\h{15}&\rat(11/4){3}{\spec{2}}&\rat(3/4){3}{}&\rat(3/4){3}{}&\root\int(0){\0}{}&\int(0){\0}{}&\int(1){2}{}&\rat(3/4){1}{}&\int(0){\0}{}&&\cnt{496}&\bounds{138}{148}{fail}{false}&&\cr
\l{1}&\int(3){\4}{}&\int(0){\0}{}&\int(0){\0}{}&\int(0){\0}{}&\int(0){\0}{}&\int(0){\0}{}&\int(0){\0}{}&\int(0){\0}{}&\ln{1}&\cnt{1}&\bounds{1}{1}{1}{true}&&\cr
\l{3}&\int(2){\2}{}&\int(1){\2}{}&\int(0){\0}{}&\int(0){\0}{}&\int(0){\0}{}&\int(0){\0}{}&\int(0){\0}{}&\int(0){\0}{}&\ln{2}&\cnt{6}&\bounds{2}{2}{1}{true}&&\cr
\l{5}&\int(2){\2}{}&\int(0){\0}{}&\int(0){\0}{}&\int(0){\0}{}&\int(0){\0}{}&\int(1){\2}{}&\int(0){\0}{}&\int(0){\0}{}&\ln{1}&\cnt{8}&\bounds{3}{4}{5}{true}&&\cr
\l{7}&\int(2){\2}{}&\int(0){\0}{}&\int(0){\0}{}&\int(0){\0}{}&\int(0){\0}{}&\int(0){\0}{}&\int(1){\2}{}&\int(0){\0}{}&\ln{1}&\cnt{6}&\bounds{2}{2}{1}{true}&&\cr
\l{9}&\rat(3/2){2}{}&\rat(1/2){2}{}&\rat(1/2){2}{}&\int(0){\0}{}&\int(0){\0}{}&\int(0){\0}{}&\rat(1/2){2}{}&\int(0){\0}{}&\ln{1}&\cnt{27}&\bounds{9}{9}{1}{true}&&\cr
\l{11}&\rat(3/2){2}{}&\rat(1/2){2}{}&\int(0){\0}{}&\int(0){\0}{}&\int(0){\0}{}&\int(1){2}{}&\int(0){\0}{}&\int(0){2}{}&\ln{2}&\cnt{18}&\bounds{4}{6}{5}{true}&&\cr
\l{13}&\rat(3/2){2}{}&\rat(3/4){3}{}&\rat(3/4){3}{}&\int(0){\0}{}&\int(0){3}{}&\int(0){\0}{}&\int(0){\0}{}&\int(0){1}{}&\ln{1}&\cnt{16}&\bounds{4}{4}{1}{true}&&\cr
\l{15}&\rat(3/2){2}{}&\rat(3/4){3}{}&\rat(1/4){1}{}&\int(0){1}{}&\int(0){\0}{}&\rat(1/2){1}{}&\int(0){\0}{}&\int(0){\0}{}&\ln{2}&\cnt{12}&\bounds{4}{4}{1}{true}&&\cr
\l{17}&\rat(3/2){2}{}&\rat(3/4){3}{}&\int(0){\0}{}&\int(0){3}{}&\int(0){\0}{}&\int(0){\0}{}&\rat(3/4){1}{}&\int(0){3}{}&\ln{2}&\cnt{8}&\bounds{2}{2}{1}{true}&&\cr
\l{19}&\rat(3/2){2}{}&\rat(3/4){3}{}&\int(0){\0}{}&\int(0){\0}{}&\int(0){1}{}&\rat(1/2){3}{}&\rat(1/4){3}{}&\int(0){\0}{}&\ln{2}&\cnt{24}&\bounds{6}{6}{1}{true}&&\cr
\l{21}&\rat(3/2){2}{}&\rat(1/4){1}{}&\int(0){\0}{}&\int(0){\0}{}&\int(0){3}{}&\rat(1/2){1}{}&\rat(3/4){1}{}&\int(0){\0}{}&\ln{2}&\cnt{24}&\bounds{6}{6}{1}{true}&&\cr
\l{23}&\rat(3/2){2}{}&\int(0){\0}{}&\int(0){\0}{}&\int(0){2}{\spec{1}}&\int(0){\0}{}&\int(1){2}{}&\rat(1/2){2}{}&\int(0){\0}{}&\ln{2}&\cnt{3}&\bounds{1}{1}{1}{true}&&\cr
\endorbit
\orbit{8}{16}{124}
\h{16}&\rat(11/4){3}{\spec{2}}&\root\rat(3/4){3}{}&\rat(3/4){3}{}&\int(0){\0}{}&\int(0){\0}{}&\int(1){2}{}&\rat(3/4){1}{}&\int(0){\0}{}&&\cnt{468}&\bounds{124}{134}{fail}{false}&&\cr
\l{1}&\int(3){\4}{}&\int(0){\0}{}&\int(0){\0}{}&\int(0){\0}{}&\int(0){\0}{}&\int(0){\0}{}&\int(0){\0}{}&\int(0){\0}{}&\ln{1}&\cnt{1}&\bounds{1}{1}{1}{true}&&\cr
\l{3}&\int(2){\2}{}&\int(1){\2}{}&\int(0){\0}{}&\int(0){\0}{}&\int(0){\0}{}&\int(0){\0}{}&\int(0){\0}{}&\int(0){\0}{}&\ln{1}&\cnt{2}&\bounds{1}{1}{1}{true}&&\cr
\l{5}&\int(2){\2}{}&\int(0){\0}{}&\int(1){\2}{}&\int(0){\0}{}&\int(0){\0}{}&\int(0){\0}{}&\int(0){\0}{}&\int(0){\0}{}&\ln{2}&\cnt{6}&\bounds{2}{2}{1}{true}&&\cr
\l{7}&\int(2){\2}{}&\int(0){\0}{}&\int(0){\0}{}&\int(0){\0}{}&\int(0){\0}{}&\int(1){\2}{}&\int(0){\0}{}&\int(0){\0}{}&\ln{1}&\cnt{8}&\bounds{3}{4}{5}{true}&&\cr
\l{9}&\rat(3/2){2}{}&\rat(1/2){2}{}&\rat(1/2){2}{}&\int(0){\0}{}&\int(0){\0}{}&\int(0){\0}{}&\rat(1/2){2}{}&\int(0){\0}{}&\ln{1}&\cnt{9}&\bounds{3}{3}{1}{true}&&\cr
\l{11}&\rat(3/2){2}{}&\rat(1/2){2}{}&\int(0){\0}{}&\int(0){\0}{}&\int(0){\0}{}&\int(1){2}{}&\int(0){\0}{}&\int(0){2}{}&\ln{1}&\cnt{6}&\bounds{2}{2}{1}{true}&&\cr
\l{13}&\rat(3/2){2}{}&\rat(3/4){3}{}&\rat(3/4){3}{}&\int(0){\0}{}&\int(0){3}{}&\int(0){\0}{}&\int(0){\0}{}&\int(0){1}{}&\ln{2}&\cnt{16}&\bounds{4}{4}{1}{true}&&\cr
\l{15}&\rat(3/2){2}{}&\rat(3/4){3}{}&\rat(1/4){1}{}&\int(0){1}{}&\int(0){\0}{}&\rat(1/2){1}{}&\int(0){\0}{}&\int(0){\0}{}&\ln{2}&\cnt{24}&\bounds{6}{6}{1}{true}&&\cr
\l{17}&\rat(3/2){2}{}&\rat(1/4){1}{}&\rat(3/4){3}{}&\int(0){3}{}&\int(0){\0}{}&\rat(1/2){3}{}&\int(0){\0}{}&\int(0){\0}{}&\ln{2}&\cnt{8}&\bounds{2}{2}{1}{true}&&\cr
\l{19}&\rat(3/2){2}{}&\int(0){\0}{}&\rat(1/2){2}{}&\int(0){\0}{}&\int(0){2}{}&\int(1){2}{}&\int(0){\0}{}&\int(0){\0}{}&\ln{2}&\cnt{18}&\bounds{4}{6}{5}{true}&&\cr
\l{21}&\rat(3/2){2}{}&\int(0){\0}{}&\rat(3/4){3}{}&\int(0){1}{}&\int(0){1}{}&\int(0){\0}{}&\rat(3/4){1}{}&\int(0){\0}{}&\ln{1}&\cnt{16}&\bounds{4}{4}{1}{true}&&\cr
\l{23}&\rat(3/2){2}{}&\int(0){\0}{}&\rat(3/4){3}{}&\int(0){\0}{}&\int(0){\0}{}&\rat(1/2){1}{}&\rat(1/4){3}{}&\int(0){3}{}&\ln{2}&\cnt{24}&\bounds{6}{6}{1}{true}&&\cr
\endorbit
\orbit{8}{17}{116}
\h{17}&\rat(11/4){3}{\spec{2}}&\rat(3/4){3}{}&\rat(3/4){3}{}&\int(0){\0}{}&\int(0){\0}{}&\root\int(1){2}{\spec{1}}&\rat(3/4){1}{}&\int(0){\0}{}&&\cnt{440}&\bounds{116}{128}{fail}{false}&&\cr
\l{1}&\int(3){\4}{}&\int(0){\0}{}&\int(0){\0}{}&\int(0){\0}{}&\int(0){\0}{}&\int(0){\0}{}&\int(0){\0}{}&\int(0){\0}{}&\ln{1}&\cnt{1}&\bounds{1}{1}{1}{true}&&\cr
\l{3}&\int(2){\2}{}&\int(1){\2}{}&\int(0){\0}{}&\int(0){\0}{}&\int(0){\0}{}&\int(0){\0}{}&\int(0){\0}{}&\int(0){\0}{}&\ln{3}&\cnt{6}&\bounds{2}{2}{1}{true}&&\cr
\l{5}&\rat(3/2){2}{}&\rat(1/2){2}{}&\rat(1/2){2}{}&\int(0){\0}{}&\int(0){\0}{}&\int(0){\0}{}&\rat(1/2){2}{}&\int(0){\0}{}&\ln{1}&\cnt{27}&\bounds{9}{9}{1}{true}&&\cr
\l{7}&\rat(3/2){2}{}&\rat(1/2){2}{}&\int(0){\0}{}&\int(0){\0}{}&\int(0){\0}{}&\int(1){2}{}&\int(0){\0}{}&\int(0){2}{}&\ln{3}&\cnt{18}&\bounds{4}{6}{5}{true}&&\cr
\l{9}&\rat(3/2){2}{}&\rat(3/4){3}{}&\rat(3/4){3}{}&\int(0){\0}{}&\int(0){3}{}&\int(0){\0}{}&\int(0){\0}{}&\int(0){1}{}&\ln{3}&\cnt{16}&\bounds{4}{4}{1}{true}&&\cr
\l{11}&\rat(3/2){2}{}&\rat(3/4){3}{}&\int(0){\0}{}&\int(0){\0}{}&\int(0){1}{}&\rat(1/2){3}{}&\rat(1/4){3}{}&\int(0){\0}{}&\ln{3}&\cnt{24}&\bounds{6}{6}{1}{true}&&\cr
\endorbit
\orbit{8}{18}{146}
\h{18}&\rat(3/4){3}{}&\rat(3/4){3}{}&\rat(3/4){3}{}&\int(0){\0}{}&\int(0){\0}{}&\root\int(3){2}{\spec{1}}&\rat(3/4){1}{}&\int(0){\0}{}&&\cnt{600}&\bounds{146}{152}{fail}{false}&&\cr
\l{1}&\rat(1/2){2}{}&\rat(3/4){3}{}&\int(0){\0}{}&\int(0){\0}{}&\int(0){1}{}&\rat(3/2){3}{}&\rat(1/4){3}{}&\int(0){\0}{}&\ln{12\sd}&\cnt{36}&\bounds{9}{9}{1}{true}&&\cr
\l{2}&\rat(3/4){3}{}&\rat(3/4){3}{}&\int(0){\0}{}&\int(0){2}{}&\int(0){\0}{}&\rat(3/2){3}{}&\int(0){\0}{}&\int(0){1}{}&\ln{6\sd}&\cnt{24}&\bounds{5}{6}{5}{true}&&\cr
\l{3}&\int(1){\2}{}&\int(0){\0}{}&\int(0){\0}{}&\int(0){\0}{}&\int(0){\0}{}&\int(2){\2}{}&\int(0){\0}{}&\int(0){\0}{}&\ln{4}&\cnt{3}&\bounds{1}{1}{1}{true}&&\cr
\endorbit
\orbit{8}{19}{140}
\h{19}&\rat(3/4){3}{}&\rat(3/4){3}{}&\rat(3/4){3}{}&\root\int(0){\0}{}&\int(0){\0}{}&\int(3){2}{\spec{1}}&\rat(3/4){1}{}&\int(0){\0}{}&&\cnt{520}&\bounds{140}{144}{fail}{false}&&\cr
\l{1}&\rat(1/2){2}{}&\rat(3/4){3}{}&\int(0){\0}{}&\int(0){\0}{}&\int(0){1}{}&\rat(3/2){3}{}&\rat(1/4){3}{}&\int(0){\0}{}&\ln{4\sd}&\cnt{36}&\bounds{9}{9}{1}{true}&&\cr
\l{2}&\rat(1/2){2}{}&\rat(1/4){1}{}&\rat(3/4){3}{}&\int(0){3}{}&\int(0){\0}{}&\rat(3/2){3}{}&\int(0){\0}{}&\int(0){\0}{}&\ln{4\sd}&\cnt{18}&\bounds{6}{6}{1}{true}&&\cr
\l{3}&\rat(1/2){2}{}&\int(0){\0}{}&\rat(1/4){1}{}&\int(0){\0}{}&\int(0){\0}{}&\rat(3/2){3}{}&\rat(3/4){1}{}&\int(0){1}{}&\ln{4\sd}&\cnt{36}&\bounds{9}{9}{1}{true}&&\cr
\l{4}&\rat(3/4){3}{}&\rat(3/4){3}{}&\int(0){\0}{}&\int(0){2}{\spec{1}}&\int(0){\0}{}&\rat(3/2){3}{}&\int(0){\0}{}&\int(0){1}{}&\ln{4\sd}&\cnt{4}&\bounds{1}{1}{1}{true}&&\cr
\l{5}&\rat(3/4){3}{}&\int(0){\0}{}&\rat(3/4){3}{}&\int(0){\0}{}&\int(0){1}{}&\rat(3/2){3}{}&\int(0){\0}{}&\int(0){2}{}&\ln{2\sd}&\cnt{24}&\bounds{5}{6}{5}{true}&&\cr
\l{6}&\rat(3/4){3}{}&\int(0){\0}{}&\int(0){\0}{}&\int(0){3}{}&\int(0){2}{}&\rat(3/2){3}{}&\rat(3/4){1}{}&\int(0){\0}{}&\ln{2\sd}&\cnt{12}&\bounds{3}{4}{5}{true}&&\cr
\l{7}&\int(1){\2}{}&\int(0){\0}{}&\int(0){\0}{}&\int(0){\0}{}&\int(0){\0}{}&\int(2){\2}{}&\int(0){\0}{}&\int(0){\0}{}&\ln{4}&\cnt{9}&\bounds{3}{3}{1}{true}&&\cr
\endorbit
\orbit{8}{20}{122}
\h{20}&\root\rat(3/4){3}{}&\rat(3/4){3}{}&\rat(3/4){3}{}&\int(0){\0}{}&\int(0){\0}{}&\int(3){2}{\spec{1}}&\rat(3/4){1}{}&\int(0){\0}{}&&\cnt{492}&\bounds{122}{128}{fail}{false}&&\cr
\l{1}&\int(1){\2}{}&\int(0){\0}{}&\int(0){\0}{}&\int(0){\0}{}&\int(0){\0}{}&\int(2){\2}{}&\int(0){\0}{}&\int(0){\0}{}&\ln{1}&\cnt{3}&\bounds{1}{1}{1}{true}&&\cr
\l{3}&\rat(1/2){2}{}&\rat(3/4){3}{}&\int(0){\0}{}&\int(0){\0}{}&\int(0){1}{}&\rat(3/2){3}{}&\rat(1/4){3}{}&\int(0){\0}{}&\ln{3}&\cnt{12}&\bounds{3}{3}{1}{true}&&\cr
\l{5}&\rat(3/4){3}{}&\rat(1/2){2}{}&\rat(1/4){1}{}&\int(0){\0}{}&\int(0){3}{}&\rat(3/2){3}{}&\int(0){\0}{}&\int(0){\0}{}&\ln{3}&\cnt{36}&\bounds{9}{9}{1}{true}&&\cr
\l{7}&\rat(3/4){3}{}&\rat(3/4){3}{}&\rat(3/4){3}{}&\int(0){\0}{}&\int(0){\0}{}&\int(1){2}{}&\rat(-1/4){1}{}&\int(0){\0}{}&\ln{3}&\cnt{9}&\bounds{3}{3}{1}{true}&&\cr
\l{9}&\rat(3/4){3}{}&\rat(3/4){3}{}&\int(0){\0}{}&\int(0){2}{}&\int(0){\0}{}&\rat(3/2){3}{}&\int(0){\0}{}&\int(0){1}{}&\ln{3}&\cnt{24}&\bounds{5}{6}{5}{true}&&\cr
\endorbit
\orbit{8}{21}{140}
\h{21}&\int(1){2}{}&\int(1){2}{}&\int(1){2}{}&\int(2){\2}{}&\root\int(0){\0}{}&\int(0){\0}{}&\int(1){2}{}&\int(0){\0}{}&&\cnt{472}&\bounds{140}{164}{fail}{false}&&\cr
\l{1}&\int(1){2}{}&\int(1){2}{}&\int(0){\0}{}&\int(1){2}{}&\int(0){2}{\spec{1}}&\int(0){\0}{}&\int(0){\0}{}&\int(0){\0}{}&\ln{4\sd}&\cnt{2}&\bounds{1}{1}{1}{true}&&\cr
\l{2}&\int(1){2}{}&\rat(1/2){3}{}&\rat(1/2){1}{}&\int(1){1}{}&\int(0){\0}{}&\int(0){1}{}&\int(0){\0}{}&\int(0){\0}{}&\ln{16\sd}&\cnt{16}&\bounds{4}{4}{1}{true}&&\cr
\l{3}&\int(1){2}{}&\int(0){\0}{}&\int(1){2}{}&\int(1){2}{}&\int(0){\0}{}&\int(0){\0}{}&\int(0){\0}{}&\int(0){2}{}&\ln{4\sd}&\cnt{12}&\bounds{3}{4}{5}{true}&&\cr
\l{4}&\int(1){2}{}&\int(0){\0}{}&\rat(1/2){3}{}&\int(1){1}{}&\int(0){1}{}&\int(0){\0}{}&\rat(1/2){1}{}&\int(0){\0}{}&\ln{8\sd}&\cnt{8}&\bounds{3}{4}{5}{true}&&\cr
\l{5}&\rat(1/2){3}{}&\rat(1/2){1}{}&\rat(1/2){1}{}&\int(1){2}{}&\int(0){\0}{}&\int(0){\0}{}&\rat(1/2){1}{}&\int(0){\0}{}&\ln{2\sdd}&\cnt{32}&\bounds{10}{16}{2}{true}&&\cr
\l{6}&\int(1){\2}{}&\int(0){\0}{}&\int(0){\0}{}&\int(2){\2}{}&\int(0){\0}{}&\int(0){\0}{}&\int(0){\0}{}&\int(0){\0}{}&\ln{4}&\cnt{4}&\bounds{2}{2}{1}{true}&&\cr
\endorbit
\orbit{8}{22}{112}
\h{22}&\root\int(1){2}{\spec{1}}&\int(1){2}{}&\int(1){2}{}&\int(2){\2}{}&\int(0){\0}{}&\int(0){\0}{}&\int(1){2}{}&\int(0){\0}{}&&\cnt{416}&\bounds{112}{124}{fail}{false}&&\cr
\l{1}&\rat(1/2){3}{}&\int(1){2}{}&\rat(1/2){3}{}&\int(1){1}{}&\int(0){\0}{}&\int(0){\0}{}&\int(0){\0}{}&\int(0){3}{}&\ln{6\sd}&\cnt{16}&\bounds{4}{4}{1}{true}&&\cr
\l{2}&\rat(1/2){3}{}&\rat(1/2){1}{}&\rat(1/2){1}{}&\int(1){2}{}&\int(0){\0}{}&\int(0){\0}{}&\rat(1/2){1}{}&\int(0){\0}{}&\ln{1\sdd}&\cnt{32}&\bounds{10}{16}{2}{true}&&\cr
\l{3}&\int(1){2}{}&\int(1){2}{}&\int(1){2}{}&\int(0){\0}{}&\int(0){\0}{}&\int(0){\0}{}&\int(0){2}{}&\int(0){\0}{}&\ln{3}&\cnt{4}&\bounds{2}{2}{1}{true}&&\cr
\l{5}&\int(1){2}{}&\int(1){2}{}&\int(0){\0}{}&\int(1){2}{}&\int(0){2}{}&\int(0){\0}{}&\int(0){\0}{}&\int(0){\0}{}&\ln{3}&\cnt{12}&\bounds{3}{4}{5}{true}&&\cr
\l{7}&\int(1){2}{}&\rat(1/2){3}{}&\rat(1/2){1}{}&\int(1){1}{}&\int(0){\0}{}&\int(0){1}{}&\int(0){\0}{}&\int(0){\0}{}&\ln{6}&\cnt{16}&\bounds{4}{4}{1}{true}&&\cr
\endorbit
\orbit{8}{23}{112}
\h{23}&\int(1){2}{}&\int(1){2}{}&\int(1){2}{}&\root\int(2){\2}{}&\int(0){\0}{}&\int(0){\0}{}&\int(1){2}{}&\int(0){\0}{}&&\cnt{416}&\bounds{112}{112}{1}{false}&&\cr
\l{1}&\int(1){2}{}&\rat(1/2){3}{}&\rat(1/2){1}{}&\int(1){1}{}&\int(0){\0}{}&\int(0){1}{}&\int(0){\0}{}&\int(0){\0}{}&\ln{24\sd}&\cnt{16}&\bounds{4}{4}{1}{true}&&\cr
\l{2}&\int(1){\2}{}&\int(0){\0}{}&\int(0){\0}{}&\int(2){\2}{}&\int(0){\0}{}&\int(0){\0}{}&\int(0){\0}{}&\int(0){\0}{}&\ln{4}&\cnt{4}&\bounds{2}{2}{1}{true}&&\cr
\endorbit
\orbit{8}{24}{126}
\h{24}&\int(2){\2}{}&\int(2){\2}{}&\int(2){\2}{}&\int(0){\0}{}&\int(0){\0}{}&\int(0){\0}{}&\root\int(0){\0}{}&\int(0){\0}{}&&\cnt{520}&\bounds{126}{148}{fail}{false}&&\cr
\l{1}&\int(2){\2}{}&\int(1){\2}{\spec{1}}&\int(0){\0}{}&\int(0){\0}{}&\int(0){\0}{}&\int(0){\0}{}&\int(0){\0}{}&\int(0){\0}{}&\ln{12\sd}&\cnt{2}&\bounds{1}{1}{1}{true}&&\cr
\l{2}&\int(1){2}{}&\int(1){2}{}&\int(1){2}{}&\int(0){\0}{}&\int(0){\0}{}&\int(0){\0}{}&\int(0){2}{\spec{1}}&\int(0){\0}{}&\ln{2\sd}&\cnt{8}&\bounds{3}{4}{5}{true}&&\cr
\l{3}&\int(1){2}{}&\int(1){3}{}&\int(1){3}{}&\int(0){\0}{}&\int(0){3}{}&\int(0){\0}{}&\int(0){\0}{}&\int(0){1}{}&\ln{12\sd}&\cnt{32}&\bounds{7}{8}{5}{true}&&\cr
\l{4}&\int(1){3}{}&\int(1){3}{}&\int(1){3}{}&\int(0){\0}{}&\int(0){\0}{}&\int(0){2}{}&\int(0){1}{}&\int(0){\0}{}&\ln{8\sd}&\cnt{12}&\bounds{3}{4}{5}{true}&&\cr
\endorbit
\orbit{8}{25}{122}
\h{25}&\int(2){\2}{}&\int(2){\2}{}&\int(2){\2}{}&\root\int(0){\0}{}&\int(0){\0}{}&\int(0){\0}{}&\int(0){\0}{}&\int(0){\0}{}&&\cnt{520}&\bounds{122}{140}{fail}{false}&&\cr
\l{1}&\int(2){\2}{}&\int(1){\2}{\spec{1}}&\int(0){\0}{}&\int(0){\0}{}&\int(0){\0}{}&\int(0){\0}{}&\int(0){\0}{}&\int(0){\0}{}&\ln{12\sd}&\cnt{2}&\bounds{1}{1}{1}{true}&&\cr
\l{2}&\int(1){2}{}&\int(1){2}{}&\int(1){2}{}&\int(0){\0}{}&\int(0){\0}{}&\int(0){\0}{}&\int(0){2}{}&\int(0){\0}{}&\ln{1\sdd}&\cnt{48}&\bounds{10}{16}{5}{true}&&\cr
\l{3}&\int(1){2}{}&\int(1){3}{}&\int(1){3}{}&\int(0){\0}{}&\int(0){3}{}&\int(0){\0}{}&\int(0){\0}{}&\int(0){1}{}&\ln{6\sd}&\cnt{32}&\bounds{7}{8}{5}{true}&&\cr
\l{4}&\int(1){2}{}&\int(1){3}{}&\int(1){1}{}&\int(0){1}{}&\int(0){\0}{}&\int(0){1}{}&\int(0){\0}{}&\int(0){\0}{}&\ln{6\sd}&\cnt{16}&\bounds{4}{4}{1}{true}&&\cr
\l{5}&\int(1){3}{}&\int(1){3}{}&\int(1){3}{}&\int(0){\0}{}&\int(0){\0}{}&\int(0){2}{}&\int(0){1}{}&\int(0){\0}{}&\ln{6\sd}&\cnt{24}&\bounds{5}{6}{5}{true}&&\cr
\l{6}&\int(1){3}{}&\int(1){1}{}&\int(1){1}{}&\int(0){2}{\spec{1}}&\int(0){\0}{}&\int(0){\0}{}&\int(0){1}{}&\int(0){\0}{}&\ln{4\sd}&\cnt{4}&\bounds{1}{1}{1}{true}&&\cr
\endorbit
\orbit{8}{26}{104}
\h{26}&\root\int(2){\2}{}&\int(2){\2}{}&\int(2){\2}{}&\int(0){\0}{}&\int(0){\0}{}&\int(0){\0}{}&\int(0){\0}{}&\int(0){\0}{}&&\cnt{464}&\bounds{104}{120}{fail}{false}&&\cr
\l{1}&\int(1){3}{}&\int(1){2}{}&\int(1){3}{}&\int(0){1}{}&\int(0){\0}{}&\int(0){\0}{}&\int(0){\0}{}&\int(0){3}{}&\ln{8\sd}&\cnt{32}&\bounds{7}{8}{5}{true}&&\cr
\l{2}&\int(1){3}{}&\int(1){3}{}&\int(1){3}{}&\int(0){\0}{}&\int(0){\0}{}&\int(0){2}{}&\int(0){1}{}&\int(0){\0}{}&\ln{8\sd}&\cnt{24}&\bounds{5}{6}{5}{true}&&\cr
\l{3}&\int(2){\2}{}&\int(1){\2}{\spec{1}}&\int(0){\0}{}&\int(0){\0}{}&\int(0){\0}{}&\int(0){\0}{}&\int(0){\0}{}&\int(0){\0}{}&\ln{4}&\cnt{2}&\bounds{1}{1}{1}{true}&&\cr
\endorbit
\orbit{8}{27}{140}
\h{27}&\root\rat(3/4){3}{}&\rat(3/4){3}{}&\rat(3/4){3}{}&\rat(3/4){3}{}&\rat(3/4){3}{}&\rat(3/4){3}{}&\rat(3/4){3}{}&\rat(3/4){3}{}&&\cnt{420}&\bounds{140}{140}{1}{false}&&\cr
\l{1}&\rat(1/2){2}{}&\rat(3/4){3}{}&\rat(3/4){3}{}&\int(0){\0}{}&\rat(3/4){3}{}&\int(0){\0}{}&\int(0){\0}{}&\rat(1/4){1}{}&\ln{7}&\cnt{3}&\bounds{1}{1}{1}{true}&&\cr
\l{3}&\rat(3/4){3}{}&\rat(1/2){2}{}&\rat(3/4){3}{}&\rat(1/4){1}{}&\int(0){\0}{}&\int(0){\0}{}&\int(0){\0}{}&\rat(3/4){3}{}&\ln{21}&\cnt{9}&\bounds{3}{3}{1}{true}&&\cr
\endorbit
\orbit{8}{28}{42}
\h{28}&\int(6){\6}{\spec{1}}&\int(0){\0}{}&\int(0){\0}{}&\int(0){\0}{}&\int(0){\0}{}&\int(0){\0}{}&\int(0){\0}{}&\int(0){\0}{}&&\cnt{168}&\bounds{42}{56}{fail}{false}&&\cr
\l{1}&\int(3){\2}{}&\int(0){\2}{}&\int(0){\0}{}&\int(0){\0}{}&\int(0){\0}{}&\int(0){\0}{}&\int(0){\0}{}&\int(0){\0}{}&\ln{7\sdd}&\cnt{24}&\bounds{6}{8}{5}{true}&&\cr
\endorbit
\endorbits

\subsection{Configuration~\hlink{8A3}{1}}\label{s.8A3-1}
The only maximal geometric set, denoted $\Lsub1$, does have $132$ lines. It
is characterized by the constant pattern
(see \autoref{rem.pattern.unique})
\[
\pat=\pattern{12,12}.
\sublabel1
\]


\subsection{Configuration~\hlink{8A3}{4}}\label{s.8A3-4}
There are $162$ sets $\fL\in\Bnd_{10}(\oorb4)$. Most
are maximal; one of them,
denoted $\Lsub2$, has $132$ lines and is determined by the pattern
\[
\pat=\pattern{0,3,*,4,*,0,*}.
\sublabel2
\]
Eleven sets are of rank $19$;
extending
these sets by a maximal orbit
(see \autoref{ss.ext.max}), we arrive at the bound $\ls|\fL|\le112$.


On the other hand, there is a unique set
$\fL\in\Bnd_6(\oorb1\cup\oorb2\cup\oorbs2\cup\oorb6\cup\oorbs6)$.
One has $\ls|\fL|=92$, and this set is maximal, see \autoref{ss.maximal}.
\endproof


\section{The lattice $\N(12\bA\sb2)$}\label{S.12A2}
The ultimate result of this section is the following theorem.

\newcs{check-12A2-1}{}
\newcs{check-12A2-2}{}
\newcs{check-12A2-3}{}
\newcs{check-12A2-4}{}
\newcs{check-12A2-5}{\noexpand\checkdone}
\newcs{check-12A2-6}{}
\newcs{check-12A2-7}{}
\newcs{check-12A2-8}{}

\theorem\label{th.12A2}
Let $\fL\subset\N(12\bA_2)$ be a geometric set. Then, unless
\roster*
\item
$\ls|\fL|=144$ and $\fL$ is conjugate to $\Lmax2$,
see~\eqref{eq.Lmax.2}, or
\item
$\ls|\fL|=132$ and $\fL$ is conjugate to $\Lsub3$,
see~\eqref{eq.Lsub.3},
\endroster
one has $\ls|\fL|\le130$.
\endtheorem

In the course of the proof of this theorem we also discover
and describe (by means of their patterns,
see \autoref{rem.pattern.unique})
several geometric sets~$\fL$ of size $\ls|\fL|\ge124$.

\proof[Proof of \autoref{th.12A2}]
We proceed as in \autoref{S.few}, analyzing pairs
$(\hbar,\rbar)$
one by one. (The notation in the table is explained in \autoref{S.few}.)
{\def\vectorbox#1{#1\,\,}\def\countskip{\ }
\newcs{max-12A2}{190}%
\newcs{count-12A2}{9}%
\newcs{counts-12A2}{\DO{1}{190}\DO{2}{170}\DO{3}{160}\DO{4}{162}\DO{5}{144}\DO{6}{180}\DO{7}{164}\DO{8}{150}\DO{9}{66}}%
\orbits
\orbit{12}{1}{190}
\h{1}&\root\rat(8/3){2}{\spec{2}}&\rat(2/3){2}{}&\rat(2/3){2}{}&\rat(2/3){2}{}&\int(0){\0}{}&\rat(2/3){2}{}&\int(0){\0}{}&\int(0){\0}{}&\rat(2/3){2}{}&\int(0){\0}{}&\int(0){\0}{}&\int(0){\0}{}&&\cnt{572}&\bounds{190}{196}{fail}{false}&&\cr
\l{1}&\rat(4/3){1}{}&\rat(2/3){2}{}&\rat(2/3){2}{}&\rat(1/3){1}{}&\int(0){\0}{}&\int(0){\0}{}&\int(0){2}{}&\int(0){1}{}&\int(0){\0}{}&\int(0){\0}{}&\int(0){\0}{}&\int(0){\0}{}&\ln{30\sd}&\cnt{18}&\bounds{6}{6}{1}{true}&&\cr
\l{2}&\rat(4/3){1}{}&\rat(1/3){1}{}&\rat(1/3){1}{}&\rat(1/3){1}{}&\int(0){\0}{}&\rat(1/3){1}{}&\int(0){\0}{}&\int(0){\0}{}&\rat(1/3){1}{}&\int(0){\0}{}&\int(0){\0}{}&\int(0){\0}{}&\ln{1\sdd}&\cnt{32}&\bounds{10}{16}{2}{true}&&\cr
\endorbit
\orbit{12}{2}{170}
\h{2}&\rat(8/3){2}{\spec{2}}&\rat(2/3){2}{}&\rat(2/3){2}{}&\rat(2/3){2}{}&\root\int(0){\0}{}&\rat(2/3){2}{}&\int(0){\0}{}&\int(0){\0}{}&\rat(2/3){2}{}&\int(0){\0}{}&\int(0){\0}{}&\int(0){\0}{}&&\cnt{492}&\bounds{170}{176}{fail}{false}&&\cr
\l{1}&\rat(4/3){1}{}&\rat(2/3){2}{}&\rat(2/3){2}{}&\rat(1/3){1}{}&\int(0){\0}{}&\int(0){\0}{}&\int(0){2}{}&\int(0){1}{}&\int(0){\0}{}&\int(0){\0}{}&\int(0){\0}{}&\int(0){\0}{}&\ln{20\sd}&\cnt{18}&\bounds{6}{6}{1}{true}&&\cr
\l{2}&\rat(4/3){1}{}&\rat(2/3){2}{}&\rat(2/3){2}{}&\int(0){\0}{}&\int(0){1}{}&\rat(1/3){1}{}&\int(0){\0}{}&\int(0){\0}{}&\int(0){\0}{}&\int(0){2}{}&\int(0){\0}{}&\int(0){\0}{}&\ln{10\sd}&\cnt{6}&\bounds{2}{2}{1}{true}&&\cr
\l{3}&\rat(4/3){1}{}&\rat(1/3){1}{}&\rat(1/3){1}{}&\rat(1/3){1}{}&\int(0){\0}{}&\rat(1/3){1}{}&\int(0){\0}{}&\int(0){\0}{}&\rat(1/3){1}{}&\int(0){\0}{}&\int(0){\0}{}&\int(0){\0}{}&\ln{1\sdd}&\cnt{32}&\bounds{10}{16}{2}{true}&&\cr
\l{4}&\int(2){\2}{}&\int(1){\2}{}&\int(0){\0}{}&\int(0){\0}{}&\int(0){\0}{}&\int(0){\0}{}&\int(0){\0}{}&\int(0){\0}{}&\int(0){\0}{}&\int(0){\0}{}&\int(0){\0}{}&\int(0){\0}{}&\ln{5}&\cnt{4}&\bounds{2}{2}{1}{true}&&\cr
\endorbit
\orbit{12}{3}{160}
\h{3}&\rat(8/3){2}{\spec{2}}&\root\rat(2/3){2}{}&\rat(2/3){2}{}&\rat(2/3){2}{}&\int(0){\0}{}&\rat(2/3){2}{}&\int(0){\0}{}&\int(0){\0}{}&\rat(2/3){2}{}&\int(0){\0}{}&\int(0){\0}{}&\int(0){\0}{}&&\cnt{464}&\bounds{160}{160}{1}{false}&&\cr
\l{1}&\int(2){\2}{}&\int(0){\0}{}&\int(1){\2}{}&\int(0){\0}{}&\int(0){\0}{}&\int(0){\0}{}&\int(0){\0}{}&\int(0){\0}{}&\int(0){\0}{}&\int(0){\0}{}&\int(0){\0}{}&\int(0){\0}{}&\ln{4}&\cnt{4}&\bounds{2}{2}{1}{true}&&\cr
\l{3}&\rat(4/3){1}{}&\rat(2/3){2}{}&\rat(2/3){2}{}&\rat(1/3){1}{}&\int(0){\0}{}&\int(0){\0}{}&\int(0){2}{}&\int(0){1}{}&\int(0){\0}{}&\int(0){\0}{}&\int(0){\0}{}&\int(0){\0}{}&\ln{12}&\cnt{18}&\bounds{6}{6}{1}{true}&&\cr
\endorbit
\orbit{12}{4}{162}
\h{4}&\rat(2/3){2}{}&\rat(2/3){2}{}&\rat(2/3){2}{}&\rat(2/3){2}{}&\rat(2/3){2}{}&\rat(2/3){1}{}&\root\int(0){\0}{}&\int(0){\0}{}&\int(0){\0}{}&\rat(2/3){1}{}&\rat(2/3){2}{}&\rat(2/3){1}{}&&\cnt{444}&\bounds{162}{186}{fail}{false}&&\cr
\l{1}&\rat(2/3){2}{}&\rat(2/3){2}{}&\rat(2/3){2}{}&\rat(2/3){2}{}&\int(0){\0}{}&\rat(1/3){2}{}&\int(0){\0}{}&\int(0){\0}{}&\int(0){2}{}&\int(0){\0}{}&\int(0){\0}{}&\int(0){\0}{}&\ln{36\sd}&\cnt{6}&\bounds{2}{2}{1}{true}&&\cr
\l{2}&\rat(2/3){2}{}&\rat(2/3){2}{}&\rat(2/3){2}{}&\int(0){\0}{}&\int(0){\0}{}&\int(0){\0}{}&\int(0){1}{}&\int(0){\0}{}&\int(0){\0}{}&\rat(2/3){1}{}&\rat(1/3){1}{}&\int(0){\0}{}&\ln{18\sd}&\cnt{2}&\bounds{1}{1}{1}{true}&&\cr
\l{3}&\rat(2/3){2}{}&\rat(2/3){2}{}&\rat(1/3){1}{}&\rat(1/3){1}{}&\int(0){\0}{}&\int(0){\0}{}&\int(0){\0}{}&\int(0){\0}{}&\int(0){\0}{}&\rat(1/3){2}{}&\int(0){\0}{}&\rat(2/3){1}{}&\ln{24\sd}&\cnt{8}&\bounds{3}{4}{5}{true}&&\cr
\endorbit
\orbit{12}{5}{144}
\h{5}&\root\rat(2/3){2}{}&\rat(2/3){2}{}&\rat(2/3){2}{}&\rat(2/3){2}{}&\rat(2/3){2}{}&\rat(2/3){1}{}&\int(0){\0}{}&\int(0){\0}{}&\int(0){\0}{}&\rat(2/3){1}{}&\rat(2/3){2}{}&\rat(2/3){1}{}&&\cnt{416}&\bounds{144}{160}{fail}{false}&&\cr
\l{1}&\rat(2/3){2}{}&\rat(2/3){2}{}&\rat(2/3){2}{}&\rat(2/3){2}{}&\int(0){\0}{}&\rat(1/3){2}{}&\int(0){\0}{}&\int(0){\0}{}&\int(0){2}{}&\int(0){\0}{}&\int(0){\0}{}&\int(0){\0}{}&\ln{24}&\cnt{6}&\bounds{2}{2}{1}{true}&&\cr
\l{3}&\rat(2/3){2}{}&\rat(2/3){2}{}&\rat(1/3){1}{}&\rat(1/3){1}{}&\int(0){\0}{}&\int(0){\0}{}&\int(0){\0}{}&\int(0){\0}{}&\int(0){\0}{}&\rat(1/3){2}{}&\int(0){\0}{}&\rat(2/3){1}{}&\ln{8}&\cnt{8}&\bounds{3}{4}{5}{true}&&\cr
\endorbit
\orbit{12}{6}{180}
\h{6}&\int(2){\2}{}&\int(2){\2}{}&\int(2){\2}{}&\root\int(0){\0}{}&\int(0){\0}{}&\int(0){\0}{}&\int(0){\0}{}&\int(0){\0}{}&\int(0){\0}{}&\int(0){\0}{}&\int(0){\0}{}&\int(0){\0}{}&&\cnt{516}&\bounds{180}{180}{1}{false}&&\cr
\l{1}&\int(2){\2}{}&\int(1){\2}{\spec{1}}&\int(0){\0}{}&\int(0){\0}{}&\int(0){\0}{}&\int(0){\0}{}&\int(0){\0}{}&\int(0){\0}{}&\int(0){\0}{}&\int(0){\0}{}&\int(0){\0}{}&\int(0){\0}{}&\ln{12\sd}&\cnt{1}&\bounds{1}{1}{1}{true}&&\cr
\l{2}&\int(1){2}{}&\int(1){2}{}&\int(1){2}{}&\int(0){2}{}&\int(0){\0}{}&\int(0){2}{}&\int(0){\0}{}&\int(0){\0}{}&\int(0){2}{}&\int(0){\0}{}&\int(0){\0}{}&\int(0){\0}{}&\ln{8\sd}&\cnt{9}&\bounds{3}{3}{1}{true}&&\cr
\l{3}&\int(1){2}{}&\int(1){2}{}&\int(1){2}{}&\int(0){\0}{}&\int(0){2}{}&\int(0){\0}{}&\int(0){\0}{}&\int(0){2}{}&\int(0){\0}{}&\int(0){\0}{}&\int(0){\0}{}&\int(0){2}{}&\ln{16\sd}&\cnt{27}&\bounds{9}{9}{1}{true}&&\cr
\endorbit
\orbit{12}{7}{164}
\h{7}&\rat(2/3){2}{}&\rat(2/3){2}{}&\rat(2/3){2}{}&\rat(2/3){2}{}&\int(2){\2}{}&\rat(2/3){2}{}&\root\int(0){\0}{}&\int(0){\0}{}&\rat(2/3){2}{}&\int(0){\0}{}&\int(0){\0}{}&\int(0){\0}{}&&\cnt{468}&\bounds{164}{164}{1}{false}&&\cr
\l{1}&\rat(2/3){2}{}&\rat(2/3){2}{}&\rat(2/3){2}{}&\int(0){\0}{}&\int(1){2}{}&\int(0){\0}{}&\int(0){\0}{}&\int(0){2}{}&\int(0){\0}{}&\int(0){\0}{}&\int(0){\0}{}&\int(0){2}{}&\ln{12\sd}&\cnt{9}&\bounds{3}{3}{1}{true}&&\cr
\l{2}&\rat(2/3){2}{}&\rat(2/3){2}{}&\rat(1/3){1}{}&\int(0){\0}{}&\int(1){1}{}&\int(0){\0}{}&\int(0){2}{}&\int(0){\0}{}&\rat(1/3){1}{}&\int(0){\0}{}&\int(0){\0}{}&\int(0){\0}{}&\ln{6\sd}&\cnt{4}&\bounds{2}{2}{1}{true}&&\cr
\l{3}&\rat(2/3){2}{}&\rat(1/3){1}{}&\rat(2/3){2}{}&\rat(1/3){1}{}&\int(1){1}{}&\int(0){\0}{}&\int(0){\0}{}&\int(0){\0}{}&\int(0){\0}{}&\int(0){\0}{}&\int(0){2}{}&\int(0){\0}{}&\ln{24\sd}&\cnt{12}&\bounds{4}{4}{1}{true}&&\cr
\l{4}&\rat(2/3){2}{}&\int(0){\0}{}&\rat(2/3){2}{}&\rat(2/3){2}{}&\int(1){2}{}&\int(0){\0}{}&\int(0){2}{}&\int(0){\0}{}&\int(0){\0}{}&\int(0){2}{}&\int(0){\0}{}&\int(0){\0}{}&\ln{8\sd}&\cnt{3}&\bounds{1}{1}{1}{true}&&\cr
\l{5}&\int(1){\2}{}&\int(0){\0}{}&\int(0){\0}{}&\int(0){\0}{}&\int(2){\2}{}&\int(0){\0}{}&\int(0){\0}{}&\int(0){\0}{}&\int(0){\0}{}&\int(0){\0}{}&\int(0){\0}{}&\int(0){\0}{}&\ln{6}&\cnt{2}&\bounds{1}{1}{1}{true}&&\cr
\endorbit
\orbit{12}{8}{150}
\h{8}&\root\rat(2/3){2}{}&\rat(2/3){2}{}&\rat(2/3){2}{}&\rat(2/3){2}{}&\int(2){\2}{}&\rat(2/3){2}{}&\int(0){\0}{}&\int(0){\0}{}&\rat(2/3){2}{}&\int(0){\0}{}&\int(0){\0}{}&\int(0){\0}{}&&\cnt{440}&\bounds{150}{150}{1}{false}&&\cr
\l{1}&\rat(2/3){2}{}&\rat(2/3){2}{}&\rat(2/3){2}{}&\rat(2/3){2}{}&\int(0){\0}{}&\rat(2/3){2}{}&\int(0){\0}{}&\int(0){\0}{}&\rat(-1/3){2}{}&\int(0){\0}{}&\int(0){\0}{}&\int(0){\0}{}&\ln{5}&\cnt{2}&\bounds{1}{1}{1}{true}&&\cr
\l{3}&\rat(2/3){2}{}&\rat(2/3){2}{}&\rat(2/3){2}{}&\int(0){\0}{}&\int(1){2}{}&\int(0){\0}{}&\int(0){\0}{}&\int(0){2}{}&\int(0){\0}{}&\int(0){\0}{}&\int(0){\0}{}&\int(0){2}{}&\ln{10}&\cnt{9}&\bounds{3}{3}{1}{true}&&\cr
\l{5}&\rat(2/3){2}{}&\rat(2/3){2}{}&\rat(1/3){1}{}&\int(0){\0}{}&\int(1){1}{}&\int(0){\0}{}&\int(0){2}{}&\int(0){\0}{}&\rat(1/3){1}{}&\int(0){\0}{}&\int(0){\0}{}&\int(0){\0}{}&\ln{10}&\cnt{12}&\bounds{4}{4}{1}{true}&&\cr
\endorbit
\orbit{12}{9}{66}
\h{9}&\int(6){\6}{\spec{1}}&\int(0){\0}{}&\int(0){\0}{}&\int(0){\0}{}&\int(0){\0}{}&\int(0){\0}{}&\int(0){\0}{}&\int(0){\0}{}&\int(0){\0}{}&\int(0){\0}{}&\int(0){\0}{}&\int(0){\0}{}&&\cnt{132}&\bounds{66}{66}{1}{false}&&\cr
\l{1}&\int(3){\2}{}&\int(0){\2}{}&\int(0){\0}{}&\int(0){\0}{}&\int(0){\0}{}&\int(0){\0}{}&\int(0){\0}{}&\int(0){\0}{}&\int(0){\0}{}&\int(0){\0}{}&\int(0){\0}{}&\int(0){\0}{}&\ln{11\sdd}&\cnt{12}&\bounds{6}{6}{1}{false}&&\cr
\endorbit
\endorbits
 
}

\def\porb#1{\obj1{\cluster_{#1}}}



\subsection{Configuration~\hlink{12A2}{1}}\label{s.12A2-1}
We subdivide the orbit~$\oorb1$ into five pairwise disjoint
clusters $\porb1,\ldots,\porb5$ constituting
$\Orbit_5(\oorb1,6)$. Explicitly,
\[*
\porb{k}:=\bigl\{l\in\oorb1\bigm|l_k\cdot\hh_k=\tfrac13\bigr\},
\quad k\in\CK:=\bigl\{k\in\Omega\bigm|\mbox{$\hh_k=\frac23$}\bigr\}.
\]
Then, arguing as in~\autoref{ss.clusters},
we compute $\Bnd_{58}(\oorb1)=\varnothing$.

\subsection{Configuration~\hlink{12A2}{2}}\label{s.12A2-2}
There are four sets $\fL\in\Bnd_{30}(\oorb1)$. One,
denoted $\Lmisc2$, is maximal and has $130$ lines;
it is characterized by the pattern
\[
\pat:=\pattern{6,0,10,0,*}\quad\mbox{(not $\Q$-complete)}.
\misclabel2{130}
\]
The
three other sets are of rank $19$; extending them by an extra orbit
(see \autoref{ss.ext.other}), we
arrive at a number of sets of size $\ls|\fL|\le118$ and one,
up to $\OG_{\hh}(\bN)$, maximal set $\Lmisc3$ of
size~$124$.
The latter is characterized by any of the five patterns
\[
\pat_\cluster=\pattern{5|4,2,8,0,*},
\qquad \cluster:=\cluster_1\in\Orbit_{5}(\oorb1,16);
\misclabel3{124}
\]
explicitly, $\cluster=\bigl\{l\in\oorb1\bigm|l_k\cdot\hh_k\ne\frac13\bigr\}$
for some $k\in\CK$ (see \autoref{s.12A2-1}).

On the other hand,
there are $13$ sets $\fL\in\Bnd_6(\oorb2\cup\oorb3\cup\oorb4\cup\oorbs4)$, all
saturated and with $\ls|\fL|\le94$. One set is of rank~$19$; extending it by
an extra orbit (see \autoref{ss.ext.other}), we obtain a number of sets
with at most $92$ lines.




\subsection{Configuration~\hlink{12A2}{3}}\label{s.12A2-3}
There is a single set $\fL\in\Bnd_{14}(\oorb3)$;
this set is maximal, and one has $\ls|\fL|=120$.


\subsection{Configuration~\hlink{12A2}{4}}\label{s.12A2-4}
There are $733$ sets $\fL\in\Bnd_{14}(\oorb1)$, which are all saturated;
we have $\ls|\fL|\le126$. Extending the $32$ sets of rank $19$ by a maximal
orbit (see \autoref{ss.ext.max}), we arrive at $\ls|\fL|\le112$.
The only, up to $\OG_{\hh}(\bN)$,
set $\Lmisc4$ with $126$ lines constitutes $\Bnd_{0}(\oorb1)$; it
is characterized by any of the four
patterns
\[
\pat_\cluster=\pattern{2,0,3|2},
\qquad \cluster:=\cluster_3\in\Orbit_4(\oorb3,6).
\misclabel4{126}
\]

On the other hand, there are $105$ sets $\fL\in\Bnd_{14}(\oorb2\cup\oorb3)$,
which are all of rank $18$ or~$19$.
Extending them by a maximal orbit (see \autoref{ss.ext.max}), we arrive
at $\ls|\fL|\le118$.


\subsection{Configuration~\hlink{12A2}{6}}\label{s.12A2-6}
There are $16$ sets $\fL\in\Bnd_{48}(\oorb3)$. One of them, denoted $\Lmax2$,
is maximal and contains $144$ lines. It is determined by
the pattern
\[
\pat=\pattern{0,0,9}\quad\mbox{(not $\Q$-complete)}.
\maxlabel2
\]
Extending the remaining $15$ sets by one or two extra orbits
(see \autoref{ss.ext.other}),
we obtain, among others,
a set with $124$ lines
and one with $128$ lines. The latter,
denoted by $\Lmisc5$, is characterized by the pattern
\[
\pat=\pattern{0,2,7}.
\misclabel5{128}
\]
The $124$-element set $\Lmisc6$ is characterized by any of the six patterns
\[
\pat_\cluster=\pattern{0,3|2,7|6},
\misclabel6{124}
\]
where $\cluster:=\cluster_3\in\Orbit_6(\oorb3,8)$ and
$\cluster_2\in\Orbit_1(\oorb2,4,\stab\cluster)$
is determined by~$\cluster$.


\subsection{Configuration~\hlink{12A2}{7}}\label{s.12A2-7}
There are $244$ sets $\fL\in\Bnd_{22}(\oorb3)$, all saturated. One of these
sets, denoted $\Lsub3$, has $132$ lines; it is determined by the pattern
\[
\pat=\pattern{3,0,4,0,0,*}.
\sublabel3
\]
For the other sets, one has $\ls|\fL|\le116$.
Twenty sets are of rank $19$;
extending them by a maximal orbit (see \autoref{ss.ext.max}), we obtain
at most $120$ lines.

On the other hand, there are nine sets
$\fL\in\Bnd_8(\oorb1\cup\oorb2\cup\oorb4\cup\oorb5\cup\oorbs5)$, which are
all saturated and have $\ls|\fL|\le104$.
Extensions of the two sets of rank $19$ by an extra orbit
(see \autoref{ss.ext.other}) have at most $108$ lines.


\subsection{Configuration~\hlink{12A2}{8}}\label{s.12A2-8}
There are $94$ sets $\fL\in\Bnd_9(\oorb5)$.
Most are maximal, and one has $\ls|\fL|\le116$.
The extensions of the two sets of rank $19$ by a maximal orbit
(see \autoref{ss.ext.max}) are maximal sets with at most $110$ lines.
\endproof


\makeatletter
\let\@int\@@int
\let\@rat\@@rat
\makeatother

\section{The lattice $\N(24\bA\sb1)$}\label{S.24A1}
The results of this section are summarized by the following theorem.

\newcs{check-24A1-1}{\noexpand\checkshow}
\newcs{check-24A1-2}{}
\newcs{check-24A1-3}{}

\theorem\label{th.24A1}
Let $\fL\subset\N(24\bA_1)$ be a geometric set. Then, unless
\roster*
\item
$\ls|\fL|=144$
and $\fL$ is conjugate to $\Lmax3$, see~\eqref{eq.Lmax.3}, or
\item
$\ls|\fL|=132$
and $\fL$ is conjugate to one of the sets
$\Lsub4$,
$\Lsub5$,
$\Lsub6$, or
$\Lsub7$,
see~\eqref{eq.Lsub.4},
\eqref{eq.Lsub.5},
\eqref{eq.Lsub.6}, or
\eqref{eq.Lsub.7}, respectively,
\endroster
one has $\ls|\fL|\le130$.
\endtheorem

\proof
We proceed as in the previous sections.
Each component $v_k\in\bR_k\dual$,
$k\in\Omega=[1,\ldots,24]$, of a vector $v\in\bN$
is a multiple of the generator $r_k\in\bR_k$.
To save space, we use the notation
\[*
\mbox{$\Azero$ (if $v_k=0$),\quad
$\Ahalf$ (if $v_k=\pm\frac12r_k$),\quad
$\Aroot$ (if $v_k=\pm r_k$),\quad
$\bullet$ (the position of $\rbar$)}.
\]
The signs always agree, so that
$\hbar_k\cdot l_k\ge0$ for any line $l\in\fF(\hbar)$ and $k\in\Omega$.

There are three configurations%
\iffullversion.\else\ (see \autoref{tab.24A1}).\fi

{\def\latticetext{}
\newcs{max-24A1}{220}%
\newcs{count-24A1}{3}%
\newcs{counts-24A1}{\DO{1}{160}\DO{2}{184}\DO{3}{220}}%
\orbits
\orbit{24}{1}{160}
\h{1}&\int(2){\2}{}&\int(2){\2}{}&\int(2){\2}{}&\root\int(0){\0}{}&\int(0){\0}{}&\int(0){\0}{}&\int(0){\0}{}&\int(0){\0}{}&\int(0){\0}{}&\int(0){\0}{}&\int(0){\0}{}&\int(0){\0}{}&\int(0){\0}{}&\int(0){\0}{}&\int(0){\0}{}&\int(0){\0}{}&\int(0){\0}{}&\int(0){\0}{}&\int(0){\0}{}&\int(0){\0}{}&\int(0){\0}{}&\int(0){\0}{}&\int(0){\0}{}&\int(0){\0}{}&&\cnt{512}&\bounds{160}{256}{fail}{false}&&\cr
\l{1}&\int(1){1}{}&\int(1){1}{}&\int(1){1}{}&\int(0){\0}{}&\int(0){1}{}&\int(0){1}{}&\int(0){\0}{}&\int(0){\0}{}&\int(0){\0}{}&\int(0){1}{}&\int(0){\0}{}&\int(0){\0}{}&\int(0){\0}{}&\int(0){\0}{}&\int(0){\0}{}&\int(0){\0}{}&\int(0){\0}{}&\int(0){\0}{}&\int(0){1}{}&\int(0){1}{}&\int(0){\0}{}&\int(0){\0}{}&\int(0){\0}{}&\int(0){\0}{}&\ln{16\sdd}&\cnt{32}&\bounds{10}{16}{5}{true}&&\cr
\endorbit
\orbit{24}{2}{184}
\h{2}&\rat(1/2){1}{}&\rat(1/2){1}{}&\rat(1/2){1}{}&\rat(1/2){1}{}&\rat(1/2){1}{}&\int(2){\2}{}&\root\int(0){\0}{}&\rat(1/2){1}{}&\int(0){\0}{}&\int(0){\0}{}&\rat(1/2){1}{}&\int(0){\0}{}&\rat(1/2){1}{}&\int(0){\0}{}&\int(0){\0}{}&\int(0){\0}{}&\int(0){\0}{}&\int(0){\0}{}&\int(0){\0}{}&\int(0){\0}{}&\int(0){\0}{}&\int(0){\0}{}&\int(0){\0}{}&\int(0){\0}{}&&\cnt{464}&\bounds{184}{240}{fail}{false}&&\cr
\l{1}&\rat(1/2){1}{}&\rat(1/2){1}{}&\rat(1/2){1}{}&\rat(1/2){1}{}&\int(0){\0}{}&\int(1){1}{}&\int(0){\0}{}&\int(0){\0}{}&\int(0){\0}{}&\int(0){\0}{}&\int(0){\0}{}&\int(0){\0}{}&\int(0){\0}{}&\int(0){\0}{}&\int(0){1}{}&\int(0){\0}{}&\int(0){\0}{}&\int(0){1}{}&\int(0){\0}{}&\int(0){\0}{}&\int(0){\0}{}&\int(0){\0}{}&\int(0){\0}{}&\int(0){1}{}&\ln{56\sd}&\cnt{8}&\bounds{3}{4}{5}{true}&&\cr
\l{2}&\int(1){\2}{}&\int(0){\0}{}&\int(0){\0}{}&\int(0){\0}{}&\int(0){\0}{}&\int(2){\2}{}&\int(0){\0}{}&\int(0){\0}{}&\int(0){\0}{}&\int(0){\0}{}&\int(0){\0}{}&\int(0){\0}{}&\int(0){\0}{}&\int(0){\0}{}&\int(0){\0}{}&\int(0){\0}{}&\int(0){\0}{}&\int(0){\0}{}&\int(0){\0}{}&\int(0){\0}{}&\int(0){\0}{}&\int(0){\0}{}&\int(0){\0}{}&\int(0){\0}{}&\ln{8}&\cnt{1}&\bounds{1}{1}{1}{true}&&\cr
\endorbit
\orbit{24}{3}{220}
\h{3}&\rat(1/2){1}{}&\rat(1/2){1}{}&\rat(1/2){1}{}&\rat(1/2){1}{}&\rat(1/2){1}{}&\rat(1/2){1}{}&\rat(1/2){1}{}&\rat(1/2){1}{}&\root\int(0){\0}{}&\rat(1/2){1}{}&\int(0){\0}{}&\int(0){\0}{}&\int(0){\0}{}&\rat(1/2){1}{}&\int(0){\0}{}&\int(0){\0}{}&\int(0){\0}{}&\int(0){\0}{}&\int(0){\0}{}&\int(0){\0}{}&\rat(1/2){1}{}&\int(0){\0}{}&\int(0){\0}{}&\rat(1/2){1}{}&&\cnt{440}&\bounds{220}{220}{1}{false}&&\cr
\l{1}&\rat(1/2){1}{}&\rat(1/2){1}{}&\rat(1/2){1}{}&\rat(1/2){1}{}&\rat(1/2){1}{}&\int(0){\0}{}&\int(0){\0}{}&\rat(1/2){1}{}&\int(0){\0}{}&\int(0){\0}{}&\int(0){1}{}&\int(0){\0}{}&\int(0){1}{}&\int(0){\0}{}&\int(0){\0}{}&\int(0){\0}{}&\int(0){\0}{}&\int(0){\0}{}&\int(0){\0}{}&\int(0){\0}{}&\int(0){\0}{}&\int(0){\0}{}&\int(0){\0}{}&\int(0){\0}{}&\ln{110\sd}&\cnt{4}&\bounds{2}{2}{1}{true}&&\cr
\endorbit
\endorbits
 
}

Fix a basis $\{r_k\}$, $k\in\Omega$, for $24\bA_1$ consisting of roots.
The kernel
\[*
\bN\bmod24\bA_2\subset\discr24\bA_1\cong(\Z/2)^{24}
\]
of the extension is the Golay code $\CC_{24}$
(see~\cite{Conway.Sloane}).
The map $\supp$ identifies codewords with
subsets of $\Omega$; then, $\CC_{24}$ is invariant
under complement and, in addition to~$\varnothing$ and~$\Omega$,
it consists of $759$ octads, $759$ complements thereof,
and $2576$ dodecads.

To simplify the notation, we identify the basis vectors $r_k$
(assumed fixed) with their indices $k\in\Omega$. For a subset
$\CS\subset\Omega$, we let $\vv\CS:=\sum r$, $r\in\CS$,
and abbreviate $\cw\CS:=\frac12\vv\CS\in\bN$ if $\CS\in\CC_{24}$ is a
codeword.

\subsection{Configuration~\hlink{24A1}{1}}\label{s.24A1-1}
We have $\ls|\stab\hh|=5760$ and $\hh$ is the sum of three roots. Using
patterns (see \autoref{s.patterns}), we compute
$\Bnd_{38}(\oorb1)=\varnothing$,
arriving at
$\ls|\fL|\le120$.

\def\rh{r_{\hh}}
\def\porb#1{\obj1{\cluster_{#1}}}

\subsection{Configuration~\hlink{24A1}{2}}\label{s.24A1-2}
We have $\ls|\stab\hh|=1344$ and $\hh=\cw\CO+\rh$,
where $\CO\in\CC_{24}$
is an octad and $\rh\notin\CO$.
Let
$\CK:=\Omega\sminus(\CO\cup\{\rh,\rbar\})$
and break $\oorb1$ into eight clusters
\[*
\porb{o}:=\bigl\{\orb\subset\oorb1\bigm|(\CK\cap\supp\orb)\subset o\bigr\},
\quad o\in\CC_{24},\ \ls|o|=8,\ o\cap\CO=\varnothing,\
 \rbar\notin o,\ \rh\in o.
\]
They constitute the orbit $\Orbit_8(\oorb1,14)$.
Each combinatorial orbit
$\orb\subset\oorb1$ belongs to two clusters, and each pair of clusters
intersects in a single pair of dual orbits.

The set $\Bnd_{52}(\oorb1)$ is computed
cluster by cluster, as explained in \autoref{ss.clusters}.
We arrive at a number of sets~$\fL$ of size $\ls|\fL|\le120$
and a few those with
$124\le\ls|\fL|\le132$.
All sets are maximal.
The large sets found can be described as
\[*
\fL=\oorb1\cap\spn(\rbar,\bh,\ortu_s,v)^\perp,
\]
where
\roster*
\item
$\rbar$ and $\bh:=\cw\CO-2\rh=\hh-3\rh$ generate the
subspace $\oorb1^\perp\subset N$,
\item
$\ortu_s:=\vv\CK-2s$ for a certain fixed point $s\in\CO$,
\endroster
and the extra vector $v$ is specified below, using
the \latin{ad hoc} notation
\roster*
\item
$\ortv_o:=\cw{o\sminus\CO}-\cw{o\cap\CO}$ for a codeword $o\in\CC_{24}$.
\endroster
Then, the large sets are as follows:
\bgroup\allowdisplaybreaks
\def\boxi#1#2#3{\mbox{$\rbar#1o$, $\rh#2o$, $s#3o$},}
\def\boxii#1#2{\mbox{$\ls|o\cap\CO|=#1$, $\ls|o|=#2$}}
\begin{alignat}{3}
\Lsub4:\quad
 &v=k,&\qquad&k\in\CK;&\quad&\sublabel4\\
\Lsub5:\quad
 &v=\ortv_o,&&\boxi\in\in\in&&\boxii28;\sublabel5\\
\Lmisc7:\quad
 &v=\ortv_o,&&\boxi\notin\in\in&&\boxii28;\misclabel7{126}\\
\Lmisc8:\quad
 &v=\ortv_o+\rh,&&\boxi\in\notin\in&&\boxii48;\misclabel8{126}\\
\Lmisc9:\quad
 &v=\cw{o},&&\boxi\in\notin\notin&&\boxii08;\misclabel9{126}\\
\Lmisc{10}:\quad
 &v=\ortv_o,&&\boxi\notin\in\in&&\boxii2{12};\misclabel{10}{126}\\
\Lmisc{11}:\quad
 &v=\ortv_o,&&\boxi\in\in\in&&\boxii2{12};\misclabel{11}{124}\\
\Lmisc{12}:\quad
 &v=\cw{o}-\rh,&&\boxi\in\in\notin&&\boxii28.\misclabel{12}{124}
\end{alignat}\egroup
In each case, it is straightforward that the set of data required for the
description is unique up to $\OG_{\hh}(\bN)$.

\def\pporb#1{\obj1{\cluster#1}}
\def\porb#1{\obj1{\cluster_{#1}}}

\subsection{Configuration~\hlink{24A1}{3}}\label{s.24A1-3}
We have $\ls|\stab\hh|=7920$ and $\hh=\cw\CO$, where $\CO\in\CC_{24}$
ia a dodecad. Let $\CK:=\Omega\sminus(\CO\cup\rbar)$. Each support
$o:=\supp\orb$, $\orb\in\oorb1$, is an octad,
so that $\ls|o\cap\CO|=6$ and
$\ls|o\cap\CK|=2$; conversely,
each $2$-element set $s\subset\CK$ extends to a unique
pair of such octads, representing a pair of dual orbits
$\orb,\orb^*\subset\oorb1$.
We break $\oorb1$ into eleven clusters
\[*
\porb{k}:=\{\orb\subset\oorb1\,|\,\supp\orb\ni k\},
\quad k\in\CK.
\]
Each orbit belongs to two
clusters, and each pair of clusters intersects in a single pair of dual
orbits.
We compute the set $\Bnd_{88}(\oorb1)$ cluster by cluster, as explained in
\autoref{ss.clusters}.
Note that the first
cluster~$\pporb{}$ has $\ls|\fL\cap\pporb{}|\ge24$ (and, hence,
$\dif_0(\pporb{})\ge4$) and, in case of equality,
also $\ls|\fL\cap\porb{k}|=24$ for each $k\in\CK$.
In this latter case, we reduce overcounting by using the following
observations:
\roster*
\item
if $\dif_0(\pporb{})=6$, there is exactly one other cluster $\pporb'$ with
$\dif_0(\pporb')=6$, so that $\fL\cap\orb=\fL\cap\orb^*=\varnothing$
for the two orbits $\orb,\orb^*\subset\pporb{}\cap\pporb'$;
\item
if $\dif_0(\pporb{})=4$, then $\dif_0(\porb{k})=4$ for each $k\in\CK$
(thus, no preferred order), and we can choose $\pporb'$ so that
$\ls|\fL\cap\orb|=\ls|\fL\cap\orb^*|=1$ for
$\orb,\orb^*\subset\pporb{}\cap\pporb'$.
\endroster
In these two cases, we start with the pair $\pporb{}$, $\pporb'$ and employ
the extra symmetry.

The result is one maximal set $\Lmax3$ and two submaximal
sets $\Lsub6$, $\Lsub7$. As a by-product, we
have found six sets $\fL$ with $124\le\ls|\fL|\le130$ and a
number of sets of size $\ls|\fL|\le120$.
Most large sets can be described as
\[*
\fL=\oorb1\cap\spn(\rbar,\vv\CK,r,v)^\perp,
\]
where $r\in\CK$ is a certain fixed point and the extra vector $v$ is
described below. This description depends on a codeword $o\in\CC_{24}$
(we use the shortcut $\ortw_o:=3\ortv_o+\cw\CO$, where $\ortv_o$ is as
in \autoref{s.24A1-2}) and,
occasionally, an extra point $s\in o\cap\CO$ or $t\in\CK\sminus r$.
Then, the large sets are as follows:
\bgroup\allowdisplaybreaks
\def\boxi#1#2#3{\mbox{$\rbar#1o$, $r#2o$,\ifx#3\relax\else\ $t#3o$,\fi}}
\def\boxii#1#2{\mbox{$\ls|o\cap\CO|=#1$, $\ls|o|=#2$}}
\begin{alignat}{3}
\Lmax3:\quad
 &v=t;&\qquad&&\quad&\maxlabel3\\
\Lsub6:\quad
 &v=\cw{o}-s,&&\boxi\notin\notin\relax&&\boxii28;\sublabel6\\
\Lsub7:\quad
 &v=\ortw_o,&&\boxi\in\in\relax&&\boxii48;\sublabel7\\
\Lmisc{13}:\quad
 &v=\cw{o}-s,&&\boxi\in\notin\relax&&\boxii28;\misclabel{13}{130}\\
\Lmisc{14}:\quad
 &v=\ortw_o-3t,&&\boxi\in\in\in&&\boxii48;\misclabel{14}{128}\\
\Lmisc{15}:\quad
 &v=\ortw_o,&&\boxi\notin\in\relax&&\boxii48;\misclabel{15}{126}\\
\Lmisc{16}:\quad
 &v=\ortw_o,&&\boxi\notin\in\relax&&\boxii4{12};\misclabel{16}{126}\\
\Lmisc{17}:\quad
 &v=\ortw_o,&&\boxi\notin\notin\relax&&\boxii48.\misclabel{17}{124}
\end{alignat}\egroup
In~\eqref{eq.Lmisc.17} we require, in addition, that the $6$-element set
$(o\cap\CO)\cup\{\rbar,r\}$ should not be contained in an octad. Under this
extra assumption, the set of data needed for the description is unique up to
$\OG_{\hh}(\bN)$.
\endproof

\remark
In \autoref{s.24A1-3},
there is one more $126$-element set $\fL$. However, since $\fL$ is graph
isomorphic to $\Lmisc{15}\cong\Lmisc{16}$ and, on the other hand, we do not
assert the completeness in this range, we omit its description, which is more
complicated.
\endremark

\remark\label{rem.counts}
In the course of this computation, we have also observed all even line counts
$38\le\ls|\fL|\le120$ realized by geometric sets of rank~$20$.
\endremark

\section{Proofs of the main results}\label{S.proofs}

In this concluding section,
we fill in a few missing links to complete the proofs of the
principal results of the paper stated in the introduction.

\subsection{Proof of \autoref{th.main} and \autoref{ad.main}}\label{proof.main}
As explained in \autoref{s.K3}, instead of counting tritangents to smooth
sextics one can study (doubling the numbers) the Fano graphs of smooth
$2$-polarized $K3$-surfaces.
By \autoref{th.Fano.graph}, the latter task is equivalent to the study of the
Fano graphs of certain $2$-polarized lattices $\NS\ni h$, and
\autoref{prop.Niemeier} and subsequent definitions reduce it further to the
study of geometric subsets $\fL\subset\fF(\hbar)$ in $6$-polarized Niemeier
lattices $\bN\ni\hbar$ \emph{other than the Leech lattice~$\Lambda$} (as we
can always assume that there is a root $\rbar\in\hbar^\perp$).
This is done in Theorems~\ref{th.<=6}, \ref{th.8A3}, \ref{th.12A2},
and~\ref{th.24A1}, and there remains to observe that all sets of size~$144$
are isomorphic as abstract graphs,
\roster
\item\label{144}
size $144$:
$\Lmax1\cong\Lmax2\cong\Lmax3$,
$T=[12,6,12]\realcurve$,
\endroster
and there are two isomorphism classes of sets of size~$132$:
\roster[\lastitem]
\item\label{132.1}
size $132$:
$\Lsub1\cong\Lsub2\cong\Lsub3\cong\Lsub4\cong\Lsub5\cong\Lsub6$,
$T=[2,0,66]\realcurve$,
\item\label{132.2}
size $132$:
$\Lsub7$,
$T=[4,0,32]\realcurve$.
\endroster
The graphs are compared by means of the \texttt{GRAPE} package
\cite{GRAPE:nauty,GRAPE:paper,GRAPE} in \GAP~\cite{GAP4}.
\latin{A posteriori}, the large graphs $\fL$
found in the paper
are distinguished by their size $\ls|\fL|$, discriminant
form $\discr(\spn_\Z\fL)$, and, in a
few cases below, the size $\ls|\Aut\fL|$ of the group of abstract graph
automorphisms (also computed by \texttt{GRAPE}).
Instead of $\discr(\spn_\Z\fL)$, we give a list of representatives of the genus of
the transcendental lattice $T:=\NS^\perp$ of the corresponding
$2$-polarized $K3$-surface, using the inline notation $[2a,b,2c]$ for the
even rank~$2$
form
$T=\Z u+\Z v$, $u^2=2a$, $u\cdot v=b$, $v^2=2c$.
The meaning of the superscript $\realcurve$ is explained in
\autoref{proof.main.2}\iref{item.real.curve} below.

This observation establishes the bounds stated in \autoref{th.main}; the
uniqueness is proved in \autoref{proof.main.2} below.
For the record, we give a similar classification for the other large
geometric sets found in the course of the computation and
described elsewhere in the paper:
\roster[\lastitem]
\item\label{130}
size $130$:
$\Lmisc{1}\cong\Lmisc{2}\cong\Lmisc{13}$,
$T=[12,3,12]\realcurve$;
\item\label{128}
size $128$:
$\Lmisc{5}\cong\Lmisc{14}$,
$T=[12,2,12]\realcurve$;
\item\label{126.1}
size $126$:
$\Lmisc{7}\cong\Lmisc{8}\cong\Lmisc{15}\cong\Lmisc{16}$,
$T= [ 2, 1, 72 ]\realcurve$, $[ 6, 1, 24 ]$, or $[ 8, 1, 18 ]$;
\item\label{126.2}
size $126$:
$\Lmisc{4}\cong\Lmisc{10}$,
$\ls|\Aut\fL|=144$,
$T=[ 14, 7, 14 ]\realcurve$;
\item\label{126.3}
size $126$:
$\Lmisc{9}$,
$\ls|\Aut\fL|=504$,
$T=[ 14, 7, 14 ]\realcurve$;
\item\label{124}
size $124$:
$\Lmisc{3}\cong\Lmisc{6}\cong\Lmisc{11}\cong\Lmisc{12}\cong\Lmisc{17}$,
$T=[ 4, 0, 38 ]\realcurve$ or $[ 6, 2, 26 ]$.
\endroster
Besides, we have found nine isomorphism classes of geometric sets of
size~$120$. Note that,
unlike~\iref{144}--\iref{132.2},
we do not assert the completeness of these lists.

The statement of \autoref{ad.main} is essentially given by
\autoref{rem.counts} and the above list, as geometric sets with fewer than
$38$ lines are easily constructed directly, mostly within an
appropriate single
combinatorial orbit~$\orb$.

\subsection{Proof of \autoref{th.main} (uniqueness)}\label{proof.main.2}
In full agreement with \autoref{cor.Elkies}, all large geometric sets
listed in \autoref{proof.main} are of
the maximal rank $\rank\fL=20$.
Therefore, the isomorphism classes of the smooth $2$-polarized $K3$-surfaces
(or, equivalently, the projective equivalence classes of sextics $C\subset\Cp2$) with
the given Fano graph $\Fn(X,h)\cong\fL$ are given by the global Torelli
theorem~\cite{Pjatecki-Shapiro.Shafarevich} (\cf. also
\cite[Theorem~3.11]{DIS}) as the classes of primitive embeddings
$\NS\into L:=-2\bE_8\oplus3\bU$ up to the group $\OG_h^+(L)$ of
auto-isometries of~$L$ preserving~$h$ and the \emph{positive sign structure}
(\ie, orientation of maximal positive definite subspaces of $L\otimes\R$;
here, $\NS\ni h$ is the $2$-polarized lattice obtained from a mild extension
$S\supset\spn\fL\ni\hbar$ by the
inverse construction of \autoref{s.S}).

The classification of embeddings is done using
Nikulin~\cite{Nikulin:forms}.
With an extension~$S$
(and, hence, lattice~$\NS$) fixed, the genus of the transcendental
lattice $T:=\NS^\perp$ is
determined by the discriminant $\discr\NS\cong-\discr T$.
Then, for each representative $T$ of this genus,
the isomorphism classes of the embeddings with $\NS^\perp\cong T$ are in a
one-to-one correspondence with the double cosets
\[*
\OG_h(\NS)\backslash\Aut(\discr\NS)/O^+(T).
\]
There are obvious identities and inclusions
\[*
\OG_h(\NS)=\OG_\hbar(S)\subset\OG_\hbar(\spn\fL)\subset\OG_\hbar(\spn_\Z\fL)=\Aut\fL.
\]
Besides, we have the following lemma.

\lemma\label{lem.detT}
Let $\fL\subset\bN$ be a geometric set, and assume that
\[*
\rank\fL=20,
\qquad
\det(\spn_\Z\fL)<1296.
\]
Then the only mild extension $S\supset\spn\fL$ is
$S=\spn\fL=\spn_\Z\fL$.
\endlemma

\proof
For the N\'{e}ron--Severi lattice~$\NS(X)$ of a $K3$-surface~$X$ corresponding
to $S$ we
have
\[*
\ls|\det\NS(X)|=\det(\spn_\Z\fL)/3i^2,
\quad i:=[S:\spn_\Z\fL].
\]
On the other hand, since $X$ is smooth and of the maximal Picard rank~$20$,
we have $\ls|\det\NS(X)|\ge108$ by
\cite[Theorem 1.5]{degt:singular.K3}. This implies $i<2$, hence $i=1$.
\endproof

\autoref{lem.detT} applies to all geometric sets listed
in \autoref{proof.main}\iref{144}--\iref{124}
(and to the nine sets of size~$120$ mentioned
thereafter).
Then, a direct computation shows that the natural map
$\Aut\fL\to\Aut(\discr\NS)$
is surjective. This fact renders the other group $\OG^+(T)$ redundant and
proves that
each pair $(\fL,T)$ listed is realized by either
\roster
\item\label{item.real.curve}
a single curve $C\cong\bar C$ (marked with a $^*$ in the list) or
\item
a pair $C$, $\bar C$ of complex conjugate curves.
\endroster
(The former is the case whenever $T$ admits an orientation reversing
auto-isometry; note that
we do not assert that the
curve admits a real structure, although most likely it does.)
In particular, each of the three configurations listed in items~\iref{144},
\iref{132.1}, \iref{132.2}
is realized by a single curve, as stated in \autoref{th.main}.
\qed

\subsection{Proof of \autoref{th.real}}\label{proof.real}
Let $C\subset\Cp2$ be a real sextic. By definition, this means that $C$ is
invariant under a certain fixed \emph{real structure}
(anti-holomorphic involution) $c\:\Cp2\to\Cp2$.
This involution lifts to two commuting anti-holomorphic automorphisms $c_\pm$
of the covering $K3$-surface $\Gf\:X\to\Cp2$, so that $c_+\circ c_-=\tau$ is
the deck translation. \latin{A priori}, $c_\pm$ are either involutions or of
order~$4$, with $c_\pm^2=\tau$.
However, if we assume that $C$ has a real tritangent~$L$, then at
least one of the three tangency points (possibly, infinitely near) must be
real; thus, the ramification locus of~$\Gf$ has a real point and both
lifts~$c_\pm$
are involutions, \ie, real structures on~$X$.
Furthermore, we can select~$c_+$ so that both pull-backs
$L_1,L_2\subset\Gf\1(L)$ are real (and then they are complex conjugate with
respect to~$c_-$); then, the pull-backs $L_1',L_2'\subset\Gf\1(L')$ of any
other real tritangent~$L'$ are also real, as each of $L_i'$ intersects
exactly one of~$L_1$, $L_2$ at a single point, which must be real.

Thus, we reduce the problem to counting real lines in a real $2$-polarized
$K3$-surface $(X,h)$. More precisely, this means that we fix a real structure
$c\:X\to X$, $c_*(h)=-h$, and count lines $L\subset X$ satisfying
$c_*[L]=-[L]$.
(Recall that each line is unique in its homology class and that
anti-holomorphic maps reverse the orientation of complex curves.)
Arguing as in \cite{DIS}, we can perturb the period of~$X$ to change
$\NS(X)$ to the sublattice rationally generated by the real lines; then,
\emph{all} lines in
the new surface~$X$ are real and $\Ker(1-c_*)\subset T(X)=\NS(X)^\perp$.
Using the classification of real structures found in~\cite{Nikulin:forms}, we
obtain the following statement.

\lemma[\cf. {\cite[Lemma 3.10]{DIS}}]\label{lem.real}
A smooth $2$-polarized $K3$-surface~$X$ is equilinear deformation equivalent
to a real surface~$Y$ in which all lines are real if and only if
the orthogonal complement $\Fn(X,h)^\perp$ contains $\bA_1$ or $\bU(2)$ as a
sublattice.
\done
\endlemma

If $\rank\Fn(X,h)=20$,
the Picard rank $\rank\NS(X)=20$ is also maximal, the moduli space is
finite, and the statement can be made more precise.

\lemma\label{lem.real.max}
Let $X$ be a $2$-polarized $K3$-surface, $\rank\Fn(X,h)=20$. Then, the real
structures $c\:X\to X$ with respect to which all lines are real are in a
one-to-one correspondence with pairs of roots $\pm r\in T(X)$. Under this
correspondence, $-c_*$ is the reflection against the hyperplane
$r^\perp\subset H_2(X;\Z)$.
\done
\endlemma

There remains to examine the list found in \autoref{proof.main} and observe
that the maximum, which is $132$ lines, is realized by a unique graph, \viz.
the one in item~\iref{132.1}, and the corresponding transcendental lattice~$T$
contains a single pair of roots. The next \emph{known} examples are
item~\iref{126.1} with
$126$ lines and two of the nine graphs with $120$ lines,
but we do not assert the completeness
of our lists in this range.
\qed


{
\let\.\DOTaccent
\def\cprime{$'$}
\bibliographystyle{amsplain}
\bibliography{degt}

\providecommand{\bysame}{\leavevmode\hbox to3em{\hrulefill}\thinspace}
\providecommand{\MR}{\relax\ifhmode\unskip\space\fi MR }
\providecommand{\MRhref}[2]{%
  \href{http://www.ams.org/mathscinet-getitem?mr=#1}{#2}
}
\providecommand{\href}[2]{#2}
\begin{thebibliography}{10}

\bibitem{Bourbaki:Lie}
Nicolas Bourbaki, \emph{Lie groups and {L}ie algebras. {C}hapters 4--6},
  Elements of Mathematics (Berlin), Springer-Verlag, Berlin, 2002, Translated
  from the 1968 French original by Andrew Pressley. \MR{1890629 (2003a:17001)}

\bibitem{Clebsch}
A.~Clebsch, \emph{Zur {T}heorie der algebraischen {F}l\"{a}chen}, J. Reine
  Angew. Math. \textbf{58} (1861), 93--108. \MR{1579142}

\bibitem{Conway.Sloane}
J.~H. Conway and N.~J.~A. Sloane, \emph{Sphere packings, lattices and groups},
  Grundlehren der Mathematischen Wissenschaften [Fundamental Principles of
  Mathematical Sciences], vol. 290, Springer-Verlag, New York, 1988, With
  contributions by E. Bannai, J. Leech, S. P. Norton, A. M. Odlyzko, R. A.
  Parker, L. Queen and B. B. Venkov. \MR{920369 (89a:11067)}

\bibitem{degt:supersingular}
Alex Degtyarev, \emph{Lines in supersingular quartics}, To appear,
  \arXiv{1604.05836}, 2016.

\bibitem{degt:lines}
\bysame, \emph{Lines on {S}mooth {P}olarized {K}3-{S}urfaces}, Discrete Comput.
  Geom. \textbf{62} (2019), no.~3, 601--648. \MR{3996938}

\bibitem{degt:singular.K3}
\bysame, \emph{Smooth models of singular {$K3$}-surfaces}, Rev. Mat. Iberoam.
  \textbf{35} (2019), no.~1, 125--172. \MR{3914542}

\bibitem{degt:sextics.full}
\bysame, \emph{Tritangents to smooth sextic curves}, Preprint MPIM 19-51, 2019.

\bibitem{DIS}
Alex Degtyarev, Ilia Itenberg, and Ali~Sinan Sert\"oz, \emph{Lines on quartic
  surfaces}, Math. Ann. \textbf{368} (2017), no.~1-2, 753--809. \MR{3651588}

\bibitem{Elkies}
Noam~D. Elkies, \emph{Upper bounds on the number of lines on a surface}, New
  Trends in Arithmetic and Geometry of Algebraic Surfaces, BIRS conference
  17w5146, 2017.

\bibitem{GAP4}
\emph{{GAP} {\textendash} {G}roups, {A}lgorithms, and {P}rogramming, {V}ersion
  4.10.1}, \href {https://www.gap-system.org} {https://www.gap-system.org}, Feb
  2019.

\bibitem{Rams.Gonzales}
V{\'\i}ctor Gonz{\'a}lez-Alonso and S{\l}awomir Rams, \emph{Counting lines on
  quartic surfaces}, Taiwanese J. Math. \textbf{20} (2016), no.~4, 769--785.
  \MR{3535673}

\bibitem{Kondo}
Shigeyuki Kond{\=o}, \emph{Niemeier lattices, {M}athieu groups, and finite
  groups of symplectic automorphisms of {$K3$} surfaces}, Duke Math. J.
  \textbf{92} (1998), no.~3, 593--603, With an appendix by Shigeru Mukai.
  \MR{1620514}

\bibitem{Kulikov:periods}
Vik.~S. Kulikov, \emph{Surjectivity of the period mapping for {$K3$} surfaces},
  Uspehi Mat. Nauk \textbf{32} (1977), no.~4(196), 257--258. \MR{0480528 (58
  \#688)}

\bibitem{GRAPE:nauty}
Brendan~D. McKay, \emph{Nauty user's guide (version 1.5)}, Tech. Report
  TR-CS-90-0, Australian National University, Computer Science Department,
  1990.

\bibitem{GRAPE:paper}
Brendan~D. McKay and Adolfo Piperno, \emph{Practical graph isomorphism, {II}},
  J. Symbolic Comput. \textbf{60} (2014), 94--112. \MR{3131381}

\bibitem{Mukai}
Shigeru Mukai, \emph{Finite groups of automorphisms of {$K3$} surfaces and the
  {M}athieu group}, Invent. Math. \textbf{94} (1988), no.~1, 183--221.
  \MR{958597}

\bibitem{Niemeier}
Hans-Volker Niemeier, \emph{Definite quadratische {F}ormen der {D}imension
  {$24$} und {D}iskriminante {$1$}}, J. Number Theory \textbf{5} (1973),
  142--178. \MR{0316384}

\bibitem{Nikulin:forms}
V.~V. Nikulin, \emph{Integer symmetric bilinear forms and some of their
  geometric applications}, Izv. Akad. Nauk SSSR Ser. Mat. \textbf{43} (1979),
  no.~1, 111--177, 238, English translation: Math USSR-Izv. 14 (1979), no. 1,
  103--167 (1980). \MR{525944 (80j:10031)}

\bibitem{Nikulin:degenerations}
\bysame, \emph{Degenerations of {K}\"{a}hlerian {K}3 surfaces with finite
  symplectic automorphism groups}, Izv. Ross. Akad. Nauk Ser. Mat. \textbf{79}
  (2015), no.~4, 103--158. \MR{3397421}

\bibitem{Nishiyama}
Ken-ichi Nishiyama, \emph{The {J}acobian fibrations on some {$K3$} surfaces and
  their {M}ordell-{W}eil groups}, Japan. J. Math. (N.S.) \textbf{22} (1996),
  no.~2, 293--347. \MR{1432379}

\bibitem{Persson:sextics}
Ulf Persson, \emph{Double sextics and singular {$K$}-{$3$} surfaces}, Algebraic
  geometry, {S}itges ({B}arcelona), 1983, Lecture Notes in Math., vol. 1124,
  Springer, Berlin, 1985, pp.~262--328. \MR{805337 (87i:14036)}

\bibitem{Pjatecki-Shapiro.Shafarevich}
I.~I. Pjatecki{\u\i}-{\v{S}}apiro and I.~R. {\v{S}}afarevi{\v{c}},
  \emph{Torelli's theorem for algebraic surfaces of type {${\rm K}3$}}, Izv.
  Akad. Nauk SSSR Ser. Mat. \textbf{35} (1971), 530--572, English translation:
  Math. USSR-Izv. 5, 547--588. \MR{0284440 (44 \#1666)}

\bibitem{rams.schuett:char3}
S{\l}awomir Rams and Matthias Sch{\"u}tt, \emph{112 lines on smooth quartic
  surfaces (characteristic 3)}, Q. J. Math. \textbf{66} (2015), no.~3,
  941--951. \MR{3396099}

\bibitem{rams.schuett}
\bysame, \emph{64 lines on smooth quartic surfaces}, Math. Ann. \textbf{362}
  (2015), no.~1-2, 679--698. \MR{3343894}

\bibitem{rams.schuett:char2}
\bysame, \emph{At most 64 lines on smooth quartic surfaces (characteristic 2)},
  Nagoya Math. J. \textbf{232} (2018), 76--95. \MR{3866501}

\bibitem{Saint-Donat}
B.~Saint-Donat, \emph{Projective models of {$K$-$3$} surfaces}, Amer. J. Math.
  \textbf{96} (1974), 602--639. \MR{0364263 (51 \#518)}

\bibitem{Schur:quartics}
Friedrich Schur, \emph{Ueber eine besondre {C}lasse von {F}l\"achen vierter
  {O}rdnung}, Math. Ann. \textbf{20} (1882), no.~2, 254--296. \MR{1510168}

\bibitem{Segre}
B.~Segre, \emph{The maximum number of lines lying on a quartic surface}, Quart.
  J. Math., Oxford Ser. \textbf{14} (1943), 86--96. \MR{0010431 (6,16g)}

\bibitem{Shimada:X56}
Ichiro Shimada and Tetsuji Shioda, \emph{On a smooth quartic surface containing
  56 lines which is isomorphic as a {K}3 surface to the {F}ermat quartic},
  Manuscripta Math. \textbf{153} (2017), no.~1-2, 279--297. \MR{3635983}

\bibitem{GRAPE}
Leonard~H. Soicher, \emph{{GRAPE}, {GR}aph {A}lgorithms using {PE}rmutation
  groups, {V}ersion 4.8.1}, \href {https://gap-packages.github.io/grape}
  {https://gap-packages.github.io/grape}, Oct 2018, Refereed GAP package.

\bibitem{Veniani:char2}
Davide~Cesare Veniani, \emph{Lines on {K}3 quartic surfaces in characteristic
  2}, Q. J. Math. \textbf{68} (2017), no.~2, 551--581. \MR{3667213}

\bibitem{Veniani}
\bysame, \emph{The maximum number of lines lying on a {K}3 quartic surface},
  Math. Z. \textbf{285} (2017), no.~3-4, 1141--1166. \MR{3623744}

\bibitem{Veniani:equations}
\bysame, \emph{Symmetries and equations of smooth quartic surfaces with many
  lines}, To appear, \arXiv{1708.01219}, 2017.

\end{thebibliography}
}

\endgroup
\end{document}